\documentclass[12pt,a4paper]{article}
\usepackage{fontspec}
\usepackage[margin=25mm,headheight=27pt]{geometry}
\usepackage{amsmath,amssymb,amsthm,mathtools,bm,mathrsfs}
\usepackage{booktabs,longtable,array,enumitem,microtype}
\usepackage[dvipsnames]{xcolor}
\usepackage{sectsty,needspace}
\allsectionsfont{\raggedright}
\usepackage{hyperref,fancyhdr}
\hypersetup{colorlinks=true,linkcolor=blue,citecolor=blue,urlcolor=blue,
 pdftitle={Classification and C2-cofiniteness of rational vertex operator algebras generated in low conformal weights},pdfauthor={Jiancai Sun and Ruochen Yang}}
\renewcommand{\refname}{References}
\newtheorem{theorem}{Theorem}[section]
\newtheorem{proposition}[theorem]{Proposition}
\newtheorem{lemma}[theorem]{Lemma}
\newtheorem{corollary}[theorem]{Corollary}
\theoremstyle{definition}\newtheorem{definition}[theorem]{Definition}
\newtheorem{example}[theorem]{Example}
\theoremstyle{remark}\newtheorem{remark}[theorem]{Remark}

\newcommand{\C}{\mathbb C}
\newcommand{\R}{\mathbb R}
\newcommand{\Z}{\mathbb Z}
\newcommand{\N}{\mathbb Z_{\geq0}}
\newcommand{\one}{\mathbf1}
\newcommand{\wt}{\operatorname{wt}}
\newcommand{\Res}{\operatorname{Res}}
\newcommand{\Span}{\operatorname{span}_{\C}}
\newcommand{\id}{\operatorname{id}}
\newcommand{\End}{\operatorname{End}}
\newcommand{\Aut}{\operatorname{Aut}}
\newcommand{\Sym}{\operatorname{Sym}}
\newcommand{\Tr}{\operatorname{Tr}}
\newcommand{\Soc}{\operatorname{Soc}}
\newcommand{\Ann}{\operatorname{Ann}}
\newcommand{\rad}{\operatorname{rad}}

\newcommand{\g}{\mathfrak g}
\newcommand{\q}{\mathfrak q}
\newcommand{\s}{\mathfrak s}
\newcommand{\m}{\mathfrak m}
\newcommand{\vir}{\mathfrak{vir}}
\numberwithin{equation}{section}
\setlist{itemsep=3pt,topsep=5pt}
\setlist[enumerate]{label=(\roman*)}
\begin{document}
\begin{center}
{\Large\bfseries Classification and $C_2$-cofiniteness of rational vertex\\[3pt]
operator algebras generated in low conformal weights\par}
\end{center}
\vspace{6pt}
\begingroup
\renewcommand{\thefootnote}{}
\endgroup
\begin{center}
Jiancai Sun\\
Department of Mathematics and Newtouch Center for Mathematics\\
Shanghai University, Shanghai 200444, China\\
E-mail address: \href{mailto:jcsun@shu.edu.cn}{jcsun@shu.edu.cn}
\end{center}
\vspace{0pt}
\begin{center}
Ruochen Yang\\
Department of Mathematics\\
Shanghai University, Shanghai 200444, China\\
E-mail address: \href{mailto:yangruochen@shu.edu.cn}{yangruochen@shu.edu.cn}
\end{center}
\thispagestyle{plain}
\begin{abstract}
We obtain classification results for rational vertex operator algebras generated in low conformal weights.
More precisely, we prove that every rational vertex operator algebra of OZ type strongly generated by its conformal weight two subspace $V_2$ is isomorphic to a finite tensor product of Virasoro minimal models.
We also prove that every rational vertex operator algebra of CFT type generated as a vertex algebra by its conformal weight one subspace $V_1$ is isomorphic to a finite tensor product of simple affine vertex operator algebras associated with finite-dimensional simple Lie algebras at positive integral levels.
Consequently, Zhu's conjecture that rationality implies $C_2$-cofiniteness holds for both classes, and every such rational vertex operator algebra is regular.
\end{abstract}
\clearpage
{\linespread{0.92}\small\tableofcontents}
\clearpage

\section{Introduction}\label{sec:intro}

\Needspace{3\baselineskip}
\subsection{Background}\label{subsec:background}
Complete reducibility of modules is one of the fundamental problems in the representation theory of vertex operator algebras \cite{Zhu,DLMreg,DLMtw}.
The notion of rationality was introduced in \cite[Definition~1.2.4]{Zhu}; we adopt the definition of \cite[Section~2]{DLMreg},
namely, that every admissible module is completely reducible.
Regularity was introduced in \cite[Definition~3.1]{DLMreg} and requires every weak module to be a direct sum of irreducible ordinary modules.
It is also proved in \cite{DLMreg} that basic examples, including simple affine vertex operator algebras at positive integral levels, Virasoro minimal models,
and vertex operator algebras associated with positive-definite even lattices, are regular.

Another condition closely related to complete reducibility is $C_2$-cofiniteness:
\[
 C_2(V)=\Span\{u_{-2}v:u,v\in V\},\qquad
 \dim V/C_2(V)<\infty.
\]
This condition was introduced in \cite[Section~4.4]{Zhu} in the study of modular invariance of characters of vertex operator algebras.
By \cite[Theorem~3.8]{Li}, regularity implies $C_2$-cofiniteness,
while \cite[Theorem~4.5]{ABD} shows that, for vertex operator algebras of CFT type,
rationality and $C_2$-cofiniteness together imply regularity. Hence, for vertex operator algebras of CFT type,
\begin{equation}\label{eq:introregular}
 \text{regularity}\quad\Longleftrightarrow\quad
 \text{rationality and }C_2\text{-cofiniteness}.
\end{equation}

A fundamental question is whether rationality implies $C_2$-cofiniteness.
This question goes back to \cite[Section~4.4]{Zhu}, where it is conjectured that every rational vertex operator algebra satisfies the finiteness condition C,
which includes $C_2$-cofiniteness; see \cite[Section~1]{DLMmod}.
It is conjectured in \cite[Conjecture~3.17]{DLMreg} that every rational vertex operator algebra is regular;
by \eqref{eq:introregular}, for vertex operator algebras of CFT type this conjecture is precisely equivalent to the statement that rationality implies $C_2$-cofiniteness; see \cite[Section~1]{Li}.
To the best of our knowledge, this question remains open. The converse does not hold: the even part $SF^+$ of the symplectic fermionic vertex operator superalgebra $SF$
is $C_2$-cofinite \cite[Theorem~3.13]{AbeSF} but not rational \cite[Proposition~5.3]{AbeSF}; see \cite{DY}.

Twisted modules arise naturally in orbifold theory.
Orbifold constructions and twisted sectors originated in the physics literature on conformal field theory \cite{DHVW,DVVV},
and their mathematical theory can be found in \cite{DLMtw,DLMmod,Miy,CM,DRX}.
For a finite-order automorphism $g$, weak, admissible, and ordinary modules have corresponding $g$-twisted versions,
leading to the notions of $g$-rationality \cite{DLMtw} and $g$-regularity \cite[Definition~2.5]{Yamauchi}.
On the other hand, unitary vertex operator algebras were introduced in \cite[Definition~2.2]{DL}, where it is also proved that most of the familiar rational vertex operator algebras are unitary.
For unitary vertex operator algebras of CFT type, \cite[Theorem~1.1]{SY} proves, for every finite-order automorphism $g$, that
\[
 g\text{-regularity}\quad\Longleftrightarrow\quad
 g\text{-rationality and }C_2\text{-cofiniteness}.
\]
All of the above results treat $C_2$-cofiniteness as a condition alongside rationality.

In this paper we study the relation between rationality and $C_2$-cofiniteness from a different point of view.
For a vertex operator algebra $V$ generated by a subspace of low conformal weight,
the algebraic structure on the generating subspace largely determines $V$ itself \cite{Lian,Lam,KL,LYu}.
If these algebraic structures can be determined from rationality, then $V$ can be identified with a tensor product of vertex operator algebras
already known to be rational and $C_2$-cofinite, so that $C_2$-cofiniteness and regularity follow directly from rationality.

\Needspace{3\baselineskip}
\subsection{The problem and the main results}
For a vertex operator algebra $V=\bigoplus_{n\in\Z}V_n$, we denote by $V_n$ the subspace of conformal weight $n$.
We study two classes of vertex operator algebras generated by subspaces of low conformal weight;
the main classification theorems are based on the commutative algebra structure on $V_2$ and on the Lie algebra structure on $V_1$, respectively.

The first class consists of simple OZ-type vertex operator algebras strongly generated by $V_2$,
where OZ type means CFT type with $V_1=0$ \cite[Definition~2.7]{GNAVOA}.
In this case, $V_2$ is a commutative associative algebra under the multiplication $ab=\frac12a_1b$ and carries a symmetric invariant bilinear form
\cite[Propositions~3.7, 3.9]{LYu}; this is the commutative Frobenius algebra of Section~\ref{subsec:frob-weight-two}.
For a finite-order automorphism $g\in\Aut(V)$, besides rationality we also consider $g$-rationality and the finiteness of the set of isomorphism classes of irreducible ordinary $g$-twisted modules
\cite[Theorem~8.1~(b), (c)]{DLMtw}.
The first main theorem shows that, for this class of vertex operator algebras, these three conditions are equivalent to one another,
and each of them suffices to determine the commutative Frobenius algebra structure on $V_2$,
thereby identifying $V$ with a finite tensor product of Virasoro minimal models \cite[Theorem~4.2]{Wang}.

\begin{theorem}\label{thm:two}
Let $V$ be a simple OZ-type vertex operator algebra strongly generated by $V_2$,
and let $g\in\Aut(V)$ be an automorphism of finite order. The following conditions are equivalent:
\begin{enumerate}
\item there are only finitely many isomorphism classes of irreducible ordinary $g$-twisted $V$-modules;
\item $V$ is $g$-rational;
\item $V$ is rational;
\item there is a vertex operator algebra isomorphism
\[
 V\cong\bigotimes_{i=1}^r L_{\mathfrak{vir}}(c_{p_i,q_i},0),
\]
where $r\in\Z_{\ge0}$, the $p_i,q_i$ are integers satisfying
\[
 2\le p_i<q_i,\qquad \gcd(p_i,q_i)=1
 \quad(1\le i\le r),
\]
and
\[
 c_{p_i,q_i}=1-\frac{6(p_i-q_i)^2}{p_iq_i};
\]
each $L_{\mathfrak{vir}}(c_{p_i,q_i},0)$ is a Virasoro minimal model.
\end{enumerate}
For $r=0$, the right-hand side is taken to be $\C\one$.
When these equivalent conditions hold, $V$ is $C_2$-cofinite and regular.
\end{theorem}

If a nonzero vertex operator algebra $V$ is strongly generated by $V_2$,
then $V$ is of OZ type by Remark~\ref{rem:strong-generation}~(1);
if, in addition, $V$ is rational, then $V$ is simple by Lemma~\ref{lem:simple}.
Hence Theorem~\ref{thm:two} immediately yields the following corollary.

\begin{corollary}\label{cor:untwot}
Let $V$ be a nonzero vertex operator algebra. If $V$ is rational and strongly generated by $V_2$,
then $V$ is a finite tensor product of Virasoro minimal models,
and is $C_2$-cofinite and regular.
\end{corollary}

Corollary~\ref{cor:untwot} may be viewed as an extension of \cite[Theorem~5.2.21]{LamT} and \cite[Theorem~4.7]{ZD} in the following sense.
For the former, instead of assuming that $V_2$ is semisimple as a commutative associative algebra, we assume rationality and deduce the semisimplicity of $V_2$.
For the latter, we omit the assumptions that $V$ is $C_2$-cofinite, that $\dim V_2=2$, and that its central charge and effective central charge (see \cite[(1.2)]{DMec}) are both equal to one, while adding the assumption that $V$ is strongly generated by $V_2$.

The second class consists of rational vertex operator algebras of CFT type that are generated by $V_1$ as a vertex algebra.
In this case, $V_1$ is a Lie algebra under $[x,y]=x_0y$ and carries the symmetric invariant bilinear form given by $x_1y$
\cite[Lemmas~3.4--3.5, Theorem~3.7]{Lian}.
By Lemma~\ref{lem:simple}, such a $V$ is simple; hence this bilinear form is nondegenerate \cite[Corollary~5.3~(iii)]{Lian},
and $V_1$ together with this form is a quadratic Lie algebra in the sense of Section~\ref{subsec:currents} \cite[Section~5]{TW}, \cite[Section~4.1]{AleM};
moreover, \cite[Theorem~5.4~(i)]{Lian} identifies $V$ with the simple graded quotient of the corresponding affine vertex algebra.
The second main theorem shows that, for this class of vertex operator algebras, rationality forces $V_1$ to be a semisimple Lie algebra
and forces the restriction of the above bilinear form to each simple ideal to be a positive integer multiple of the normalized Killing form \cite[Remark~6.2.15]{LL};
hence the above simple graded quotient is a finite tensor product of simple affine vertex operator algebras
\cite[Theorem~3.1.3]{FZ} at positive integral levels.

\begin{theorem}\label{thm:one}
Let $V$ be a rational vertex operator algebra of CFT type that is generated by $V_1$ as a vertex algebra.
Then $V_1$ is a semisimple Lie algebra with respect to $[x,y]=x_0y$, $x,y\in V_1$.
Let $\s_1,\ldots,\s_r$ be all its simple ideals, where $r\in\N$. Then there exist $k_1,\ldots,k_r\in\Z_{>0}$
and a vertex operator algebra isomorphism
\[
 V\cong\bigotimes_{i=1}^r L_{(\s_i,(\, ,\,)_i)}(k_i,0),
\]
where, for $1\le i\le r$, $(\, ,\,)_i$ is the normalized Killing form of $\s_i$, that is, the symmetric invariant bilinear form on $\s_i$ for which long roots have squared length $2$,
and $L_{(\s_i,(\, ,\,)_i)}(k_i,0)$ is the simple affine vertex operator algebra of level $k_i$;
for $r=0$, the right-hand side is taken to be $\C\one$. Consequently, $V$ is $C_2$-cofinite and regular.
\end{theorem}

Theorem~\ref{thm:one} strengthens the classification in \cite[Corollary~5.8~(a)]{DMinteg} by removing the assumptions of $C_2$-cofiniteness and strong CFT type.
By Theorems~\ref{thm:two} and~\ref{thm:one}, the conjecture recalled in Section~\ref{subsec:background} that rationality implies $C_2$-cofiniteness
holds for these two classes of vertex operator algebras.

\Needspace{3\baselineskip}
\subsection{Related work}
In \cite{Lam}, vertex operator algebras are constructed from commutative associative algebras and their invariant bilinear forms;
in the semisimple case, the corresponding simple quotient is a tensor product of simple Virasoro vertex operator algebras,
and its automorphisms are determined by the corresponding algebraic structure \cite[Theorems~5.2.21, 5.3.1]{LamT}.
The characterization of the conformal weight two structure used in this paper is given in \cite[Propositions~3.7, 3.9]{LYu},
where the irreducible graded modules of the corresponding vertex operator algebra are also classified \cite[Theorem~5.5, Corollary~5.6]{LYu}.
These results are the starting point of this paper in the conformal weight two case.

Several concrete precedents exist in the nonsemisimple case.
For the $W$-algebra $W(2,2)$, \cite[Theorems~2.1, 2.3]{ZD} proves the irreducibility of the corresponding vacuum quotient
and gives the complete list $L(c,h_1,h_2)$, $h_1,h_2\in\C$, of its irreducible modules, so that there are infinitely many isomorphism classes;
in our notation, this corresponds to the commutative Frobenius algebra $(A,f)$ with $A=\C[\varepsilon]/(\varepsilon^2)$ and $f(1_A)=f(\varepsilon)=c\in\C\setminus\{0\}$.
For $A=\C[\varepsilon]/(\varepsilon^N)$, $N\in\Z_{\ge2}$, \cite[Theorem~3.5]{JZ} proves that the vertex operator algebra corresponding to a nondegenerate Frobenius form
is simple, and \cite[Proposition~3.7, Remark~3.8]{JZ} also determines its irreducible admissible modules
and shows that it is irrational.
We extend both aspects to an arbitrary nonsemisimple local commutative Frobenius algebra $(A,f)$:
the corresponding vertex operator algebra is always simple (Theorem~\ref{thm:frobsimple}),
and for every finite-order automorphism of $(A,f)$ there are infinitely many pairwise non-isomorphic irreducible ordinary twisted modules (Proposition~\ref{prop:infinite}).
This is the key to the proof of (i)$\Rightarrow$(iv) in Theorem~\ref{thm:two}.
When $g$ in Theorem~\ref{thm:two} permutes isomorphic tensor factors, the required permutation-twisted modules and permutation rationality
are taken directly from \cite[Theorems~4.4, 6.4]{BDM}.

In the conformal weight one case, the corresponding affine vertex algebras and their simple graded quotients are studied in \cite[Sections~4--6]{Lian},
where their conformal vectors (called Virasoro elements in \cite{Lian}) are also classified.
Under the assumption of strong rationality (that is, rationality, $C_2$-cofiniteness, and strong CFT type, where strong CFT type means CFT type with $L(1)V_1=0$), \cite[Theorem~1]{DMec} proves that $V_1$ is reductive.
For a simple, rational, $C_2$-cofinite vertex operator algebra $V$ of strong CFT type, \cite[Corollary~5.8~(a)]{DMinteg} shows that $V$ is generated by $V_1$ as a vertex algebra if and only if it is isomorphic to a finite tensor product of simple affine vertex operator algebras associated with finite-dimensional simple Lie algebras at positive integral levels.
The hypotheses of $C_2$-cofiniteness and strong CFT type in these results are not assumed in Theorem~\ref{thm:one}.
A classification in the unitary case is given in \cite[Theorem~4.10]{AL}:
a unitary vertex operator algebra of CFT type generated by $V_1$ as a vertex algebra is a tensor product of
simple affine vertex operator algebras at positive integral levels and a Heisenberg vertex operator algebra;
if rationality is also assumed, then, since Heisenberg vertex operator algebras have infinitely many pairwise non-isomorphic irreducible ordinary modules,
the Heisenberg factor is excluded, which gives the unitary special case of Theorem~\ref{thm:one}.
To treat the quadratic summands without nonzero semisimple ideals, our general proof follows the scaling idea of \cite[proof of Proposition~2]{Sch}, using scalings that preserve the bilinear form;
these summands are then excluded by the finite dimensionality of the Zhu algebra.

Many important vertex operator algebras are generated by vectors of low conformal weight.
For example, the lattice vertex operator algebra $V_{E_8}$ associated with the $E_8$ root lattice is generated by $V_1$ as a vertex algebra \cite[Section~5.6]{Kac}
and is regular \cite[Theorem~3.16]{DLMreg}; by \cite[Theorem~5.6~(b)]{Kac}, it is isomorphic to the simple affine vertex operator algebra of level $1$ associated with $\mathfrak e_8$,
which is consistent with Theorem~\ref{thm:one}.
As another example, the simple vertex operator algebras associated with the $W$-algebra $W(2,2)$ in \cite{ZD} and with the truncated Virasoro algebras in \cite{JZ}
are simple OZ-type vertex operator algebras strongly generated by $V_2$, and hence satisfy the hypotheses of Theorem~\ref{thm:two};
both are irrational: the irreducible modules $L(c,h_1,h_2)$, $h_1,h_2\in\C$, of the former fall into infinitely many isomorphism classes
\cite[Theorem~2.3~(4)]{ZD}, and so do the irreducible admissible modules of the latter
\cite[Proposition~3.7~(3), Remark~3.8]{JZ}, which is consistent with Theorem~\ref{thm:two}.
Note that Theorem~\ref{thm:two} requires strong generation: the moonshine module $V^\natural$ is generated by $V^\natural_2$,
but $V^\natural_2$ is a nonassociative algebra \cite[Section~1]{DGL}, so by Proposition~\ref{prop:wt2}, $V^\natural$ is not strongly generated by $V^\natural_2$;
similarly, the OZ-type vertex operator algebras generated by Ising vectors of $\sigma$-type in \cite{JLY} are all rational,
and by Corollary~\ref{cor:untwot}, those among them that are not finite tensor products of Virasoro minimal models are not strongly generated by $V_2$.
On the other hand, our results cover neither vertex operator algebras that require generators of higher conformal weight
nor the case where generators of conformal weights $1$ and $2$ are needed simultaneously.
An example of the former is the parafermion vertex operator algebra $K(\mathfrak{sl}_2,k)$ with $k\in\Z_{\ge3}$, which is generated by vectors of conformal weights $2$ and $3$
\cite[Theorem~4.1]{DLWY} and is $C_2$-cofinite \cite[Theorem~10.5]{ALY1} and rational \cite[Corollary~6.2]{ALY2}.
An example of the latter is the Haagerup vertex operator algebra of central charge $8$ \cite{Riesen,HuangHaag},
which is generated by vectors of conformal weights $1$ and $2$ and is strongly rational \cite[Theorem~1]{Riesen}.

\Needspace{3\baselineskip}
\subsection{Outline of the proof and organization of the paper}
The proof of Theorem~\ref{thm:two} consists of two steps.
The first step determines the tensor product decomposition of $V$: after the commutative Frobenius algebra corresponding to $V_2$ is decomposed into local summands $(A_i,f_i)$, Proposition~\ref{prop:frobfactors} shows that
$V$ is the tensor product of the simple quotients corresponding to the $(A_i,f_i)$.
For $A_i\not\cong\C$, Section~\ref{sec:frobgram} uses graph contractions to compare, in each degree, the bilinear form pulled back from the normalized invariant bilinear form
with an auxiliary bilinear form depending only on the Frobenius form, and proves that their Gram determinants are equal and nonzero;
hence the vertex operator algebra corresponding to $(A_i,f_i)$ is already simple, and the corresponding factor is this vertex operator algebra itself; see Theorem~\ref{thm:frobsimple}.
For $A_i\cong\C$, writing $c_i=f_i(1_{A_i})$, we see from \eqref{eq:vir-special} that the corresponding factor is precisely
the simple Virasoro vertex operator algebra $L_{\mathfrak{vir}}(c_i,0)$ of central charge $c_i$.
In the second step, for all factors that are not Virasoro minimal models, we use \cite[Theorem~3.6]{CDW}
to convert restricted modules into twisted modules, and construct, by highest-weight induction,
irreducible ordinary twisted modules whose lowest conformal weights can take arbitrary values; see Proposition~\ref{prop:infinite}.
Hence the finiteness of the set of isomorphism classes of irreducible ordinary $g$-twisted modules allows only Virasoro minimal models as tensor factors;
the remaining implications and regularity follow from Propositions~\ref{prop:vir}, \ref{prop:tensor}, \ref{prop:permutation}, and~\ref{prop:knownC2}, together with Remark~\ref{rem:module-relations}~(3).

The proof of Theorem~\ref{thm:one} starts from Proposition~\ref{prop:lian}. Write $\g=V_1$, and let $B_{\g}$ be the bilinear form given by $x_1y$;
then $V$ is isomorphic, as a vertex algebra, to the simple graded quotient of the corresponding affine vertex algebra $V_{(\g,B_{\g})}(1,0)$ introduced in Definition~\ref{def:affine-vacuum}.
By Lemma~\ref{lem:split}, $\g=\s\perp\q$, where $\s$ is semisimple and $\q$ has no nonzero semisimple ideal.
The successive polynomial limits of Proposition~\ref{prop:degenerate} reduce the Lie bracket of $\q$ to zero while keeping the bilinear form fixed.
The Gram determinant of the normalized contravariant form in each degree is unchanged at every step;
hence $V_{(\q,B_{\q})}(1,0)$ has no nonzero proper graded vertex ideals; see Theorem~\ref{thm:currentsimple}.
We then decompose $\s$ into simple ideals $\s_1,\ldots,\s_r$, and write the restriction of $B_{\g}$ to each $\s_i$ as $k_i$ times the normalized Killing form.
We obtain a grading-preserving tensor product decomposition of $V$ from Proposition~\ref{prop:curfactors},
and an injection $U(\q)\to A(V)$ from Proposition~\ref{prop:injection}, where $A(V)$ is the Zhu algebra.
Rationality gives $\dim A(V)<\infty$ by Proposition~\ref{prop:zhu}, so $\q=0$.
Next, using the Casimir elements of the $\s_i$, we determine the conformal vector of $V$ and prove that no level $k_i$ is the critical level
defined in Definition~\ref{def:affine-sugawara}, so that the above decomposition is a vertex operator algebra isomorphism; see Proposition~\ref{prop:sugawara}.
Finally, Propositions~\ref{prop:tensor}~(1) and~\ref{prop:affineclass} give $k_i\in\Z_{>0}$,
and hence, by Propositions~\ref{prop:tensor}~(2) and~\ref{prop:knownC2}, $V$ is $C_2$-cofinite and regular.

The paper is organized according to the logical dependence of the above proofs.
Section~\ref{sec:prelim} recalls the basic definitions and results on vertex operator algebras, twisted modules, $C_2$-cofiniteness,
invariant bilinear forms, and the Zhu algebra.
Sections~\ref{sec:frob}--\ref{sec:twists} prove Theorem~\ref{thm:two}:
Section~\ref{sec:frob} obtains from $V_2$ the commutative Frobenius algebra and the corresponding tensor product decomposition,
Section~\ref{sec:frobgram} uses graph contractions to prove that the vertex operator algebras corresponding to nonsemisimple local factors are simple,
and Section~\ref{sec:twists} constructs the required twisted modules and completes the proof.
Sections~\ref{sec:current}--\ref{sec:affine} prove Theorem~\ref{thm:one}:
Section~\ref{sec:current} studies quadratic Lie algebras and graded vertex ideals of the corresponding affine vertex algebras,
and Section~\ref{sec:affine} completes the proof using the Zhu algebra and Casimir elements.

\Needspace{3\baselineskip}
\subsection{Acknowledgment}
J.-C. Sun acknowledges support from the National Natural Science Foundation of China (Nos.~12071276, 11931009, and 12226402) and the Science and Technology Commission of Shanghai Municipality (No.~25ZR1401126).

\section{Preliminaries}\label{sec:prelim}
In this section, we recall the basic definitions and results needed later and fix our notation.
Section~\ref{subsec:voa} reviews vertex operator algebras, vertex ideals, and automorphisms following \cite{FHL,Kac,LL},
and fixes our conventions for CFT type, OZ type, generation, and strong generation.
Section~\ref{subsec:twisted-modules} recalls weak, admissible, and ordinary twisted modules and $g$-rationality as in \cite{DLMtw,SY},
and shows that a rational vertex operator algebra of CFT type is simple; see Lemma~\ref{lem:simple}.
Section~\ref{subsec:c2} recalls $C_2$-cofiniteness and its relation to regularity, as in \cite{Zhu,ABD},
and records the rationality and $C_2$-cofiniteness of finite tensor products.
In Section~\ref{subsec:invariant-forms}, we review invariant bilinear forms following \cite{LiForm,SY}, and record their radicals,
a nondegeneracy criterion in terms of Gram determinants in each conformal weight, and the corresponding results for finite tensor products;
these will be used to decide simplicity and to establish tensor product decompositions.
Section~\ref{subsec:zhu} recalls the Zhu algebra and its action on the top levels of admissible modules \cite{Zhu,DLMtw};
these results will be used in Proposition~\ref{prop:injection} of Section~\ref{sec:affine}.

\Needspace{3\baselineskip}
\subsection{Vertex operator algebras, automorphisms, and generation}\label{subsec:voa}
All vector spaces in this paper are over $\C$.
We first recall the definition of a vertex operator algebra, following \cite{FHL,Kac,LL}.

\begin{definition}\label{def:vertex-operator-algebra}
A \emph{vertex operator algebra} is a quadruple $(V,Y,\one,\omega)$, where
$V$ is a $\Z$-graded vector space such that
\[
 V=\bigoplus_{d\in\Z}V_d,\qquad
 \dim V_d<\infty,\qquad V_d=0\quad\text{for $d$ sufficiently small};
\]
$Y$ is a linear map
\[
 \begin{aligned}
 Y(\cdot,z):V&\longrightarrow(\End V)[[z,z^{-1}]],\\
 a&\longmapsto Y(a,z)=\sum_{n\in\Z}a_nz^{-n-1},
 \end{aligned}
\]
with $a_n\in\End V$ for $a\in V$ and $n\in\Z$;
$\one\in V_0$ is the \emph{vacuum vector} and $\omega\in V_2$ is the \emph{conformal vector}.
These data are subject to the following conditions:
\begin{enumerate}[label=(\arabic*)]
\item For any $a,b\in V$, one has $a_nb=0$ for $n$ sufficiently large.
\item The vacuum property and the creation property hold:
\[
 Y(\one,z)=\id_V,\qquad
 Y(a,z)\one\in V[[z]],\qquad
 \left.Y(a,z)\one\right|_{z=0}=a,\quad a\in V.
\]
\item For any $a,b\in V$, the Jacobi identity holds:
\begin{equation*}
 \begin{aligned}
 &z_0^{-1}\delta\!\left(\frac{z_1-z_2}{z_0}\right)Y(a,z_1)Y(b,z_2)
 -z_0^{-1}\delta\!\left(\frac{z_2-z_1}{-z_0}\right)Y(b,z_2)Y(a,z_1)\\
 &\hspace{12mm}=
 z_2^{-1}\delta\!\left(\frac{z_1-z_0}{z_2}\right)Y(Y(a,z_0)b,z_2).
 \end{aligned}
\end{equation*}
Here $\delta(z)=\sum_{n\in\Z}z^n$, and each binomial expression is expanded in nonnegative integral powers of its second variable.
\item Write
\[
 Y(\omega,z)=\sum_{n\in\Z}L(n)z^{-n-2},\qquad L(n)=\omega_{n+1},
\]
where $L(n)\in\End V$ for $n\in\Z$. There exists $c\in\C$ such that
\begin{equation*}
 [L(m),L(n)]=(m-n)L(m+n)
 +\frac{m^3-m}{12}\delta_{m+n,0}\,c\,\id_V,
 \qquad m,n\in\Z.
\end{equation*}
The number $c$ is called the \emph{central charge} of $V$. Moreover, $L(0)a=da$ for $d\in\Z$ and $a\in V_d$, and
\begin{equation}\label{eq:voa-derivative}
 Y(L(-1)a,z)=\partial_zY(a,z),\qquad a\in V.
\end{equation}
The operator $L(-1)$ is called the \emph{translation operator} of $V$.
\end{enumerate}
For a nonzero homogeneous vector $a\in V_d$, the integer $d$ is called the \emph{conformal weight} of $a$, and we write $\wt a=d$.
\end{definition}

More generally, let $V$ be a vector space, let $Y$ be a linear map of the above form, and let $\one\in V$.
If these data satisfy conditions (1)--(3) of Definition~\ref{def:vertex-operator-algebra},
then the triple $(V,Y,\one)$ is called a \emph{vertex algebra}; no grading or conformal vector is required here
\cite[Definition~3.1.1]{LL}. Forgetting the conformal vector and the grading,
one obtains from every vertex operator algebra its \emph{underlying vertex algebra}.
For a vertex algebra $(V,Y,\one)$, the operator $D\in\End V$ defined by $Da=a_{-2}\one$ for $a\in V$ is called the \emph{translation operator} of $V$ \cite[Proposition~3.1.18, (3.1.67)]{LL}.
When $V$ is a vertex operator algebra, $D=L(-1)$.

Unless otherwise stated, we fix for the rest of this subsection a vertex operator algebra $(V,Y,\one,\omega)$
together with its conformal weight grading $V=\bigoplus_{n\in\Z}V_n$.

\begin{remark}\label{rem:basic-identities}
The following basic formulas, which will be used repeatedly, can be found in \cite[Section~3.1]{LL}.
For homogeneous $a,b\in V$ and $n\in\Z$ with $a_nb\ne0$, we have
\begin{equation}\label{eq:weights}
 \wt(a_nb)=\wt a+\wt b-n-1.
\end{equation}
For any $a,b\in V$ and $m,n\in\Z$, we have the commutator formula
\begin{equation}\label{eq:mode-commutator}
 [a_m,b_n]=\sum_{j=0}^{\infty}\binom mj(a_jb)_{m+n-j},
\end{equation}
where, for $\alpha\in\C$ and $j\in\N$, we put $\binom{\alpha}{j}=\frac{\alpha(\alpha-1)\cdots(\alpha-j+1)}{j!}$, with the convention $\binom{\alpha}{0}=1$;
by truncation, only finitely many terms on the right-hand side are nonzero. We also have the skew-symmetry formula
\begin{equation}\label{eq:skew-symmetry}
 Y(a,z)b=e^{zL(-1)}Y(b,-z)a,\qquad a,b\in V.
\end{equation}
\end{remark}

For $a\in V$ and $n\in\Z$, the coefficient
$a_n\in\End V$ of $z^{-n-1}$ in $Y(a,z)$ is called the \emph{$n$-th mode} of $a$.
For $b\in V$, the vector $a_nb$ is called the \emph{$n$-th product} of $a$ and $b$.

For CFT type we adopt the convention of \cite[Section~2]{ABD}; see also \cite{KL} and \cite[Definition~2.3]{SY}.
OZ type (``One-Zero'' type) was introduced by Griess; see \cite[Definition~2.7]{GNAVOA} and also \cite[Section~1]{JLY}.

\begin{definition}\label{def:cft-oz}
A vertex operator algebra $V$ is said to be of \emph{CFT type} if
\[
 V_0=\C\one\ne0,\qquad V_n=0,\quad n\in\Z_{<0}.
\]
It is said to be of \emph{OZ type} if it is of CFT type and $V_1=0$.
\end{definition}

Ideals and simplicity are introduced as in \cite[Definition~3.9.7, Remark~3.9.8]{LL}.

\begin{definition}\label{def:ideal-simple}
A \emph{vertex ideal} of a vertex algebra $V$ is a subspace $I\subseteq V$
that is stable under all modes $a_n$ ($a\in V$, $n\in\Z$) and under the translation operator $D$, that is,
\[
 DI\subseteq I,\qquad a_nI\subseteq I,\quad a\in V,\ n\in\Z.
\]
We say that $V$ is \emph{simple} if $V\ne0$ and its only ideals are $0$ and $V$.
\end{definition}

For homomorphisms and automorphisms we use the conventions of \cite[Section~3.9]{LL}.

\begin{definition}\label{def:vertex-homomorphism}\label{def:operator-homomorphism}
Let $(V,Y_V,\one_V)$ and $(W,Y_W,\one_W)$ be vertex algebras.
A linear map $f\colon V\to W$ is called a \emph{vertex algebra homomorphism} if $f(\one_V)=\one_W$ and
\begin{equation}\label{eq:vertex-homomorphism}
 f(a_nb)=(f(a))_nf(b),\qquad a,b\in V,\ n\in\Z.
\end{equation}
Here the modes in \eqref{eq:vertex-homomorphism} are taken with respect to $Y_V$ and $Y_W$, respectively, and \eqref{eq:vertex-homomorphism} is equivalent to
\[
 f\bigl(Y_V(a,z)b\bigr)=Y_W(f(a),z)f(b),\qquad a,b\in V.
\]
If, in addition, $V$ and $W$ are vertex operator algebras with conformal vectors
$\omega_V$ and $\omega_W$, respectively, and $f$ satisfies
\begin{equation}\label{eq:conformal-homomorphism}
 f(\omega_V)=\omega_W,
\end{equation}
then $f$ is called a \emph{vertex operator algebra homomorphism}.
In either case, a bijective homomorphism is called an \emph{isomorphism} of the corresponding kind.
\end{definition}

\begin{definition}\label{def:automorphism}
An \emph{automorphism} of a vertex operator algebra $V$ is a
vertex operator algebra isomorphism from $V$ to itself. Under composition of maps, all such automorphisms form a group,
denoted by $\Aut(V)$.
\end{definition}

By \eqref{eq:vertex-homomorphism} and \eqref{eq:conformal-homomorphism},
every $g\in\Aut(V)$ commutes with $L(n)$ for $n\in\Z$, and hence preserves each $V_d$, $d\in\Z$.
We next treat subalgebras and generation together, distinguishing whether or not the subalgebra is required to contain the original conformal vector
\cite[Definitions~3.9.1, 3.9.2, 3.9.5]{LL}.

\begin{definition}\label{def:generation}
A \emph{vertex subalgebra} of a vertex operator algebra $V$ is a subspace $U\subseteq V$ such that
\[
 \one\in U,\qquad a_nb\in U,\qquad a,b\in U,\ n\in\Z.
\]
In this case $L(-1)U\subseteq U$, since $L(-1)a=a_{-2}\one$ for $a\in U$.
If, moreover, $\omega\in U$, then $U$ is called a \emph{vertex operator subalgebra}.

Given $S\subseteq V$, the \emph{vertex subalgebra generated by $S$}
is the smallest vertex subalgebra containing $S$, denoted by $\langle S\rangle$.
If $\langle S\rangle=V$, we say that $S$ \emph{generates} $V$.
The \emph{vertex operator subalgebra generated by $S$}
is the smallest vertex operator subalgebra containing $S$, that is, $\langle S\cup\{\omega\}\rangle$.
\end{definition}

\begin{remark}\label{rem:generation-convention}
Throughout this paper, every hypothesis on generation refers to generation of the underlying vertex algebra.
In particular, generation by $V_1$ in Theorem~\ref{thm:one} means $\langle V_1\rangle=V$;
the vector $\omega$ is not added to the generating set.
For finite generation, the two conventions lead to equivalent conditions,
since adjoining $\omega$ to a finite generating set still gives a finite set.
\end{remark}

We now recall the definition of strong generation. For the related spanning results, see
\cite[(3.2), Proposition~3.4]{KL}.

\begin{definition}\label{def:strong}
Let $V$ be a vertex operator algebra, and let $S\subseteq V$. We say that $S$ \emph{strongly generates} $V$
if $V$ is spanned by $\one$ and the vectors
\[
 a^1_{-n_1}\cdots a^\ell_{-n_\ell}\one,
 \qquad \ell\in\Z_{\ge1},\quad a^j\in S,\quad n_j\in\Z_{\ge1},\quad 1\le j\le\ell.
\]
We say that $V$ is \emph{finitely generated} (respectively, \emph{finitely strongly generated})
if one can choose a finite generating set (respectively, a finite strong generating set).
Strong generation implies generation.
\end{definition}

\begin{remark}\label{rem:strong-generation}
\begin{enumerate}[label=(\arabic*)]
\item If $S\subseteq V_2$, then by \eqref{eq:weights}, each nonzero monomial in Definition~\ref{def:strong}
has conformal weight $\sum_{j=1}^{\ell}(n_j+1)$.
Hence, if $V\ne0$ and $V_2$ strongly generates $V$, then $V$ is of OZ type, with
$V_1=0$ and $V_3=L(-1)V_2$. If $S\subseteq V_1$, the corresponding nonzero monomials
have conformal weight $\sum_{j=1}^{\ell}n_j$.

\item Suppose that $V$ is of CFT type and is generated by $V_1$ as a vertex algebra. Then $V_1$ strongly generates $V$.
Indeed, for $a,b\in V_1$, it follows from \eqref{eq:weights} that $a_0b\in V_1$, $a_1b\in\C\one$, and $a_jb=0$ for $j\in\Z_{\ge2}$.
Hence \eqref{eq:mode-commutator} reduces to
\[
 [a_m,b_n]=(a_0b)_{m+n}+m(a_1b)_{m+n-1},\qquad m,n\in\Z,
\]
where the second term is a scalar multiple of the identity operator: if $a_1b=\gamma\one$ with $\gamma\in\C$, then
\[
 m(a_1b)_{m+n-1}=m\gamma\delta_{m+n,0}\id_V.
\]
By \cite[Proposition~3.9.3]{LL}, $V$ is spanned by $\one$ and the vectors
\[
 a^1_{m_1}\cdots a^\ell_{m_\ell}\one,
 \qquad \ell\in\Z_{\ge1},\quad a^j\in V_1,\quad m_j\in\Z,\quad 1\le j\le\ell.
\]
We argue by induction on the number $\ell$ of modes in such a product; for $\ell=0$ the vector is $\one$.
If all $m_j<0$, the vector already has the form required by strong generation.
Otherwise, using the above commutation relations, we move a mode with nonnegative index step by step past the mode immediately to its right.
Each commutator term produced in this way replaces the two modes by a single mode of some vector in $V_1$
or by a scalar multiple of the identity operator; it therefore contains strictly fewer than $\ell$ modes, and the induction hypothesis applies.
In the remaining term, the nonnegative mode eventually acts directly on $\one$, and this term vanishes because
$a_m\one=0$ for $a\in V_1$ and $m\in\Z_{\ge0}$.
Consequently, each of the above vectors is a finite linear combination of $\one$ and negative-mode monomials.

\item For $a,b\in V$, the \emph{normally ordered product} is defined as in \cite[Section~3.8]{LL} by
\[
 :Y(a,z)Y(b,z):=
 \left(\sum_{m<0}a_mz^{-m-1}\right)Y(b,z)
 +Y(b,z)\left(\sum_{m\ge0}a_mz^{-m-1}\right),
\]
where both sums run over integers $m$.
From \eqref{eq:voa-derivative} and the iterate formula, one obtains the identity \cite[(3.8.8)]{LL}
\begin{equation*}
 Y(a_{-n}b,z)=\frac1{(n-1)!}
 :\bigl(\partial_z^{\,n-1}Y(a,z)\bigr)Y(b,z):,
 \qquad n\in\Z_{\ge1}.
\end{equation*}
For any untwisted weak $V$-module $M$, the corresponding identity for modules reads
\begin{equation}\label{eq:negative-mode-module}
 Y_M(a_{-n}b,z)=\frac1{(n-1)!}
 :\bigl(\partial_z^{\,n-1}Y_M(a,z)\bigr)Y_M(b,z):,
 \qquad a,b\in V,\quad n\in\Z_{\ge1},
\end{equation}
where the normally ordered product is defined in the same way.
In the proof of Proposition~\ref{prop:injection}, we shall use \eqref{eq:negative-mode-module} and the
strong generation by $V_1$ to verify that the grading of the module constructed there satisfies the conditions of an admissible grading.

\item By Definition~\ref{def:vertex-operator-algebra}, both $V_1$ and $V_2$ are finite-dimensional.
If $V_2$ strongly generates $V$, then any basis of $V_2$ is a finite strong generating set of $V$,
so $V$ is finitely strongly generated.
If $V$ is of CFT type and is generated by $V_1$ as a vertex algebra, then by (2),
$V_1$ strongly generates $V$, so $V$ is again finitely strongly generated.
Thus the hypotheses of Theorems~\ref{thm:two} and~\ref{thm:one} already imply finite strong generation,
and no separate finiteness assumption on generation is needed.
\end{enumerate}
\end{remark}

\Needspace{3\baselineskip}
\subsection{Weak, admissible, and ordinary twisted modules}\label{subsec:twisted-modules}
In this subsection we use the terminology of \cite{DLMtw} and the notation of \cite{SY}.
Throughout this subsection, we fix a vertex operator algebra $(V,Y,\one,\omega)$ and an automorphism
$g\in\Aut(V)$ of finite order $T\in\Z_{\ge1}$.
For $r\in\{0,\ldots,T-1\}$, set
\[
 V_g^{[\frac{r}{T}]}=\{a\in V:ga=e^{\frac{2\pi\mathrm i r}{T}}a\}.
\]
The subscript $g$ indicates the automorphism under consideration, and in the superscript $[\frac{r}{T}]$, the number $\frac{r}{T}$ is the representative
in $[0,1)$ of its class modulo $\Z$. Since $g^T=\id_V$, there is an eigenspace decomposition
\[
 V=\bigoplus_{r=0}^{T-1}V_g^{[\frac{r}{T}]}.
\]
More generally, for any vector space $A$ with a $g$-action such that $g^T=\id_A$, we write
\[
 A_g^{[\frac{r}{T}]}=\{a\in A:ga=e^{\frac{2\pi\mathrm i r}{T}}a\},
 \qquad r\in\{0,\ldots,T-1\};
\]
thus $A_g^{[0]}$ is the \emph{fixed-point subspace} of $g$.
If $A\subseteq V$ is a $g$-stable subspace, then $A_g^{[\frac{r}{T}]}=A\cap V_g^{[\frac{r}{T}]}$.

We begin with the definition of a weak twisted module, following \cite[Definition~3.1]{DLMtw}.

\begin{definition}\label{def:weak}
A \emph{weak $g$-twisted $V$-module} is a vector space $M$ equipped with a linear map
\[
 \begin{aligned}
 Y_M(\cdot,z):V&\longrightarrow(\End M)[[z^{\frac{1}{T}},z^{-\frac{1}{T}}]],\\
 a&\longmapsto Y_M(a,z)=\sum_{n\in\frac1T\Z}a_n^Mz^{-n-1},
 \end{aligned}
\]
with $a_n^M\in\End M$ for $a\in V$ and $n\in\frac1T\Z$,
such that, for any $r\in\{0,\ldots,T-1\}$, $a\in V_g^{[\frac{r}{T}]}$, $b\in V$, and $w\in M$,
\[
 \begin{aligned}
 Y_M(a,z)&=\sum_{n\in \frac{r}{T}+\Z}a_n^Mz^{-n-1},\\
 a_n^Mw&=0\qquad\text{for $n\in \frac{r}{T}+\Z$ sufficiently large},\\
 Y_M(\one,z)&=\id_M,
 \end{aligned}
\]
and the twisted Jacobi identity holds:
\begin{align*}
 &z_0^{-1}\delta\!\left(\frac{z_1-z_2}{z_0}\right)Y_M(a,z_1)Y_M(b,z_2)
 -z_0^{-1}\delta\!\left(\frac{z_2-z_1}{-z_0}\right)Y_M(b,z_2)Y_M(a,z_1)\notag\\
 &\hspace{1em}=
 z_2^{-1}\left(\frac{z_1-z_0}{z_2}\right)^{-\frac{r}{T}}
 \delta\!\left(\frac{z_1-z_0}{z_2}\right)Y_M(Y(a,z_0)b,z_2).
\end{align*}
When $g=\id_V$, we have $T=1$, and the definition reduces to the usual definition of a weak $V$-module.
\end{definition}

Since $g\omega=\omega$, we can write
\[
 Y_M(\omega,z)=\sum_{n\in\Z}L_M(n)z^{-n-2},\qquad L_M(n)\in\End M,\quad n\in\Z.
\]
By \cite[Section~3]{DLMtw}, these operators satisfy the Virasoro relations,
and $Y_M$ satisfies the $L(-1)$-derivative property. When no confusion can arise, we omit the superscript $M$ on modes.
A weak twisted module is not required to carry a prescribed grading, and $L_M(0)$ need not be semisimple or locally finite.
Admissible and ordinary modules are defined as follows, following \cite[Definition~3.3]{DLMtw}
and \cite[Definition~3.2]{DLMtw}, respectively.

\begin{definition}\label{def:admissible}
Let $M$ be a weak $g$-twisted $V$-module as in Definition~\ref{def:weak}.
We say that $M$ is \emph{admissible} if it admits a compatible grading
\[
 M=\bigoplus_{d\in\frac1T\N}M(d),\qquad
 a_nM(d)\subseteq M(d+\wt a-n-1),
\]
where the second condition is required for every nonzero homogeneous $a\in V$, $n\in\frac1T\Z$, and
$d\in\frac1T\N$, with the convention that components of negative degree are zero.
Here $M(d)$ is not required to be finite-dimensional, nor to be an eigenspace of $L_M(0)$.
When $g=\id_V$, such an $M$ is called an \emph{admissible $V$-module} \cite{DLMreg,Li,ABD}.

We say that $M$ is \emph{ordinary} if it has an eigenspace decomposition with respect to $L_M(0)$
\[
 M=\bigoplus_{\lambda\in\C}M_\lambda,\qquad
 L_M(0)|_{M_\lambda}=\lambda\,\id_{M_\lambda},\qquad \dim M_\lambda<\infty,
\]
and, for each fixed $\lambda\in\C$, one has $M_{\lambda+s}=0$ for $s\in\frac1T\Z$ sufficiently small.
When $g=\id_V$, such an $M$ is called an \emph{ordinary $V$-module} \cite{DLMreg,Li,ABD,DY}.
\end{definition}

Complete reducibility for these classes of modules leads to the following definitions
\cite{DLMreg,DLMtw,Li,ABD,Yamauchi,DY,YS,SY}.

\begin{definition}\label{def:rational}
With $V$ and $g$ fixed as in this subsection, we say that $V$ is \emph{$g$-rational} if every admissible $g$-twisted $V$-module
is a direct sum of irreducible admissible modules.
We say that $V$ is \emph{$g$-regular} if every weak $g$-twisted module is a direct sum of irreducible
ordinary modules. When $g=\id_V$, we say that $V$ is \emph{rational} and \emph{regular}, respectively.
\end{definition}

\begin{remark}\label{rem:module-relations}\label{prop:twisted-finiteness}
\begin{enumerate}[label=(\arabic*)]
\item By \cite[Lemma~3.4]{DLMtw}, every ordinary $g$-twisted module admits a
compatible admissible grading, and every admissible $g$-twisted module is a weak
$g$-twisted module.

\item By the definitions and (1), $g$-regularity implies $g$-rationality.

\item By \cite[Theorem~8.1~(b), (c)]{DLMtw}, if $V$ is $g$-rational,
then there are only finitely many isomorphism classes of irreducible admissible $g$-twisted $V$-modules,
and all of them are ordinary.
\end{enumerate}
\end{remark}

In the untwisted case, applying rationality to the adjoint module gives the following simplicity result.
It has already been used in deriving Corollary~\ref{cor:untwot}, and it will also be used in the proof of Theorem~\ref{thm:one}.

\begin{lemma}\label{lem:simple}
A rational vertex operator algebra of CFT type is simple.
\end{lemma}
\begin{proof}
By the CFT type assumption, the conformal weight grading of $V$ makes the adjoint module $V$ admissible.
By rationality, we can write
$V=\bigoplus_{\alpha\in\Lambda}I_\alpha$,
where $\Lambda$ is an index set and each $I_\alpha$ is an irreducible admissible $V$-submodule.
For $v\in V_0=\ker L(0)$, write $v$ as a finite sum $v=\sum_{\alpha\in\Lambda} v_\alpha$
with $v_\alpha\in I_\alpha$. Since $L(0)I_\alpha\subseteq I_\alpha$ and the sum is direct,
$0=L(0)v=\sum_{\alpha\in\Lambda} L(0)v_\alpha$ implies $L(0)v_\alpha=0$ for every $\alpha$,
so $v_\alpha\in I_\alpha\cap V_0$. Hence, by the CFT type assumption,
\[
V_0=\bigoplus_{\alpha\in\Lambda} (I_\alpha\cap V_0)=\C\one.
\]
Therefore there is a unique $\alpha_0\in\Lambda$ such that $\one\in I_{\alpha_0}$.
For every $v\in V$, we have $v=v_{-1}\one\in I_{\alpha_0}$,
and thus $V=I_{\alpha_0}$. Hence $V$ is simple.
\end{proof}

The above proof uses the untwisted adjoint module and does not apply directly to $g$-rationality for nontrivial $g$.
For this reason, the simplicity assumption is retained in Theorem~\ref{thm:two}.

\Needspace{3\baselineskip}
\subsection{\texorpdfstring{$C_2$}{C2}-cofiniteness, regularity, and tensor products}\label{subsec:c2}
This subsection recalls the definition of $C_2$-cofiniteness, its relation to regularity, and some properties of finite tensor products.
We begin with the notion of $C_2$-cofiniteness introduced by Zhu \cite{Zhu}.

\begin{definition}\label{def:C2}
Let $V$ be a vertex operator algebra. Set
\[
 C_2(V)=\Span\{a_{-2}b:a,b\in V\},\qquad R_V=V/C_2(V).
\]
We say that $V$ is \emph{$C_2$-cofinite} if $\dim_{\C}R_V<\infty$.
\end{definition}

For $a\in V$, write $\bar a=a+C_2(V)$. By \cite[Section~4.4]{Zhu},
$R_V$ becomes a Poisson algebra with respect to
\[
 \bar a\,\bar b=\overline{a_{-1}b},\qquad
 \{\bar a,\bar b\}=\overline{a_0b},\qquad a,b\in V,
\]
with identity element $\overline{\one}$.
Here a \emph{Poisson bracket} means a Lie bracket satisfying the Leibniz rule.

The following regularity criterion \cite[Theorem~4.5]{ABD} will be used in the proofs of both main theorems.

\begin{proposition}\label{prop:knownC2}
Let $V$ be of CFT type. Then $V$ is regular if and only if $V$ is rational and $C_2$-cofinite.
\end{proposition}

Our account of tensor products follows \cite[Section~2.5]{FHL} and \cite[Section~3.12]{LL}.

\begin{definition}\label{def:tensor}
Let $(U,Y_U,\one_U,\omega_U)$ and $(W,Y_W,\one_W,\omega_W)$ be vertex operator algebras.
Their \emph{tensor product} is the vertex operator algebra on the vector space $U\otimes W$
whose vertex operator is determined by
\begin{equation}\label{eq:tensor-vertex}
 Y_{U\otimes W}(u\otimes w,z)=Y_U(u,z)\otimes Y_W(w,z),
 \qquad u\in U,\quad w\in W,
\end{equation}
and linearity, and whose vacuum vector and conformal vector are
\begin{equation}\label{eq:tensor-vectors}
 \one_{U\otimes W}=\one_U\otimes\one_W,\qquad
 \omega_{U\otimes W}=\omega_U\otimes\one_W+\one_U\otimes\omega_W.
\end{equation}
The corresponding conformal weight grading is
\begin{equation}\label{eq:tensor-grading}
 (U\otimes W)_n=\bigoplus_{\substack{p,q\in\Z\\p+q=n}}U_p\otimes W_q,
 \qquad n\in\Z.
\end{equation}
Tensor products of finitely many vertex operator algebras are defined by iteration; by convention, the empty tensor product is $\C\one$.
\end{definition}

We next record the rationality and $C_2$-cofiniteness of finite tensor products.
A finite tensor product is rational if and only if each tensor factor is rational; see \cite[Proposition~2.7]{DMZ}
and also \cite[Section~2]{DLMreg}. For the implication from the rationality of the factors to that of the tensor product,
see also \cite[Theorem~5.2]{Milas} (with respect to the definition of rationality used there) and \cite[Proposition~6.2]{BDM}, the latter with all automorphisms taken to be the identity.
Similarly, a finite tensor product is $C_2$-cofinite if and only if each tensor factor is $C_2$-cofinite;
see \cite[Proposition~5.1~(1)]{AbeT}. The implication from the factors to the tensor product
is also the special case of \cite[Lemma~2.8]{AvE} for identity automorphisms.
Combining these facts with the above construction, the definition of the $C_2$ quotient, and its Poisson algebra structure, we obtain the following result.

\begin{proposition}\label{prop:tensor}
\begin{enumerate}[label=(\arabic*)]
\item A finite tensor product of vertex operator algebras is rational if and only if each tensor factor is rational.
\item For vertex operator algebras $U$ and $W$, we have
\[
 C_2(U\otimes W)=C_2(U)\otimes W+U\otimes C_2(W),
\]
and the natural map
\[
 \begin{aligned}
 R_{U\otimes W}&\longrightarrow R_U\otimes R_W,\\
 \overline{u\otimes w}&\longmapsto\bar u\otimes\bar w,
 \qquad u\in U,\quad w\in W,
 \end{aligned}
\]
is an isomorphism of Poisson algebras, where each bar denotes the image in the corresponding $C_2$ quotient.
In particular, a finite tensor product of vertex operator algebras is $C_2$-cofinite if and only if each tensor factor is $C_2$-cofinite.
\end{enumerate}
\end{proposition}

\Needspace{3\baselineskip}
\subsection{Invariant bilinear forms}\label{subsec:invariant-forms}
Following \cite[Section~3.2]{SY}, we review the basic properties of invariant bilinear forms and record the formulas needed later
for computing values of bilinear forms, deciding nondegeneracy, and studying quotient algebras.
Throughout this subsection, $(V,Y,\one,\omega)$ is a vertex operator algebra with
conformal weight grading $V=\bigoplus_{n\in\Z}V_n$,
and $B_V$ denotes a bilinear form on $V$.
The notion of invariance is recalled from \cite[Definition~2.4]{LiForm}; see also
\cite[Remark~5.3.3]{FHL}.

\begin{definition}\label{def:invariant}
A \emph{bilinear form} on $V$ is a map $B_V:V\times V\to\C$ that is $\C$-linear in each variable.
It is said to be \emph{invariant} if, for any $a,u,v\in V$,
\begin{equation}\label{eq:invariance}
 B_V(Y(a,z)u,v)=
 B_V\bigl(u,Y(e^{zL(1)}(-z^{-2})^{L(0)}a,z^{-1})v\bigr),
\end{equation}
where the equality is understood coefficientwise in the powers of $z$.
If $B_V(\one,\one)=1$, the form $B_V$ is said to be \emph{normalized}.
\end{definition}

\begin{remark}\label{rem:adjoint}
For $d\in\Z$ and $a\in V_d$, we have $(-z^{-2})^{L(0)}a=(-1)^dz^{-2d}a$.
Since $L(1)$ lowers conformal weights and $V$ is lower bounded, $e^{zL(1)}a$ is a finite sum.
A vector $a\in V_d$ is said to be \emph{quasi-primary} if $L(1)a=0$.
If $a$ is quasi-primary and $B_V$ is invariant, then comparing the coefficients of $z^{-n-1}$ on both sides of \eqref{eq:invariance}
gives the following adjoint identity \cite[Remark~5.3.5]{FHL}:
\begin{equation}\label{eq:adjoint}
 B_V(a_nu,v)=(-1)^dB_V(u,a_{2d-n-2}v),
 \qquad n\in\Z,\quad u,v\in V.
\end{equation}
\end{remark}

For a symmetric bilinear form $B_V$, its \emph{radical} is defined by
\[
 \rad B_V=\{v\in V:B_V(v,w)=0\text{ for all }w\in V\}.
\]
If $\rad B_V=0$, we say that $B_V$ is \emph{nondegenerate}.
The following properties can be found in \cite[Remark~2.5, Proposition~2.6, (2.23), Theorem~3.1, Corollary~3.2]{LiForm}
and \cite[Proposition~4.6]{CKLW}; for the uniqueness in the CFT type case and the passage from nondegeneracy
to simplicity, see also \cite[Lemma~3.6]{SY}.

\begin{proposition}\label{prop:form}
\begin{enumerate}[label=(\arabic*)]
\item Every invariant bilinear form $B_V$ is symmetric.
Conformal weight spaces of different conformal weights are orthogonal with respect to $B_V$, that is,
\[
 B_V(V_m,V_n)=0,\qquad m,n\in\Z,\quad m\ne n.
\]
The radical of $B_V$ is a vertex ideal.
\item The restriction map $B_V\mapsto B_V(\one,-)|_{V_0}$ induces a natural isomorphism between the space of invariant bilinear forms
on $V$ and $(V_0/L(1)V_1)^*$.
\item If $V$ is of CFT type, then $V$ admits a normalized invariant bilinear form
if and only if $L(1)V_1=0$. In this case the form is unique, and $V$ is simple if and only if this form is nondegenerate.
\end{enumerate}
\end{proposition}

\begin{remark}\label{rem:weight-gram}
Let $B_V$ be an invariant bilinear form, and write
\[
 (B_V)_n=B_V|_{V_n\times V_n},\qquad n\in\Z.
\]
By Proposition~\ref{prop:form}~(1),
\[
 \rad B_V=\bigoplus_{n\in\Z}\rad (B_V)_n.
\]
For each $n\in\Z$ with $V_n\ne0$, let $d_n=\dim V_n$,
and choose a basis $v_{n,1},\ldots,v_{n,d_n}$ of $V_n$.
The corresponding \emph{Gram matrix} is
\[
 \bigl((B_V)_n(v_{n,i},v_{n,j})\bigr)_{1\le i,j\le d_n}.
\]
Consequently, $B_V$ is nondegenerate if and only if
\[
 \det\!\left[\bigl((B_V)_n(v_{n,i},v_{n,j})\bigr)_{1\le i,j\le d_n}\right]\ne0
\]
for every $n\in\Z$ with $V_n\ne0$.
\end{remark}

\begin{remark}\label{rem:form-ideal}
Let $B_V$ be a normalized invariant bilinear form.
By Proposition~\ref{prop:form}~(1), $\rad B_V$ is a vertex ideal,
and it is proper since $B_V(\one,\one)=1$.
The bilinear form $\overline B_V$ on the quotient algebra $V/\rad B_V$ induced by $B_V$,
\[
 \overline B_V(u+\rad B_V,v+\rad B_V)=B_V(u,v),\qquad u,v\in V,
\]
is well defined, normalized, invariant, and nondegenerate.
\end{remark}

We now record the results on bilinear forms needed for finite tensor products; the underlying linear algebra is contained in
\cite[Lemmas~1.4, 2.1~(2)]{EKM}.

\begin{proposition}\label{prop:tensorforms}
Let $s\in\Z_{\ge1}$. For $1\le i\le s$, let
$U_i=\bigoplus_{n\in\Z}U_{i,n}$ be a graded vector space with each $U_{i,n}$ finite-dimensional,
and let $B_i$ be a symmetric bilinear form on $U_i$ such that
\[
 B_i(U_{i,m},U_{i,n})=0,\qquad m,n\in\Z,\quad m\ne n.
\]
Let $U=\bigotimes_{i=1}^s U_i$, equipped with the grading
\[
 U=\bigoplus_{n\in\Z}U_n,\qquad
 U_n=\bigoplus_{\substack{(n_1,\ldots,n_s)\in\Z^s\\n_1+\cdots+n_s=n}}
 U_{1,n_1}\otimes\cdots\otimes U_{s,n_s},\qquad n\in\Z.
\]
The \textit{tensor product} of the forms $B_i$ is the
symmetric bilinear form $B_U=\bigotimes_{i=1}^s B_i$ on $U$ determined by
\begin{equation}\label{eq:tensor-form}
 B_U(u_1\otimes\cdots\otimes u_s,v_1\otimes\cdots\otimes v_s)
 =\prod_{i=1}^s B_i(u_i,v_i)
\end{equation}
and bilinearity, where $u_i,v_i\in U_i$ for $1\le i\le s$. Then
\begin{equation}\label{eq:tensor-radical}
 \rad B_U=\sum_{i=1}^s
 U_1\otimes\cdots\otimes\rad B_i\otimes\cdots\otimes U_s.
\end{equation}
Each $\rad B_i$ is a graded subspace, so $U_i/\rad B_i$ inherits the grading of $U_i$.
The natural map
\[
 \begin{aligned}
 \mathsf q_U:U&\longrightarrow\bigotimes_{i=1}^s(U_i/\rad B_i),\\
 w_1\otimes\cdots\otimes w_s
 &\longmapsto(w_1+\rad B_1)\otimes\cdots\otimes(w_s+\rad B_s),
 \end{aligned}
\]
where $w_i\in U_i$ for $1\le i\le s$,
is a surjection preserving the total grading, and $\ker\mathsf q_U=\rad B_U$.
Hence $\mathsf q_U$ induces a grading-preserving linear isomorphism
\begin{equation}\label{eq:tensor-quotient}
 \overline{\mathsf q}_U:U/\rad B_U\xrightarrow{\ \sim\ }
 \bigotimes_{i=1}^s(U_i/\rad B_i),
\end{equation}
and the induced bilinear form on the quotient space $U/\rad B_U$ is nondegenerate.

Suppose, in addition, that for each $1\le i\le s$, $U_i$ is a vertex operator algebra,
$U_{i,n}$ is its subspace of conformal weight $n$ for $n\in\Z$,
and $B_i$ is normalized and invariant. Then, with respect to the tensor product vertex operator algebra structure on $U$ from Definition~\ref{def:tensor},
$B_U$ is a normalized invariant bilinear form on $U$. In this case, the subspaces $\rad B_i$ and $\rad B_U$ are all proper vertex ideals,
and if $U_i/\rad B_i$ and $U/\rad B_U$ are given the corresponding quotient vertex operator algebra structures,
then the map $\overline{\mathsf q}_U$ in \eqref{eq:tensor-quotient}
is a vertex operator algebra isomorphism.
\end{proposition}
\begin{proof}
For $1\le i\le s$, the orthogonality of distinct homogeneous components gives
\[
 \rad B_i=\bigoplus_{n\in\Z}\rad\bigl(B_i|_{U_{i,n}\times U_{i,n}}\bigr).
\]
Thus $\rad B_i$ is a graded subspace. For each $n\in\Z$, choose a complement $E_{i,n}$ of
$\rad(B_i|_{U_{i,n}\times U_{i,n}})$ in $U_{i,n}$. By \cite[Lemma~1.4]{EKM},
the restriction of $B_i$ to $E_{i,n}$ is nondegenerate. Setting
$E_i=\bigoplus_{n\in\Z}E_{i,n}$, we have $U_i=\rad B_i\oplus E_i$.

For any $(n_1,\ldots,n_s)\in\Z^s$, \cite[Lemma~2.1~(2)]{EKM} shows that the restriction of
$B_U$ to the finite-dimensional space $E_{1,n_1}\otimes\cdots\otimes E_{s,n_s}$ is
nondegenerate. The subspaces corresponding to different tuples $(n_1,\ldots,n_s)$ are mutually orthogonal,
and $E=\bigotimes_{i=1}^s E_i$ is their direct sum; hence the restriction of $B_U$ to $E$ is
nondegenerate. Expanding $U=\bigotimes_{i=1}^s(\rad B_i\oplus E_i)$, we obtain the direct sum decomposition
\[
 U=\left(\sum_{i=1}^s
 U_1\otimes\cdots\otimes\rad B_i\otimes\cdots\otimes U_s\right)\oplus E.
\]
By \eqref{eq:tensor-form},
\[
 \sum_{i=1}^s U_1\otimes\cdots\otimes\rad B_i\otimes\cdots\otimes U_s
 \subseteq\rad B_U.
\]
Conversely, let $v\in\rad B_U$, and write $v=u+e$ according to the above direct sum decomposition,
where $u$ lies in the subspace given by the sum and $e\in E$. For every $f\in E$,
\[
 B_U(e,f)=B_U(v,f)-B_U(u,f)=0.
\]
Since $B_U$ is nondegenerate on $E$, we get $e=0$, and hence
\[
 v=u\in\sum_{i=1}^s
 U_1\otimes\cdots\otimes\rad B_i\otimes\cdots\otimes U_s.
\]
This proves \eqref{eq:tensor-radical}.

For $1\le i\le s$, denote by $\mathsf q_i:U_i\longrightarrow U_i/\rad B_i$ the natural quotient map.
Since $U_i=\rad B_i\oplus E_i$, its restriction
\[
 \mathsf q_i|_{E_i}:E_i\longrightarrow U_i/\rad B_i
\]
is a linear isomorphism. Moreover,
\[
 \mathsf q_U|_E=\bigotimes_{i=1}^s\bigl(\mathsf q_i|_{E_i}\bigr),
\]
so $\mathsf q_U|_E$ is also a linear isomorphism. Since
\[
 \mathsf q_U\bigl(U_1\otimes\cdots\otimes\rad B_i\otimes\cdots\otimes U_s\bigr)=0,
 \qquad 1\le i\le s,
\]
the above direct sum decomposition and \eqref{eq:tensor-radical} give $\ker\mathsf q_U=\rad B_U$. As each $\rad B_i$ is a graded subspace,
$\mathsf q_U$ preserves the total grading and induces the grading-preserving linear isomorphism in \eqref{eq:tensor-quotient}.
The nondegeneracy of the induced bilinear form on $U/\rad B_U$ follows directly from the definition of the radical.

Finally, suppose that the $U_i$ and $B_i$ satisfy the hypotheses of the proposition concerning vertex operator algebras.
For $1\le i\le s$, let $Y_i$ and $L_i(n)$, $n\in\Z$, denote the vertex operator and the
Virasoro operators of $U_i$, and let $Y_U$ and $L_U(n)$ denote the corresponding operators on $U$.
Consider pure tensors
\[
 a=\bigotimes_{i=1}^s a_i,\qquad
 u=\bigotimes_{i=1}^s u_i,\qquad
 v=\bigotimes_{i=1}^s v_i,\qquad
 a_i,u_i,v_i\in U_i,\quad 1\le i\le s.
\]
By \eqref{eq:tensor-vectors}, $L_U(0)$ and $L_U(1)$ are the sums of the corresponding operators on the factors.
Since operators acting on different factors commute, the identity $e^{z(A+B)}=e^{zA}e^{zB}$ for commuting operators applies to the summands of $L_U(1)$;
together with the additivity of conformal weights, this gives
\[
 e^{zL_U(1)}(-z^{-2})^{L_U(0)}a
 =\bigotimes_{i=1}^s\left(e^{zL_i(1)}(-z^{-2})^{L_i(0)}a_i\right).
\]
Combining \eqref{eq:tensor-vertex}, \eqref{eq:tensor-form}, and the invariance of each $B_i$, we obtain
\[
 \begin{aligned}
 B_U(Y_U(a,z)u,v)
 &=\prod_{i=1}^s B_i(Y_i(a_i,z)u_i,v_i)\\
 &=\prod_{i=1}^s B_i\bigl(u_i,
 Y_i(e^{zL_i(1)}(-z^{-2})^{L_i(0)}a_i,z^{-1})v_i\bigr)\\
 &=B_U\bigl(u,Y_U(e^{zL_U(1)}(-z^{-2})^{L_U(0)}a,z^{-1})v\bigr).
 \end{aligned}
\]
Since pure tensors span $U$, it follows by linearity that $B_U$ is invariant.
Furthermore, by \eqref{eq:tensor-vectors} and \eqref{eq:tensor-form},
\[
 B_U(\one_U,\one_U)
 =\prod_{i=1}^s B_i(\one_{U_i},\one_{U_i})=1,
\]
so $B_U$ is normalized.

By Remark~\ref{rem:form-ideal}, the subspaces $\rad B_i$ and $\rad B_U$ are all proper vertex ideals,
so the corresponding quotient vertex operator algebras are defined.
The quotient maps of the factors preserve vertex operators, vacuum vectors, and conformal vectors;
hence, by \eqref{eq:tensor-vertex} and \eqref{eq:tensor-vectors}, $\mathsf q_U$ is a
vertex operator algebra homomorphism. Consequently, the induced map $\overline{\mathsf q}_U$ is a
vertex operator algebra isomorphism.
\end{proof}

\Needspace{3\baselineskip}
\subsection{Zhu algebra and top-level representations}\label{subsec:zhu}
Throughout this subsection, $V$ is a vertex operator algebra. Following \cite{Zhu,DLMtw},
we recall the Zhu algebra and its action on the degree-zero spaces of admissible modules;
these results will be used in Proposition~\ref{prop:injection}.
We start with the construction in \cite[Section~2]{Zhu}; see also the case $g=\id_V$ of
\cite[Proposition~2.3, Theorem~2.4]{DLMtw}.

\begin{definition}\label{def:zhu}
For $d\in\Z$, $u\in V_d$, and $v\in V$, set
\[
 u*v=\Res_z\frac{(1+z)^d}{z}Y(u,z)v,
 \qquad
 u\circ v=\Res_z\frac{(1+z)^d}{z^2}Y(u,z)v,
\]
and extend these products bilinearly to $V\times V$. Set
\[
 O(V)=\Span\{u\circ v:u,v\in V\},\qquad A(V)=V/O(V),
\]
and write $[u]=u+O(V)$ for $u\in V$. The multiplication
\[
 [u]*[v]=[u*v],\qquad u,v\in V,
\]
is well defined and makes $A(V)$ an associative algebra with identity element $[\one]$;
this algebra $A(V)$ is called the \emph{Zhu algebra} of $V$.
\end{definition}

\begin{definition}\label{def:toplevel}
Let $M=\bigoplus_{d\in\N}M(d)$ be an admissible $V$-module.
If $M(0)\ne0$, then $M(0)$ is called the \emph{top level} of $M$ with respect to this admissible grading.
\end{definition}

For an admissible $V$-module $M$, $d\in\Z$, and $u\in V_d$,
the \emph{zero mode} of $u$ on $M$ is defined by
\[
 o(u)=u^M_{d-1}\in\End M,
\]
and $o(u)$ is defined for general $u\in V$ by linearity.
By the admissible grading, $o(u)$ preserves $M(0)$; in particular,
\[
 o(x)=x_0^M,\qquad x\in V_1.
\]
We record two properties given by \cite[Theorems~5.3, 8.1~(a)]{DLMtw} for $g=\id_V$.

\begin{proposition}\label{prop:zhu}
Let $M=\bigoplus_{d\in\N}M(d)$ be an admissible $V$-module.
Then $M(0)$ is an $A(V)$-module with action
\[
 [u]\cdot w=o(u)w,\qquad u\in V,\quad w\in M(0).
\]
If $V$ is rational, then $\dim A(V)<\infty$.
\end{proposition}

\begin{remark}\label{rem:top-level-action}
Let $V$ be of CFT type. By \cite[Theorem~3.7]{Lian}, $V_1$ is a Lie algebra under the bracket
$[x,y]=x_0y$ for $x,y\in V_1$.
For $x,y\in V_1$, the case $g=\id_V$ of \cite[Lemma~2.2~(iii)]{DLMtw} gives
\begin{equation*}
 [x]*[y]-[y]*[x]=[x_0y].
\end{equation*}
Hence the map $x\mapsto[x]$ is a Lie algebra homomorphism with respect to the commutator bracket on $A(V)$ determined by the multiplication $*$,
and it therefore extends uniquely to a unital associative algebra homomorphism
\[
 U(V_1)\longrightarrow A(V),
\]
where, for a Lie algebra $\mathfrak a$, $U(\mathfrak a)$ denotes its universal enveloping algebra.
For any Lie subalgebra $\q\subseteq V_1$, we thus obtain a homomorphism $U(\q)\to A(V)$.
This homomorphism will be used in Proposition~\ref{prop:injection}, where it is shown to be injective under the hypotheses there.
\end{remark}

\section{Algebras in conformal weight two and tensor product decompositions}\label{sec:frob}
This section studies the structure of a simple vertex operator algebra $V$ strongly generated by $V_2$.
The Griess product on $V_2$ and the corresponding bilinear form determine a commutative Frobenius algebra.
Using the constructions of \cite{Lam,LamT,LYu}, we show that $V$ is isomorphic to the tensor product of the
simple vertex operator algebras associated with the local summands of this Frobenius algebra; see Proposition~\ref{prop:frobfactors}.
The tensor factors corresponding to semisimple local summands are simple Virasoro vertex operator algebras,
while those corresponding to nonsemisimple local summands will be studied in Section~\ref{sec:frobgram}.
This decomposition will be used in the proof of Theorem~\ref{thm:two} in Section~\ref{sec:twists}.

Section~\ref{subsec:frob-weight-two} explains, following \cite{LYu}, how the Griess product on $V_2$ and the corresponding bilinear form determine a commutative Frobenius algebra.
Section~\ref{subsec:local-euler} recalls from \cite{AKl,AbE} the local decomposition of a finite-dimensional commutative associative unital algebra
and the corresponding Frobenius structure; it also introduces the Euler element and records its basic properties (see Proposition~\ref{prop:local} and Lemma~\ref{lem:euler});
the latter will be used in Section~\ref{sec:frobgram}.
In Section~\ref{subsec:frob-voa}, we construct, as in \cite{LYu}, the vertex operator algebra associated with a commutative Frobenius algebra,
and obtain from the above local decomposition the tensor product decomposition of Proposition~\ref{prop:frobfactors};
finally, we explain that the semisimple local case corresponds to the Virasoro vertex operator algebra,
and recall its simplicity and the rationality, regularity, and $C_2$-cofiniteness of its simple quotient.

\Needspace{3\baselineskip}
\subsection{The commutative Frobenius algebra attached to conformal weight two}\label{subsec:frob-weight-two}
Throughout this subsection, $V$ is a nonzero vertex operator algebra.
When $V$ is of OZ type, $V_2$ is a commutative algebra under the multiplication $a\cdot b=a_1b$,
called the \emph{Griess algebra} of $V$ \cite[Section~1]{JLY};
this multiplication is called the \emph{Griess product}; see also \cite[Remark~3.2.3]{LL}.
The name comes from the commutative nonassociative algebra introduced by Griess in his construction of the Monster simple group; see \cite{G} and also \cite[Section~2]{GNAVOA}.
The realization of this algebra in the subspace of conformal weight two of the moonshine vertex operator algebra can be found in \cite[Section~10.3]{FLM}.
When $V$ is strongly generated by $V_2$, \cite[Proposition~3.7]{LYu} gives
$a_0b=\frac12L(-1)a_1b$ for $a,b\in V_2$.
Hence \cite[Proposition~3.9]{LYu} applies
and yields the associativity of the multiplication on $V_2$ together with the properties of the corresponding bilinear form.
These results lead to the following proposition.

\begin{proposition}\label{prop:wt2}
If $V$ is strongly generated by $V_2$, then $V$ is of CFT type, $V_1=0$,
and $A=V_2$ is a finite-dimensional
commutative associative algebra under the multiplication $ab=\frac12a_1b$, with identity element $1_A=\omega$.
Define a bilinear form $B_A:A\times A\to\C$ by $B_A(a,b)\one=a_3b$ for $a,b\in A$. Then
\begin{equation}\label{eq:wt2products}
 a_0b=L(-1)(ab),\quad a_1b=2ab,\quad a_2b=0,\quad a_3b=B_A(a,b)\one,\quad a_jb=0,\quad j\in\Z_{\ge4}.
\end{equation}
Clearly, this bilinear form is symmetric and invariant, that is,
\[
 B_A(ab,c)=B_A(a,bc),\qquad a,b,c\in A.
\]
If $V$ is simple, then $B_A$ is nondegenerate.
\end{proposition}
\begin{proof}
Apart from the final nondegeneracy statement, all assertions follow from Remark~\ref{rem:strong-generation}~(1)
and the above results \cite[Propositions~3.7, 3.9]{LYu}.
Suppose that $V$ is simple. By $V_1=0$ and Proposition~\ref{prop:form}~(3), $V$ has a
normalized nondegenerate invariant bilinear form $B_V$.
For $a,b\in V_2$, we have $L(1)a=0$, so \eqref{eq:adjoint} gives
\[
 B_V(a,b)=B_V(a_{-1}\one,b)
 =B_V(\one,a_3b)=B_A(a,b).
\]
Thus $B_A=(B_V)_2$, in the notation of Remark~\ref{rem:weight-gram}.
By that remark and the nondegeneracy of $B_V$, the form $B_A$ is nondegenerate.
\end{proof}

\begin{remark}\label{rem:griess-normalization}
Under the Griess product $a\cdot b=a_1b$, the identity element is $\frac{\omega}{2}$.
Proposition~\ref{prop:wt2} uses the normalized multiplication $ab=\frac12a_1b$,
for which the identity element of $A$ is $1_A=\omega$.
\end{remark}

To describe the above structure on $V_2$, we recall the definition of a commutative Frobenius algebra,
using the notation $(\mathcal A,f)$ and $\eta$ of \cite[Section~3]{AbE};
see also \cite[Section~2]{DMu}.

\begin{definition}\label{def:frobenius}
A \emph{commutative Frobenius algebra} is a pair $(\mathcal A,f)$, where $\mathcal A$ is a
finite-dimensional commutative associative unital $\C$-algebra and
$f:\mathcal A\to\C$ is a linear functional such that the bilinear form
\[
 \eta(a,b)=f(ab),\qquad a,b\in \mathcal A,
\]
is nondegenerate. We call $f$ the \emph{Frobenius functional} and $\eta$ the \emph{Frobenius form}, and denote the identity element of $\mathcal A$ by $1_{\mathcal A}$.
Clearly, this bilinear form is symmetric and invariant, that is,
\[
 \eta(ab,c)=\eta(a,bc),\qquad a,b,c\in \mathcal A.
\]
\end{definition}

\begin{remark}\label{rem:frobenius-normalization}
Let $V$ be a nonzero simple vertex operator algebra strongly generated by $V_2$.
By Proposition~\ref{prop:wt2}, $A=V_2$ is a
finite-dimensional commutative associative algebra under the multiplication $ab=\frac12a_1b$, with identity element $1_A=\omega$,
and the corresponding bilinear form $B_A$ is nondegenerate.
Note that $1_A\in V_2$ is the identity element of $A$, whereas $\one\in V_0$ is the vacuum vector of $V$.
Define
\[
 f(a)=2B_A(a,1_A),\qquad a\in A.
\]
Let $c_V\in\C$ denote the central charge of $V$. Using the invariance of $B_A$,
$b1_A=b$, and $\omega_3\omega=\frac{c_V}{2}\one$, we obtain, for $a,b\in A$,
\begin{align}
 \eta(a,b)&=f(ab)=2B_A(ab,1_A)
 =2B_A(a,b1_A)=2B_A(a,b),\label{eq:normalization}\\
 f(1_A)&=2B_A(1_A,1_A)=2B_A(\omega,\omega)
 =2\cdot\frac{c_V}{2}=c_V.\notag
\end{align}
Hence $\eta=2B_A=2(B_V)_2$ is nondegenerate, and $(A,f)$ is a commutative Frobenius algebra.
This normalization will be used below to construct the corresponding vertex operator algebra and to compute its invariant bilinear form.
\end{remark}

\Needspace{3\baselineskip}
\subsection{Local decomposition and the Euler element}\label{subsec:local-euler}
Throughout this subsection, $\mathcal A$ is a finite-dimensional commutative associative unital $\C$-algebra
with identity element $1_{\mathcal A}$. By finite dimensionality, every
descending chain of ideals of $\mathcal A$ stabilizes, so $\mathcal A$ is Artinian
\cite[Definition~16.21, Exercise~16.42]{AKl}.
We recall its local decomposition and the corresponding Frobenius structure, and record the basic properties of the Euler element.

Our terminology for local algebras and idempotents follows
\cite[Definition~3.4, (1.10)]{AKl};
for the socle and semisimplicity, we use the formulation of \cite[Section~3]{AbE} in the finite-dimensional commutative case.

\begin{definition}\label{def:local}
Let $\mathcal A$ be a finite-dimensional commutative associative unital $\C$-algebra.
We say that $\mathcal A$ is \emph{local} if it has exactly one maximal ideal $\m$.
In this case, $\mathcal A/\m$ is called its \emph{residue field}, and the \emph{socle} of $\mathcal A$ is defined by
\[
 \Soc(\mathcal A)=\Ann_{\mathcal A}(\m)
 =\{x\in\mathcal A:xu=0\text{ for all }u\in\m\}.
\]

An element $e\in\mathcal A$ with $e^2=e$ is called an \emph{idempotent}.
Two idempotents $e,e'$ are said to be \emph{orthogonal} if $ee'=0$.
A nonzero idempotent $e$ is said to be \emph{primitive} if it cannot be written as a sum of two nonzero orthogonal idempotents.

By the Wedderburn--Artin theorem
\cite[Theorem~3.5 and p.~39]{Lamt},
together with the commutativity of $\mathcal A$ and the algebraic closedness of $\C$,
$\mathcal A$ is semisimple if and only if
\[
 \mathcal A\cong\C^s,\qquad s\in\Z_{\ge0}.
\]
This characterization can also be found in \cite[proof of Theorem~3.4]{AbE}.
\end{definition}

\begin{remark}\label{rem:local-properties}
Let $\mathcal A$ be a local algebra as in Definition~\ref{def:local}, with maximal ideal $\m$.
\begin{enumerate}[label=(\arabic*)]
\item The residue field $\mathcal A/\m$ is a finite-dimensional field extension of $\C$;
since $\C$ is algebraically closed, $\mathcal A/\m\cong\C$;
see \cite[Theorem~15.4, proof of Corollary~15.5]{AKl}.
\item By (1) and the characterization of semisimplicity in Definition~\ref{def:local},
\[
 \mathcal A\not\cong\C
 \quad\Longleftrightarrow\quad \m\ne0
 \quad\Longleftrightarrow\quad \mathcal A\text{ is nonsemisimple}.
\]
\item The units of $\mathcal A$, that is, the elements having a multiplicative inverse, are precisely the elements of $\mathcal A\setminus\m$; see \cite[Lemma~3.5]{AKl}.
\end{enumerate}
\end{remark}

The following proposition decomposes $\mathcal A$ into local algebras \cite[Corollary~19.11, Proposition~11.18]{AKl}
and describes the restriction of the Frobenius structure to the summands \cite[Proposition~2 and its proof]{Saw}.

\begin{proposition}\label{prop:local}
Let $\mathcal A$ be a finite-dimensional commutative associative unital $\C$-algebra.
Then $\mathcal A$ has a local decomposition, unique up to the ordering of the summands,
\[
 \mathcal A=\bigoplus_{i=1}^s\mathcal A_i,\qquad
 \mathcal A_i=e_i\mathcal A,\qquad s\in\Z_{\ge0},
\]
where $e_1,\ldots,e_s$ are nonzero primitive idempotents satisfying
\[
 1_{\mathcal A}=\sum_{i=1}^s e_i,\qquad
 e_i^2=e_i\quad(1\le i\le s),\qquad
 e_ie_j=0\quad(1\le i,j\le s,\ i\ne j).
\]
Each $\mathcal A_i$ is a local algebra with identity element $e_i$.
If $(\mathcal A,f)$ is a commutative Frobenius algebra and $\eta(a,b)=f(ab)$ for $a,b\in\mathcal A$,
then distinct summands are orthogonal with respect to $\eta$, that is,
\[
 \eta(\mathcal A_i,\mathcal A_j)=0,
 \qquad 1\le i,j\le s,\quad i\ne j.
\]
Set
\[
 f_i=f|_{\mathcal A_i},\qquad
 \eta_i=\eta|_{\mathcal A_i\times\mathcal A_i},\qquad 1\le i\le s.
\]
Then $\eta_i$ is nondegenerate and $\eta_i(a,b)=f_i(ab)$ for $a,b\in\mathcal A_i$;
hence each $(\mathcal A_i,f_i)$ is a commutative Frobenius algebra.
If $\mathcal A$ itself is a local commutative Frobenius algebra, then $\dim\Soc(\mathcal A)=1$;
if, in addition, $\mathcal A\not\cong\C$, then $\Soc(\mathcal A)\subseteq\m$,
where $\m$ is the maximal ideal of $\mathcal A$.
\end{proposition}
\begin{proof}
Since $\mathcal A$ is Artinian, its local decomposition and the uniqueness of this decomposition can be found in
\cite[Corollary~19.11, Proposition~11.18]{AKl}, respectively;
for the formulation in terms of idempotents, see \cite[(1.10), Proposition~1.12]{AKl}.
For the orthogonality of the summands, the nondegeneracy of the restricted forms,
and the one-dimensionality of the socle of a local commutative Frobenius algebra, see \cite[Proposition~2 and its proof]{Saw}.
Finally, suppose that $\mathcal A$ is local with $\mathcal A\not\cong\C$, and that
$x\in\Soc(\mathcal A)\setminus\m$. Then $x$ is a unit,
and $x\m=0$ gives $\m=0$, contradicting $\mathcal A\not\cong\C$.
Therefore $\Soc(\mathcal A)\subseteq\m$.
\end{proof}

We now recall the Euler element of \cite[Definition~2.2, (2.10)]{DMu},
which corresponds to the characteristic element of \cite[Section~3]{AbE}.
We denote it by $\mathbf e$ to distinguish it from the conformal vector $\omega$.

\begin{definition}\label{def:euler}
Let $(\mathcal A,f)$ be a commutative Frobenius algebra, and write $\eta(a,b)=f(ab)$ for $a,b\in\mathcal A$.
Let $d=\dim\mathcal A$, and choose a basis $b_1,\ldots,b_d$ of $\mathcal A$.
Its $\eta$-dual basis $b^1,\ldots,b^d$ is uniquely determined by
\[
 \eta(b_i,b^j)=\delta_{ij},\qquad 1\le i,j\le d.
\]
Write
\[
 \eta_{ij}=\eta(b_i,b_j),\qquad
 (\eta^{ij})=(\eta_{ij})^{-1},\qquad 1\le i,j\le d.
\]
Then
\[
 b^i=\sum_{j=1}^d\eta^{ij}b_j,\qquad 1\le i\le d.
\]
The \emph{Euler element} is defined by
\[
 \mathbf e=\sum_{i=1}^d b_ib^i
 =\sum_{i,j=1}^d\eta^{ij}b_ib_j.
\]
The notation $\mathbf e$ should not be confused with the primitive idempotents $e_i$ in the local decomposition.
\end{definition}

The next lemma records the trace formula for the Euler element \cite[Section~3, Footnote~11]{MS}
and its relation to the socle in the local case \cite[Section~3, Proposition~3.3]{AbE}.
The property $\mathbf e^2=0$ will be used in Section~\ref{sec:frobgram} to show that certain tensor contractions vanish.

\begin{lemma}\label{lem:euler}
Let $(\mathcal A,f)$ be a commutative Frobenius algebra, and let
$\mathbf e$ be the Euler element of Definition~\ref{def:euler}.
For $a\in\mathcal A$, denote by
\[
\begin{aligned}
 m_a:\mathcal A&\longrightarrow\mathcal A,\\
 x&\longmapsto ax,\qquad x\in\mathcal A,
\end{aligned}
\]
the multiplication operator.
Then $\mathbf e$ is independent of the choice of basis and satisfies the trace formula
\[
 f(\mathbf e a)=\Tr(m_a),\qquad a\in\mathcal A.
\]
In particular, $f(\mathbf e)=\dim\mathcal A$.
If $\mathcal A$ is local, then $\Soc(\mathcal A)=\C\mathbf e$.
If $\mathcal A$ is local and $\mathcal A\not\cong\C$, with maximal ideal $\m$,
then Proposition~\ref{prop:local} gives $\mathbf e\in\Soc(\mathcal A)\subseteq\m$, and moreover $\mathbf e^2=0$.
\end{lemma}
\begin{proof}
The independence of the Euler element from the choice of basis follows from the canonical construction in \cite[Section~3]{AbE},
and the trace formula is given in \cite[Section~3, Footnote~11]{MS},
where $\theta$ and $\chi$ correspond to our $f$ and $\mathbf e$, respectively.
If $\mathcal A$ is local, then \cite[Proposition~3.3]{AbE} gives
$\mathcal A\mathbf e=\Soc(\mathcal A)$.
Since $f(\mathbf e)=\dim\mathcal A\ne0$, we have $\mathbf e\ne0$;
combined with Proposition~\ref{prop:local}, this gives $\Soc(\mathcal A)=\C\mathbf e$.
If, in addition, $\mathcal A\not\cong\C$, then Proposition~\ref{prop:local} again gives
$\mathbf e\in\Soc(\mathcal A)\subseteq\m$,
and since $\mathbf e\m=0$, it follows that $\mathbf e^2=0$.
\end{proof}

\Needspace{3\baselineskip}
\subsection{Vertex operator algebras associated with commutative Frobenius algebras}\label{subsec:frob-voa}
Let $\mathcal A$ be a finite-dimensional commutative associative unital algebra
with identity element $1_{\mathcal A}$, and let $f:\mathcal A\to\C$ be a linear functional.
In this subsection we first construct the vertex operator algebra associated with $(\mathcal A,f)$,
and then use the local decomposition of Proposition~\ref{prop:local} to establish the corresponding tensor product decomposition.
Finally, we show that the case of a semisimple local algebra corresponds to the usual Virasoro vertex operator algebra,
and recall its simplicity and the rationality, regularity, and $C_2$-cofiniteness of its simple quotient.

We use the Lie algebra construction of \cite[Theorem~4.1, (4.3)]{LYu} in the purely even case,
written in the following form adapted to our conventions for the multiplication, the bilinear form, and the indexing of modes.

\begin{definition}\label{def:froblie}
The \emph{Lie algebra $\vir(\mathcal A,f)$ associated with $(\mathcal A,f)$} has underlying vector space
\[
 \vir(\mathcal A,f)=(\mathcal A\otimes\C[t,t^{-1}])\oplus\C K,
 \qquad L_m(a)=a\otimes t^m,
\]
where $t$ is a formal variable, $a\in\mathcal A$, $m\in\Z$, and $K$ is a central element. For $a,b\in\mathcal A$ and $m,n\in\Z$, its bracket is
\begin{equation}\label{eq:frobbracket}
 [L_m(a),L_n(b)]
 =(m-n)L_{m+n}(ab)
 +\frac{m^3-m}{12}\delta_{m+n,0}f(ab)K.
\end{equation}
\end{definition}

The commutativity and associativity of $\mathcal A$ imply, by a direct computation, that this bracket satisfies the Jacobi identity.
Let
\[
 \mathfrak p=\C K\oplus\Span\{L_m(a):a\in \mathcal A,\ m\in\Z_{\ge-1}\}.
\]
By \eqref{eq:frobbracket}, $\mathfrak p$ is a Lie subalgebra, and
\[
 K\one=\one,\qquad L_m(a)\one=0,
 \qquad a\in \mathcal A,\quad m\in\Z_{\ge-1},
\]
defines a one-dimensional $\mathfrak p$-module $\C\one$.

The following construction can be found in \cite{Lam}, \cite[Sections~5.1--5.2]{LamT},
and \cite[(5.8), Proposition~5.1]{LYu}.
To compare with \cite[(4.3)]{LYu}, take the multiplication $a\cdot b=2ab$ on $\mathcal A$
and the bilinear form $\frac{f(ab)}{2}$, in the purely even case and at central level $1$;
the element $a(m+1)$ there corresponds to our $L_m(a)$.
By \eqref{eq:frobbracket}, the conformal vector is $L_{-2}(1_{\mathcal A})\one$
and the central charge is $f(1_{\mathcal A})$. For the existence of a unique maximal proper ideal, see
\cite[Proposition~5.2.18]{LamT}.

\begin{proposition}\label{prop:frobuniversal}
Let
\[
 V_{\vir(\mathcal A,f)}(1,0)
 =U(\vir(\mathcal A,f))\otimes_{U(\mathfrak p)}\C\one
\]
be the induced vacuum module, and denote $1\otimes\one$ again by $\one$. This induced vacuum module carries a unique vertex operator algebra structure with vacuum vector $\one$
satisfying the formulas below.
Equipped with this structure, $V_{\vir(\mathcal A,f)}(1,0)$ is called the \textit{vertex operator algebra associated with $(\mathcal A,f)$}.
For $a\in \mathcal A$, set $v_a=L_{-2}(a)\one$, and write $v_{\mathcal A}=\{v_a:a\in \mathcal A\}$. Then
\[
 Y(v_a,z)=\sum_{m\in\Z}L_m(a)z^{-m-2},
 \qquad \omega=v_{1_{\mathcal A}},\qquad c=f(1_{\mathcal A}).
\]
In particular, $L(n)=L_n(1_{\mathcal A})$ for $n\in\Z$. This vertex operator algebra is of CFT type, and
\[
 V_{\vir(\mathcal A,f)}(1,0)_0=\C\one,\qquad
 V_{\vir(\mathcal A,f)}(1,0)_1=0,\qquad
 V_{\vir(\mathcal A,f)}(1,0)_2=v_{\mathcal A}\cong \mathcal A,
\]
where $v_{\mathcal A}\cong \mathcal A$ is the vector space isomorphism given by $a\mapsto v_a$.
Fix an ordered basis of $\mathcal A$, and fix a total order on the corresponding modes $L_{-n}(a)$, $n\in\Z_{\ge2}$.
Then a PBW basis of $V_{\vir(\mathcal A,f)}(1,0)$ is given by the ordered monomials
\[
 L_{-n_1}(a_1)\cdots L_{-n_\ell}(a_\ell)\one,
\]
where $\ell\in\N$, $n_j\in\Z_{\ge2}$, and each $a_j$, $1\le j\le\ell$, is taken from this basis;
for $\ell=0$ the monomial is $\one$. Such a monomial has conformal weight
$\sum_{j=1}^{\ell}n_j$.
This vertex operator algebra is strongly generated by $v_{\mathcal A}\cong\mathcal A$ and has a unique maximal proper ideal;
the corresponding simple quotient is denoted by $L_{\vir(\mathcal A,f)}(1,0)$.
\end{proposition}
\begin{remark}\label{rem:frob-vacuum-pbw}
\begin{enumerate}[label=(\arabic*)]
\item Following the notation of \cite[Section~6.1]{LL} for Virasoro vertex operator algebras,
we write $V_{\vir(\mathcal A,f)}(1,0)$ for the above induced vacuum module.
Here $1$ indicates that $K$ acts as the identity, and $0$ indicates that each $L_0(a)$ acts as zero on the vacuum.

\item Let $\mathfrak n_-=\Span\{L_m(a):a\in \mathcal A,\ m\in\Z_{\le-2}\}$.
By the PBW theorem, there is a vector space isomorphism
\begin{equation}\label{eq:frob-pbw}
 V_{\vir(\mathcal A,f)}(1,0)\cong U(\mathfrak n_-).
\end{equation}

\item Arguing as in \cite[Remarks~6.1.3, 6.1.6]{LL}, one obtains the following universal property of
$V_{\vir(\mathcal A,f)}(1,0)$.
Let $W$ be a vertex operator algebra, and let $\iota:\mathcal A\to W_2$ be a linear map
such that $\iota(1_{\mathcal A})=\omega_W$ and, for $a,b\in \mathcal A$,
\[
\begin{aligned}
 \iota(a)_0\iota(b)&=L_W(-1)\iota(ab),\\
 \iota(a)_1\iota(b)&=2\iota(ab),\\
 \iota(a)_2\iota(b)&=0,\\
 \iota(a)_3\iota(b)&=\tfrac12f(ab)\one_W,\\
 \iota(a)_j\iota(b)&=0,\qquad j\in\Z_{\ge4}.
\end{aligned}
\]
Then there exists a unique vertex operator algebra homomorphism
\[
\begin{aligned}
 \widetilde\iota:V_{\vir(\mathcal A,f)}(1,0)&\longrightarrow W,\\
 v_a&\longmapsto\iota(a),\qquad a\in\mathcal A.
\end{aligned}
\]
\end{enumerate}
\end{remark}

Now let $V$ satisfy the hypotheses of Remark~\ref{rem:frobenius-normalization},
that is, $V$ is a nonzero simple vertex operator algebra strongly generated by $V_2$.
Let $A=V_2$ with the multiplication $ab=\tfrac12a_1b$ and identity element $1_A=\omega$,
and let $f(a)=2B_A(a,1_A)$ for $a\in A$.
By Proposition~\ref{prop:wt2} and \eqref{eq:normalization},
the universal property in Remark~\ref{rem:frob-vacuum-pbw}~(3) yields a vertex operator algebra homomorphism
\[
\begin{aligned}
 \mathsf q_V^{\mathfrak{vir}}:V_{\vir(A,f)}(1,0)&\longrightarrow V,\\
 v_a&\longmapsto a,\qquad a\in A.
\end{aligned}
\]
Its image contains the generating space $V_2$, so $\mathsf q_V^{\mathfrak{vir}}$ is surjective.
Hence $V$ is the simple quotient of $V_{\vir(A,f)}(1,0)$.

The Virasoro tensor decomposition given by Lam for semisimple algebras can be found in \cite[Theorem~5.2.21]{LamT}.
Using the local decomposition of Proposition~\ref{prop:local} and the above universal property,
we now formulate this decomposition in a form that applies to the local commutative Frobenius summands.

\begin{proposition}\label{prop:frobfactors}
Let $V$ be simple and strongly generated by $V_2$,
and let $(A,f)$ be as in Remark~\ref{rem:frobenius-normalization}.
Let $A=\bigoplus_{i=1}^s A_i$ be the decomposition of Proposition~\ref{prop:local},
where $A_i=e_iA$, the $e_i$ are primitive idempotents, and $f_i=f|_{A_i}$ for $1\le i\le s$.
Then there are vertex operator algebra isomorphisms
\begin{align}
 V_{\vir(A,f)}(1,0)&\cong\bigotimes_{i=1}^s V_{\vir(A_i,f_i)}(1,0),\label{eq:frobfactors-universal}\\
 V&\cong\bigotimes_{i=1}^s L_{\vir(A_i,f_i)}(1,0).\label{eq:frobfactors}
\end{align}
Every automorphism of $V$ permutes the factors $L_{\vir(A_i,f_i)}(1,0)$.
For every $g\in\Aut(V)$, the restriction $g|_A$ is a unital algebra automorphism preserving $f$, that is,
\[
 g(ab)=g(a)g(b),\qquad g(1_A)=1_A,\qquad f(g(a))=f(a),
 \qquad a,b\in A.
\]
\end{proposition}
\begin{proof}
If $V=\C\one$, then $A=0$, and the assertion holds by the convention for the empty tensor product.
Assume from now on that $V\ne\C\one$, so that $A\ne0$ and hence $s\ge1$.
Denote the vacuum of $V_{\vir(A_i,f_i)}(1,0)$ by $\one_i$, and set
$v_a^{(i)}=L_{-2}(a)\one_i$ for $a\in A_i$.
Embed each $v_a^{(i)}$ into the corresponding tensor factor.
Modes of different factors commute; together with $A_iA_j=0$ for $i\ne j$ and \eqref{eq:frobbracket},
this shows that the linear map given by these embeddings satisfies the singular product relations in Remark~\ref{rem:frob-vacuum-pbw}~(3).
Moreover, since $1_A=\sum_{i=1}^s e_i$, the image of $1_A$ is the sum of the conformal vectors of the factors,
which is the conformal vector of the tensor product; see \eqref{eq:tensor-vectors}.
Hence Remark~\ref{rem:frob-vacuum-pbw}~(3) gives a vertex operator algebra homomorphism
\[
\begin{aligned}
 \Phi:V_{\vir(A,f)}(1,0)&\longrightarrow\bigotimes_{i=1}^s V_{\vir(A_i,f_i)}(1,0),\\
 v_a&\longmapsto\one_1\otimes\cdots\otimes v_a^{(i)}
 \otimes\cdots\otimes\one_s,\qquad a\in A_i.
\end{aligned}
\]
Choose a basis of each $A_i$, order the modes $L_{-n}(a)$, $n\in\Z_{\ge2}$, according to the order of the summands,
and fix a total order within each summand.
The PBW basis vectors of the source space are then of the form $P_1\cdots P_s\one$,
where $P_i$ is an ordered monomial in the negative modes associated with the basis elements of the $i$-th summand,
and $P_i=1$ is allowed.
By the construction of $\Phi$,
\[
 \Phi(P_1\cdots P_s\one)
 =(P_1\one_1)\otimes\cdots\otimes(P_s\one_s).
\]
The right-hand side runs exactly once through the tensor product basis of the target space,
so $\Phi$ is an isomorphism.

For each $i$, the vertex operator algebra $L_{\vir(A_i,f_i)}(1,0)$ is simple and of CFT type,
and $L_{\vir(A_i,f_i)}(1,0)_1=0$.
By Proposition~\ref{prop:form}~(3), each of them has a normalized nondegenerate invariant bilinear form $B_i$.
By Proposition~\ref{prop:tensorforms}, $B=\bigotimes_{i=1}^s B_i$ is again normalized, invariant, and nondegenerate.
The corresponding tensor product
$\bigotimes_{i=1}^s L_{\vir(A_i,f_i)}(1,0)$ is a vertex operator algebra of CFT type,
and hence is simple by Proposition~\ref{prop:form}~(3).
Composing $\Phi$ with the tensor product of the natural quotient maps of the factors, we obtain a surjective vertex operator algebra homomorphism
\[
 V_{\vir(A,f)}(1,0)\longrightarrow
 \bigotimes_{i=1}^s L_{\vir(A_i,f_i)}(1,0).
\]
Therefore $\bigotimes_{i=1}^s L_{\vir(A_i,f_i)}(1,0)$ is the
simple quotient of $V_{\vir(A,f)}(1,0)$.
Since $V$ was shown above to be the simple quotient of $V_{\vir(A,f)}(1,0)$,
the uniqueness of the maximal proper ideal in Proposition~\ref{prop:frobuniversal} gives
\[
 V\cong\bigotimes_{i=1}^s L_{\vir(A_i,f_i)}(1,0),
\]
which is \eqref{eq:frobfactors}.

Finally, let $g\in\Aut(V)$. Since $g$ preserves conformal weights and the conformal vector,
we have $g(A)=A$ and $g(1_A)=1_A$.
By \eqref{eq:wt2products} and \eqref{eq:normalization}, we have, for $a,b\in A$,
\[
\begin{aligned}
 g(ab)&=\tfrac12g(a_1b)=g(a)g(b),\\
 \tfrac12f(g(a))\one&=(1_A)_3g(a)=g((1_A)_3a)=\tfrac12f(a)\one.
\end{aligned}
\]
Thus $g|_A$ is a unital algebra automorphism preserving $f$.
By the uniqueness of the local decomposition in Proposition~\ref{prop:local}, $g$ permutes the primitive
idempotents $e_i$, and hence permutes the summands $A_i=e_iA$.
By the strong generation statement of Proposition~\ref{prop:frobuniversal},
under the isomorphism \eqref{eq:frobfactors}, each factor $L_{\vir(A_i,f_i)}(1,0)$ corresponds to the vertex subalgebra of $V$ generated by $A_i$,
so $g$ permutes the factors $L_{\vir(A_i,f_i)}(1,0)$.
\end{proof}

We next consider the case of a semisimple local algebra, that is, $\mathcal A=\C1_{\mathcal A}$;
see Remark~\ref{rem:local-properties}~(2).
The case of a nonsemisimple local commutative Frobenius algebra is treated in the next section;
in that case, $V_{\vir(\mathcal A,f)}(1,0)$ is already simple (see Theorem~\ref{thm:frobsimple}).

We first recall the construction of the Virasoro vertex operator algebra; see \cite[Section~4]{FZ} and
\cite[Section~6.1]{LL}.

\begin{definition}\label{def:virasoro-voa}
Let $\mathfrak{vir}$ be the Virasoro algebra, with basis consisting of $L(n)$, $n\in\Z$, and a central element $K$, satisfying
\[
 [L(m),L(n)]=(m-n)L(m+n)
 +\frac{m^3-m}{12}\delta_{m+n,0}K,\qquad m,n\in\Z.
\]
For $c\in\C$, set
\[
 \mathfrak b=\C K\oplus\operatorname{span}_{\C}\{L(n):n\ge-1\},
\]
and let $\C_c=\C\one$ be the one-dimensional $\mathfrak b$-module with
$K\one=c\one$ and $L(n)\one=0$ for $n\ge-1$.
Following the notation of \cite[Section~6.1]{LL}, we denote the corresponding induced vacuum module by
\[
 V_{\mathfrak{vir}}(c,0)=U(\mathfrak{vir})\otimes_{U(\mathfrak b)}\C_c,
\]
and denote $1\otimes\one$ again by $\one$.
This module carries a unique vertex operator algebra structure such that
\[
 \omega=L(-2)\one,\qquad
 Y(\omega,z)=\sum_{n\in\Z}L(n)z^{-n-2}.
\]
It is called the \emph{Virasoro vertex operator algebra of central charge $c$}.
Its unique simple quotient is denoted by $L_{\mathfrak{vir}}(c,0)$ and is called the \emph{simple Virasoro vertex operator algebra of central charge $c$}.
\end{definition}

Now let $\mathcal A=\C1_{\mathcal A}$ and write $c=f(1_{\mathcal A})$.
On $V_{\vir(\mathcal A,f)}(1,0)$, the bracket \eqref{eq:frobbracket} becomes
\[
 [L_m(1_{\mathcal A}),L_n(1_{\mathcal A})]
 =(m-n)L_{m+n}(1_{\mathcal A})
 +\frac{m^3-m}{12}\delta_{m+n,0}c\,\id,
 \qquad m,n\in\Z.
\]
Combining this with $L_m(1_{\mathcal A})\one=0$ for $m\ge-1$ and Definition~\ref{def:virasoro-voa}, we obtain
vertex operator algebra isomorphisms
\begin{align}
 V_{\vir(\mathcal A,f)}(1,0)&\cong V_{\mathfrak{vir}}(c,0),\label{eq:vir-special-universal}\\
 L_{\vir(\mathcal A,f)}(1,0)&\cong L_{\mathfrak{vir}}(c,0).\label{eq:vir-special}
\end{align}
The isomorphism \eqref{eq:vir-special-universal} sends $v_{1_{\mathcal A}}$ to $\omega=L(-2)\one$,
and the isomorphism \eqref{eq:vir-special} is obtained by passing to the unique simple quotients.
Thus $V_{\vir(\mathcal A,f)}(1,0)$ is a generalization of $V_{\mathfrak{vir}}(c,0)$ to an arbitrary $\mathcal A$.

The following results will be used in the analysis of the Virasoro factors in Section~\ref{sec:twists}.
For the simplicity of $V_{\mathfrak{vir}}(c,0)$, see \cite{FZ,Wang,ZD}.
The rationality, regularity, and $C_2$-cofiniteness of the Virasoro minimal models
are established in \cite[Theorem~4.2]{Wang}, \cite[Theorem~3.13]{DLMreg},
and \cite[Lemma~12.3]{DLMmod}, respectively.
Since these vertex operator algebras are of CFT type, Proposition~\ref{prop:knownC2}
(see \cite[Theorem~4.5]{ABD}) shows that the above regularity by itself already implies their rationality and $C_2$-cofiniteness;
the latter also follows from \cite[Theorem~3.8]{Li}.

\begin{proposition}\label{prop:vir}
Let $c\in\C$. For integers $p,q$ with $2\le p<q$ and $\gcd(p,q)=1$, write
\[
 c_{p,q}=1-\frac{6(p-q)^2}{pq}.
\]
If $c=c_{p,q}$, then $L_{\mathfrak{vir}}(c,0)$ is called a \textit{Virasoro minimal model}; it is rational, regular, and $C_2$-cofinite.
In particular, $c_{2,3}=0$ and $L_{\mathfrak{vir}}(0,0)=\C\one$.
If $c$ is not of this form, then $V_{\mathfrak{vir}}(c,0)=L_{\mathfrak{vir}}(c,0)$.
\end{proposition}

\section{Graph contractions, invariant bilinear forms, and simplicity}\label{sec:frobgram}
This section studies the vertex operator algebra $V_{\vir(\mathcal A,f)}(1,0)$ associated with a
nonsemisimple local commutative Frobenius algebra $(\mathcal A,f)$, and proves that it is simple.
This result determines the factors in Proposition~\ref{prop:frobfactors} corresponding to the nonsemisimple local summands,
and will be used in the proof of Theorem~\ref{thm:two} in Section~\ref{sec:twists}.

When $\mathcal A=\C[\varepsilon]/(\varepsilon^2)$ and $f(1_{\mathcal A})=f(\varepsilon)\ne0$,
the corresponding result is \cite[Theorem~2.1]{ZD}; when $\mathcal A=\C[\varepsilon]/(\varepsilon^N)$ with $N\ge2$,
the nondegeneracy of the Frobenius form is equivalent to $f(\varepsilon^{N-1})\ne0$,
and the corresponding result is \cite[Theorem~3.5, Proposition~3.7~(1)]{JZ}.
We treat here a general nonsemisimple local commutative Frobenius algebra.

By Propositions~\ref{prop:frobuniversal} and~\ref{prop:form}~(3),
this vertex operator algebra has a unique normalized invariant bilinear form.
We obtain the required simplicity by proving that the restriction of this form to each conformal weight space is nondegenerate.
Section~\ref{subsec:graphs} first represents the tensor contractions determined by the Frobenius form by finite graphs,
and uses the dual basis identity to reduce the evaluation of a connected graph to the evaluation of a power of the Euler element.
In the nonsemisimple local case, the identity $\mathbf e^2=0$ of Lemma~\ref{lem:euler}
then gives the vanishing statement of Lemma~\ref{lem:graph}.
Section~\ref{subsec:det-criterion} uses this vanishing to establish a criterion for Gram determinants.
In Section~\ref{subsec:pbw-simplicity}, by means of PBW symmetrization, we compare the invariant bilinear form
with an auxiliary bilinear form determined only by the central terms of the commutation relations,
show that the Gram determinants of the two forms with respect to the same basis are equal and nonzero,
and finally deduce simplicity from Remark~\ref{rem:weight-gram} and Proposition~\ref{prop:form}~(3).

\Needspace{3\baselineskip}
\subsection{Graphs and tensor contractions}\label{subsec:graphs}
In this subsection we fix a commutative Frobenius algebra $(\mathcal A,f)$.
For $a,b\in\mathcal A$, write $\eta(a,b)=f(ab)$,
and let $\mathbf e$ be the Euler element of Definition~\ref{def:euler}.
We represent the tensor contractions determined by the Frobenius form by finite graphs, and compute their values by edge contraction.
For a connected graph in which every vertex has degree at least $2$,
we express its evaluation as the value under $f$ of a power of the Euler element;
moreover, when $\mathcal A$ is local and $\mathcal A\not\cong\C$,
we show that the corresponding contraction vanishes as soon as some vertex has degree greater than $2$.
This vanishing will be used in the determinant criterion of the next subsection.

Let $d=\dim\mathcal A$, choose a basis $b_1,\ldots,b_d$ of $\mathcal A$, and write
\[
 \eta_{ij}=\eta(b_i,b_j),\qquad
 (\eta^{ij})=(\eta_{ij})^{-1},\qquad 1\le i,j\le d.
\]
For each integer $k\ge2$ and $a_1,\ldots,a_k\in\mathcal A$, define the $k$-linear form
\[
 \Theta_k:\mathcal A^k\longrightarrow\C,\qquad
 \Theta_k(a_1,\ldots,a_k)=f(a_1\cdots a_k).
\]
Since $\mathcal A$ is commutative, $\Theta_k$ is symmetric, and $\Theta_2=\eta$.
The notation $\Theta_k$ is ours; the product formula above appears in \cite[(3.11)--(3.12)]{DMu}.
In coordinates, a tensor contraction is the operation of setting specified pairs of indices equal and summing over them;
when several tensors are involved, one first takes their tensor product.
For the general definition, see \cite[Section~3.3, (117)--(118), p.~35]{KKEG}.
Below we use graphs to record how the inputs of these multilinear forms are paired,
and use the entries $\eta^{ij}$ of the inverse matrix to define the specific sums needed in this paper.

\begin{definition}\label{def:graph}
A \emph{finite graph} $\Gamma$ has a finite vertex set $\operatorname{Vert}(\Gamma)$
and a finite edge set $\operatorname{Edge}(\Gamma)$.
Each edge has two distinct edge ends; each end belongs to a unique edge
and is attached to a unique vertex. If the two ends of an edge are attached to the same vertex,
the edge is called a \emph{self-loop}. Distinct edges may join the same two vertices (which are allowed to coincide);
such edges are called \emph{multiple edges}.
Denote by $\mathcal H$ the set of all edge ends;
for $v\in\operatorname{Vert}(\Gamma)$, denote by $\mathcal H_v$ the set of ends attached to $v$.
The \emph{degree} of a vertex $v$ is defined by $\deg v=|\mathcal H_v|$,
so each self-loop contributes $2$ to the degree.
We write each edge as the unordered pair $\{h,h'\}$ of its two ends, where $h,h'\in\mathcal H$.

Suppose that every vertex has degree at least $2$. For each map
$\kappa:\mathcal H\to\{1,\ldots,d\}$, we regard $\kappa(h)$ as the basis index assigned to the end $h$.
The \emph{tensor contraction} corresponding to $\Gamma$ is defined by
\begin{equation}\label{eq:graphsum}
 \mathcal Z(\Gamma)=
 \sum_{\kappa:\mathcal H\to\{1,\ldots,d\}}
 \left\{\prod_{v\in\operatorname{Vert}(\Gamma)}
 f\!\left(\prod_{h\in\mathcal H_v}b_{\kappa(h)}\right)\right\}
 \left\{\prod_{\{h,h'\}\in\operatorname{Edge}(\Gamma)}
 \eta^{\kappa(h)\kappa(h')}\right\}.
\end{equation}
For this construction of graph evaluations, see \cite[Theorem~2.28, Section~2.7]{FM}.
By the commutativity of $\mathcal A$ and the symmetry $\eta^{ij}=\eta^{ji}$, the right-hand side does not depend on the ordering of the ends at each vertex,
nor on the ordering of the two ends of each edge.
The index of each end occurs in exactly one vertex factor and in exactly one edge factor,
and all indices are summed over, so $\mathcal Z(\Gamma)\in\C$.
If the vertex set is empty, $\Gamma$ is called the \emph{empty graph}; in this case the edge set is also empty.
With the convention that an empty product equals $1$, the value of the empty graph is $1$.

A \emph{path} is a finite sequence, alternating between vertices and edges and beginning and ending with a vertex, in which each edge
joins the two vertices immediately before and after it, and the vertices are pairwise distinct; a sequence consisting of a single vertex is also regarded as a path.
A \emph{cycle} is a finite sequence, alternating between vertices and edges and beginning and ending with the same vertex, in which each edge joins the two vertices immediately before and after it,
at least one edge occurs, the edges are pairwise distinct, and no vertex is repeated except that the first and last vertices coincide.
This definition includes self-loops and cycles formed by two distinct edges joining the same two vertices.
If $\Gamma$ is nonempty and any two vertices can be joined by a path, $\Gamma$ is said to be \emph{connected}.
If the vertices are grouped according to whether they can be joined by a path, then each group of vertices, together with all edges joining them,
forms a \emph{connected component}.
A \emph{spanning tree} of a connected graph is a connected graph without cycles obtained by keeping all vertices
and selecting some of the edges.
\end{definition}

\begin{remark}\label{rem:graph-contraction}
Denote the tensor corresponding to $\eta$ by
\[
 C_\eta=\sum_{i,j=1}^{d}\eta^{ij}b_i\otimes b_j
 \in\mathcal A\otimes\mathcal A.
\]
The notation $C_\eta$ is ours; its coordinate expression can be found in \cite[Section~2.7]{FM}.
Formula \eqref{eq:graphsum} is the coordinate expression obtained by using $\Theta_{\deg v}$ at each vertex $v$
and pairing along each edge by means of $C_\eta$.
Since $\eta$ is nondegenerate and symmetric and each $\Theta_k$ is symmetric,
the graph evaluation construction of \cite[Theorem~2.28, Section~2.7]{FM} shows that $\mathcal Z(\Gamma)$ is independent of the choice of basis.
\end{remark}

In our notation, the dual basis identity \cite[(2.7)]{DMu} reads as follows.
\begin{lemma}\label{lem:dual-basis}
For every $a\in\mathcal A$,
\[
 a=\sum_{i,j=1}^d\eta(a,b_i)\eta^{ij}b_j
  =\sum_{i,j=1}^df(ab_i)\eta^{ij}b_j.
\]
\end{lemma}

Using Lemma~\ref{lem:dual-basis} and the formula for the Euler element, we now compute two examples involving a self-loop or multiple edges.

\begin{example}\label{ex:contractions}
If $\Gamma$ has a single vertex and a single self-loop, then Definition~\ref{def:euler} and Lemma~\ref{lem:euler} give
\[
 \mathcal Z(\Gamma)=\sum_{i,j=1}^d\eta^{ij}f(b_ib_j)
 =f(\mathbf e)=d.
\]
If $\Gamma$ consists of two vertices and three edges joining them, then both vertices have degree $3$, and
\[
 \mathcal Z(\Gamma)=\sum_{i,j,k,p,q,r=1}^d
 f(b_ib_jb_k)f(b_pb_qb_r)
 \eta^{ip}\eta^{jq}\eta^{kr}.
\]
Fix $j,k,q,r$, and take $a=b_jb_k$ in the expansion of Lemma~\ref{lem:dual-basis}. This gives
\[
 b_jb_k=\sum_{i,p=1}^d\eta(b_jb_k,b_i)\eta^{ip}b_p
       =\sum_{i,p=1}^df(b_ib_jb_k)\eta^{ip}b_p.
\]
Multiplying both sides by $b_qb_r$ and applying $f$, we get
\[
 \sum_{i,p=1}^df(b_ib_jb_k)\eta^{ip}f(b_pb_qb_r)
 =f(b_jb_kb_qb_r).
\]
Hence, summing first over $i,p$ and then using the formula for the Euler element in Definition~\ref{def:euler}, we obtain
\[
\mathcal Z(\Gamma)=\sum_{j,k,q,r=1}^d\eta^{jq}\eta^{kr}f(b_jb_kb_qb_r)
=f\!\left(
\left(\sum_{j,q=1}^d\eta^{jq}b_jb_q\right)
\left(\sum_{k,r=1}^d\eta^{kr}b_kb_r\right)
\right)=f(\mathbf e^2).
\]
If $\mathcal A$ is local and $\mathcal A\not\cong\C$, then $\mathbf e^2=0$ by Lemma~\ref{lem:euler},
so the second tensor contraction vanishes, while the first equals $\dim\mathcal A$.
\end{example}

We now compute $\mathcal Z(\Gamma)$ for a general connected graph by edge contraction.
Here an \emph{edge contraction} means deleting an edge joining two distinct vertices and merging these two vertices,
with the remaining ends previously attached to the two vertices reattached to the merged vertex.
Using Lemma~\ref{lem:dual-basis}, we show that this operation leaves $\mathcal Z(\Gamma)$ unchanged;
we then contract successively along a spanning tree and treat the remaining self-loops by means of Definition~\ref{def:euler}.

\begin{lemma}\label{lem:graph}
Let $\Gamma$ be a connected finite graph with $v_\Gamma$ vertices and $e_\Gamma$ edges,
in which every vertex has degree at least $2$. Then
\begin{equation}\label{eq:graphEuler}
 \mathcal Z(\Gamma)=f(\mathbf e^{\,e_\Gamma-v_\Gamma+1}).
\end{equation}
If, in addition, $\mathcal A$ is local, $\mathcal A\not\cong\C$,
and at least one vertex has degree greater than $2$, then $\mathcal Z(\Gamma)=0$.
For a finite graph $\Gamma$, not necessarily connected, in which every vertex has degree at least $2$,
let $\ell$ be the number of its connected components, and denote them by $\Gamma_1,\ldots,\Gamma_\ell$. Then
\[
 \mathcal Z(\Gamma)=\prod_{j=1}^\ell \mathcal Z(\Gamma_j).
\]
\end{lemma}
\begin{proof}
For $a,b\in\mathcal A$, multiplying both sides of the expansion of $a$ in Lemma~\ref{lem:dual-basis} by $b$ and applying $f$, we obtain
\begin{equation}\label{eq:collapse}
 \sum_{i,j=1}^df(ab_i)\eta^{ij}f(b_jb)=f(ab).
\end{equation}
Now consider an edge joining distinct vertices $v,w$.
Denote by $i,j$ the indices of its two ends attached to $v,w$, respectively, and fix the indices of all other ends.
Let $a$ and $b$ be the products of the basis vectors corresponding to the remaining ends at $v$ and at $w$, respectively.
Then the factors given by these two vertices are $f(ab_i)$ and $f(b_jb)$, and the coefficient of the chosen edge is $\eta^{ij}$.
By \eqref{eq:collapse}, summing over $i,j$ replaces $f(ab_i)\eta^{ij}f(b_jb)$ by $f(ab)$,
while the remaining factors in \eqref{eq:graphsum} are unchanged.
Here $ab$ is precisely the product of the basis vectors corresponding to all ends at the merged vertex.
Therefore, after summing over the remaining indices, we obtain the sum for the contracted graph, so this edge contraction leaves $\mathcal Z(\Gamma)$ unchanged.

The contraction deletes the two ends of the chosen edge, so the merged vertex has degree
\[
 \deg v+\deg w-2\ge2+2-2=2.
\]
The degrees of the other vertices are unchanged, so the contracted graph still satisfies the degree condition of the lemma.

We construct a spanning tree of $\Gamma$ starting from a single vertex.
As long as some vertex is not yet included, connectedness allows us to choose an edge
with one end at a vertex already chosen and the other end at a vertex not yet chosen; we add this edge and the new vertex.
Each step adds exactly one vertex and one edge, and creates no cycle.
After finitely many steps we obtain a spanning tree containing all $v_\Gamma$ vertices,
with $v_\Gamma-1$ edges selected.
We contract these edges in the reverse of the order in which they were selected.
At each step, the edge to be contracted still joins two distinct vertices;
otherwise, the tree edges contracted earlier would contain a path joining the two endpoints of this edge,
and this path together with the edge would form a cycle, contradicting the fact that a spanning tree contains no cycles.
Hence the above edge contraction can be applied at each step, and $\mathcal Z(\Gamma)$ remains unchanged.
Each contraction deletes only the chosen edge; all other edges are kept.
After $v_\Gamma-1$ contractions, only one vertex remains, and all remaining edges are self-loops,
including the self-loops originally present in $\Gamma$.
Denote the number of these self-loops by $q$. Then
\[
 q=e_\Gamma-(v_\Gamma-1)=e_\Gamma-v_\Gamma+1.
\]
Since each edge has two ends and each vertex has degree at least $2$,
\begin{equation}\label{eq:graph-degree-sum}
 2e_\Gamma=\sum_{v\in\operatorname{Vert}(\Gamma)}\deg v\ge2v_\Gamma.
\end{equation}
Hence $q=e_\Gamma-v_\Gamma+1\ge1$.
For every $a\in\mathcal A$, Definition~\ref{def:euler} gives
\[
 \sum_{i,j=1}^d\eta^{ij}f(ab_ib_j)=f(a\mathbf e).
\]
In the contracted graph, choose a self-loop, fix the indices of all other ends,
and let $a$ be the product of the basis vectors corresponding to these ends, the empty product being $1_{\mathcal A}$.
The above identity shows that summing over the indices of the two ends of the chosen self-loop amounts to inserting a factor $\mathbf e$ into the product.
Since the indices of the different self-loops are summed independently, the linearity of $f$ and the commutativity of $\mathcal A$ give
\[
 \mathcal Z(\Gamma)=f\!\left(\left(\sum_{i,j=1}^d\eta^{ij}b_ib_j\right)^q\right)
 =f(\mathbf e^q),
\]
which is \eqref{eq:graphEuler}.

Now suppose that $\mathcal A$ is local, $\mathcal A\not\cong\C$, and at least one vertex has degree greater than $2$.
Then \eqref{eq:graph-degree-sum} gives $2e_\Gamma>2v_\Gamma$.
Since $e_\Gamma$ and $v_\Gamma$ are integers, $e_\Gamma\ge v_\Gamma+1$, and hence $q=e_\Gamma-v_\Gamma+1\ge2$.
By Lemma~\ref{lem:euler}, $\mathbf e^2=0$, so $\mathcal Z(\Gamma)=0$.
Finally, for a general $\Gamma$, we keep the notation $\Gamma_1,\ldots,\Gamma_\ell$ of the lemma for its connected components.
In \eqref{eq:graphsum}, the summation indices of different components are independent of each other,
so $\mathcal Z(\Gamma)=\prod_{j=1}^\ell \mathcal Z(\Gamma_j)$.
\end{proof}
\Needspace{3\baselineskip}
\subsection{Symmetric powers and a determinant criterion}\label{subsec:det-criterion}
In this subsection we fix a commutative Frobenius algebra $(\mathcal A,f)$.
For $a,b\in\mathcal A$, write $\eta(a,b)=f(ab)$.
We first define the vector space $W_n^{\mathfrak{vir}}$ and an auxiliary bilinear form $H_n^{\mathfrak{vir}}$ on it, and prove that $H_n^{\mathfrak{vir}}$ is nondegenerate.
Then, in the case where $\mathcal A$ is local and $\mathcal A\not\cong\C$, we give a criterion for a bilinear form $G_n^{\mathfrak{vir}}$:
expanding the values of $G_n^{\mathfrak{vir}}-H_n^{\mathfrak{vir}}$ according to groupings of the input vectors
and applying Lemma~\ref{lem:graph}, we prove that $G_n^{\mathfrak{vir}}$ and $H_n^{\mathfrak{vir}}$ have the same Gram determinant.
In the next subsection, by means of PBW symmetrization, this result will be applied to
the normalized invariant bilinear form on $V_{\vir(\mathcal A,f)}(1,0)$.

We next recall, as in \cite[Problem~2.11.3~(d), (e)]{EGH}, symmetric powers and their identification with permutation-invariant tensors.
\begin{definition}\label{def:symmetric-weight-space}
Let $r\in\Z_{\ge1}$, and denote by $S_r$ the symmetric group on the set $\{1,\ldots,r\}$.
For $\sigma\in S_r$ and $a_1,\ldots,a_r\in\mathcal A$, set
\[
 \sigma(a_1\otimes\cdots\otimes a_r)
 =a_{\sigma^{-1}(1)}\otimes\cdots\otimes a_{\sigma^{-1}(r)},
\]
and extend linearly; this gives an action of $S_r$ on $\mathcal A^{\otimes r}$.
The \emph{$r$-th symmetric power} is defined by
\[
 \Sym^r\mathcal A
 :=\mathcal A^{\otimes r}\Big/
 \Span\{u-\sigma u:\ u\in\mathcal A^{\otimes r},\ \sigma\in S_r\}.
\]
The image of $a_1\otimes\cdots\otimes a_r$ in this quotient is denoted by $a_1\cdots a_r$
and is called a \emph{symmetric product}; it is different from the product in $\mathcal A$.
Define the linear map
\[
 \begin{aligned}
 P_r:\mathcal A^{\otimes r}&\longrightarrow\mathcal A^{\otimes r},\\
 u&\longmapsto\frac1{r!}\sum_{\sigma\in S_r}\sigma u.
 \end{aligned}
\]
It is the symmetrization projection: $P_r^2=P_r$, and its image is
\[
 \operatorname{im}P_r
 =\{u\in\mathcal A^{\otimes r}:\ \sigma u=u\text{ for all }\sigma\in S_r\}.
\]
In characteristic zero, the following map is a vector space isomorphism:
\[
 \begin{aligned}
 \Sym^r\mathcal A&\xrightarrow{\ \sim\ }\operatorname{im}P_r,\\
 a_1\cdots a_r&\longmapsto P_r(a_1\otimes\cdots\otimes a_r).
 \end{aligned}
\]
In what follows, we identify these two spaces via this isomorphism.
For $r=0$, we use the conventions $\mathcal A^{\otimes0}=\Sym^0\mathcal A=\C$ and $P_0=\id_\C$.

\end{definition}

Following the conformal weight grading in Proposition~\ref{prop:frobuniversal}, we define the following spaces.
\begin{definition}\label{def:weight-space}
For $m\in\Z_{\ge2}$, let $\mathcal A_m$ be a copy of the vector space $\mathcal A$, assigned weight $m$.
Fix $n\in\N$. A sequence $(r_m)_{m\ge2}$ with $r_m\in\N$ and $\sum_{m\ge2}mr_m=n$
is called \emph{occupancy data} of weight $n$.
Define
\begin{equation}\label{eq:Wn}
 W_n^{\mathfrak{vir}}=\bigoplus_{\sum_{m=2}^{\infty}mr_m=n}
       \bigotimes_{m=2}^{\infty}\Sym^{r_m}\mathcal A_m.
\end{equation}
The summand corresponding to given occupancy data is called an \emph{occupancy summand}.
Since $r_m=0$ for $m>n$ and $0\le r_m\le\frac{n}{m}$, there are only finitely many summands,
and each summand has only finitely many nontrivial tensor factors.
Omitting the factors $\Sym^0\mathcal A_m=\C$ and taking the empty tensor product to be $\C$, we have
$W_0^{\mathfrak{vir}}=\C$ and $W_1^{\mathfrak{vir}}=0$. From now on we use the notation
\[
 \prod_{m\ge2}\prod_{i=1}^{r_m}a_{m,i}
 :=\bigotimes_{m\ge2}(a_{m,1}\cdots a_{m,r_m}),\qquad a_{m,i}\in\mathcal A_m,
\]
where each parenthesized expression is a symmetric product; such a vector is called a \emph{symmetric monomial}.
\end{definition}

\begin{example}\label{ex:symmetric-weight-space}
By the definition of weight, for the weights $2,3,4,5$ we have
\[
 W_2^{\mathfrak{vir}}=\mathcal A_2,\qquad W_3^{\mathfrak{vir}}=\mathcal A_3,\qquad
 W_4^{\mathfrak{vir}}=\mathcal A_4\oplus\Sym^2\mathcal A_2,\qquad
 W_5^{\mathfrak{vir}}=\mathcal A_5\oplus(\mathcal A_2\otimes\mathcal A_3).
\]
For example, there are exactly two occupancy data of weight $4$: $r_4=1$, or $r_2=2$, with all other $r_m$ equal to zero in each case.
\end{example}

We now define the auxiliary bilinear form $H_n^{\mathfrak{vir}}$ on $W_n^{\mathfrak{vir}}$. For $m\in\Z_{\ge2}$, write
\[
 \alpha_m=\frac{m^3-m}{12}.
\]
We first define a bilinear form on each $\Sym^r\mathcal A_m$, and then construct $H_n^{\mathfrak{vir}}$ according to the decomposition \eqref{eq:Wn}.

\begin{definition}\label{def:auxiliary-pairing}
For $m\ge2$ and $r\in\N$, we define a bilinear form $B_{\Sym^{r}\mathcal A_m}$ on $\Sym^r\mathcal A_m$.
For $x_i,y_i\in\mathcal A_m$, $1\le i\le r$, set
\[
 B_{\Sym^{r}\mathcal A_m}(x_1\cdots x_r,y_1\cdots y_r)
 =\alpha_m^r\sum_{\sigma\in S_r}\prod_{i=1}^r\eta(x_i,y_{\sigma(i)}),
\]
where $\eta$ is computed via $\mathcal A_m\simeq\mathcal A$.
For $r=0$, define $B_{\Sym^{0}\mathcal A_m}(1,1)=1$.
For each $n\in\Z_{\ge0}$, we define the \emph{auxiliary bilinear form} $H_n^{\mathfrak{vir}}$ on $W_n^{\mathfrak{vir}}$ as follows: distinct occupancy summands are mutually orthogonal,
and on a given occupancy summand we take the tensor product of the forms $B_{\Sym^{r_m}\mathcal A_m}$. That is,
for $x_{m,i},y_{m,i}\in\mathcal A_m$,
\begin{equation}\label{eq:osc}
 H_n^{\mathfrak{vir}}\!\left(\prod_{m\ge2}\prod_{i=1}^{r_m}x_{m,i},
             \prod_{m\ge2}\prod_{i=1}^{r_m}y_{m,i}\right)
 =\prod_{m\ge2}\left(\alpha_m^{r_m}\sum_{\sigma\in S_{r_m}}
       \prod_{i=1}^{r_m}\eta(x_{m,i},y_{m,\sigma(i)})\right).
\end{equation}
This formula is extended bilinearly, and a factor with $r_m=0$ is taken to be $1$; in particular, $H_0^{\mathfrak{vir}}(1,1)=1$.
\end{definition}

In order to write the inverse Gram matrix of $H_n^{\mathfrak{vir}}$ as a tensor, we recall the definition of the inverse tensor from \cite[Sections~2.1.2 and~2.7]{FM}.
\begin{definition}\label{def:inverse-tensor}
Let $E$ be a finite-dimensional space, and let $B$ be a nondegenerate symmetric bilinear form on $E$.
Let $N=\dim E$, and choose a basis $w_1,\ldots,w_N$ of $E$. The tensor
\[
 \sum_{i,j=1}^N
 \left[\left(B(w_k,w_l)\right)_{1\le k,l\le N}^{-1}\right]_{ij}
 w_i\otimes w_j\in E\otimes E
\]
is called the \emph{inverse tensor} of $B$.
If we number the $N^2$ terms of this sum and write it as $\sum_{q=1}^{N^2}u_q\otimes v_q$ with $u_q,v_q\in E$, then the tensor is uniquely determined by the condition
\[
 \sum_{q=1}^{N^2} B(x,u_q)v_q=x,\qquad x\in E,
\]
and is therefore independent of the basis. The tensor $C_\eta$ of Remark~\ref{rem:graph-contraction} is the inverse tensor of $\eta$.
\end{definition}

\begin{lemma}\label{lem:symmetric-pairing}
For each $n\in\N$, the form $H_n^{\mathfrak{vir}}$ is symmetric and nondegenerate.
For $m\ge2$ and $r\in\N$, the inverse tensor of $B_{\Sym^{r}\mathcal A_m}$ is
\begin{equation}\label{eq:symmetric-inverse}
 \frac1{r!\alpha_m^r}
 \sum_{\substack{1\le i_1,\ldots,i_r\le d\\1\le j_1,\ldots,j_r\le d}}
 \left(\prod_{\nu=1}^r\eta^{i_\nu j_\nu}\right)
 (b_{i_1}\cdots b_{i_r})\otimes(b_{j_1}\cdots b_{j_r}).
\end{equation}
Here we keep the basis $b_1,\ldots,b_d$ and the coefficients $\eta^{ij}$ of the previous subsection,
each $b_i$ is identified via $\mathcal A_m\simeq\mathcal A$, and the products on both sides are symmetric products.
For $r=0$, the expression \eqref{eq:symmetric-inverse} is taken to be $1\otimes1$.
\end{lemma}
\begin{proof}
By the symmetry of $\eta$, both $B_{\Sym^{r}\mathcal A_m}$ and $H_n^{\mathfrak{vir}}$ are symmetric.
For $r\in\N$, define the tensor pairing on $\mathcal A^{\otimes r}$ by
\[
 \eta^{\otimes r}(x_1\otimes\cdots\otimes x_r,y_1\otimes\cdots\otimes y_r)
 =\prod_{i=1}^r\eta(x_i,y_i),\qquad x_i,y_i\in\mathcal A.
\]
Since tensor products of $\eta$-dual bases are again dual bases, $\eta^{\otimes r}$ is nondegenerate.
For $u,v\in\mathcal A^{\otimes r}$ and $\sigma\in S_r$, we have
\[
 \eta^{\otimes r}(\sigma u,v)=\eta^{\otimes r}(u,\sigma^{-1}v),\qquad
 P_r^2=P_r.
\]
Averaging the first identity gives
$\eta^{\otimes r}(P_ru,v)=\eta^{\otimes r}(u,P_rv)$.
Hence $P_r$ is a self-adjoint projection with respect to this bilinear form.
If $u\in\operatorname{im}P_r$ and $\eta^{\otimes r}(u,\operatorname{im}P_r)=0$,
then for every $v\in\mathcal A^{\otimes r}$,
\[
 \eta^{\otimes r}(u,v)=\eta^{\otimes r}(P_ru,v)
 =\eta^{\otimes r}(u,P_rv)=0.
\]
By nondegeneracy, $u=0$, so $\eta^{\otimes r}\big|_{\operatorname{im}P_r\times\operatorname{im}P_r}$ is nondegenerate.

Let $u=x_1\otimes\cdots\otimes x_r$ and $v=y_1\otimes\cdots\otimes y_r$. Then
\[
 \eta^{\otimes r}(P_ru,P_rv)=\eta^{\otimes r}(u,P_rv)
 =\frac1{r!}\sum_{\sigma\in S_r}\prod_{i=1}^r\eta(x_i,y_{\sigma(i)}).
\]
Under the identification of Definition~\ref{def:symmetric-weight-space}, we have
\[
 B_{\Sym^{r}\mathcal A_m}
 =r!\alpha_m^r\,\eta^{\otimes r}\big|_{\operatorname{im}P_r\times\operatorname{im}P_r}.
\]
Since $\alpha_m\ne0$ for $m\ge2$, the form $B_{\Sym^{r}\mathcal A_m}$ is nondegenerate.
Taking finite tensor products and then the orthogonal direct sum over the occupancy summands, we see that $H_n^{\mathfrak{vir}}$ is nondegenerate.

We now compute the inverse tensor. Placing the first factor of each copy of $C_\eta$ in $C_\eta^{\otimes r}$ on the left
and the second factor on the right, we obtain the tensor in $\mathcal A^{\otimes r}\otimes\mathcal A^{\otimes r}$
\[
 T_r=\sum_{\substack{1\le i_1,\ldots,i_r\le d\\1\le j_1,\ldots,j_r\le d}}
 \left(\prod_{\nu=1}^r\eta^{i_\nu j_\nu}\right)
 (b_{i_1}\otimes\cdots\otimes b_{i_r})\otimes
 (b_{j_1}\otimes\cdots\otimes b_{j_r}).
\]
Applying Lemma~\ref{lem:dual-basis} to each tensor factor, we see that $T_r$ is the inverse tensor of $\eta^{\otimes r}$.
Numbering the $d^{2r}$ terms of this sum, we write $T_r=\sum_{q=1}^{d^{2r}}u_q\otimes v_q$ for short,
where $u_q,v_q\in\mathcal A^{\otimes r}$.
For $x\in\operatorname{im}P_r$, since $P_r$ is self-adjoint and $P_rx=x$, we have
\[
 \sum_{q=1}^{d^{2r}}\eta^{\otimes r}(x,P_ru_q)P_rv_q
 =P_r\!\left(\sum_{q=1}^{d^{2r}}\eta^{\otimes r}(x,u_q)v_q\right)=P_rx=x.
\]
Hence $(P_r\otimes P_r)T_r$ is the inverse tensor of $\eta^{\otimes r}\big|_{\operatorname{im}P_r\times\operatorname{im}P_r}$.
Dividing by $r!\alpha_m^r$ and using the identifications $\Sym^r\mathcal A\simeq\operatorname{im}P_r$
and $\mathcal A_m\simeq\mathcal A$, we obtain \eqref{eq:symmetric-inverse}.
\end{proof}

By grouping the inputs, we now apply Lemma~\ref{lem:graph} to compare the Gram determinants of two bilinear forms.
\begin{definition}\label{def:partition-blocks}
A \emph{partition} of a finite set is a family of pairwise disjoint nonempty subsets whose union is the whole set;
each of these subsets is called a \emph{block}; below, blocks are denoted by $\mathcal I$.
The partition of the empty set is taken to be the empty family.
\end{definition}

\begin{lemma}\label{lem:detcriterion}
Let $(\mathcal A,f)$ be a local commutative Frobenius algebra with $\mathcal A\not\cong\C$.
Fix $n\in\N$, let $G_n^{\mathfrak{vir}}$ be a bilinear form on $W_n^{\mathfrak{vir}}$, and write $B_n^{\Delta}=G_n^{\mathfrak{vir}}-H_n^{\mathfrak{vir}}$.

For any two occupancy data $(r_m)_{m\ge2}$ and $(r'_m)_{m\ge2}$ of weight $n$,
take symmetric monomials in the corresponding two summands,
\[
 x=\prod_{m\ge2}\prod_{i=1}^{r_m}x_{m,i},\qquad
 y=\prod_{m\ge2}\prod_{j=1}^{r'_m}y_{m,j},
\]
where $x_{m,i}\in\mathcal A_m$ for $1\le i\le r_m$ and
$y_{m,j}\in\mathcal A_m$ for $1\le j\le r'_m$.
Write $\ell=\sum_{m\ge2}r_m$ and $k=\sum_{m\ge2}r'_m$.
Via $\mathcal A_m\simeq\mathcal A$, denote the $\ell$ inputs of $x$ successively by
$\xi_1,\ldots,\xi_\ell\in\mathcal A$ and the $k$ inputs of $y$ successively by
$\xi_{\ell+1},\ldots,\xi_{\ell+k}\in\mathcal A$, and record, for each $\xi_j$, the space $\mathcal A_m$ to which the corresponding original input belongs.
Let $\mathcal P_{\ell+k}$ be the set of all partitions of $\{1,\ldots,\ell+k\}$ such that
every block has size at least $2$ and at least one block has size at least $3$.

Suppose that, for every such pair of occupancy data, one can choose scalars $\mu_\pi\in\C$,
$\pi\in\mathcal P_{\ell+k}$, independent of the particular values of the input vectors, such that for all corresponding input vectors
\begin{equation}\label{eq:partitions}
 B_n^{\Delta}(x,y)
 =\sum_{\pi\in\mathcal P_{\ell+k}}\mu_\pi
   \prod_{\mathcal I\in\pi}\Theta_{|\mathcal I|}\bigl((\xi_j)_{j\in\mathcal I}\bigr)
 =\sum_{\pi\in\mathcal P_{\ell+k}}\mu_\pi
   \prod_{\mathcal I\in\pi}f\!\left(\prod_{j\in\mathcal I}\xi_j\right),
\end{equation}
where the products inside $f$ are taken in $\mathcal A$.

Fix an arbitrary basis of $W_n^{\mathfrak{vir}}$; the determinant of each bilinear form below means its Gram determinant with respect to this basis. Then
\[
 \det(H_n^{\mathfrak{vir}}+tB_n^{\Delta})=\det H_n^{\mathfrak{vir}}\ne0,\qquad t\in\C.
\]
In particular, taking $t=1$, we obtain $\det G_n^{\mathfrak{vir}}=\det H_n^{\mathfrak{vir}}\ne0$.
\end{lemma}
\begin{proof}
If $W_n^{\mathfrak{vir}}=0$, the assertion holds by the convention $\det(0\times0)=1$.
If $n=0$, then $W_0^{\mathfrak{vir}}=\C$ and $\mathcal P_0$ is empty, so \eqref{eq:partitions} gives $B_0^{\Delta}=0$.
Assume from now on that $n\ge2$ and $q_n:=\dim W_n^{\mathfrak{vir}}>0$. We proceed in five steps.

\smallskip\noindent\emph{(1) The auxiliary operator and the trace expansion.}
By Lemma~\ref{lem:symmetric-pairing}, $H_n^{\mathfrak{vir}}$ is nondegenerate. By \cite[Theorem~3.16~(3), (4)]{Con},
there exists a unique endomorphism $N_n:W_n^{\mathfrak{vir}}\to W_n^{\mathfrak{vir}}$ such that
\[
 B_n^{\Delta}(x,y)=H_n^{\mathfrak{vir}}(x,N_ny),\qquad x,y\in W_n^{\mathfrak{vir}}.
\]
Then for every $t\in\C$,
\[
 (H_n^{\mathfrak{vir}}+tB_n^{\Delta})(x,y)=H_n^{\mathfrak{vir}}\bigl(x,({\id_{W_n^{\mathfrak{vir}}}}+tN_n)y\bigr).
\]
Taking Gram determinants with respect to the same basis, we obtain
\[
 \det(H_n^{\mathfrak{vir}}+tB_n^{\Delta})=\det H_n^{\mathfrak{vir}}\det({\id_{W_n^{\mathfrak{vir}}}}+tN_n).
\]
We first prove that $\Tr(N_n^s)=0$ for every integer $s\ge1$, and then use this to compute the second determinant on the right-hand side.

Choose a basis of each occupancy summand and combine them into a basis $w_1,\ldots,w_{q_n}$ of $W_n^{\mathfrak{vir}}$. Set
\[
 \mathsf h_{ij}=H_n^{\mathfrak{vir}}(w_i,w_j),\qquad \mathsf b_{ij}^{\Delta}=B_n^{\Delta}(w_i,w_j),\qquad
 (\mathsf h^{ij})=(\mathsf h_{ij})^{-1},\qquad 1\le i,j\le q_n.
\]
From $B_n^{\Delta}(x,y)=H_n^{\mathfrak{vir}}(x,N_ny)$, the matrix of $N_n$ is $(\mathsf h_{ij})^{-1}(\mathsf b_{ij}^{\Delta})$.
Fix $s\in\Z_{\ge1}$. Matrix multiplication gives
\begin{equation}\label{eq:traceindices}
 \Tr(N_n^s)=
 \sum_{\substack{1\le i_1,\ldots,i_s\le q_n\\1\le j_1,\ldots,j_s\le q_n}}
 \prod_{\nu=1}^s \mathsf h^{i_\nu j_\nu}\mathsf b_{j_\nu i_{\nu+1}}^{\Delta},
 \qquad i_{s+1}=i_1.
\end{equation}

\smallskip\noindent\emph{(2) Rewriting the trace by means of the inverse tensor.}
By Definition~\ref{def:inverse-tensor}, we write the inverse tensor of $H_n^{\mathfrak{vir}}$ as a finite sum
\[
 \sum_{i,j=1}^{q_n}\mathsf h^{ij}w_i\otimes w_j
 =\sum_{\lambda\in\Lambda}u_\lambda\otimes v_\lambda,
\]
where $\Lambda$ is a finite index set and $u_\lambda,v_\lambda\in W_n^{\mathfrak{vir}}$ for $\lambda\in\Lambda$.
Substituting $\mathsf b_{ij}^{\Delta}=B_n^{\Delta}(w_i,w_j)$ into \eqref{eq:traceindices}
and using the above expansion of the inverse tensor together with the bilinearity of $B_n^{\Delta}$, we obtain
\[
 \Tr(N_n^s)
 =\sum_{\lambda_1,\ldots,\lambda_s\in\Lambda}
   \prod_{\nu=1}^{s}B_n^{\Delta}(v_{\lambda_\nu},u_{\lambda_{\nu+1}}),
 \qquad \lambda_{s+1}=\lambda_1.
\]

We now choose the vectors $u_\lambda,v_\lambda$ in this finite sum using the explicit expansion given by \eqref{eq:symmetric-inverse}.
By Definition~\ref{def:auxiliary-pairing}, distinct occupancy summands are orthogonal with respect to $H_n^{\mathfrak{vir}}$,
so its inverse tensor is the sum of the inverse tensors of the summands.
Within each summand, the inverse tensor of each $B_{\Sym^{r_m}\mathcal A_m}$ is given by \eqref{eq:symmetric-inverse}.
Explicitly, for fixed $m\ge2$, a term of that expression is
\[
 \frac1{r_m!\alpha_m^{r_m}}
 \left(\prod_{\beta=1}^{r_m}\eta^{i_{m,\beta}j_{m,\beta}}\right)
 (b_{i_{m,1}}\cdots b_{i_{m,r_m}})\otimes
 (b_{j_{m,1}}\cdots b_{j_{m,r_m}}),
\]
where $1\le i_{m,\beta},j_{m,\beta}\le d$ for $1\le \beta\le r_m$; summing over these indices gives the inverse tensor of this factor.
We take the tensor product of these expansions over all $m$, placing the first factor of each copy on the left and the second factor on the right;
for $r_m=0$ we take $1\otimes1$.
In each term of the resulting expansion, absorbing the scalar coefficient into the left vector, we obtain the form $u_\lambda\otimes v_\lambda$ above.
Each of $u_\lambda,v_\lambda$ is then a scalar multiple of a symmetric monomial formed from $b_1,\ldots,b_d$ in the various $\mathcal A_m$.

\smallskip\noindent\emph{(3) Writing the expansion as graph evaluations.}
In the explicit expansion of step (2), each pair $u_\lambda,v_\lambda$ lies in the same occupancy summand.
For each $1\le\nu\le s$, we fix the occupancy data of the summand containing $u_{\lambda_\nu},v_{\lambda_\nu}$,
denote it by $(r_{\nu,m})_{m\ge2}$, and keep the summation over the basis indices.
After pulling out the scalar coefficients by the bilinearity of $B_n^{\Delta}$, we apply \eqref{eq:partitions}
to the remaining symmetric monomials in the $\nu$-th factor $B_n^{\Delta}(v_{\lambda_\nu},u_{\lambda_{\nu+1}})$ of the trace formula,
and choose one corresponding partition $\pi_\nu$.

We now construct a graph: the blocks of the partitions $\pi_\nu$ determine the vertices, and the pairings in \eqref{eq:symmetric-inverse} determine the edges.
We take the positions at which basis vectors occur in these $s$ factors $B_n^{\Delta}$ as ends, and denote the set of all of them by $\mathcal H$.
Different positions receive different labels, even when they correspond to the same basis vector.
Via the correspondence between input positions and ends, we regard each $\pi_\nu$ as a partition of the corresponding set of ends, and keep the notation $\mathcal I$ for its blocks.
For each block $\mathcal I\in\pi_\nu$ we introduce a vertex $v_{\mathcal I}$
and attach the ends in $\mathcal I$ to it. Thus
\[
 \mathcal H_{v_{\mathcal I}}=\mathcal I,\qquad
 \deg v_{\mathcal I}=|\mathcal I|,
\]
and its vertex factor is the corresponding $\Theta_{|\mathcal I|}$ in \eqref{eq:partitions}.

For the $\beta$-th input from $\mathcal A_m$ in the $\nu$-th pair $u_{\lambda_\nu},v_{\lambda_\nu}$,
where $1\le\nu\le s$, $m\ge2$, and $1\le \beta\le r_{\nu,m}$,
denote the basis indices on the two sides by $i_{\nu,m,\beta},j_{\nu,m,\beta}\in\{1,\ldots,d\}$, respectively.
By \eqref{eq:symmetric-inverse}, the coefficient corresponding to these two inputs is $\eta^{i_{\nu,m,\beta}j_{\nu,m,\beta}}$.
We join by an edge the two ends corresponding to $b_{i_{\nu,m,\beta}}$ in $u_{\lambda_\nu}$ and to $b_{j_{\nu,m,\beta}}$ in $v_{\lambda_\nu}$.
Since each $\pi_\nu$ is a partition of input positions, each end is attached to exactly one vertex;
moreover, by the way the pairing in \eqref{eq:symmetric-inverse} is made, each end belongs to exactly one edge.
We thus obtain a finite graph in the sense of Definition~\ref{def:graph}, which we denote by $\Gamma$.
An edge whose two ends belong to the same vertex is a self-loop; distinct edges may also join the same pair of vertices.

Once the above occupancy data and partitions are fixed, the graph $\Gamma$ is determined,
while the basis indices at the input positions still range independently over $1,\ldots,d$.
Denote an assignment of these indices by $\kappa:\mathcal H\to\{1,\ldots,d\}$.
Each vertex $v$ contributes $\Theta_{\deg v}((b_{\kappa(h)})_{h\in\mathcal H_v})$,
and each edge $\{h,h'\}$ contributes $\eta^{\kappa(h)\kappa(h')}$.
The product of the scalars $\mu_{\pi_\nu}$ and $(r_{\nu,m}!\alpha_m^{r_{\nu,m}})^{-1}$ does not depend on $\kappa$,
and can therefore be taken outside the summation. The remaining sum is precisely the graph evaluation defined by \eqref{eq:graphsum}:
\[
 \sum_{\kappa:\mathcal H\to\{1,\ldots,d\}}
 \left\{\prod_{v\in\operatorname{Vert}(\Gamma)}
 \Theta_{\deg v}\bigl((b_{\kappa(h)})_{h\in\mathcal H_v}\bigr)\right\}
 \left\{\prod_{\{h,h'\}\in\operatorname{Edge}(\Gamma)}
 \eta^{\kappa(h)\kappa(h')}\right\}
 =\mathcal Z(\Gamma).
\]
For fixed $n,s$, there are only finitely many choices of the above occupancy summands and partitions.
Denote the set of resulting graphs by $\mathcal G_{n,s}$, and let $\psi_\Gamma\in\C$ be the sum of the coefficients corresponding to the same graph $\Gamma$. Then
\[
 \Tr(N_n^s)=\sum_{\Gamma\in\mathcal G_{n,s}}\psi_\Gamma\mathcal Z(\Gamma).
\]

\smallskip\noindent\emph{(4) Vanishing of the traces.}
All blocks of each $\pi_\nu$ have size at least $2$, and at least one block has size at least $3$.
Hence, for each $\Gamma\in\mathcal G_{n,s}$, all vertex degrees are at least $2$,
and at least one vertex has degree greater than $2$.
Let $\Gamma_0$ be the connected component containing this vertex.
Since $\mathcal A$ is local and $\mathcal A\not\cong\C$, Lemma~\ref{lem:graph} gives $\mathcal Z(\Gamma_0)=0$.
By the multiplicativity of graph evaluations over connected components, also given in Lemma~\ref{lem:graph}, $\mathcal Z(\Gamma)=0$.
Substituting this into the finite sum of step (3), we obtain
\begin{equation}\label{eq:traces}
 \Tr(N_n^s)=0,\qquad s\in\Z_{\ge1}.
\end{equation}

\smallskip\noindent\emph{(5) Equality of the Gram determinants.}
Let $t$ be a formal variable, let $p(t)=\det({\id_{W_n^{\mathfrak{vir}}}}+tN_n)\in\C[t]$, and denote by $p'(t)=\frac{d}{dt}p(t)$ the derivative of $p(t)$ with respect to $t$.
Since $p(0)=1$, $p(t)$ is invertible in the formal power series ring $\C[[t]]$;
moreover, in $\operatorname{End}(W_n^{\mathfrak{vir}})[[t]]$ we have
\[
 ({\id_{W_n^{\mathfrak{vir}}}}+tN_n)^{-1}=\sum_{j=0}^{\infty}(-t)^jN_n^j.
\]
Applying Jacobi's formula \cite[Section~2.1.1, (46)]{PP} to ${\id_{W_n^{\mathfrak{vir}}}}+tN_n$
and using \eqref{eq:traces} in $\C[[t]]$, we obtain
\[
 \frac{p'(t)}{p(t)}
 =\Tr\bigl(({\id_{W_n^{\mathfrak{vir}}}}+tN_n)^{-1}N_n\bigr)
 =\sum_{j=0}^{\infty}(-t)^j\Tr(N_n^{j+1})=0.
\]
In characteristic zero, $p'(t)=0$ and $p(0)=1$ imply $p(t)=1$.
Substituting this into the determinant identity of step (1), we obtain
\[
 \det(H_n^{\mathfrak{vir}}+tB_n^{\Delta})=\det H_n^{\mathfrak{vir}}\ne0.
\]
This is a polynomial identity in $t$, and therefore holds for every $t\in\C$.
\end{proof}

\Needspace{3\baselineskip}
\subsection{PBW symmetrization and simplicity}\label{subsec:pbw-simplicity}
In this subsection we fix a commutative Frobenius algebra $(\mathcal A,f)$, and keep the spaces $W_n^{\mathfrak{vir}}$ and the forms $H_n^{\mathfrak{vir}}$ of the previous subsection.
By means of PBW symmetrization, we pull back the normalized invariant bilinear form on $V_{\vir(\mathcal A,f)}(1,0)$
to $W_n^{\mathfrak{vir}}$, and show that its difference from $H_n^{\mathfrak{vir}}$ satisfies \eqref{eq:partitions}.
Finally, in the case where $\mathcal A$ is local and $\mathcal A\not\cong\C$, we apply Lemma~\ref{lem:detcriterion} and obtain simplicity.

We first record the existence and uniqueness of this bilinear form, together with the adjoint relations needed to compute its values.

\begin{lemma}\label{lem:frob-invariant-form}
The vertex operator algebra $V_{\vir(\mathcal A,f)}(1,0)$ has a unique normalized invariant bilinear form,
which we denote by $B_{V_{\vir(\mathcal A,f)}(1,0)}$.
For $a\in\mathcal A$, $m\in\Z$, and $u,v\in V_{\vir(\mathcal A,f)}(1,0)$,
\[
 B_{V_{\vir(\mathcal A,f)}(1,0)}(L_m(a)u,v)
 =B_{V_{\vir(\mathcal A,f)}(1,0)}(u,L_{-m}(a)v).
\]
\end{lemma}
\begin{proof}
By Proposition~\ref{prop:frobuniversal}, this vertex operator algebra is of CFT type, and its conformal weight one subspace is zero.
Hence Proposition~\ref{prop:form}~(3) gives the existence and uniqueness of the bilinear form in question.
Moreover, $v_a\in V_{\vir(\mathcal A,f)}(1,0)_2$, so $L(1)v_a=0$, that is, $v_a$ is quasi-primary.
Since $(v_a)_{m+1}=L_m(a)$, \eqref{eq:adjoint} gives the stated identity.
\end{proof}

For $n\in\N$, in the notation of Remark~\ref{rem:weight-gram},
we denote by $\bigl(B_{V_{\vir(\mathcal A,f)}(1,0)}\bigr)_n$ the restriction of this form to
$V_{\vir(\mathcal A,f)}(1,0)_n\times V_{\vir(\mathcal A,f)}(1,0)_n$.
For a complex vector space $E$, denote its symmetric algebra by
\[
 \Sym(E)=\bigoplus_{r\in\N}\Sym^rE.
\]
For $m\in\Z_{\ge2}$, we keep the space $\mathcal A_m$ of Definition~\ref{def:weight-space} and its weight $m$.
By \cite[Chapter~III, Section~6, Proposition~9]{B},
we have a natural vector space isomorphism
\begin{equation}\label{eq:symmetric-algebra-decomposition}
 \Sym\!\left(\bigoplus_{m\ge2}\mathcal A_m\right)
 \cong\bigotimes_{m\ge2}\Sym(\mathcal A_m),
\end{equation}
where each pure tensor on the right-hand side has only finitely many factors different from $1$.
For each $n\in\N$, this vector space isomorphism maps the part of total weight $n$ on the left-hand side onto $W_n^{\mathfrak{vir}}$.
By means of PBW symmetrization, we now define a map from $W_n^{\mathfrak{vir}}$ to $V_{\vir(\mathcal A,f)}(1,0)_n$, and use it to transfer the invariant bilinear form to $W_n^{\mathfrak{vir}}$.

\begin{definition}\label{def:frob-pbw-pairing}
Fix $n\in\N$. Let $\ell\in\N$, $m_j\in\Z_{\ge2}$, and $x_j\in\mathcal A_{m_j}$ for $1\le j\le\ell$,
with $\sum_{j=1}^{\ell}m_j=n$. Via the vector space isomorphism $\mathcal A_{m_j}\simeq\mathcal A$,
denote the vector corresponding to $x_j$ by $a_j\in\mathcal A$, $1\le j\le\ell$. Let
\[
 r_m=\bigl|\{j:1\le j\le\ell,\ m_j=m\}\bigr|,\qquad m\in\Z_{\ge2}.
\]
Then $(r_m)_{m\ge2}$ is occupancy data in the sense of Definition~\ref{def:weight-space}.
Via the vector space isomorphism \eqref{eq:symmetric-algebra-decomposition}, we again denote by $x_1\cdots x_\ell$ the vector in $W_n^{\mathfrak{vir}}$ corresponding to the symmetric product $x_1\cdots x_\ell$. Explicitly,
\[
 x_1\cdots x_\ell
 :=\bigotimes_{m\ge2}\left(\prod_{\substack{1\le j\le\ell\\m_j=m}}x_j\right)\in W_n^{\mathfrak{vir}},
\]
where each parenthesized expression is a symmetric product in $\Sym^{r_m}\mathcal A_m$, and the factors with $r_m=0$ are omitted as in Definition~\ref{def:weight-space}.
The right-hand side groups the inputs according to the spaces $\mathcal A_m$ to which they belong, and the resulting vector does not depend on the order in which the inputs are written.
Define the linear map
\[
\begin{aligned}
 \Phi_n^{\mathfrak{vir}}:W_n^{\mathfrak{vir}}&\longrightarrow V_{\vir(\mathcal A,f)}(1,0)_n,\\
 x_1\cdots x_\ell&\longmapsto\frac1{\ell!}\sum_{\sigma\in S_\ell}
 L_{-m_{\sigma(1)}}(a_{\sigma(1)})\cdots L_{-m_{\sigma(\ell)}}(a_{\sigma(\ell)})\one.
\end{aligned}
\]
The average is taken over all $\ell$ inputs; a permutation only changes the order of the operators, and for each $m\in\Z_{\ge2}$, every term contains exactly $r_m$ operators of the form $L_{-m}(a)$ with $a\in\mathcal A$.
In particular, $\Phi_0^{\mathfrak{vir}}(1)=\one$, and $\Phi_1^{\mathfrak{vir}}$ is the unique map between two zero spaces.
Using $\Phi_n^{\mathfrak{vir}}$, we define a bilinear form on $W_n^{\mathfrak{vir}}$ by
\[
 G_n^{\mathfrak{vir}}(u,v)=B_{V_{\vir(\mathcal A,f)}(1,0)}\bigl(\Phi_n^{\mathfrak{vir}}(u),\Phi_n^{\mathfrak{vir}}(v)\bigr),\qquad u,v\in W_n^{\mathfrak{vir}}.
\]
\end{definition}

The average defining $\Phi_n^{\mathfrak{vir}}$ is multilinear in the inputs and invariant under all permutations of the inputs, so $\Phi_n^{\mathfrak{vir}}$ is well defined by \eqref{eq:symmetric-algebra-decomposition}.
For $m\in\Z_{\ge2}$ and $a\in\mathcal A$, assigning to the copy of $a$ in $\mathcal A_m$ the operator $L_{-m}(a)$, we obtain a vector space isomorphism $\bigoplus_{m\ge2}\mathcal A_m\cong\mathfrak n_-$.
Combining the PBW symmetrization of \cite[Corollary~13.7]{Et} with \eqref{eq:frob-pbw}, and noting that all these maps preserve total weight, we see that each $\Phi_n^{\mathfrak{vir}}$ is a vector space isomorphism.
Next, using \eqref{eq:frobbracket} to expand the values of $G_n^{\mathfrak{vir}}$, we show that $G_n^{\mathfrak{vir}}-H_n^{\mathfrak{vir}}$ satisfies the expansion condition \eqref{eq:partitions}.

\begin{lemma}\label{lem:partition}
For each $n\in\N$, $G_n^{\mathfrak{vir}}-H_n^{\mathfrak{vir}}$ satisfies the expansion condition \eqref{eq:partitions} required in Lemma~\ref{lem:detcriterion}.
\end{lemma}
\begin{proof}
For $n=0$, the normalization and Definition~\ref{def:auxiliary-pairing} give
$G_0^{\mathfrak{vir}}(1,1)=H_0^{\mathfrak{vir}}(1,1)=1$, so $G_0^{\mathfrak{vir}}-H_0^{\mathfrak{vir}}=0$. For $n=1$, we have $W_1^{\mathfrak{vir}}=0$, and the assertion also holds.
Assume from now on that $n\ge2$. We proceed in four steps.

\smallskip\noindent\emph{(1) Values of the bilinear form and vacuum coefficients.}
Let $r,t\in\Z_{\ge1}$ and $m_i,n_j\in\Z_{\ge2}$ for $1\le i\le r$ and $1\le j\le t$, with
$\sum_{i=1}^r m_i=\sum_{j=1}^t n_j=n$.
Let $\xi_1,\ldots,\xi_{r+t}\in\mathcal A$, and set, in $U(\mathfrak n_-)$,
\[
\begin{aligned}
 P&=L_{-m_1}(\xi_1)\cdots L_{-m_r}(\xi_r),\\
 Q&=L_{-n_1}(\xi_{r+1})\cdots L_{-n_t}(\xi_{r+t}).
\end{aligned}
\]
Let $\mathcal J_{r+t}=\{1,\ldots,r+t\}$, and use $j\in \mathcal J_{r+t}$ to label the position at which the initial input $\xi_j$ occurs.
Even if two input vectors are equal, their position labels are different.
For a vector $w\in V_{\vir(\mathcal A,f)}(1,0)$, the coefficient of $\one$ in its conformal weight $0$ component is called its \emph{vacuum coefficient}.
By the normalization and the orthogonality of distinct conformal weights in Proposition~\ref{prop:form}~(1),
$B_{V_{\vir(\mathcal A,f)}(1,0)}(\one,w)$ equals the vacuum coefficient of $w$.
Hence, by Lemma~\ref{lem:frob-invariant-form}, $B_{V_{\vir(\mathcal A,f)}(1,0)}(P\one,Q\one)$ equals the vacuum coefficient of the vector
\begin{equation}\label{eq:vacuum-coefficient-product}
 L_{m_r}(\xi_r)\cdots L_{m_1}(\xi_1)
 L_{-n_1}(\xi_{r+1})\cdots L_{-n_t}(\xi_{r+t})\one.
\end{equation}

\smallskip\noindent\emph{(2) Commutation relations and the partition expansion.}
We now apply the commutation relations repeatedly to the product of operators in \eqref{eq:vacuum-coefficient-product}.
Each product term in the resulting sum, together with its scalar coefficient, is called an \emph{expansion term}.
In each expansion term, we select, from left to right, the first pair of adjacent operator factors in which the integer subscript of the left factor is greater than that of the right one,
and use \eqref{eq:frobbracket} to write their product as the sum of the product in the reversed order and the commutator.
We repeat this operation on the resulting terms until, in every term, the subscripts are in nondecreasing order from left to right.
At each stage, the operator factors not yet eliminated are of the form
\[
 L_p\!\left(\prod_{j\in\mathcal I}\xi_j\right),
 \qquad p\in\Z,\quad\varnothing\ne\mathcal I\subseteq \mathcal J_{r+t},
\]
where the product is taken in $\mathcal A$. For $p,q\in\Z$ and disjoint nonempty subsets $\mathcal I_1,\mathcal I_2\subseteq \mathcal J_{r+t}$,
\[
 \left[L_p\!\left(\prod_{j\in\mathcal I_1}\xi_j\right),
          L_q\!\left(\prod_{j\in\mathcal I_2}\xi_j\right)\right]
 =(p-q)L_{p+q}\!\left(\prod_{j\in\mathcal I_1\cup\mathcal I_2}\xi_j\right)
 +\frac{p^3-p}{12}\delta_{p+q,0}
       f\!\left(\prod_{j\in\mathcal I_1\cup\mathcal I_2}\xi_j\right)K.
\]
The first term on the right-hand side is called the \emph{noncentral term}; selecting the second term, the \emph{central term}, is called a \emph{central contraction}.
The noncentral term merges the two sets of input positions and keeps one operator factor.
Since $K$ acts as the identity, a central contraction eliminates the two factors and produces
$\Theta_{|\mathcal I_1\cup\mathcal I_2|}((\xi_j)_{j\in\mathcal I_1\cup\mathcal I_2})$
together with a numerical coefficient depending only on $p,q$.
In each expansion term, the sets of input positions corresponding to the operator factors not yet eliminated and to the $\Theta$ factors already produced
are pairwise disjoint, and their union is always $\mathcal J_{r+t}$.

A pair of positions such that the factor in the left position has a larger integer subscript than the factor in the right position is called an \emph{inversion}.
In each commutator term the number of operator factors strictly decreases; in the term with the two factors interchanged, where the number of factors is unchanged, the number of inversions strictly decreases.
Therefore the above reordering terminates after finitely many steps, and the resulting expansion has only finitely many terms.
After the reordering is complete, if some factor has subscript at least $-1$, then the rightmost factor also has subscript at least $-1$, so it annihilates $\one$;
if only a nonempty product of factors with subscripts at most $-2$ remains, then applying it to $\one$ gives a vector of positive conformal weight.
Consequently, the terms that contribute nontrivially to the vacuum coefficient contain no operator factors; only a numerical coefficient remains, multiplied by a product of scalar factors
$\Theta_k(\xi_{j_1},\ldots,\xi_{j_k})=f(\xi_{j_1}\cdots\xi_{j_k})$,
where $k$ is the number of inputs contained in the factor, $2\le k\le r+t$, and $j_1,\ldots,j_k$ are pairwise distinct original input labels in $\mathcal J_{r+t}$.
Each such factor determines a set of input positions, denoted by $\mathcal I=\{j_1,\ldots,j_k\}$.
In each such expansion term, these sets are pairwise disjoint and their union is $\mathcal J_{r+t}$, so they form a partition of this set, which we denote by $\pi$.
Each block $\mathcal I\in\pi$ arises from a central contraction and hence contains at least two positions.
Thus the vacuum coefficient of the vector in \eqref{eq:vacuum-coefficient-product} can be written as a finite linear combination of products
\[
 \prod_{\mathcal I\in\pi}\Theta_{|\mathcal I|}\bigl((\xi_j)_{j\in\mathcal I}\bigr),
\]
where $\pi$ runs over the partitions of $\mathcal J_{r+t}$ all of whose blocks have size at least $2$, and the corresponding coefficients are allowed to be zero.
The coefficients of this linear combination depend only on the initial integer subscripts,
and not on the particular values of the input vectors $\xi_j$.

\smallskip\noindent\emph{(3) The contribution of the terms all of whose blocks have size $2$.}
If a noncentral term was selected at some stage in producing an expansion term, then at least two initial positions have been merged into the same operator factor,
and from then on these positions always remain together in some factor.
If this term contributes nontrivially to the vacuum coefficient, these positions must eventually form a central term together with another nonempty set of positions.
Hence the corresponding block contains at least three positions.
Therefore the terms all of whose blocks have size $2$ are precisely those obtained by central contractions of the initial operators in pairs until all of them are eliminated.
In this process, the order of the operators may be interchanged, but no noncentral term is selected.

We keep the notation $\alpha_m=\frac{m^3-m}{12}$, $m\in\Z_{\ge2}$, of the previous subsection. For $a,b\in\mathcal A$ and $m,k\in\Z_{\ge2}$,
the central term of $[L_m(a),L_{-k}(b)]$ is $\delta_{m,k}\alpha_m\eta(a,b)K$, and the corresponding scalar coefficient is
\[
 \delta_{m,k}\alpha_m\eta(a,b).
\]
Hence every nonzero pairing must join an initial factor $L_m(a)$ with an initial factor $L_{-m}(b)$.
When the subscripts of two factors are both positive or both negative, the central term vanishes, so in this part of the computation we only need to interchange their order.
We may therefore first move the rightmost positive-mode factor, one step at a time, past the negative-mode factors to its right.
After selecting its central term with one of these factors, we repeat the process with the remaining factors;
if no central term is ever selected, this positive-mode factor eventually acts on $\one$ and the resulting term is zero.
Therefore each of the above pairings is counted exactly once, and it suffices to sum directly over all pairings.

Let $r_m=|\{i:1\le i\le r,\ m_i=m\}|$ and
$r'_m=|\{j:1\le j\le t,\ n_j=m\}|$ for $m\ge2$.
These two sequences are, respectively, the occupancy data corresponding to the operator expressions of $P\one$ on the left
and $Q\one$ on the right in $B_{V_{\vir(\mathcal A,f)}(1,0)}(P\one,Q\one)$ of step (1). If the occupancy data of the two sides differ, then there is no pairing that eliminates all initial operators in pairs by nonzero central contractions, so the total contribution of the terms all of whose blocks have size $2$ is zero.
If the occupancy data of the two sides coincide, fix $m\in\Z_{\ge2}$, and denote the inputs $\xi_j$ on the left with $m_j=m$, $1\le j\le r$,
and the inputs $\xi_{r+j}$ on the right with $n_j=m$, $1\le j\le t$,
in increasing order of $j$ by $\xi_{m,i}$ and $\xi'_{m,i}\in\mathcal A$, $1\le i\le r_m$, respectively.
For each permutation $\sigma\in S_{r_m}$, we pair $\xi_{m,i}$ with $\xi'_{m,\sigma(i)}$ for $1\le i\le r_m$.
As $\sigma$ ranges over $S_{r_m}$, this runs through each bijective pairing of the input positions on the two sides exactly once, and the corresponding total contribution is
\[
 \alpha_m^{r_m}\sum_{\sigma\in S_{r_m}}
 \prod_{i=1}^{r_m}\eta(\xi_{m,i},\xi'_{m,\sigma(i)}).
\]
Taking the product of these contributions over all $m\ge2$ gives \eqref{eq:osc}.
Together with the case of different occupancy data on the two sides, this shows that the total contribution of the terms all of whose blocks have size $2$ to the vacuum coefficient of the vector in \eqref{eq:vacuum-coefficient-product}
is exactly the value of $H_n^{\mathfrak{vir}}$ on the corresponding two symmetric monomials.
This sum over pairings does not depend on the order of the operator factors on either side.

\smallskip\noindent\emph{(4) PBW symmetrization and \eqref{eq:partitions}.}
Symmetrizing the operator factors in $P$ and in $Q$ as in Definition~\ref{def:frob-pbw-pairing}, we obtain
\[
\frac{1}{r!\,t!}\!\sum_{\mathclap{\substack{\sigma\in S_r\\\tau\in S_t}}}
B_{V_{\vir(\mathcal A,f)}(1,0)}\Bigl(
L_{-m_{\sigma(1)}}(\xi_{\sigma(1)})\cdots L_{-m_{\sigma(r)}}(\xi_{\sigma(r)})\one,
L_{-n_{\tau(1)}}(\xi_{r+\tau(1)})\cdots L_{-n_{\tau(t)}}(\xi_{r+\tau(t)})\one\Bigr).
\]
By bilinearity, this is precisely the value of $G_n^{\mathfrak{vir}}$ on the corresponding two symmetric monomials.
In the expansion after each permutation, we always label the input $\xi_j$ by its original label $j\in\mathcal J_{r+t}$.
By step (2), the set $\mathcal I$ of input positions of each scalar factor $\Theta_{|\mathcal I|}((\xi_j)_{j\in\mathcal I})$
is a block of the corresponding partition $\pi$.
Keeping the original labels does not change the size of any block.
By step (3), the total contribution of the terms all of whose blocks have size $2$ is the same for every permutation, so after symmetrization it is still the value of $H_n^{\mathfrak{vir}}$ on the corresponding two symmetric monomials.
Therefore, subtracting $H_n^{\mathfrak{vir}}$ from $G_n^{\mathfrak{vir}}$ removes exactly the total contribution of the terms all of whose blocks have size $2$.
For each remaining term, every block $\mathcal I$ of the corresponding partition $\pi$ contains at least two positions, and at least one block contains at least three positions; this is exactly the summation condition in \eqref{eq:partitions}.
Combining the coefficients corresponding to the same partition and denoting the result by $\mu_\pi\in\C$, we obtain the expansion \eqref{eq:partitions} of $G_n^{\mathfrak{vir}}-H_n^{\mathfrak{vir}}$.
Here $\mathcal J_{r+t}=\{1,\ldots,r+t\}$ corresponds to the set of input positions in that formula according to the order of the inputs on the two sides.
None of the above numerical coefficients depends on the particular values of the input vectors, and hence neither do the combined coefficients $\mu_\pi$.
\end{proof}

Combining the expansion of Lemma~\ref{lem:partition} with the determinant criterion of Lemma~\ref{lem:detcriterion},
we obtain the following simplicity result.

\begin{theorem}\label{thm:frobsimple}
If $(\mathcal A,f)$ is a local commutative Frobenius algebra
with $\mathcal A\not\cong\C$, then $V_{\vir(\mathcal A,f)}(1,0)$ is simple.
\end{theorem}
\begin{proof}
By Lemmas~\ref{lem:detcriterion} and~\ref{lem:partition}, we have
$\det G_n^{\mathfrak{vir}}=\det H_n^{\mathfrak{vir}}\ne0$ for every $n\in\N$.
By the definition of $G_n^{\mathfrak{vir}}$ in Definition~\ref{def:frob-pbw-pairing} and the fact that $\Phi_n^{\mathfrak{vir}}$ is a vector space isomorphism,
$\bigl(B_{V_{\vir(\mathcal A,f)}(1,0)}\bigr)_n$ is nondegenerate.
Hence, by Remark~\ref{rem:weight-gram},
\[
 \rad B_{V_{\vir(\mathcal A,f)}(1,0)}
 =\bigoplus_{n=0}^{\infty}
 \rad\bigl(B_{V_{\vir(\mathcal A,f)}(1,0)}\bigr)_n=0.
\]
By Proposition~\ref{prop:form}~(3), $V_{\vir(\mathcal A,f)}(1,0)$ is simple.
\end{proof}

Theorem~\ref{thm:frobsimple} treats the case of a nonsemisimple local commutative Frobenius algebra.
The semisimple local case was discussed at the end of Section~\ref{subsec:frob-voa}: there $\mathcal A=\C1_{\mathcal A}$, and, with $c=f(1_{\mathcal A})$,
\eqref{eq:vir-special-universal} and \eqref{eq:vir-special} show that the corresponding vertex operator algebra and its unique simple quotient
are isomorphic to $V_{\mathfrak{vir}}(c,0)$ and $L_{\mathfrak{vir}}(c,0)$, respectively; for their properties, see Proposition~\ref{prop:vir}.

\section{Twisted modules and the classification in conformal weight two}\label{sec:twists}
In this section we prove Theorem~\ref{thm:two}. In view of the tensor decomposition of Proposition~\ref{prop:frobfactors},
the key point of the proof is the following: if there are only finitely many isomorphism classes of irreducible ordinary $g$-twisted modules,
then every tensor factor in Proposition~\ref{prop:frobfactors} is a Virasoro minimal model.

Section~\ref{subsec:virasoro-permutation} first determines the automorphisms of Virasoro tensor products, and recalls from \cite{BDM} the
permutation-twisted modules and tensor product constructions that we need.
Section~\ref{subsec:twisted-lie-modules} uses \cite[Theorem~3.6]{CDW}
to convert restricted modules for the corresponding Lie algebra into weak twisted modules.
Section~\ref{subsec:highest-weight-twisted} then constructs, by highest-weight induction, irreducible ordinary twisted modules
with different lowest conformal weights.

Section~\ref{subsec:proof-weight-two} combines Theorem~\ref{thm:frobsimple} with Proposition~\ref{prop:vir}
and applies the above construction to the tensor product of the factors that are not Virasoro minimal models, taking tensor products with twisted modules for the remaining factors.
If such a factor exists, we obtain infinitely many pairwise non-isomorphic irreducible ordinary $g$-twisted modules, contradicting the finiteness assumption.
Finally, using the permutation rationality of \cite{BDM} and the relevant results of Section~\ref{sec:prelim},
we complete the equivalences in Theorem~\ref{thm:two} and obtain $C_2$-cofiniteness and regularity.

\Needspace{3\baselineskip}
\subsection{Virasoro tensor products and permutation-twisted modules}\label{subsec:virasoro-permutation}
By \eqref{eq:vir-special}, the simple vertex operator algebras corresponding to one-dimensional local summands are simple Virasoro vertex operator algebras.
In this subsection we first show that every automorphism of a Virasoro tensor product is given by a permutation of factors with the same central charge (Lemma~\ref{lem:virasoro-permutations});
this refines the statement in Proposition~\ref{prop:frobfactors} that automorphisms permute the factors.
We then recall the relevant results of \cite{BDM}, in preparation for the construction of irreducible ordinary twisted modules and the proof of $g$-rationality in Section~\ref{subsec:proof-weight-two}.

In what follows, we omit the trivial factors $L_{\mathfrak{vir}}(0,0)=\C\one$;
the only automorphism of the empty tensor product $\C\one$ is the identity.

\begin{lemma}\label{lem:virasoro-permutations}
Let $s\in\Z_{\ge1}$ and $c_i\in\C\setminus\{0\}$ for $1\le i\le s$, and let
\[
 V=\bigotimes_{i=1}^sL_{\mathfrak{vir}}(c_i,0).
\]
For indices $i,j\in\{1,\ldots,s\}$ with $c_i=c_j$,
we regard the $i$-th and $j$-th factors as two copies of the same $L_{\mathfrak{vir}}(c_i,0)$,
identified via the unique vertex operator algebra isomorphism
sending the conformal vector of the $i$-th factor to that of the $j$-th factor.
Then for every $g\in\Aut(V)$, there exists a unique $\sigma\in S_s$ such that
\[
 c_{\sigma(i)}=c_i,\qquad 1\le i\le s,
\]
and
\begin{equation}\label{eq:virasoro-factor-permutation}
 g(u_1\otimes\cdots\otimes u_s)
 =u_{\sigma^{-1}(1)}\otimes\cdots\otimes u_{\sigma^{-1}(s)},
 \quad u_i\in L_{\mathfrak{vir}}(c_i,0),\quad 1\le i\le s.
\end{equation}
\end{lemma}
\begin{proof}
Denote the conformal vector of the $i$-th factor by $\omega_i$, and set
$e_i=\one\otimes\cdots\otimes\omega_i\otimes\cdots\otimes\one\in V_2$ for $1\le i\le s$.
On $V_2$ we use the multiplication of Proposition~\ref{prop:wt2}, given by $ab=\tfrac12a_1b$ for $a,b\in V_2$.
By \eqref{eq:tensor-vertex}--\eqref{eq:tensor-grading} and the Virasoro relations,
\[
 V_2=\bigoplus_{i=1}^s\C e_i,\qquad
 \omega=\sum_{i=1}^s e_i,\qquad
 e_ie_j=\delta_{ij}e_i,\quad 1\le i,j\le s.
\]
Hence each $\C e_i$ is a one-dimensional local summand, and by Proposition~\ref{prop:local},
$e_1,\ldots,e_s$ are precisely the primitive idempotents of $V_2$.

Let $g\in\Aut(V)$. Since $g(V_2)=V_2$ and $g(ab)=g(a)g(b)$,
$g$ preserves the set of primitive idempotents, so there exists a unique $\sigma\in S_s$ such that
$g(e_i)=e_{\sigma(i)}$ for $1\le i\le s$. Furthermore,
\[
 \frac{c_i}{2}\one
 =g\bigl((e_i)_3e_i\bigr)
 =(e_{\sigma(i)})_3e_{\sigma(i)}
 =\frac{c_{\sigma(i)}}{2}\one,
 \qquad 1\le i\le s.
\]
Hence $c_i=c_{\sigma(i)}$.

By \cite[Section~6.1]{LL}, each Virasoro factor $L_{\mathfrak{vir}}(c_i,0)$ is generated by its conformal vector $\omega_i$,
so the restriction of $g$ to the $i$-th factor is the unique isomorphism onto the $\sigma(i)$-th factor
sending $\omega_i$ to $\omega_{\sigma(i)}$.
Once these factors are identified via the isomorphisms in the statement, $g$ moves the vector in the $i$-th position to the $\sigma(i)$-th position,
so the vector in the $j$-th position is $u_{\sigma^{-1}(j)}$, $1\le j\le s$. This gives \eqref{eq:virasoro-factor-permutation}.
\end{proof}

\begin{remark}\label{rem:virasoro-permutation-example}
Lemma~\ref{lem:virasoro-permutations} shows that the automorphisms of $V$ can only permute factors with the same central charge.
For example, let $c,d\in\C\setminus\{0\}$ with $c\ne d$, and let
\[
 V=L_{\mathfrak{vir}}(c,0)^{\otimes3}\otimes L_{\mathfrak{vir}}(d,0)^{\otimes2}.
\]
The first three factors and the last two factors can be permuted independently. For instance, take
\[
 g(u_1\otimes u_2\otimes u_3\otimes v_1\otimes v_2)
 =u_3\otimes u_1\otimes u_2\otimes v_2\otimes v_1,
\]
where $u_1,u_2,u_3\in L_{\mathfrak{vir}}(c,0)$ and $v_1,v_2\in L_{\mathfrak{vir}}(d,0)$.
By \eqref{eq:tensor-vertex}, this map is an automorphism of $V$.
\end{remark}

We now record the permutation-twisted module constructions of \cite{BDM},
together with the tensor product properties of general twisted modules.

\begin{proposition}\label{prop:permutation}
\begin{enumerate}[label=(\arabic*)]
\item Let $U$ be a vertex operator algebra, let $k\in\Z_{\ge1}$,
and let $\sigma$ be a $k$-cycle permutation of $U^{\otimes k}$, for example
\[
 \sigma(u_1\otimes\cdots\otimes u_k)
 =u_k\otimes u_1\otimes\cdots\otimes u_{k-1},
 \qquad u_1,\ldots,u_k\in U.
\]
Every irreducible ordinary $U$-module gives rise to an irreducible ordinary
$\sigma$-twisted $U^{\otimes k}$-module.

\item Let $s\in\Z_{\ge1}$, and for $1\le i\le s$ let $U_i$ be a vertex operator algebra,
$g_i\in\Aut(U_i)$ an automorphism of finite order,
and $M_i$ an irreducible ordinary $g_i$-twisted $U_i$-module.
Then $M=\bigotimes_{i=1}^sM_i$ is an irreducible ordinary
$\bigotimes_{i=1}^s g_i$-twisted $\bigotimes_{i=1}^s U_i$-module.

\item If $U$ is a rational vertex operator algebra, then for every $k\in\Z_{\ge1}$
and every permutation automorphism $\tau\in S_k$, $U^{\otimes k}$ is $\tau$-rational.

\item Let $s\in\Z_{\ge1}$, and for $1\le i\le s$ let $U_i$ be a vertex operator algebra and
$g_i\in\Aut(U_i)$ an automorphism of finite order.
If each $U_i$ is $g_i$-rational, then $\bigotimes_{i=1}^sU_i$ is
$\bigotimes_{i=1}^sg_i$-rational.
\end{enumerate}
\end{proposition}
\begin{proof}
Assertion (1) follows from \cite[Theorem~4.4, Corollary~5.1]{BDM},
and (3) and (4) follow from \cite[Theorem~6.4~(3), Proposition~6.2]{BDM}, respectively.
We now prove (2). The weak twisted module structure on a general tensor product is given at the beginning of \cite[Section~6]{BDM},
and its irreducibility is guaranteed by \cite[Lemma~6.1]{BDM}.

It remains to show that this tensor product is ordinary.
The verification of the grading and of finite dimensionality follows the idea of \cite[Section~4.6, Proposition~4.6.2]{LL},
with the integral grading replaced by the corresponding fractional grading.
Denote the order of $g_i$ by $T_i\in\Z_{\ge1}$, $1\le i\le s$,
and let $T$ be a positive common multiple of these orders.
Each $M_i$ is an irreducible ordinary $g_i$-twisted $U_i$-module.
By irreducibility and the definition of an ordinary module (Definition~\ref{def:admissible}), we can choose a lowest conformal weight $h_i\in\C$ such that
\[
 M_i=\bigoplus_{n\in\N}(M_i)_{h_i+\frac nT},\qquad
 (M_i)_{h_i}\ne0,\qquad 1\le i\le s.
\]
On $M$ we have
\[
 L_M(0)=\sum_{i=1}^s
 \id_{M_1}\otimes\cdots\otimes L_{M_i}(0)\otimes\cdots\otimes\id_{M_s}.
\]
Hence $L_M(0)$ is diagonalizable, and for each $n\in\N$,
\[
 M_{\sum_{i=1}^s h_i+\frac{n}{T}}
 =\bigoplus_{\substack{n_1,\ldots,n_s\in\N\\n_1+\cdots+n_s=n}}
   \bigotimes_{i=1}^s(M_i)_{h_i+\frac{n_i}{T}}.
\]
The right-hand side has only finitely many summands, each of which is finite-dimensional,
and $\operatorname{Spec}L_M(0)\subset\sum_{i=1}^s h_i+\frac1T\N$, so $M$ is ordinary.
\end{proof}

\begin{remark}\label{rem:permutation-cycle-construction}
By \cite[Theorem~1.1]{Rotman}, every permutation is a product of disjoint cycles.
After reordering the tensor factors as in \cite[Section~6, paragraph preceding Theorem~6.4]{BDM},
one applies Proposition~\ref{prop:permutation}~(1) to each cycle and then Proposition~\ref{prop:permutation}~(2)
to obtain the corresponding irreducible ordinary twisted module.
This construction will be used in the proof of Theorem~\ref{thm:two} in Section~\ref{subsec:proof-weight-two}.
\end{remark}

\Needspace{3\baselineskip}
\subsection{Twisted Lie algebras and weak twisted modules}\label{subsec:twisted-lie-modules}
This subsection prepares the twisted module construction needed for the highest-weight induction in Section~\ref{subsec:highest-weight-twisted}.
We first recall from \cite[Section~3]{CDW} Lie conformal algebras, the corresponding twisted Lie algebras,
and their restricted modules; then, combining these with the construction of Section~\ref{sec:frob}, we give the Lie algebra corresponding to $V_{\vir(A,f)}(1,0)$
and its commutation relations. Finally, using \cite[Theorem~3.6]{CDW}, we explain how weak twisted modules for this vertex operator algebra
are obtained from the corresponding restricted modules.

Lie conformal algebras are also called vertex Lie algebras; for their definition, see \cite{Kac,Primc}. We use the formal series formulation of \cite[Definition~3.1]{CDW}.
For a vertex algebra $(V,Y,\one)$, $Y_-(a,z)b=\sum_{j\ge0}(a_jb)z^{-j-1}$
is the singular part of $Y(a,z)b$; a Lie conformal algebra retains these singular products and the translation operator.

\begin{definition}\label{def:lca}
A \emph{Lie conformal algebra} is a vector space $\mathcal C$ equipped with an endomorphism
$\partial\in\End\mathcal C$ and a linear map
\[
\begin{aligned}
 Y_-(\,\cdot\,,z):\mathcal C&\longrightarrow
 \operatorname{Hom}(\mathcal C,z^{-1}\mathcal C[z^{-1}]),\\
 a&\longmapsto Y_-(a,z)=\sum_{j\ge0}a_jz^{-j-1},
\end{aligned}
\]
where $a_j\in\End\mathcal C$, such that for $a,b\in\mathcal C$,
\begin{align*}
 [\partial,Y_-(a,z)]&=\frac{d}{dz}Y_-(a,z)=Y_-(\partial a,z),\\
 Y_-(a,z)b&=\operatorname{Sing}\bigl(e^{z\partial}Y_-(b,-z)a\bigr),\\
 [Y_-(a,z_1),Y_-(b,z_2)]
 &=\operatorname{Sing}Y_-(Y_-(a,z_1-z_2)b,z_2).
\end{align*}
Here $\operatorname{Sing}$ takes the singular part in each of the formal variables,
that is, it keeps the terms with strictly negative powers of each corresponding variable.
Equivalently, for all $a,b\in\mathcal C$ and $m,n\in\Z_{\ge0}$,
\[
\begin{aligned}
 (\partial a)_n&=-n\,a_{n-1},\\
 a_n b&=\sum_{j\ge0}\frac{(-1)^{n+j+1}}{j!}\,\partial^j(b_{n+j}a),\\
 [a_m,b_n]&=\sum_{j=0}^{m}\binom mj(a_jb)_{m+n-j}.
\end{aligned}
\]
For $n=0$, the first identity is understood as $(\partial a)_0=0$.
These are precisely the three conditions of \cite[Definition~2.7, (C1)--(C3)]{Kac},
corresponding to translation, skew symmetry, and the commutator identity, respectively.
\end{definition}

We next recall the construction of the Lie algebra and the vertex algebra corresponding to a Lie conformal algebra;
see \cite[Lemma~3.2, (3.5)--(3.8) and the construction following them]{CDW}.

\begin{proposition}\label{prop:lca-vertex-construction}
Let $\mathcal C$ be a Lie conformal algebra.
\begin{enumerate}[label=(\arabic*)]
\item Let
\[
 \widehat{\mathcal C}
 =\bigl(\mathcal C\otimes\C[t,t^{-1}]\bigr)
 /\operatorname{im}\!\left(\partial\otimes1+1\otimes\frac{d}{dt}\right),
\]
where $t$ is a formal variable; denote by $a(m)$ the class of $a\otimes t^m$.
Then $\widehat{\mathcal C}$ is a Lie algebra under the bracket
\begin{equation}\label{eq:lcamodes}
 [a(m),b(n)]=\sum_{j\ge0}\binom mj(a_jb)(m+n-j),
 \qquad a,b\in\mathcal C,\quad m,n\in\Z.
\end{equation}
The relation imposed by the quotient is
\[
 (\partial a)(m)=-m\,a(m-1),\qquad a\in\mathcal C,\quad m\in\Z.
\]
\item Let
\[
 \widehat{\mathcal C}_+=\Span\{a(m):a\in\mathcal C,\ m\in\Z_{\ge0}\}.
\]
Then $\widehat{\mathcal C}_+$ is a Lie subalgebra. Take its one-dimensional trivial module $\C$, and write
\[
 V_{\widehat{\mathcal C}}=U(\widehat{\mathcal C})
 \otimes_{U(\widehat{\mathcal C}_+)}\C,\qquad \one=1\otimes1.
\]
This induced vacuum module carries a unique vertex algebra structure with vacuum vector $\one$ such that
\[
 Y(a(-1)\one,z)=\sum_{m\in\Z}a(m)z^{-m-1},\qquad a\in\mathcal C.
\]
The map $a\mapsto a(-1)\one$ is injective, and we identify $\mathcal C$ with its image.
By the PBW theorem,
\[
 V_{\widehat{\mathcal C}}
 =\Span\bigl\{a_1(-n_1)\cdots a_\ell(-n_\ell)\one:
       \ell\in\N,\ a_i\in\mathcal C,\ n_i\in\Z_{\ge1}\bigr\},
\]
where $1\le i\le\ell$, and $\ell=0$ corresponds to $\one$.
In particular, $\mathcal C$ generates $V_{\widehat{\mathcal C}}$.
\end{enumerate}
\end{proposition}

To introduce the corresponding twisted Lie algebra, we first recall the definition of an automorphism of a Lie conformal algebra;
see \cite[Definition~3.2]{CS} and \cite[Section~3, the definition of homomorphisms preceding Lemma~3.3]{CDW}.

\begin{definition}\label{def:lca-automorphism}
Let $(\mathcal C,\partial,Y_-)$ be a Lie conformal algebra.
A linear isomorphism $g:\mathcal C\to\mathcal C$ is called an \emph{automorphism} of $\mathcal C$ if
\[
 g\partial=\partial g,\qquad
 g\bigl(Y_-(a,z)b\bigr)=Y_-(g(a),z)g(b),\qquad a,b\in\mathcal C.
\]
Equivalently, the second condition can be written as
\[
 g(a_jb)=(g(a))_j\,g(b),\qquad a,b\in\mathcal C,\quad j\in\Z_{\ge0}.
\]
\end{definition}

From now on, let $g$ be an automorphism of $\mathcal C$ of finite order $T\in\Z_{\ge1}$.
As in Section~\ref{sec:prelim}, for $r\in\{0,\ldots,T-1\}$ we write
\[
 \mathcal C_g^{[\frac{r}{T}]}
 =\{a\in\mathcal C:g(a)=e^{\frac{2\pi\mathrm i r}{T}}a\}.
\]
Since $g^T=\id_{\mathcal C}$,
\[
 \mathcal C=\bigoplus_{r=0}^{T-1}\mathcal C_g^{[\frac{r}{T}]}.
\]
We use the construction and definition of \cite[(3.9), Definition~3.4]{CDW}.

\begin{definition}\label{def:twisted-lie-restricted}
The \emph{twisted Lie algebra} corresponding to $(\mathcal C,g)$ is
\[
 \widehat{\mathcal C}[g]
 =\left(\bigoplus_{r=0}^{T-1}
 \mathcal C_g^{[\frac{r}{T}]}\otimes t^{\frac{r}{T}}\C[t,t^{-1}]\right)
 /\operatorname{im}\!\left(\partial\otimes1+1\otimes\frac{d}{dt}\right).
\]
For $a\in\mathcal C_g^{[\frac{r}{T}]}$ and $p\in \frac{r}{T}+\Z$, we again denote by $a(p)$ the class of $a\otimes t^p$.
The translation relation and the bracket are
\[
 (\partial a)(p)=-p\,a(p-1),\qquad
 [a(p),b(q)]=\sum_{j\ge0}\binom pj(a_jb)(p+q-j),
\]
where $b\in\mathcal C_g^{[\frac{s}{T}]}$, $q\in \frac{s}{T}+\Z$, and $r,s\in\{0,\ldots,T-1\}$.
When $g=\id$, we have $\widehat{\mathcal C}[g]=\widehat{\mathcal C}$.

For a $\widehat{\mathcal C}[g]$-module $W$, write
\[
 a^g(z)=\sum_{p\in \frac{r}{T}+\Z}a(p)z^{-p-1},\qquad a\in\mathcal C_g^{[\frac{r}{T}]},
\]
and extend linearly according to the eigenspace decomposition. We say that $W$ is \emph{restricted} if
\[
 a^g(z)\in\operatorname{Hom}\bigl(W,W((z^{\frac{1}{T}}))\bigr),\qquad a\in\mathcal C.
\]
Equivalently, for each $a\in\mathcal C_g^{[\frac{r}{T}]}$ and $w\in W$, we have $a(p)w=0$ for all sufficiently large $p\in \frac{r}{T}+\Z$.
\end{definition}

By \cite[(3.11)]{CDW}, $g$ extends uniquely to an automorphism of $V_{\widehat{\mathcal C}}$, which we again denote by $g$.
We now quote \cite[Theorem~3.6]{CDW}, which gives the correspondence between restricted $\widehat{\mathcal C}[g]$-modules
and weak $g$-twisted $V_{\widehat{\mathcal C}}$-modules.

\begin{proposition}\label{prop:restricted-twisted-modules}
Let $\mathcal C$ be a Lie conformal algebra, and let $g$ be an automorphism of $\mathcal C$ of finite order.
Every restricted $\widehat{\mathcal C}[g]$-module $W$ carries a natural
weak $g$-twisted $V_{\widehat{\mathcal C}}$-module structure such that
\[
 Y_W(a,z)=a^g(z),\qquad a\in\mathcal C.
\]
Conversely, every weak $g$-twisted $V_{\widehat{\mathcal C}}$-module
becomes naturally a restricted $\widehat{\mathcal C}[g]$-module via the coefficients of the above vertex operators.
\end{proposition}

We now apply this result to the construction of Section~\ref{sec:frob}. Let $A$ be a finite-dimensional commutative associative
unital algebra with identity element $1_A$, and let $f:A\to\C$ be a linear functional.
Let $g$ be an algebra automorphism of $A$ of order $T\in\Z_{\ge1}$ preserving $f$. Let
\[
 \mathcal C=(\C[\partial]\otimes A)\oplus\C\mathbf k,
 \qquad \partial\mathbf k=0,
 \qquad v_a=1\otimes a\quad(a\in A).
\]
On $\C[\partial]\otimes A$, set
\[
 \partial(\partial^j\otimes a)=\partial^{j+1}\otimes a,
 \qquad a\in A,\quad j\in\Z_{\ge0},
\]
and extend linearly. For $a,b\in A$, define the singular products by
\begin{equation}\label{eq:frob-lca-products}
 Y_-(v_a,z)v_b=z^{-1}\partial v_{ab}
             +2z^{-2}v_{ab}+\frac{f(ab)}2z^{-4}\mathbf k,
 \qquad a,b\in A;
\end{equation}
all singular products involving $\mathbf k$ are zero, and the remaining products are obtained by extending via the translation identities
in Definition~\ref{def:lca}. Explicitly, for $a,b\in A$ and $i,j\in\Z_{\ge0}$, set
\[
 Y_-(\partial^iv_a,z)\partial^jv_b
 =\left(\frac{d}{dz}\right)^i
  \left(\partial-\frac{d}{dz}\right)^j\bigl(Y_-(v_a,z)v_b\bigr),
\]
where $\partial$ on the right-hand side acts on the coefficients of each term, and extend bilinearly in the two input vectors.
We again denote by $g$ the linear automorphism determined by
\[
 g(\partial^jv_a)=\partial^jv_{g(a)},\qquad
 g\mathbf k=\mathbf k,\qquad a\in A,\ j\in\N.
\]
The following lemma shows that this construction is compatible with the vertex operator algebra of Section~\ref{sec:frob},
and gives the twisted Lie algebra needed later.

\begin{lemma}\label{lem:frob-lca-construction}
In the above setting, the following hold.
\begin{enumerate}[label=(\arabic*)]
\item $\mathcal C$ is a Lie conformal algebra, and $g$ is an automorphism of $\mathcal C$.
\item Let $\mathcal J_{\mathbf k}$ be the vertex ideal of $V_{\widehat{\mathcal C}}$ generated by $\mathbf k-\one$.
Then $\mathcal J_{\mathbf k}$ is $g$-stable, and there is a vertex operator algebra isomorphism
\begin{equation}\label{eq:frob-central-quotient}
 V_{\widehat{\mathcal C}}/\mathcal J_{\mathbf k}
 \cong V_{\vir(A,f)}(1,0),\qquad
 v_a+\mathcal J_{\mathbf k}\longmapsto L_{-2}(a)\one,
\end{equation}
where the conformal vector on the left-hand side is $v_{1_A}+\mathcal J_{\mathbf k}$.
This isomorphism satisfies
\[
 \partial^jv_a+\mathcal J_{\mathbf k}
 \longmapsto L(-1)^jL_{-2}(a)\one,
 \qquad a\in A,\quad j\in\Z_{\ge0}.
\]
\item In what follows, $r,s\in\{0,\ldots,T-1\}$. Write $\vir_g(A,f)=\widehat{\mathcal C}[g]$, and set
\[
 L_p(a)=v_a(p+1),\qquad K=\mathbf k(-1),
 \qquad a\in A_g^{[\frac{r}{T}]},\quad p\in \frac{r}{T}+\Z.
\]
The underlying vector space of this Lie algebra can be identified with
\[
 \vir_g(A,f)=
 \left(\bigoplus_{r=0}^{T-1}A_g^{[\frac{r}{T}]}\otimes t^{\frac{r}{T}}\C[t,t^{-1}]\right)\oplus\C K,
\]
where $a\otimes t^p$ corresponds to $L_p(a)$, $K$ is central, and
\begin{equation}\label{eq:frobtwist}
 [L_p(a),L_q(b)]=(p-q)L_{p+q}(ab)
 +\frac{p^3-p}{12}\delta_{p+q,0}f(ab)K
\end{equation}
for $a\in A_g^{[\frac{r}{T}]}$, $b\in A_g^{[\frac{s}{T}]}$, $p\in \frac{r}{T}+\Z$, and $q\in \frac{s}{T}+\Z$.
\end{enumerate}
\end{lemma}
\begin{proof}
We first prove (1). We use \cite[Lemma~3.2]{CDW} to show that $\mathcal C$ is a Lie conformal algebra.
By the above extension formulas, the singular products satisfy the translation identities of Definition~\ref{def:lca}.
Hence it suffices to verify that the quotient vector space constructed as in Proposition~\ref{prop:lca-vertex-construction}~(1)
is a Lie algebra under the bracket \eqref{eq:lcamodes}.

Set
\[
 L_m(a)=v_a(m+1),\qquad K=\mathbf k(-1),
 \qquad a\in A,\quad m\in\Z.
\]
By the translation relation
\[
 (\partial u)(n)=-n\,u(n-1),
 \qquad u\in\mathcal C,\quad n\in\Z,
\]
the modes of $\partial^jv_a$ can be written in terms of the modes of $v_a$; moreover, $\partial\mathbf k=0$ and the translation relation
give $\mathbf k(n)=0$ for $n\ne-1$. Hence this quotient vector space can be identified with
\[
 (A\otimes\C[t,t^{-1}])\oplus\C K.
\]
Substituting \eqref{eq:frob-lca-products} into \eqref{eq:lcamodes}, we obtain precisely the bracket \eqref{eq:frobbracket}, with $K$ central.
The argument following Definition~\ref{def:froblie} shows that this bracket satisfies the Jacobi identity, so the above quotient is a Lie algebra.
By \cite[Lemma~3.2]{CDW}, $\mathcal C$ is a Lie conformal algebra.

We next show that $g$ is an automorphism of $\mathcal C$. Since $g(ab)=g(a)g(b)$ and $f\circ g=f$,
\eqref{eq:frob-lca-products} gives, for all $a,b\in A$,
\[
 g\bigl(Y_-(v_a,z)v_b\bigr)=Y_-(v_{g(a)},z)v_{g(b)}.
\]
Together with $g\partial=\partial g$, $g\mathbf k=\mathbf k$, and the above extension formulas, this shows that $g$ preserves all singular products.
Hence $g$ is an automorphism of $\mathcal C$ by Definition~\ref{def:lca-automorphism}.

For (2), note first that $g$ fixes $\mathbf k-\one$, so $\mathcal J_{\mathbf k}$ is $g$-stable.
Since $K=\mathbf k(-1)$, we have
\[
 \mathcal J_{\mathbf k}=(K-\id)V_{\widehat{\mathcal C}}.
\]

We now construct the vertex operator algebra isomorphism \eqref{eq:frob-central-quotient}.
The vacuum relations $v_a(n)\one=0$, $n\ge0$, in $V_{\widehat{\mathcal C}}$
correspond exactly to $L_m(a)\one=0$, $m\ge-1$.
By the definitions of the induced vacuum modules in Propositions~\ref{prop:lca-vertex-construction} and~\ref{prop:frobuniversal},
sending the vacuum vector to the vacuum vector gives a surjective $\vir(A,f)$-module homomorphism
\[
 V_{\widehat{\mathcal C}}\longrightarrow V_{\vir(A,f)}(1,0).
\]
Since $K$ acts as the identity on the latter, this map induces a surjection
\[
 V_{\widehat{\mathcal C}}/\mathcal J_{\mathbf k}\longrightarrow V_{\vir(A,f)}(1,0),
 \qquad v_a+\mathcal J_{\mathbf k}\longmapsto L_{-2}(a)\one.
\]

Fix an ordered basis of $A$, and fix a total order on the corresponding modes $L_{-n}(a)$, $n\ge2$.
By the PBW theorem, $V_{\widehat{\mathcal C}}$ has a basis
\[
 K^jL_{-n_1}(a_1)\cdots L_{-n_\ell}(a_\ell)\one,
 \qquad j,\ell\in\Z_{\ge0},\quad n_i\ge2,
\]
where each $a_i$ is taken from this basis and the modes $L_{-n_i}(a_i)$ are arranged according to the above total order.
In the quotient, $K$ acts as the identity, so the images of these monomials with $j=0$ form a basis of the quotient.
Under the above surjection they are mapped to the corresponding PBW basis vectors of Proposition~\ref{prop:frobuniversal}, so the map is bijective.

This map is a $\vir(A,f)$-module homomorphism, and hence commutes with the action of each $L_m(a)$.
Since it also preserves the vacuum vector, and $v_a+\mathcal J_{\mathbf k}$ and $L_{-2}(a)\one$ generate the vertex algebras on the two sides, respectively,
the map preserves the vertex algebra structure.
By Proposition~\ref{prop:frobuniversal}, it maps $v_{1_A}+\mathcal J_{\mathbf k}$ to the conformal vector $L_{-2}(1_A)\one$.
Therefore, taking $v_{1_A}+\mathcal J_{\mathbf k}$ as the conformal vector of the quotient,
we obtain the vertex operator algebra isomorphism \eqref{eq:frob-central-quotient}.

Moreover, the restriction of the translation operator of $V_{\widehat{\mathcal C}}$ to $\mathcal C$ is $\partial$;
see \cite[Section~3, the remark following (3.8)]{CDW}.
Since vertex operator algebra homomorphisms preserve translation operators, and the translation operator of the target space is $L(-1)$, we have
\[
 \partial^jv_a+\mathcal J_{\mathbf k}\longmapsto L(-1)^jL_{-2}(a)\one,
 \qquad a\in A,\quad j\in\Z_{\ge0}.
\]

Finally, the above isomorphism transfers the action of $g$ on the quotient to an automorphism of $V_{\vir(A,f)}(1,0)$,
again denoted by $g$, which satisfies
\[
 g\bigl(L_{-2}(a)\one\bigr)=L_{-2}(g(a))\one,
 \qquad a\in A.
\]

For (3), we first determine the underlying vector space of $\widehat{\mathcal C}[g]$.
As in (1), the translation relation in Definition~\ref{def:twisted-lie-restricted}
allows us to write the modes of $\partial^jv_a$ in terms of the modes of $v_a$.
Since $g\mathbf k=\mathbf k$, the mode indices of $\mathbf k$ are integers;
by $\partial\mathbf k=0$ and the translation relation,
\[
 0=(\partial\mathbf k)(n+1)=-(n+1)\mathbf k(n),
 \qquad n\in\Z,
\]
so $\mathbf k(n)=0$ for $n\ne-1$, and we write $K=\mathbf k(-1)$.
By the definition of $\mathcal C$, for each $r\in\{0,\ldots,T-1\}$ and $m\in\frac{r}{T}+\Z$,
the subspace spanned by the vectors $v_a(m)$, $a\in A_g^{[\frac{r}{T}]}$, is isomorphic to $A_g^{[\frac{r}{T}]}$.
Relabeling by
\[
 L_p(a)=v_a(p+1)
\]
and letting $L_p(a)$ correspond to $a\otimes t^p$, we obtain the vector space decomposition in the statement.
Since all singular products involving $\mathbf k$ are zero,
the bracket formula in Definition~\ref{def:twisted-lie-restricted} shows that $K$ is central.

We now compute the bracket. Let $a\in A_g^{[\frac{r}{T}]}$, $b\in A_g^{[\frac{s}{T}]}$,
$p\in\frac{r}{T}+\Z$, and $q\in\frac{s}{T}+\Z$.
Since $g$ preserves both the multiplication of $A$ and $f$,
\[
 ab\in A_g^{[\frac{r+s}{T}]},\qquad
 f(ab)=0\quad\text{if }r+s\not\equiv0\pmod T.
\]
Substituting \eqref{eq:frob-lca-products} into the bracket formula of Definition~\ref{def:twisted-lie-restricted}, we obtain
\begin{align*}
 [L_p(a),L_q(b)]
 &=(\partial v_{ab})(p+q+2)+2(p+1)v_{ab}(p+q+1)
   +\binom{p+1}{3}\frac{f(ab)}2\delta_{p+q,0}K\\
 &=(p-q)L_{p+q}(ab)+\frac{p^3-p}{12}\delta_{p+q,0}f(ab)K,
\end{align*}
where the second equality uses the translation relation
$(\partial v_{ab})(p+q+2)=-(p+q+2)v_{ab}(p+q+1)$
and $\frac12\binom{p+1}{3}=\frac{p^3-p}{12}$. This is \eqref{eq:frobtwist}.
\end{proof}

Let $W$ be any restricted $\vir_g(A,f)$-module on which $K$ acts as the identity.
By Proposition~\ref{prop:restricted-twisted-modules} and Lemma~\ref{lem:frob-lca-construction},
we have $Y_W(\mathbf k-\one,z)=K-\id_W=0$,
so the corresponding weak $g$-twisted $V_{\widehat{\mathcal C}}$-module structure
descends to a weak $g$-twisted $V_{\widehat{\mathcal C}}/\mathcal J_{\mathbf k}$-module structure.
Via the isomorphism \eqref{eq:frob-central-quotient},
we then obtain the corresponding weak $g$-twisted $V_{\vir(A,f)}(1,0)$-module structure;
relabeling the modes by $L_p(a)=v_a(p+1)$ gives the following corollary.

\begin{corollary}\label{cor:frob-twisted-modules}
Keep the notation $A,f,g,T$ of Lemma~\ref{lem:frob-lca-construction}.
Every restricted $\vir_g(A,f)$-module $W$ on which $K$ acts as the identity
carries a natural weak $g$-twisted $V_{\vir(A,f)}(1,0)$-module structure such that
\[
 Y_W(v_a,z)=\sum_{p\in \frac{r}{T}+\Z}L_p(a)z^{-p-2},
 \qquad a\in A_g^{[\frac{r}{T}]},\quad r\in\{0,\ldots,T-1\}.
\]
\end{corollary}

In Section~\ref{subsec:highest-weight-twisted} we construct such restricted modules by highest-weight induction and verify the ordinary condition.

\Needspace{3\baselineskip}
\subsection{Highest-weight induction and ordinary twisted modules}\label{subsec:highest-weight-twisted}
Throughout this subsection we keep the notation $A,f,g,T$ and the Lie algebra
$\vir_g(A,f)$ of Lemma~\ref{lem:frob-lca-construction}. We construct restricted $\vir_g(A,f)$-modules by highest-weight induction,
and use Corollary~\ref{cor:frob-twisted-modules} to obtain the corresponding weak $g$-twisted modules.
We then prove that their irreducible quotients are ordinary; when $A\ne0$, varying the lowest conformal weight yields
infinitely many pairwise non-isomorphic irreducible ordinary $g$-twisted $V_{\vir(A,f)}(1,0)$-modules.

In what follows, every $L_p(a)$ is taken with $a\in A_g^{[\frac rT]}$ and $p\in\frac rT+\Z$,
where $r\in\{0,\ldots,T-1\}$. We follow the construction of the untwisted induction in
\cite[(5.17)--(5.19)]{LYu}.

\begin{definition}\label{def:frob-highest-weight-module}
Let
\[
 \mathfrak p_g=\C K\oplus L_0(A_g^{[0]})
 \oplus\Span\{L_p(a):p>0\},\qquad
 \mathfrak n_{g,-}=\Span\{L_p(a):p<0\}.
\]
Denote by
\[
 (A_g^{[0]})^*=\operatorname{Hom}_{\C}(A_g^{[0]},\C)
\]
the dual space of $A_g^{[0]}$.
For $\chi\in(A_g^{[0]})^*$, let $\C u_\chi$ be the one-dimensional $\mathfrak p_g$-module with action
\[
 K u_\chi=u_\chi,\qquad L_p(a)u_\chi=0\quad(p>0),\qquad
 L_0(a)u_\chi=\chi(a)u_\chi\quad(a\in A_g^{[0]}).
\]
Define the \emph{induced highest-weight module}
\[
 \mathcal M_{\vir_g(A,f)}(1,\chi)
 =U(\vir_g(A,f))\otimes_{U(\mathfrak p_g)}\C u_\chi,
\]
and denote $1\otimes u_\chi$ again by $u_\chi$.
\end{definition}

By \eqref{eq:frobtwist}, $\mathfrak p_g$ and $\mathfrak n_{g,-}$ are Lie subalgebras,
and the above one-dimensional action respects the commutation relations of $\mathfrak p_g$, so the definition makes sense.
This induction differs from the vacuum induction of Section~\ref{sec:frob}. When $g=\id$,
the subalgebra $\mathfrak p_g$ of this subsection contains only the modes $L_p(a)$ with nonnegative indices, whereas the subalgebra $\mathfrak p$ of Section~\ref{sec:frob} also contains the modes $L_{-1}(a)$.
In this case, $\mathfrak n_{g,-}$ contains the modes $L_{-1}(a)$, and $L_{-1}(a)u_\chi=0$ is not required.

For the classification of irreducible graded modules in the untwisted case, see
\cite[Theorem~5.5, Corollary~5.6]{LYu}. We now prove the results needed for a general automorphism $g$ of finite order.

\begin{proposition}\label{prop:infinite}
For every $\chi\in(A_g^{[0]})^*$, the module $\mathcal M_{\vir_g(A,f)}(1,\chi)$ has a unique maximal proper
$\vir_g(A,f)$-submodule. Denote the corresponding irreducible quotient by $L_{\vir_g(A,f)}(1,\chi)$.
Then it is an ordinary $g$-twisted $V_{\vir(A,f)}(1,0)$-module
with lowest conformal weight $\chi(1_A)$.
If $A\ne0$, then for each such $g$ of finite order there exist infinitely many pairwise non-isomorphic
irreducible ordinary $g$-twisted $V_{\vir(A,f)}(1,0)$-modules.
\end{proposition}
\begin{proof}
By the vector space decomposition
$\vir_g(A,f)=\mathfrak n_{g,-}\oplus\mathfrak p_g$ and the PBW theorem,
\[
 \mathcal M_{\vir_g(A,f)}(1,\chi)\cong U(\mathfrak n_{g,-})\otimes\C u_\chi.
\]
Fix an ordered basis of $\mathfrak n_{g,-}$ consisting of negative modes.
Then the corresponding ordered PBW monomials applied to $u_\chi$ form a basis of $\mathcal M_{\vir_g(A,f)}(1,\chi)$.
Set $\deg K=\deg u_\chi=0$ and $\deg L_p(a)=-p$. The relations \eqref{eq:frobtwist}
are homogeneous with respect to this degree, so the above induction yields a grading
\[
 \mathcal M_{\vir_g(A,f)}(1,\chi)=\bigoplus_{d\in\frac1T\N}M(d),\qquad
 M(0)=\C u_\chi,\qquad L_p(a)M(d)\subseteq M(d-p),
\]
where $M(d)$ is the subspace of degree $d$, and we set $M(d)=0$ for $d<0$.
For a PBW monomial $L_{p_1}(a_1)\cdots L_{p_\ell}(a_\ell)u_\chi$ in $M(d)$, we have
\[
 p_j<0,\qquad \sum_{j=1}^{\ell}(-p_j)=d,\qquad
 \ell\le Td,\qquad -d\le p_j\le-\tfrac1T.
\]
Since $\dim A<\infty$, there are only finitely many such PBW monomials, so $\dim M(d)<\infty$.

Let $h=\chi(1_A)$. By \eqref{eq:frobtwist} and $L_0(1_A)u_\chi=hu_\chi$,
\[
 [L_0(1_A),L_p(a)]=-pL_p(a),\qquad
 L_0(1_A)|_{M(d)}=(h+d)\id.
\]
On the other hand, for $w\in M(d)$ and $p>d$, we have $L_p(a)w=0$.
Since every vector has only finitely many homogeneous components, the description of the modes in Lemma~\ref{lem:frob-lca-construction}~(3)
shows that $\mathcal M_{\vir_g(A,f)}(1,\chi)$ is restricted.
By Corollary~\ref{cor:frob-twisted-modules},
$\mathcal M_{\vir_g(A,f)}(1,\chi)$ carries a weak $g$-twisted $V_{\vir(A,f)}(1,0)$-module structure
with $L(0)=L_0(1_A)$. Hence $M(d)$ has conformal weight $h+d$.

We now construct the irreducible quotient. Let $\mathcal W\subseteq\mathcal M_{\vir_g(A,f)}(1,\chi)$ be any
$\vir_g(A,f)$-submodule. For $w\in \mathcal W$, write its finite homogeneous decomposition as
$w=\sum_{i=1}^k w_{d_i}$, where $w_{d_i}\in M(d_i)$ and the $d_i$ are pairwise distinct.
Since $\mathcal W$ is stable under $L_0(1_A)$,
Lagrange interpolation gives
\[
 w_{d_i}
 =\left(\prod_{\substack{1\le j\le k\\j\ne i}}
 \frac{L_0(1_A)-(h+d_j)\id}{d_i-d_j}\right)w\in \mathcal W.
\]
Hence every submodule is a graded subspace with respect to the above grading. If $\mathcal W$ is proper, then $\mathcal W\cap M(0)=0$;
otherwise $\mathcal W$ contains the generator $u_\chi$ and hence equals $\mathcal M_{\vir_g(A,f)}(1,\chi)$.
Therefore the sum of all proper submodules,
\[
 R=\sum_{\substack{\mathcal W\subsetneq\mathcal M_{\vir_g(A,f)}(1,\chi)\\\mathcal W\text{ a }\vir_g(A,f)\text{-submodule}}}\mathcal W
 \subseteq\bigoplus_{d\in\frac1T\Z_{>0}}M(d),
\]
is again a proper submodule and contains all proper submodules, so it is the unique maximal proper submodule.
Thus $L_{\vir_g(A,f)}(1,\chi)=\mathcal M_{\vir_g(A,f)}(1,\chi)/R$ is an irreducible $\vir_g(A,f)$-module.

This quotient is still restricted, and $K$ acts on it as the identity, so by
Corollary~\ref{cor:frob-twisted-modules} it carries a weak $g$-twisted
$V_{\vir(A,f)}(1,0)$-module structure. Its conformal weight spaces satisfy
\[
 (L_{\vir_g(A,f)}(1,\chi))_{h+d}\cong M(d)/(R\cap M(d)),\qquad d\in\tfrac1T\N.
\]
These spaces are all finite-dimensional, all conformal weights lie in $h+\frac1T\N$, and
\[
 (L_{\vir_g(A,f)}(1,\chi))_h=\C(u_\chi+R)\ne0.
\]
Therefore $L_{\vir_g(A,f)}(1,\chi)$ is ordinary with lowest conformal weight $h$.
Every $g$-twisted $V_{\vir(A,f)}(1,0)$-submodule is stable under all $L_p(a)$,
and $K$ acts on it as the identity, so it is also a $\vir_g(A,f)$-submodule.
By irreducibility as a Lie algebra module, $L_{\vir_g(A,f)}(1,\chi)$ is also irreducible as a twisted module.

If $A\ne0$, then $0\ne1_A\in A_g^{[0]}$, so $\chi(1_A)$ can be chosen to be any $h\in\R_{\ge0}$.
All conformal weights of the resulting modules are real, and the lowest one is exactly $h$,
so different values of $h$ give pairwise non-isomorphic irreducible ordinary $g$-twisted modules.
\end{proof}

In the next subsection, combining the simplicity results of Theorem~\ref{thm:frobsimple} and Proposition~\ref{prop:vir},
we apply Proposition~\ref{prop:infinite} to complete the proof of Theorem~\ref{thm:two}.

\Needspace{3\baselineskip}
\subsection{Proof of the theorem in conformal weight two}\label{subsec:proof-weight-two}
In this subsection we prove Theorem~\ref{thm:two} by combining the tensor decomposition of Proposition~\ref{prop:frobfactors},
the simplicity result of Theorem~\ref{thm:frobsimple}, and the twisted module constructions of this section.
The key point of the proof is that, if the decomposition contains a factor that is not a Virasoro minimal model,
then one can construct infinitely many pairwise non-isomorphic irreducible ordinary $g$-twisted $V$-modules.

\begin{proof}[Proof of Theorem~\ref{thm:two}]
Let $V$ and $g$ satisfy the hypotheses of Theorem~\ref{thm:two}.
If $V=\C\one$, then $g=\id$, every module is a direct sum of one-dimensional modules,
and (iv) is given by the empty tensor product, so the assertion holds.
Assume from now on that $V\ne\C\one$, and choose $T\in\Z_{>0}$ with $g^T=1$.

We first prove (i)$\Rightarrow$(iv). Assume (i). We begin by determining the tensor decomposition of $V$ and the action of $g$ on it.
In the notation of Remark~\ref{rem:frobenius-normalization},
let $(A,f)$ be the commutative Frobenius algebra obtained from $A=V_2$.
Applying Propositions~\ref{prop:local} and~\ref{prop:frobfactors}, we obtain
\[
 A=\bigoplus_{i\in J}A_i,\qquad
 V\cong\bigotimes_{i\in J}L_{\vir(A_i,f_i)}(1,0),\qquad f_i=f|_{A_i},\quad i\in J,
\]
where $J$ is finite and each $(A_i,f_i)$ is a local commutative Frobenius algebra, whose identity element is denoted by $e_i$.
If $A_i=\C e_i$, write $c_i=f_i(e_i)$.
Since the bilinear form defined by $f_i$ is nondegenerate, $c_i\ne0$. By \eqref{eq:vir-special-universal} and \eqref{eq:vir-special},
the corresponding Virasoro vertex operator algebra and its simple quotient are
$V_{\mathfrak{vir}}(c_i,0)$ and $L_{\mathfrak{vir}}(c_i,0)$, respectively.

Accordingly, let
\[
 \begin{aligned}
 J_{\mathrm{min}}
 &=\{i\in J:A_i=\C e_i,\ f_i(e_i)=c_{p,q}
       \text{ for some coprime }2\le p<q\},\\
 J_{\mathrm u}&=J\setminus J_{\mathrm{min}},
 \end{aligned}
\]
and write
\[
 A_{\mathrm u}=\bigoplus_{i\in J_{\mathrm u}}A_i,\qquad
 f_{\mathrm u}=f|_{A_{\mathrm u}},\qquad
 V_{\mathrm{min}}=\bigotimes_{i\in J_{\mathrm{min}}}L_{\vir(A_i,f_i)}(1,0).
\]
For $i\in J_{\mathrm u}$, if $\dim A_i>1$, then Theorem~\ref{thm:frobsimple}
shows that $V_{\vir(A_i,f_i)}(1,0)$ is simple; if $A_i=\C e_i$, then $c_i$ is not equal to any of the above $c_{p,q}$,
and the same conclusion follows from Proposition~\ref{prop:vir}. Hence, for all $i\in J_{\mathrm u}$,
$L_{\vir(A_i,f_i)}(1,0)\cong V_{\vir(A_i,f_i)}(1,0)$.
The proof of \eqref{eq:frobfactors-universal} in Proposition~\ref{prop:frobfactors} uses only the local decomposition, together with the PBW basis and the universal property in Remark~\ref{rem:frob-vacuum-pbw}~(2), (3);
hence it applies to an arbitrary commutative Frobenius algebra. Applying it to $(A_{\mathrm u},f_{\mathrm u})$ and its local decomposition
$A_{\mathrm u}=\bigoplus_{i\in J_{\mathrm u}}A_i$, we obtain
\[
 \bigotimes_{i\in J_{\mathrm u}}L_{\vir(A_i,f_i)}(1,0)
 \cong\bigotimes_{i\in J_{\mathrm u}}V_{\vir(A_i,f_i)}(1,0)
 \cong V_{\vir(A_{\mathrm u},f_{\mathrm u})}(1,0).
\]
Grouping the original tensor decomposition according to $J=J_{\mathrm{min}}\sqcup J_{\mathrm u}$, we get
\begin{equation}\label{eq:minu}
 V\cong V_{\mathrm{min}}\otimes V_{\mathrm u},\qquad
 V_{\mathrm u}:=V_{\vir(A_{\mathrm u},f_{\mathrm u})}(1,0).
\end{equation}
When $A_{\mathrm u}=0$, $V_{\mathrm u}$ is understood to be $\C\one$.

By Proposition~\ref{prop:frobfactors} and the uniqueness of the local decomposition in Proposition~\ref{prop:local}, $g$ induces a permutation
$\sigma$ of $J$ such that
\[
 g(A_i)=A_{\sigma(i)},\qquad g(e_i)=e_{\sigma(i)},\qquad
 f_{\sigma(i)}\circ g|_{A_i}=f_i,\qquad i\in J.
\]
In particular,
\[
 \dim A_{\sigma(i)}=\dim A_i,\qquad
 f_{\sigma(i)}(e_{\sigma(i)})=f_i(e_i),\qquad i\in J.
\]
Hence both $J_{\mathrm{min}}$ and $J_{\mathrm u}$ are $\sigma$-stable,
and $g$ preserves each of the two factors in \eqref{eq:minu}. Denoting the corresponding restrictions by
$g_{\mathrm{min}}$ and $g_{\mathrm u}$, we have
\[
 g=g_{\mathrm{min}}\otimes g_{\mathrm u},\qquad
 g_{\mathrm{min}}^T=g_{\mathrm u}^T=1.
\]
When $J_{\mathrm{min}}\ne\varnothing$, Lemma~\ref{lem:virasoro-permutations} shows that
$g_{\mathrm{min}}$ is the tensor factor permutation corresponding to $\sigma|_{J_{\mathrm{min}}}$.

If $J_{\mathrm u}\ne\varnothing$, then $A_{\mathrm u}\ne0$.
By Proposition~\ref{prop:infinite}, for each $h\in\R_{\ge0}$ there exists an irreducible ordinary
$g_{\mathrm u}$-twisted $V_{\mathrm u}$-module $M^{(h)}$ such that
\[
 (M^{(h)})_h\ne0,\qquad
 \operatorname{Spec}L_{M^{(h)}}(0)\subset h+\tfrac1T\N.
\]

In order to use these modules for $V$, we first construct a fixed irreducible ordinary
$g_{\mathrm{min}}$-twisted $V_{\mathrm{min}}$-module $N$.
If $V_{\mathrm{min}}=\C\one$, take $N=\C$.
Otherwise, decompose $\sigma|_{J_{\mathrm{min}}}$ into disjoint cycles
\[
 \sigma|_{J_{\mathrm{min}}}=\tau_1\cdots\tau_s,
\]
where $s\in\Z_{\ge1}$ and cycles of length $1$ are included.
For $1\le j\le s$, denote by $k_j$ the length of $\tau_j$
and by $c^{(j)}$ the common central charge of the factors it permutes, and let
\[
 U_j=L_{\mathfrak{vir}}(c^{(j)},0),\qquad 1\le j\le s.
\]
Identifying the corresponding factors via the isomorphisms of Lemma~\ref{lem:virasoro-permutations}
and reordering as in Remark~\ref{rem:permutation-cycle-construction}, we have
\[
 V_{\mathrm{min}}\cong\bigotimes_{j=1}^sU_j^{\otimes k_j},\qquad
 g_{\mathrm{min}}=\bigotimes_{j=1}^s\tau_j,
\]
where $\tau_j$ on the right-hand side denotes the permutation automorphism of $U_j^{\otimes k_j}$ induced by $\tau_j$.

The adjoint module of each $U_j$ is irreducible and ordinary,
so Proposition~\ref{prop:permutation}~(1) gives an irreducible ordinary
$\tau_j$-twisted $U_j^{\otimes k_j}$-module $N^{(j)}$.
Then, by Proposition~\ref{prop:permutation}~(2),
\[
 N=\bigotimes_{j=1}^sN^{(j)}
\]
is the required irreducible ordinary $g_{\mathrm{min}}$-twisted $V_{\mathrm{min}}$-module.

For this $N$, as in the proof of Proposition~\ref{prop:permutation}~(2),
we can choose a lowest conformal weight $h_N\in\C$ such that
\[
 N=\bigoplus_{n\in\N}N_{h_N+\frac nT},\qquad N_{h_N}\ne0.
\]

By Proposition~\ref{prop:permutation}~(2), for each $h\in\R_{\ge0}$,
$N\otimes M^{(h)}$ is an irreducible ordinary $g$-twisted $V$-module, and
\[
 (N\otimes M^{(h)})_{h_N+h}\ne0,\qquad
 \operatorname{Spec}L_{N\otimes M^{(h)}}(0)\subset h_N+h+\tfrac1T\N.
\]
If $h<h'$, then $h_N+h$ is a conformal weight of $N\otimes M^{(h)}$
but not of $N\otimes M^{(h')}$.
Since module isomorphisms preserve the eigenvalues of $L(0)$, these modules are pairwise non-isomorphic, contradicting (i).
Therefore $J_{\mathrm u}=\varnothing$, and hence $V\cong V_{\mathrm{min}}$, which is (iv).

By Remark~\ref{rem:module-relations}~(3), we have (ii)$\Rightarrow$(i).

If (iv) holds, then by Proposition~\ref{prop:vir} every tensor factor is a rational Virasoro minimal model.
By Proposition~\ref{prop:tensor}~(1), their finite tensor product is rational, so (iv)$\Rightarrow$(iii).
On the other hand, omit the trivial factors and group the factors according to their central charges;
by Lemma~\ref{lem:virasoro-permutations}, $g$ acts on each group as a permutation.
Proposition~\ref{prop:permutation}~(3) gives the corresponding twisted rationality for each group.
Then Proposition~\ref{prop:permutation}~(4), applied to the tensor product over the different groups,
shows that $V$ is $g$-rational, so (iv)$\Rightarrow$(ii).

Finally, if (iii) holds, then Remark~\ref{rem:module-relations}~(3) with $g=\id$
gives the finiteness of the set of isomorphism classes of irreducible ordinary untwisted modules.
Applying the implication (i)$\Rightarrow$(iv), already proved, to the identity automorphism, we obtain (iii)$\Rightarrow$(iv).

In summary,
\[
 \text{(ii)}\Longrightarrow\text{(i)}\Longrightarrow\text{(iv)}
 \Longrightarrow\text{(ii)},\qquad
 \text{(iii)}\Longleftrightarrow\text{(iv)}.
\]

Under these equivalent conditions, the $C_2$-cofiniteness statement of Proposition~\ref{prop:vir},
together with Proposition~\ref{prop:tensor}, gives
\[
 \dim R_V=\prod_{i=1}^r\dim R_{L_{\mathfrak{vir}}(c_{p_i,q_i},0)}<\infty.
\]
Hence $V$ is $C_2$-cofinite; together with the rationality already proved,
Proposition~\ref{prop:knownC2} shows that $V$ is regular.
\end{proof}

\section{Quadratic Lie algebras and graded vertex ideals of affine vertex algebras}\label{sec:current}
In this section we turn to conformal weight one. Let $V$ be a simple vertex operator algebra of CFT type that is generated by $V_1$ as a vertex algebra.
By Lian's structure theorem (Proposition~\ref{prop:lian}), $V_1$ carries a quadratic Lie algebra structure $(\g,B_{\g})$,
and $V$ is isomorphic, as a vertex algebra, to the simple graded quotient of $V_{(\g,B_{\g})}(1,0)$.
By Lemma~\ref{lem:split}, $\g$ decomposes as an orthogonal direct sum $\s\perp\q$, where $\s$ is a semisimple ideal and $\q$ has no nonzero semisimple ideal.
In this section we prove that $V_{(\q,B_{\q})}(1,0)$ has no nonzero proper graded vertex ideals, that is, $L_{(\q,B_{\q})}(1,0)=V_{(\q,B_{\q})}(1,0)$;
see Theorem~\ref{thm:currentsimple}.
In Section~\ref{sec:affine} we deduce from this the tensor product decomposition of $V$ (Proposition~\ref{prop:curfactors}) and complete the proof of Theorem~\ref{thm:one}.

Section~\ref{subsec:currents} recalls the Lie bracket and the bilinear form on $V_1$, introduces quadratic Lie algebras and the corresponding affine vertex algebras,
and recalls the rationality of simple affine vertex operator algebras and Lian's structure theorem.
Section~\ref{subsec:polynomial-limits} gives the above orthogonal decomposition following \cite[Corollary~5.6]{FS}, and, adapting the scaling idea of \cite[proof of Proposition~2]{Sch},
proves that a quadratic Lie algebra without nonzero semisimple ideals degenerates to an abelian Lie algebra through finitely many polynomial limits preserving the bilinear form;
see Proposition~\ref{prop:degenerate}.
Section~\ref{subsec:contravariant} constructs the normalized contravariant form on the induced vacuum module
and proves that its radical is precisely the maximal proper graded vertex ideal (Lemma~\ref{lem:radical}); hence it suffices to prove that this form is nondegenerate.
In Section~\ref{subsec:gram-affine}, as in Sections~\ref{subsec:det-criterion}--\ref{subsec:pbw-simplicity}, we compare this form, by means of PBW symmetrization, with an auxiliary bilinear form depending only on $B_{\q}$:
the Gram determinant in each degree remains unchanged along each polynomial limit, and when the Lie bracket is zero the two forms coincide and the Gram determinant is nonzero.
This gives Theorem~\ref{thm:currentsimple}.

\Needspace{3\baselineskip}
\subsection{Lie algebras in conformal weight one and affine vertex algebras}\label{subsec:currents}
In this subsection we fix a vertex operator algebra $V$ of CFT type and write $\g=V_1$;
this notation is kept in this section and the next, and $\g$ also denotes the Lie algebra structure on it.
We first recall the Lie bracket and the bilinear form on $\g$ given by vertex products, and introduce quadratic Lie algebras.
We then construct the affine vertex algebra associated with a quadratic Lie algebra, and describe its simple graded quotient and its weak modules.
For a finite-dimensional simple Lie algebra, we further recall the Sugawara construction and the rationality, $C_2$-cofiniteness, and regularity
of the corresponding simple affine vertex operator algebras.
Finally, we return to $V$: when $V$ is simple and generated by $\g$ as a vertex algebra, Lian's structure theorem realizes $V$ as
the simple graded quotient of the affine vertex algebra associated with $\g$.

Following \cite[Lemmas~3.4--3.5, Theorem~3.7]{Lian}, we recall the Lie bracket and the bilinear form on $\g$;
here the assumption of CFT type ensures that $x_1y\in V_0=\C\one$.

\begin{proposition}\label{prop:weight-one-lie}
For $x,y\in\g=V_1$, define
\begin{equation}\label{eq:currentproducts}
 [x,y]=x_0y,\qquad x_1y=B_{\g}(x,y)\one.
\end{equation}
Then $\g$ is a Lie algebra under this bracket, and $B_{\g}:\g\times\g\to\C$
is a symmetric invariant bilinear form, that is,
\[
 B_{\g}(x,y)=B_{\g}(y,x),\qquad
 B_{\g}([x,y],z)=B_{\g}(x,[y,z]),\qquad x,y,z\in\g.
\]
\end{proposition}

The form $B_{\g}$ here is defined by $x_1y$ in \eqref{eq:currentproducts}; it is related to the form $B_V$ of Section~\ref{sec:prelim} as follows.

\begin{remark}\label{rem:current-form-sign}
The above construction does not require $B_{\g}$ to be nondegenerate, nor does it require $L(1)\g=0$.
If $V$ has a normalized invariant bilinear form $B_V$,
then Proposition~\ref{prop:form}~(3) gives $L(1)\g=0$, and \eqref{eq:adjoint} then gives
\[
 B_V(x,y)=-B_{\g}(x,y),\qquad x,y\in\g.
\]
\end{remark}

We next recall the definition of a quadratic Lie algebra; see \cite[Section~4.1]{AleM}.
Such Lie algebras were first introduced in \cite[Section~5]{TW} under the name of metrisable Lie algebras, and their structure theory can be found in \cite[Theorems~I--III]{MR};
the same class of Lie algebras with a specified bilinear form is called self-dual in \cite[Section~5]{Lian} and
symmetric self-dual in \cite[Definition~2.1]{FS}.

\begin{definition}\label{def:quadratic}
A \emph{quadratic Lie algebra} is a pair $(\mathfrak a,B_{\mathfrak a})$, where
$\mathfrak a$ is a finite-dimensional complex Lie algebra and
$B_{\mathfrak a}:\mathfrak a\times\mathfrak a\to\C$ is a nondegenerate symmetric invariant bilinear form.
Here symmetry and invariance mean
\[
 B_{\mathfrak a}(x,y)=B_{\mathfrak a}(y,x),\qquad
 B_{\mathfrak a}([x,y],z)=B_{\mathfrak a}(x,[y,z]),\qquad x,y,z\in\mathfrak a,
\]
respectively, and nondegeneracy means
\[
 \rad B_{\mathfrak a}=\{x\in\mathfrak a:B_{\mathfrak a}(x,y)=0\text{ for all }y\in\mathfrak a\}=0.
\]
\end{definition}

From here through Proposition~\ref{prop:affine-restricted}, we fix a quadratic Lie algebra $(\mathfrak a,B_{\mathfrak a})$.
We construct the corresponding affine Lie algebra and induced vacuum module, and describe the vertex algebra structure on the latter;
see \cite{FZ,LiLocal,LL}.

\begin{definition}\label{def:affine-vacuum}
The \emph{affine Lie algebra} corresponding to $(\mathfrak a,B_{\mathfrak a})$ is
\[
 \widehat{\mathfrak a}_{B_{\mathfrak a}}
 =(\mathfrak a\otimes\C[t,t^{-1}])\oplus\C K,
\]
where $t$ is a formal variable and $K$ is a central element. For $x\in\mathfrak a$ and $m\in\Z$,
write $x(m)=x\otimes t^m$. The bracket is
\begin{equation}\label{eq:currentbracket}
 [x(m),y(n)]=[x,y](m+n)+m\delta_{m+n,0}B_{\mathfrak a}(x,y)K,
 \quad x,y\in\mathfrak a,\ m,n\in\Z.
\end{equation}
For $k\in\C$, if $K$ acts on a module as $k\,\id$, the module is said to be of \emph{level} $k$.

Write $\mathfrak a[t]=\mathfrak a\otimes\C[t]$. Let $\C_k=\C\one$ be the
one-dimensional $\mathfrak a[t]\oplus\C K$-module whose action is determined by
\[
 x(n)\one=0,\quad x\in\mathfrak a,\ n\in\Z_{\ge0},
 \qquad K\one=k\one.
\]
Denote the corresponding induced vacuum module by
\begin{equation}\label{eq:currentvac}
 V_{(\mathfrak a,B_{\mathfrak a})}(k,0)
 =U(\widehat{\mathfrak a}_{B_{\mathfrak a}})
 \otimes_{U(\mathfrak a[t]\oplus\C K)}\C_k,
\end{equation}
and denote $1\otimes\one$ again by $\one$. The first parameter $k$ indicates the level, and the second parameter $0$
indicates the trivial action of $\mathfrak a\otimes1$ on $\C_k$.
\end{definition}

The induced vacuum module above carries a vertex algebra structure and a natural grading;
see \cite[Section~6.2, Theorem~6.2.11]{LL} and \cite[Section~3.4, Theorem~4.8]{Lian}.

\begin{proposition}\label{prop:affine-vacuum-structure}
For each $k\in\C$, $V_{(\mathfrak a,B_{\mathfrak a})}(k,0)$ carries a unique vertex algebra structure with vacuum vector $\one$
satisfying
\begin{equation*}
 Y(x(-1)\one,z)=\sum_{m\in\Z}x(m)z^{-m-1},\qquad x\in\mathfrak a.
\end{equation*}
It is called the \textit{level-$k$ affine vertex algebra} associated with $(\mathfrak a,B_{\mathfrak a})$.
We identify $\mathfrak a$ with a subspace of $V_{(\mathfrak a,B_{\mathfrak a})}(k,0)$ via $x\mapsto x(-1)\one$; under this identification,
$x(m)$ and the vertex mode $x_m$ denote the same operator.

For $n\in\N$, denote by $V_{(\mathfrak a,B_{\mathfrak a})}(k,0)_n$ the subspace spanned by the vectors
\[
 x_1(-m_1)\cdots x_\ell(-m_\ell)\one,
 \qquad \ell\in\N,\quad x_j\in\mathfrak a,\quad m_j\in\Z_{\ge1},\quad \sum_{j=1}^\ell m_j=n,
\]
where $1\le j\le\ell$, and the vector is taken to be $\one$ when $\ell=0$. Then
\begin{equation}\label{eq:affine-grading}
 V_{(\mathfrak a,B_{\mathfrak a})}(k,0)=\bigoplus_{n\in\N}V_{(\mathfrak a,B_{\mathfrak a})}(k,0)_n,
\end{equation}
each $V_{(\mathfrak a,B_{\mathfrak a})}(k,0)_n$ is finite-dimensional, and
\[
 V_{(\mathfrak a,B_{\mathfrak a})}(k,0)_0=\C\one,\qquad
 V_{(\mathfrak a,B_{\mathfrak a})}(k,0)_1=\mathfrak a.
\]
The decomposition \eqref{eq:affine-grading} is called the \textit{natural grading} of $V_{(\mathfrak a,B_{\mathfrak a})}(k,0)$,
and nonzero vectors in $V_{(\mathfrak a,B_{\mathfrak a})}(k,0)_n$ are said to have \textit{degree} $n$.
Fix an ordered basis of $\mathfrak a$, and fix a total order on the modes $x(-m)$ with $x$ in this basis and $m\in\Z_{\ge1}$;
then the ordered monomials among the above vectors form a PBW basis of $V_{(\mathfrak a,B_{\mathfrak a})}(k,0)$.
The translation operator $D$ of $V_{(\mathfrak a,B_{\mathfrak a})}(k,0)$ satisfies
\[
 D\one=0,\qquad [D,x(m)]=-mx(m-1),\qquad x\in\mathfrak a,\ m\in\Z.
\]
\end{proposition}

Using this natural grading, we now introduce graded vertex ideals and the corresponding simple graded quotient.

\begin{definition}\label{def:affine-simple-quotient}
Let $k\in\C$.
A vertex ideal $I$ of $V_{(\mathfrak a,B_{\mathfrak a})}(k,0)$ is called a \emph{graded vertex ideal} of $V_{(\mathfrak a,B_{\mathfrak a})}(k,0)$
if, with respect to the natural grading \eqref{eq:affine-grading},
\[
 I=\bigoplus_{n\in\N}\bigl(I\cap V_{(\mathfrak a,B_{\mathfrak a})}(k,0)_n\bigr).
\]
If $I\ne V_{(\mathfrak a,B_{\mathfrak a})}(k,0)$, then $I$ is said to be \emph{proper}.
The degree-zero component of every proper graded vertex ideal of $V_{(\mathfrak a,B_{\mathfrak a})}(k,0)$ is zero.
Hence the sum of all of them is again a proper graded vertex ideal, and is therefore the unique maximal proper graded vertex ideal;
see \cite[Remarks~3.9.11, 6.2.24]{LL}.
Denote this ideal by $J_{(\mathfrak a,B_{\mathfrak a})}(k,0)$. The quotient vertex algebra
\begin{equation}\label{eq:affine-simple-quotient}
 L_{(\mathfrak a,B_{\mathfrak a})}(k,0)
 =V_{(\mathfrak a,B_{\mathfrak a})}(k,0)/J_{(\mathfrak a,B_{\mathfrak a})}(k,0)
\end{equation}
is called the \emph{simple graded quotient} of $V_{(\mathfrak a,B_{\mathfrak a})}(k,0)$.
\end{definition}

In what follows, we often absorb a nonzero level into the bilinear form; that is, we regard the induced vacuum module of level $k$ as
the induced vacuum module of level one for the bilinear form $kB_{\mathfrak a}$; see \cite[Remark~4.3.2]{LiLocal}.
Here $k\ne0$ ensures that $(\mathfrak a,kB_{\mathfrak a})$ is again a quadratic Lie algebra.

\begin{remark}\label{rem:affine-level-scaling}
For every $k\in\C$ with $k\ne0$,
\[
 V_{(\mathfrak a,kB_{\mathfrak a})}(1,0)=V_{(\mathfrak a,B_{\mathfrak a})}(k,0),\qquad
 L_{(\mathfrak a,kB_{\mathfrak a})}(1,0)=L_{(\mathfrak a,B_{\mathfrak a})}(k,0).
\]
\end{remark}

Following \cite[Theorems~6.2.12--6.2.13]{LL}, we describe the correspondence between restricted $\widehat{\mathfrak a}_{B_{\mathfrak a}}$-modules
and weak $V_{(\mathfrak a,B_{\mathfrak a})}(k,0)$-modules.

\begin{proposition}\label{prop:affine-restricted}
Let $k\in\C$. An $\widehat{\mathfrak a}_{B_{\mathfrak a}}$-module $M$ is said to be \textit{restricted} if for every
$x\in\mathfrak a$ and $w\in M$, we have $x(n)w=0$ for all sufficiently large integers $n$.
Then every restricted $\widehat{\mathfrak a}_{B_{\mathfrak a}}$-module $M$ of level $k$
carries a unique weak $V_{(\mathfrak a,B_{\mathfrak a})}(k,0)$-module structure such that
\[
 Y_M(x(-1)\one,z)=\sum_{n\in\Z}x(n)z^{-n-1},\qquad x\in\mathfrak a.
\]
Conversely, every weak $V_{(\mathfrak a,B_{\mathfrak a})}(k,0)$-module becomes, via this formula,
a restricted $\widehat{\mathfrak a}_{B_{\mathfrak a}}$-module of level $k$.
\end{proposition}

We now turn to vertex operator algebra structures. Following \cite[Remarks~6.2.15, 6.2.20, Theorem~6.2.18]{LL},
we recall the Sugawara construction for a finite-dimensional simple Lie algebra,
and, following \cite[Remarks~3.9.11, 6.2.24]{LL}, describe the corresponding simple quotient; see also \cite[Section~2.2]{JS}.

\begin{definition}\label{def:affine-sugawara}
Let $\s$ be a finite-dimensional simple Lie algebra, let $h^\vee$ be its dual Coxeter number, and let $d=\dim\s$.
Write
\[
 (x,y)=\frac1{2h^\vee}\Tr_{\s}(\operatorname{ad}x\operatorname{ad}y),\qquad x,y\in\s.
\]
The form $(\, ,\,)$ is called the \emph{normalized Killing form} of $\s$; with respect to this form, long roots have squared length $2$.
Choose a basis $u_1,\ldots,u_d$ of $\s$ and its dual basis $u^1,\ldots,u^d$ with respect to $(\, ,\,)$.
For $k\in\C$ with $k\ne-h^\vee$, set
\[
 \omega_{\s,k}=\frac1{2(k+h^\vee)}\sum_{i=1}^d u_i(-1)u^i(-1)\one.
\]
Then $\omega_{\s,k}$ is a conformal vector of $V_{(\s,(\, ,\,))}(k,0)$ with central charge $\frac{kd}{k+h^\vee}$,
and the corresponding $L(0)$-eigenspace decomposition is the natural grading \eqref{eq:affine-grading}.
The vector $\omega_{\s,k}$ is called the \emph{Sugawara conformal vector};
$V_{(\s,(\, ,\,))}(k,0)$ equipped with $\omega_{\s,k}$ is called the \emph{level-$k$ affine vertex operator algebra},
and $k=-h^\vee$ is called the \emph{critical level}.
Since every vertex ideal is stable under $L(0)$, $J_{(\s,(\, ,\,))}(k,0)$ is the unique maximal proper vertex ideal;
hence $L_{(\s,(\, ,\,))}(k,0)$, equipped with the image of $\omega_{\s,k}$, is the unique simple quotient of $V_{(\s,(\, ,\,))}(k,0)$,
and is called the \emph{level-$k$ simple affine vertex operator algebra}.
\end{definition}

We now recall the rationality, $C_2$-cofiniteness, and regularity of simple affine vertex operator algebras.
For rationality and $C_2$-cofiniteness, see \cite[Theorem~11.1]{AM2}, and, for $C_2$-cofiniteness, see also \cite[Proposition~12.6]{DLMmod};
for regularity, see \cite[Theorem~3.7]{DLMreg}.
When $k=0$, we have $L_{(\s,(\, ,\,))}(0,0)\cong\C\one$.

\begin{proposition}\label{prop:affineclass}
Let $\s$ be a finite-dimensional simple Lie algebra, let $h^\vee$ be its dual Coxeter number,
and let $(\, ,\,)$ be the normalized Killing form of $\s$ as in Definition~\ref{def:affine-sugawara}.
Let $k\in\C$ with $k\ne-h^\vee$. Then the
simple affine vertex operator algebra $L_{(\s,(\, ,\,))}(k,0)$, equipped with the Sugawara conformal vector $\omega_{\s,k}$, is rational if and only if $k\in\N$.
In this case, $L_{(\s,(\, ,\,))}(k,0)$ is $C_2$-cofinite and regular.
\end{proposition}

Finally, we return to the vertex operator algebra $V$ fixed in this subsection. When $V$ is simple and generated by $\g$ as a vertex algebra, Lian obtained the following structure theorem;
for (1), (2), and (3), see \cite[Corollary~5.3~(iii)]{Lian}, \cite[Theorem~4.11 and its proof]{Lian},
and \cite[Theorem~5.4~(i)]{Lian}, respectively. Assertion (1) ensures that
$V_{(\g,B_{\g})}(1,0)$ in Definition~\ref{def:affine-vacuum} is defined.

\begin{proposition}\label{prop:lian}
Let $V$ be a simple vertex operator algebra of CFT type that is generated by $\g=V_1$ as a vertex algebra.
\begin{enumerate}[label=(\arabic*)]
\item $B_{\g}$ is nondegenerate, and hence $(\g,B_{\g})$ is a quadratic Lie algebra.
\item There exists a unique surjective vertex algebra homomorphism
\[
 \begin{aligned}
 \mathsf q_V^{\mathfrak{Lie}}:V_{(\g,B_{\g})}(1,0)&\longrightarrow V,\\
 x(-1)\one&\longmapsto x,\qquad x\in\g.
 \end{aligned}
\]
This surjection preserves the natural gradings.
\item The above surjection induces a grading-preserving vertex algebra isomorphism
\[
 L_{(\g,B_{\g})}(1,0)\cong V.
\]
\end{enumerate}
\end{proposition}

This isomorphism is only a vertex algebra isomorphism, and need not preserve conformal vectors.
Under the hypotheses of Theorem~\ref{thm:one}, Proposition~\ref{prop:sugawara} will show that the conformal vector of $V$ is the sum of the
Sugawara conformal vectors (Definition~\ref{def:affine-sugawara}) of the simple ideals, which yields a vertex operator algebra isomorphism.

\Needspace{3\baselineskip}
\subsection{Orthogonal decompositions and polynomial limits}\label{subsec:polynomial-limits}
In this subsection, $(\mathfrak a,[\,,\,],B_{\mathfrak a})$ denotes an arbitrary quadratic Lie algebra, and the following definitions and results apply to every quadratic Lie algebra;
since we shall change the Lie bracket while keeping the bilinear form, we include the bracket $[\,,\,]$ in the notation.
We first recall the solvable radical and the Levi decomposition, and, following \cite[Section~2]{FS}, the orthogonal structure of quadratic Lie algebras;
then, by the argument of \cite[Corollary~5.6]{FS}, we decompose $\mathfrak a$ into a semisimple part and a part without nonzero semisimple ideals; see Lemma~\ref{lem:split}.
Next, following the scaling idea of \cite[proof of Proposition~2]{Sch}, which preserves the bilinear form, we introduce polynomial limits,
and prove that a quadratic Lie algebra without nonzero semisimple ideals degenerates to an abelian Lie algebra through finitely many polynomial limits;
see Proposition~\ref{prop:degenerate}.

We recall solvable Lie algebras and the solvable radical following \cite[Section~3.1]{Hum},
and Levi subalgebras following \cite[Chapter~III, Section~9]{Jacobson}.

\begin{definition}\label{def:solvable-radical}
Let $\mathfrak u$ be a finite-dimensional complex Lie algebra. We say that $\mathfrak u$ is \emph{abelian} if $[\mathfrak u,\mathfrak u]=0$.
We say that $\mathfrak u$ is \emph{solvable} if its \emph{derived series}
\[
 \mathfrak u\supseteq[\mathfrak u,\mathfrak u]\supseteq\bigl[[\mathfrak u,\mathfrak u],[\mathfrak u,\mathfrak u]\bigr]\supseteq\cdots
\]
reaches $0$ after finitely many steps, where each term is the bracket of the preceding term with itself.
The Lie algebra $\mathfrak u$ has a unique maximal solvable ideal, called the \emph{solvable radical} of $\mathfrak u$ and denoted by $\mathfrak r_{\mathfrak u}$.
We say that $\mathfrak u$ is \emph{semisimple} if $\mathfrak r_{\mathfrak u}=0$.
If a semisimple subalgebra $\mathfrak l\subseteq\mathfrak u$ satisfies $\mathfrak u=\mathfrak l\oplus\mathfrak r_{\mathfrak u}$ as vector spaces,
then $\mathfrak l$ is called a \emph{Levi subalgebra} of $\mathfrak u$, and $\mathfrak u=\mathfrak l\ltimes\mathfrak r_{\mathfrak u}$ is called a \emph{Levi decomposition} of $\mathfrak u$.
\end{definition}

We record the standard facts on Lie algebras needed later. For (1) and (2), see \cite[Section~3.1]{Hum} and \cite[Chapter~III, Section~9]{Jacobson}, respectively;
for (3), see \cite[Theorem~6.3]{Hum}; and for (4), see \cite[Sections~5.1--5.2]{Hum}.

\begin{proposition}\label{prop:lie-basic}
Let $\mathfrak u$ be a finite-dimensional complex Lie algebra.
\begin{enumerate}[label=(\arabic*)]
\item Subalgebras and quotients of solvable Lie algebras are solvable, and $\mathfrak u/\mathfrak r_{\mathfrak u}$ is semisimple.
\item $\mathfrak u$ has a Levi subalgebra.
\item If $\mathfrak u$ is semisimple, then every finite-dimensional $\mathfrak u$-module is completely reducible.
\item $\mathfrak u$ is semisimple if and only if it is a direct sum of simple ideals; in this case, every simple ideal of $\mathfrak u$ is one of these, and every ideal of $\mathfrak u$ is semisimple.
\end{enumerate}
\end{proposition}

Following \cite[Definition~2.3]{FS}, we introduce the subspaces of a quadratic Lie algebra related to the bilinear form, and recall centralizers and centers.

\begin{definition}\label{def:orthogonal}
For a linear subspace $X\subseteq\mathfrak a$, the subspace
\[
 X^\perp=\{x\in\mathfrak a:B_{\mathfrak a}(x,w)=0\text{ for all }w\in X\}
\]
is called the \emph{orthogonal complement} of $X$. We say that $X$ is \emph{isotropic} if $X\subseteq X^\perp$, and \emph{nondegenerate} if $X\cap X^\perp=0$.
The subspace $X$ is nondegenerate if and only if $B_{\mathfrak a}|_{X\times X}$ is nondegenerate.
If $\mathfrak a=X_1\oplus X_2$, where $X_1,X_2$ are ideals with $B_{\mathfrak a}(X_1,X_2)=0$,
we call this an \emph{orthogonal direct sum}, write $\mathfrak a=X_1\perp X_2$, and call $X_1$ and $X_2$ its \emph{orthogonal summands}.
The subspace
\[
 Z_{\mathfrak a}(X)=\{x\in\mathfrak a:[x,w]=0\text{ for all }w\in X\}
\]
is called the \emph{centralizer} of $X$, and $Z(\mathfrak a)=Z_{\mathfrak a}(\mathfrak a)$ is called the \emph{center} of $\mathfrak a$.
\end{definition}

The following basic properties can be found in the literature: for (1), see \cite[Propositions~1.5, 1.6]{EKM} and \cite[Definition~2.3]{FS};
for (2), see \cite[Lemma~2.4 and the remark following Lemma~2.5]{FS}; (3) is \cite[Lemma~2.5]{FS};
for (4), see the proof of \cite[Proposition~2.6~(1)]{FS} and \cite[Proposition~1.6]{EKM};
for (5), see the proof of \cite[Proposition~2.6~(2)]{FS} and the discussion around \cite[Lemma~3.2]{FS}.
Some of the arguments in \cite{FS} are stated for minimal ideals, but they apply equally to general ideals.

\begin{lemma}\label{lem:orthogonal-basic}
Let $X\subseteq\mathfrak a$ be a linear subspace.
\begin{enumerate}[label=(\arabic*)]
\item $(X^\perp)^\perp=X$ and $\dim X+\dim X^\perp=\dim\mathfrak a$. If $X$ is nondegenerate, then $\mathfrak a=X\oplus X^\perp$,
and the restrictions of $B_{\mathfrak a}$ to $X$ and to $X^\perp$ are both nondegenerate.
\item If $X$ is an ideal, then $X^\perp$ is an ideal and $[X,X^\perp]=0$.
\item $Z(\mathfrak a)=[\mathfrak a,\mathfrak a]^\perp$.
\item If $X$ is a nondegenerate ideal, then $\mathfrak a=X\perp X^\perp$, and $X$ and $X^\perp$, equipped with the restrictions of $B_{\mathfrak a}$, are quadratic Lie algebras.
\item If $X$ is an isotropic ideal, then $X$ is abelian, $X$ is an ideal of $X^\perp$, and $X^\perp/X$, equipped with
\[
 [x+X,\,y+X]_{X^\perp/X}=[x,y]+X,\qquad
 B_{X^\perp/X}(x+X,\,y+X)=B_{\mathfrak a}(x,y),\qquad x,y\in X^\perp,
\]
is a quadratic Lie algebra, denoted by $(X^\perp/X,[\,,\,]_{X^\perp/X},B_{X^\perp/X})$.
\end{enumerate}
\end{lemma}

We now decompose a quadratic Lie algebra as an orthogonal direct sum of a semisimple ideal and another ideal without nonzero semisimple ideals; this decomposition will be used in Section~\ref{sec:affine}.
By the argument of \cite[Section~3, Corollary~5.6]{FS}, one obtains the following result; the results of \cite{FS} are stated for real Lie algebras, but the argument applies equally to complex Lie algebras.

\begin{lemma}\label{lem:split}
Let $(\mathfrak a,[\,,\,],B_{\mathfrak a})$ be a quadratic Lie algebra. Then there is an orthogonal direct sum $\mathfrak a=\s\perp\q$,
where $\s$ is a semisimple ideal and $\q$ has no nonzero semisimple ideal;
moreover, $\s$ and $\q$, equipped with the restrictions of $B_{\mathfrak a}$, are quadratic Lie algebras.
\end{lemma}

Note that a quadratic Lie algebra without nonzero semisimple ideals need not be solvable.
For example, the classical double $\mathfrak{sl}_2\ltimes\mathfrak{sl}_2^*$ of $\mathfrak{sl}_2$ (see \cite[Definition~4.1, Remark~4.4, (4.3)]{FS})
has $\mathfrak{sl}_2$ as a nonzero Levi subalgebra, but one can verify that it has no nonzero semisimple ideal.

To introduce polynomial limits, we first define isometries following \cite[Section~1]{EKM}, and define isomorphisms of quadratic Lie algebras.

\begin{definition}\label{def:isometry}
Let $(\mathfrak a,[\,,\,]_{\mathfrak a},B_{\mathfrak a})$ and $(\mathfrak a',[\,,\,]_{\mathfrak a'},B_{\mathfrak a'})$ be quadratic Lie algebras.
A linear isomorphism $\varphi:\mathfrak a\to\mathfrak a'$ satisfying
\[
 B_{\mathfrak a'}(\varphi x,\varphi y)=B_{\mathfrak a}(x,y),\qquad x,y\in\mathfrak a,
\]
is called an \emph{isometry}; if $\varphi$ is also a Lie algebra isomorphism, then $\varphi$ is called an \emph{isomorphism of quadratic Lie algebras},
and we write $(\mathfrak a,[\,,\,]_{\mathfrak a},B_{\mathfrak a})\cong(\mathfrak a',[\,,\,]_{\mathfrak a'},B_{\mathfrak a'})$.
An isometry from $\mathfrak a$ to itself with respect to $B_{\mathfrak a}$ is called a \emph{$B_{\mathfrak a}$-isometry}.
\end{definition}

Following the scaling idea of \cite[proof of Proposition~2]{Sch}, which preserves the bilinear form, we now introduce polynomial limits;
the change of the bracket is described by the structure constants of the Lie bracket with respect to a basis; see \cite[Section~1.4]{Hum}.

\begin{definition}\label{def:polynomial-limit}
Let $d=\dim\mathfrak a$, and choose a basis $e_1,\ldots,e_d$ of $\mathfrak a$.
For each $\lambda\in\C$ with $\lambda\ne0$, let $\psi_\lambda$ be a $B_{\mathfrak a}$-isometry (see Definition~\ref{def:isometry}). Define
\begin{equation}\label{eq:pullback}
 [x,y]_\lambda=\psi_\lambda^{-1}[\psi_\lambda x,\psi_\lambda y],\qquad x,y\in\mathfrak a,
\end{equation}
and write $[e_i,e_j]_\lambda=\sum_{k=1}^d \mathsf c_{ij}^{k}(\lambda)e_k$ for $1\le i,j\le d$;
the coefficients $\mathsf c_{ij}^{k}(\lambda)$, $1\le i,j,k\le d$, are called the \emph{structure constants} of $[\,,\,]_\lambda$ with respect to the basis $e_1,\ldots,e_d$.
If all $\mathsf c_{ij}^{k}(\lambda)$ are polynomial functions of $\lambda$, we say that the \emph{polynomial limit} of $(\mathfrak a,[\,,\,])$ with respect to $(\psi_\lambda)_{\lambda\in\C\setminus\{0\}}$ exists.
In this case, define a bilinear map $[\,,\,]_0$ on $\mathfrak a$ by
\[
 [e_i,e_j]_0=\sum_{k=1}^d \mathsf c_{ij}^{k}(0)e_k,\qquad 1\le i,j\le d,
\]
extended by bilinearity; $(\mathfrak a,[\,,\,]_0)$ is called the \emph{polynomial limit} of $(\mathfrak a,[\,,\,])$ with respect to $(\psi_\lambda)_{\lambda\in\C\setminus\{0\}}$,
and we write $[\,,\,]\rightsquigarrow[\,,\,]_0$. This condition and $[\,,\,]_0$ do not depend on the choice of basis.
In particular, if $\psi_\lambda=\id_{\mathfrak a}$, then $[\,,\,]_0=[\,,\,]$.
\end{definition}

Polynomial limits preserve the structure of a quadratic Lie algebra.

\begin{lemma}\label{lem:polynomial-limit}
Keep the notation of Definition~\ref{def:polynomial-limit}. For $\lambda\in\C$ with $\lambda\ne0$, $\psi_\lambda$ gives an isomorphism of quadratic Lie algebras
$(\mathfrak a,[\,,\,]_\lambda,B_{\mathfrak a})\to(\mathfrak a,[\,,\,],B_{\mathfrak a})$.
If the polynomial limit of $(\mathfrak a,[\,,\,])$ with respect to $(\psi_\lambda)_{\lambda\in\C\setminus\{0\}}$ exists, then $(\mathfrak a,[\,,\,]_0,B_{\mathfrak a})$ is again a quadratic Lie algebra.
\end{lemma}
\begin{proof}
The first assertion follows directly from \eqref{eq:pullback} and the fact that $\psi_\lambda$ is a $B_{\mathfrak a}$-isometry.
With respect to the basis $e_1,\ldots,e_d$, skew-symmetry, the Jacobi identity, and the invariance of $B_{\mathfrak a}$ are polynomial identities in the $\mathsf c_{ij}^{k}(\lambda)$.
Since $\psi_\lambda$ is such an isomorphism for every nonzero $\lambda\in\C$,
these identities hold for all nonzero $\lambda\in\C$, and hence also at $\lambda=0$.
\end{proof}

To construct explicit isometries, we need to complete an isotropic subspace to a nondegenerate subspace; for the one-dimensional case, see \cite[Lemma~1.23]{EKM}.

\begin{lemma}\label{lem:isotropic-complement}
Let $\mathfrak a_{(-1)}\subseteq\mathfrak a$ be an isotropic subspace. Then there exists an isotropic subspace $\mathfrak a_{(1)}\subseteq\mathfrak a$ with $\dim \mathfrak a_{(1)}=\dim \mathfrak a_{(-1)}$
such that $\mathfrak a_{(1)}\oplus \mathfrak a_{(-1)}$ is nondegenerate. Set $\mathfrak a_{(0)}=(\mathfrak a_{(1)}\oplus \mathfrak a_{(-1)})^\perp$. Then $\mathfrak a_{(0)}$ is nondegenerate, and
\[
 \mathfrak a=\mathfrak a_{(1)}\oplus \mathfrak a_{(0)}\oplus \mathfrak a_{(-1)},\qquad \mathfrak a_{(-1)}^\perp=\mathfrak a_{(0)}\oplus \mathfrak a_{(-1)}.
\]
\end{lemma}
\begin{proof}
Let $p=\dim \mathfrak a_{(-1)}$, and choose a basis $y_1,\ldots,y_p$ of $\mathfrak a_{(-1)}$. Since $B_{\mathfrak a}$ is nondegenerate,
we can choose $x_1,\ldots,x_p\in\mathfrak a$ with $B_{\mathfrak a}(x_i,y_j)=\delta_{ij}$ for $1\le i,j\le p$. Let
\[
 \tilde x_i=x_i-\frac12\sum_{j=1}^pB_{\mathfrak a}(x_i,x_j)y_j,\qquad 1\le i\le p,
\]
and let $\mathfrak a_{(1)}=\Span\{\tilde x_1,\ldots,\tilde x_p\}$. Since $B_{\mathfrak a}(\mathfrak a_{(-1)},\mathfrak a_{(-1)})=0$, we have
$B_{\mathfrak a}(\tilde x_i,y_j)=\delta_{ij}$ and $B_{\mathfrak a}(\tilde x_i,\tilde x_j)=0$ for $1\le i,j\le p$,
so $\mathfrak a_{(1)}$ is isotropic, $\dim \mathfrak a_{(1)}=p$, and $\mathfrak a_{(1)}\oplus \mathfrak a_{(-1)}$ is nondegenerate.
By Lemma~\ref{lem:orthogonal-basic}~(1), $\mathfrak a=(\mathfrak a_{(1)}\oplus \mathfrak a_{(-1)})\oplus \mathfrak a_{(0)}$, and $\mathfrak a_{(0)}$ is nondegenerate.
Finally, $\mathfrak a_{(0)}\oplus \mathfrak a_{(-1)}\subseteq \mathfrak a_{(-1)}^\perp$, and both have dimension $\dim\mathfrak a-p$, so $\mathfrak a_{(-1)}^\perp=\mathfrak a_{(0)}\oplus \mathfrak a_{(-1)}$.
\end{proof}

The following lemma is the central computation of this subsection: in the polynomial limit obtained by scaling along an isotropic ideal $\mathfrak a_{(-1)}$,
both $\mathfrak a_{(-1)}$ and the isotropic subspace $\mathfrak a_{(1)}$ paired with it fall into the center,
while the bracket on $\mathfrak a_{(0)}$ is isomorphic to the bracket of the quadratic Lie algebra $\mathfrak a_{(-1)}^\perp/\mathfrak a_{(-1)}$.

\begin{lemma}\label{lem:isotropic-contraction}
Let $\mathfrak a_{(-1)}\subseteq\mathfrak a$ be an isotropic ideal, and take $\mathfrak a_{(1)}$ and $\mathfrak a_{(0)}$ as in Lemma~\ref{lem:isotropic-complement}.
For $\lambda\in\C$ with $\lambda\ne0$, let $\psi_\lambda|_{\mathfrak a_{(j)}}=\lambda^j\id_{\mathfrak a_{(j)}}$ for $j=1,0,-1$.
Then each $\psi_\lambda$ is a $B_{\mathfrak a}$-isometry, the polynomial limit of $(\mathfrak a,[\,,\,])$ with respect to $(\psi_\lambda)_{\lambda\in\C\setminus\{0\}}$ exists,
and its bracket $[\,,\,]_0$ satisfies
\[
 [\mathfrak a_{(1)}\oplus\mathfrak a_{(-1)},\mathfrak a]_0=0,\qquad [\mathfrak a_{(0)},\mathfrak a_{(0)}]_0\subseteq\mathfrak a_{(0)}.
\]
Hence there is an orthogonal direct sum $\mathfrak a=\mathfrak a_{(0)}\perp(\mathfrak a_{(1)}\oplus\mathfrak a_{(-1)})$ with respect to the quadratic Lie algebra $(\mathfrak a,[\,,\,]_0,B_{\mathfrak a})$,
and $\mathfrak a_{(1)}\oplus\mathfrak a_{(-1)}$ is contained in the center of $(\mathfrak a,[\,,\,]_0)$.
Moreover, the natural projection $\mathfrak a_{(0)}\to \mathfrak a_{(-1)}^\perp/\mathfrak a_{(-1)}$ gives an isomorphism of quadratic Lie algebras
\[
 (\mathfrak a_{(0)},[\,,\,]_0,B_{\mathfrak a}|_{\mathfrak a_{(0)}\times\mathfrak a_{(0)}})\cong
 \bigl(\mathfrak a_{(-1)}^\perp/\mathfrak a_{(-1)},[\,,\,]_{\mathfrak a_{(-1)}^\perp/\mathfrak a_{(-1)}},B_{\mathfrak a_{(-1)}^\perp/\mathfrak a_{(-1)}}\bigr),
\]
where the right-hand side is as in Lemma~\ref{lem:orthogonal-basic}~(5).
\end{lemma}
\begin{proof}
By Lemma~\ref{lem:isotropic-complement}, $B_{\mathfrak a}$ can be nonzero only on $\mathfrak a_{(1)}\times\mathfrak a_{(-1)}$, $\mathfrak a_{(-1)}\times\mathfrak a_{(1)}$,
and $\mathfrak a_{(0)}\times\mathfrak a_{(0)}$, so each $\psi_\lambda$ is a $B_{\mathfrak a}$-isometry.
For $x_1,x_2,x_3\in\mathfrak a$, let $\theta(x_1,x_2,x_3)=B_{\mathfrak a}([x_1,x_2],x_3)$.
By skew-symmetry and invariance, $\theta$ is an alternating trilinear form, that is, it is linear in each variable and changes sign when any two variables are interchanged.
By \eqref{eq:pullback} and the fact that $\psi_\lambda$ is an isometry,
\[
 B_{\mathfrak a}([x_1,x_2]_\lambda,x_3)=\theta(\psi_\lambda x_1,\psi_\lambda x_2,\psi_\lambda x_3),
\]
so $B_{\mathfrak a}([\,,\,]_\lambda,\,\cdot\,)$ equals $\lambda^{i+j+k}\theta$ on $\mathfrak a_{(i)}\times\mathfrak a_{(j)}\times\mathfrak a_{(k)}$,
where $i,j,k\in\{1,0,-1\}$.
By Lemma~\ref{lem:orthogonal-basic}~(2), $[\mathfrak a_{(-1)},\mathfrak a_{(-1)}^\perp]=0$; moreover, by Lemma~\ref{lem:isotropic-complement}, $\mathfrak a_{(-1)}^\perp=\mathfrak a_{(0)}\oplus\mathfrak a_{(-1)}$.
Hence, if one variable lies in $\mathfrak a_{(-1)}$ and another lies in $\mathfrak a_{(0)}\oplus\mathfrak a_{(-1)}$, then $\theta$ vanishes on the corresponding component.
Up to permutation of the variables, the components that can be nonzero correspond to
\[
 (i,j,k)=(0,0,0),\ (1,0,0),\ (1,1,0),\ (1,1,1),\ (1,1,-1),
\]
whose exponents $i+j+k$ are $0,1,2,3,1$, respectively, all nonnegative.
Choose a basis $e_1,\ldots,e_d$ compatible with $\mathfrak a=\mathfrak a_{(1)}\oplus\mathfrak a_{(0)}\oplus\mathfrak a_{(-1)}$, that is, each $e_\alpha$ lies in some $\mathfrak a_{(j)}$,
and write $\mathsf b_{\alpha\beta}=B_{\mathfrak a}(e_\alpha,e_\beta)$ and $(\mathsf b^{\alpha\beta})=(\mathsf b_{\alpha\beta})^{-1}$, where $1\le\alpha,\beta\le d$.
Since $B_{\mathfrak a}$ is nondegenerate and independent of $\lambda$, the structure constants of Definition~\ref{def:polynomial-limit} with respect to this basis satisfy
\[
 \mathsf c_{\alpha\beta}^{\gamma}(\lambda)=\sum_{\delta=1}^dB_{\mathfrak a}([e_\alpha,e_\beta]_\lambda,e_\delta)\,\mathsf b^{\delta\gamma},
 \qquad 1\le\alpha,\beta,\gamma\le d,
\]
where each $B_{\mathfrak a}([e_\alpha,e_\beta]_\lambda,e_\delta)$ equals the value of $\theta$ on one of the above components multiplied by a nonnegative integral power of $\lambda$.
Hence each $\mathsf c_{\alpha\beta}^{\gamma}(\lambda)$ is a polynomial function of $\lambda$, that is, the polynomial limit exists.

At $\lambda=0$, only the $(0,0,0)$ component survives: $B_{\mathfrak a}([x_1,x_2]_0,x_3)$ equals $\theta(x_1,x_2,x_3)$ when $x_1,x_2,x_3\in\mathfrak a_{(0)}$,
and vanishes as soon as one of $x_1,x_2,x_3$ lies in $\mathfrak a_{(1)}\oplus\mathfrak a_{(-1)}$.
Since $B_{\mathfrak a}$ is nondegenerate, $[\mathfrak a_{(1)}\oplus\mathfrak a_{(-1)},\mathfrak a]_0\subseteq\mathfrak a^\perp=0$;
moreover, by the definition of $\mathfrak a_{(0)}$, $[\mathfrak a_{(0)},\mathfrak a_{(0)}]_0\subseteq(\mathfrak a_{(1)}\oplus\mathfrak a_{(-1)})^\perp=\mathfrak a_{(0)}$.
Hence $\mathfrak a_{(0)}$ and $\mathfrak a_{(1)}\oplus\mathfrak a_{(-1)}$ are both ideals of $(\mathfrak a,[\,,\,]_0)$.
By Lemma~\ref{lem:polynomial-limit}, $(\mathfrak a,[\,,\,]_0,B_{\mathfrak a})$ is a quadratic Lie algebra;
since $\mathfrak a_{(0)}$ and $\mathfrak a_{(1)}\oplus\mathfrak a_{(-1)}$ are orthogonal with respect to $B_{\mathfrak a}$, we obtain the stated orthogonal direct sum.

Finally, let $h_1,h_2,h_3\in\mathfrak a_{(0)}$. Since $\mathfrak a_{(-1)}$ is an ideal, Lemma~\ref{lem:orthogonal-basic}~(2) shows that
$\mathfrak a_{(-1)}^\perp=\mathfrak a_{(0)}\oplus\mathfrak a_{(-1)}$ is also an ideal; since $h_1,h_2\in\mathfrak a_{(0)}\subseteq\mathfrak a_{(-1)}^\perp$,
we can write $[h_1,h_2]=h+y$ with $h\in\mathfrak a_{(0)}$ and $y\in\mathfrak a_{(-1)}$. Then
\[
 B_{\mathfrak a}([h_1,h_2]_0,h_3)=\theta(h_1,h_2,h_3)=B_{\mathfrak a}(h+y,h_3)=B_{\mathfrak a}(h,h_3).
\]
By Lemma~\ref{lem:isotropic-complement} and Definition~\ref{def:orthogonal}, $B_{\mathfrak a}|_{\mathfrak a_{(0)}\times\mathfrak a_{(0)}}$ is nondegenerate, so $[h_1,h_2]_0=h$.
Hence, for $h_1,h_2\in\mathfrak a_{(0)}$,
\[
 [h_1,h_2]_0+\mathfrak a_{(-1)}=h+\mathfrak a_{(-1)}=[h_1,h_2]+\mathfrak a_{(-1)}
 =[h_1+\mathfrak a_{(-1)},\,h_2+\mathfrak a_{(-1)}]_{\mathfrak a_{(-1)}^\perp/\mathfrak a_{(-1)}},
\]
that is, the natural projection $\mathfrak a_{(0)}\to\mathfrak a_{(-1)}^\perp/\mathfrak a_{(-1)}$ preserves brackets.
Since $\mathfrak a_{(-1)}^\perp=\mathfrak a_{(0)}\oplus\mathfrak a_{(-1)}$, it is a linear isomorphism; moreover, for $h_1,h_2\in\mathfrak a_{(0)}$,
the definition of $B_{\mathfrak a_{(-1)}^\perp/\mathfrak a_{(-1)}}$ gives
\[
 B_{\mathfrak a_{(-1)}^\perp/\mathfrak a_{(-1)}}(h_1+\mathfrak a_{(-1)},h_2+\mathfrak a_{(-1)})=B_{\mathfrak a}(h_1,h_2),
\]
that is, the projection preserves the bilinear forms.
\end{proof}

The next three lemmas constitute the three steps of the proof of Proposition~\ref{prop:degenerate}:
first, one polynomial limit reduces a quadratic Lie algebra without nonzero semisimple ideals to the solvable case;
second, in the solvable case, at most one polynomial limit splits off a nonzero orthogonal summand contained in the center;
third, polynomial limits on the orthogonal complement of this summand extend to the whole Lie algebra, which allows induction on the dimension.

Step 1. Suppose that $\mathfrak a$ has no nonzero semisimple ideal. As noted above, $\mathfrak a$ need not be solvable;
the following lemma reduces it to the solvable case by one polynomial limit.

\begin{lemma}\label{lem:to-solvable}
Suppose that the quadratic Lie algebra $(\mathfrak a,[\,,\,],B_{\mathfrak a})$ has no nonzero semisimple ideal, and let $\mathfrak r=\mathfrak r_{\mathfrak a}$ be the solvable radical of $\mathfrak a$.
Then $\mathfrak r^\perp$ is an isotropic ideal. If Lemma~\ref{lem:isotropic-contraction} is applied to $\mathfrak a_{(-1)}=\mathfrak r^\perp$,
then the resulting polynomial limit $(\mathfrak a,[\,,\,]_0,B_{\mathfrak a})$ is a solvable quadratic Lie algebra.
\end{lemma}
\begin{proof}
Write $\mathfrak a_{(-1)}=\mathfrak r^\perp$. By Lemma~\ref{lem:orthogonal-basic}~(2), $\mathfrak a_{(-1)}$ is an ideal and $[\mathfrak r,\mathfrak a_{(-1)}]=0$.
By Proposition~\ref{prop:lie-basic}~(2), choose a Levi subalgebra $\mathfrak l$.
Since $\mathfrak a_{(-1)}$ and $\mathfrak r$ are both ideals of $\mathfrak a$, under the adjoint action $x\cdot u=[x,u]$, $x\in\mathfrak l$,
$\mathfrak a_{(-1)}$ is a finite-dimensional $\mathfrak l$-module and $\mathfrak a_{(-1)}\cap\mathfrak r$ is a submodule;
by Proposition~\ref{prop:lie-basic}~(3), there is an $\mathfrak l$-stable complement $\mathfrak k$ with $\mathfrak a_{(-1)}=(\mathfrak a_{(-1)}\cap\mathfrak r)\oplus\mathfrak k$.
By $\mathfrak a=\mathfrak l\oplus\mathfrak r$, $[\mathfrak r,\mathfrak a_{(-1)}]=0$, and $\mathfrak k\subseteq\mathfrak a_{(-1)}$, we have
$[\mathfrak a,\mathfrak k]=[\mathfrak l,\mathfrak k]+[\mathfrak r,\mathfrak k]=[\mathfrak l,\mathfrak k]\subseteq\mathfrak k$, so $\mathfrak k$ is an ideal.
Since $\mathfrak k\cap\mathfrak r=0$, the natural projection $\mathfrak a\to\mathfrak a/\mathfrak r$ maps $\mathfrak k$ isomorphically onto an ideal of $\mathfrak a/\mathfrak r$;
by Proposition~\ref{prop:lie-basic}~(1), (4), $\mathfrak k$ is semisimple, and hence $\mathfrak k=0$ by assumption.
Therefore $\mathfrak a_{(-1)}\subseteq\mathfrak r=\mathfrak a_{(-1)}^\perp$, that is, $\mathfrak a_{(-1)}$ is an isotropic ideal.

By Lemma~\ref{lem:isotropic-contraction}, there is an orthogonal direct sum $\mathfrak a=\mathfrak a_{(0)}\perp(\mathfrak a_{(1)}\oplus\mathfrak a_{(-1)})$ with respect to the quadratic Lie algebra $(\mathfrak a,[\,,\,]_0,B_{\mathfrak a})$,
$\mathfrak a_{(1)}\oplus\mathfrak a_{(-1)}$ is contained in the center of $(\mathfrak a,[\,,\,]_0)$, and $(\mathfrak a_{(0)},[\,,\,]_0)\cong(\mathfrak a_{(-1)}^\perp/\mathfrak a_{(-1)},[\,,\,]_{\mathfrak a_{(-1)}^\perp/\mathfrak a_{(-1)}})
=(\mathfrak r/\mathfrak r^\perp,[\,,\,]_{\mathfrak r/\mathfrak r^\perp})$.
By Proposition~\ref{prop:lie-basic}~(1), $\mathfrak r/\mathfrak r^\perp$ is solvable, and hence so is $(\mathfrak a_{(0)},[\,,\,]_0)$;
moreover, $\mathfrak a_{(1)}\oplus\mathfrak a_{(-1)}$ is contained in the center, so $[\mathfrak a,\mathfrak a]_0=[\mathfrak a_{(0)},\mathfrak a_{(0)}]_0$, and therefore $(\mathfrak a,[\,,\,]_0)$ is solvable.
By Lemma~\ref{lem:polynomial-limit}, it is a quadratic Lie algebra when equipped with $B_{\mathfrak a}$.
\end{proof}

Step 2. For solvable quadratic Lie algebras, the following lemma splits off, by at most one polynomial limit, a nonzero orthogonal summand contained in the center,
which allows induction on the dimension; for the related structure theorem, see \cite[Theorem~3.9]{FS}.

\begin{lemma}\label{lem:solvable-step}
Let $(\mathfrak a,[\,,\,],B_{\mathfrak a})$ be a solvable quadratic Lie algebra with $\mathfrak a\ne0$.
Then there exist a polynomial limit $(\mathfrak a,[\,,\,]_0)$ of $(\mathfrak a,[\,,\,])$ and an orthogonal direct sum
$\mathfrak a=\mathfrak c\perp\mathfrak z$ with respect to the quadratic Lie algebra $(\mathfrak a,[\,,\,]_0,B_{\mathfrak a})$ such that $\mathfrak z\ne0$ is contained in the center of $(\mathfrak a,[\,,\,]_0)$
and $\mathfrak c$, equipped with the restrictions of $[\,,\,]_0$ and $B_{\mathfrak a}$, is a solvable quadratic Lie algebra.
\end{lemma}
\begin{proof}
By solvability and $\mathfrak a\ne0$, we have $[\mathfrak a,\mathfrak a]\ne\mathfrak a$;
by Lemma~\ref{lem:orthogonal-basic}~(1), (3), $Z(\mathfrak a)=[\mathfrak a,\mathfrak a]^\perp\ne0$.

If there exists $\zeta\in Z(\mathfrak a)$ with $B_{\mathfrak a}(\zeta,\zeta)\ne0$, take $\psi_\lambda=\id_{\mathfrak a}$,
so that $[\,,\,]_0=[\,,\,]$. In this case $\C\zeta$ is a nondegenerate ideal, and by Lemma~\ref{lem:orthogonal-basic}~(4)
we may take $\mathfrak z=\C\zeta$ and $\mathfrak c=(\C\zeta)^\perp$; as an ideal of a solvable Lie algebra, $\mathfrak c$ is again solvable.

Otherwise, every $\zeta\in Z(\mathfrak a)$ satisfies $B_{\mathfrak a}(\zeta,\zeta)=0$. Take $0\ne\zeta\in Z(\mathfrak a)$;
then $\mathfrak a_{(-1)}=\C\zeta$ is an isotropic ideal. Apply Lemma~\ref{lem:isotropic-contraction} to $\mathfrak a_{(-1)}$,
and take $\mathfrak c=\mathfrak a_{(0)}$ and $\mathfrak z=\mathfrak a_{(1)}\oplus\mathfrak a_{(-1)}\ne0$.
Then $(\mathfrak c,[\,,\,]_0)\cong((\C\zeta)^\perp/\C\zeta,[\,,\,]_{(\C\zeta)^\perp/\C\zeta})$, which is a quotient of an ideal of a solvable Lie algebra,
and hence is solvable by Proposition~\ref{prop:lie-basic}~(1); by Lemma~\ref{lem:isotropic-contraction}, it is a quadratic Lie algebra when equipped with the restriction of $B_{\mathfrak a}$.
\end{proof}

Step 3. For the induction on the dimension, we also need to extend polynomial limits on an orthogonal summand to the whole Lie algebra.

\begin{lemma}\label{lem:central-extension}
Let $\mathfrak a=\mathfrak c\perp\mathfrak z$ be an orthogonal direct sum with respect to the quadratic Lie algebra $(\mathfrak a,[\,,\,],B_{\mathfrak a})$, with $\mathfrak z$ contained in $Z(\mathfrak a)$.
Let $(\psi_\lambda)_{\lambda\in\C\setminus\{0\}}$ be a family of $B_{\mathfrak a}|_{\mathfrak c\times\mathfrak c}$-isometries
such that the polynomial limit $[\,,\,]_{\mathfrak c,0}$ of the bracket of $\mathfrak c$ with respect to them exists.
Then the maps $\Psi_\lambda=\psi_\lambda\oplus\id_{\mathfrak z}$ are $B_{\mathfrak a}$-isometries, and the polynomial limit $[\,,\,]_0$ of the bracket of $\mathfrak a$ with respect to $(\Psi_\lambda)_{\lambda\in\C\setminus\{0\}}$
exists. Furthermore, $[x,y]_0=[x,y]_{\mathfrak c,0}$ for $x,y\in\mathfrak c$, $\mathfrak z$ is contained in the center of $(\mathfrak a,[\,,\,]_0)$,
and $\mathfrak a=\mathfrak c\perp\mathfrak z$ is an orthogonal direct sum with respect to the quadratic Lie algebra $(\mathfrak a,[\,,\,]_0,B_{\mathfrak a})$.
\end{lemma}
\begin{proof}
Since $\psi_\lambda$ is a $B_{\mathfrak a}|_{\mathfrak c\times\mathfrak c}$-isometry and since $\mathfrak c$ and $\mathfrak z$ are orthogonal with respect to $B_{\mathfrak a}$, each $\Psi_\lambda$ is a $B_{\mathfrak a}$-isometry.
For $x,y\in\mathfrak c$ and $z_1,z_2\in\mathfrak z$, since $\mathfrak z$ is contained in the center and $\mathfrak c$ is an ideal,
\[
 [x+z_1,\,y+z_2]_\lambda=\Psi_\lambda^{-1}[\Psi_\lambda(x+z_1),\Psi_\lambda(y+z_2)]=\psi_\lambda^{-1}[\psi_\lambda x,\psi_\lambda y].
\]
Let $d=\dim\mathfrak a$ and $d'=\dim\mathfrak c$. Choose a basis $e_1,\ldots,e_{d'}$ of $\mathfrak c$ and a basis $e_{d'+1},\ldots,e_d$ of $\mathfrak z$; together they form a basis of $\mathfrak a$.
By the above identity, for $1\le i,j\le d'$, the structure constants of $[e_i,e_j]_\lambda$ are those of the bracket of $\mathfrak c$ with respect to $(\psi_\lambda)_{\lambda\in\C\setminus\{0\}}$,
which are polynomial functions of $\lambda$ by assumption; if $i>d'$ or $j>d'$, then $[e_i,e_j]_\lambda=0$.
Hence the polynomial limit exists, and setting $\lambda=0$ gives $[x+z_1,\,y+z_2]_0=[x,y]_{\mathfrak c,0}$.
Taking $x=0$ gives $[\mathfrak z,\mathfrak a]_0=0$, that is, $\mathfrak z$ is contained in the center of $(\mathfrak a,[\,,\,]_0)$.
By Lemma~\ref{lem:polynomial-limit}, $(\mathfrak a,[\,,\,]_0,B_{\mathfrak a})$ is a quadratic Lie algebra;
since $[\mathfrak c,\mathfrak c]_0\subseteq\mathfrak c$ and $[\mathfrak z,\mathfrak a]_0=0$, both $\mathfrak c$ and $\mathfrak z$ are ideals of it, so $\mathfrak a=\mathfrak c\perp\mathfrak z$ is an orthogonal direct sum with respect to the quadratic Lie algebra $(\mathfrak a,[\,,\,]_0,B_{\mathfrak a})$.
\end{proof}

Combining the above three steps, we obtain the main result of this subsection; it will be used in the proof of Theorem~\ref{thm:currentsimple}.

\begin{proposition}\label{prop:degenerate}
Suppose that the quadratic Lie algebra $(\q,[\,,\,],B_{\q})$ has no nonzero semisimple ideal. Then there exist $s\in\N$ and Lie brackets
\[
 [\,,\,]^{(0)}=[\,,\,],\quad [\,,\,]^{(1)},\quad\ldots,\quad [\,,\,]^{(s)}=0
\]
on $\q$ such that, for each $0\le j<s$, $[\,,\,]^{(j+1)}$ is the polynomial limit of $(\q,[\,,\,]^{(j)})$ with respect to some family of $B_{\q}$-isometries.
In particular, each $(\q,[\,,\,]^{(j)},B_{\q})$ is a quadratic Lie algebra, and $(\q,[\,,\,]^{(s)},B_{\q})$ is a quadratic abelian Lie algebra.
\end{proposition}
\begin{proof}
By Lemma~\ref{lem:to-solvable}, one polynomial limit reduces $(\q,[\,,\,],B_{\q})$ to a solvable quadratic Lie algebra.
Hence it suffices to prove that every solvable quadratic Lie algebra $(\mathfrak a,[\,,\,],B_{\mathfrak a})$ admits a sequence of brackets satisfying the same conditions and ending with the zero bracket;
the required sequence is then obtained by placing this sequence after the polynomial limit of Lemma~\ref{lem:to-solvable}.
We proceed by induction on $\dim\mathfrak a$; if $\mathfrak a=0$, take $s=0$.

Suppose that $\mathfrak a\ne0$. By Lemma~\ref{lem:solvable-step}, there are a polynomial limit $[\,,\,]\rightsquigarrow[\,,\,]_0$
and an orthogonal direct sum $\mathfrak a=\mathfrak c\perp\mathfrak z$ with respect to the quadratic Lie algebra $(\mathfrak a,[\,,\,]_0,B_{\mathfrak a})$,
where $\mathfrak z\ne0$ is contained in the center of $(\mathfrak a,[\,,\,]_0)$ and $\mathfrak c$ is a solvable quadratic Lie algebra of smaller dimension.
By the induction hypothesis, there exist $m\in\Z_{\ge1}$ and a sequence of polynomial limits on $\mathfrak c$
\[
 [\,,\,]_0|_{\mathfrak c\times\mathfrak c}=[\,,\,]_{\mathfrak c}^{(1)}\rightsquigarrow[\,,\,]_{\mathfrak c}^{(2)}\rightsquigarrow\cdots\rightsquigarrow[\,,\,]_{\mathfrak c}^{(m)}=0.
\]
Write $[\,,\,]_{\mathfrak a}^{(0)}=[\,,\,]$ and $[\,,\,]_{\mathfrak a}^{(1)}=[\,,\,]_0$. Applying Lemma~\ref{lem:central-extension} successively, we obtain a sequence of polynomial limits on $\mathfrak a$
$[\,,\,]_{\mathfrak a}^{(1)}\rightsquigarrow[\,,\,]_{\mathfrak a}^{(2)}\rightsquigarrow\cdots\rightsquigarrow[\,,\,]_{\mathfrak a}^{(m)}$,
where, for $1\le i\le m$, $[\,,\,]_{\mathfrak a}^{(i)}$ coincides with $[\,,\,]_{\mathfrak c}^{(i)}$ on $\mathfrak c$, and $\mathfrak z$ is contained in the center of $(\mathfrak a,[\,,\,]_{\mathfrak a}^{(i)})$.
Since $[\,,\,]_{\mathfrak c}^{(m)}=0$ and $\mathfrak z$ is contained in the center of $(\mathfrak a,[\,,\,]_{\mathfrak a}^{(m)})$, we get $[\,,\,]_{\mathfrak a}^{(m)}=0$. Thus
\[
 [\,,\,]=[\,,\,]_{\mathfrak a}^{(0)}\rightsquigarrow[\,,\,]_{\mathfrak a}^{(1)}\rightsquigarrow\cdots\rightsquigarrow[\,,\,]_{\mathfrak a}^{(m)}=0
\]
is the required sequence of brackets.
By induction on $j$, Lemma~\ref{lem:polynomial-limit} shows that each $(\q,[\,,\,]^{(j)},B_{\q})$ is a quadratic Lie algebra.
\end{proof}

\Needspace{3\baselineskip}
\subsection{Contravariant forms and their radicals}\label{subsec:contravariant}
Throughout this subsection, $(\mathfrak a,B_{\mathfrak a})$ is a quadratic Lie algebra.
In order to apply Proposition~\ref{prop:degenerate} to the induced vacuum module, we construct the normalized contravariant form on ${V_{(\mathfrak a,B_{\mathfrak a})}(1,0)}$,
and prove that its radical is precisely the maximal proper graded vertex ideal $J_{(\mathfrak a,B_{\mathfrak a})}(1,0)$; see Lemma~\ref{lem:radical}.
Hence ${V_{(\mathfrak a,B_{\mathfrak a})}(1,0)}$ has no nonzero proper graded vertex ideals if and only if this form is nondegenerate;
the latter will be verified degree by degree in Section~\ref{subsec:gram-affine}.
Since no conformal vector is specified here, the invariant bilinear forms of Section~\ref{sec:prelim} cannot be defined;
instead, we use the contravariant form determined by an anti-involution of $\widehat{\mathfrak a}_{B_{\mathfrak a}}$.

Analogously to the contravariant form on Verma modules for Kac--Moody algebras (see \cite[Section~9.4]{KacIDLA} and \cite[Section~3]{Roi}),
we introduce an anti-involution of $\widehat{\mathfrak a}_{B_{\mathfrak a}}$ and the corresponding contravariant forms.

\begin{definition}\label{def:contravariant}
Define the $\C$-linear map
\[
\begin{aligned}
 \tau:\widehat{\mathfrak a}_{B_{\mathfrak a}}&\longrightarrow\widehat{\mathfrak a}_{B_{\mathfrak a}},\\
 x(m)&\longmapsto -x(-m),\\
 K&\longmapsto K,
\end{aligned}
\]
where $x\in\mathfrak a$ and $m\in\Z$. By \eqref{eq:currentbracket} and the symmetry of $B_{\mathfrak a}$,
\[
 \tau([\xi_1,\xi_2])=[\tau(\xi_2),\tau(\xi_1)],\qquad \tau^2=\id,\qquad \xi_1,\xi_2\in\widehat{\mathfrak a}_{B_{\mathfrak a}},
\]
that is, $\tau$ is an anti-involution of $\widehat{\mathfrak a}_{B_{\mathfrak a}}$. Hence $\tau$ extends to $U(\widehat{\mathfrak a}_{B_{\mathfrak a}})$ by
\[
 \tau(1)=1,\qquad \tau(A_1A_2)=\tau(A_2)\tau(A_1),\qquad A_1,A_2\in U(\widehat{\mathfrak a}_{B_{\mathfrak a}}).
\]
A bilinear form $B$ on ${V_{(\mathfrak a,B_{\mathfrak a})}(1,0)}$ is said to be \emph{contravariant} if
\[
 B(Au,v)=B(u,\tau(A)v),\qquad A\in U(\widehat{\mathfrak a}_{B_{\mathfrak a}}),\ u,v\in{V_{(\mathfrak a,B_{\mathfrak a})}(1,0)};
\]
if, in addition, $B(\one,\one)=1$, then $B$ is called a \emph{normalized contravariant form}.
Since $\tau$ reverses the order of products, contravariance needs to be verified only for the generators
$x(m)$, $x\in\mathfrak a$, $m\in\Z$, and $K$ of $U(\widehat{\mathfrak a}_{B_{\mathfrak a}})$ as an associative algebra.
The condition for $K$ holds automatically, since $\tau(K)=K$ and $K$ acts on ${V_{(\mathfrak a,B_{\mathfrak a})}(1,0)}$ as $\id$;
moreover, since $\tau(x(m))=-x(-m)$, $B$ is a normalized contravariant form if and only if, for all $u,v\in{V_{(\mathfrak a,B_{\mathfrak a})}(1,0)}$,
\begin{equation}\label{eq:contracurrent}
 B(\one,\one)=1,\qquad
 B(x(m)u,v)=-B(u,x(-m)v),\qquad x\in\mathfrak a,\ m\in\Z.
\end{equation}
\end{definition}

Since $\widehat{\mathfrak a}_{B_{\mathfrak a}}$ is in general not a Kac--Moody algebra, \cite[Proposition~9.4]{KacIDLA} does not apply directly;
following the explicit construction of \cite[Section~9.4]{KacIDLA}, we prove that a normalized contravariant form exists and is unique.
To this end, for $u\in{V_{(\mathfrak a,B_{\mathfrak a})}(1,0)}$, denote by $\varepsilon(u)\in\C$ the scalar such that
\[
 u-\varepsilon(u)\one\in\bigoplus_{n\ge1}{V_{(\mathfrak a,B_{\mathfrak a})}(1,0)}_n,
\]
that is, the coefficient of the component of $u$ in ${V_{(\mathfrak a,B_{\mathfrak a})}(1,0)}_0=\C\one$; see \eqref{eq:affine-grading}.
It corresponds to the expectation value in \cite[Section~9.4]{KacIDLA}.

\begin{proposition}\label{prop:currentform}
There exists a unique normalized contravariant form on ${V_{(\mathfrak a,B_{\mathfrak a})}(1,0)}$, which we denote by $B_{V_{(\mathfrak a,B_{\mathfrak a})}(1,0)}$. It is symmetric, and
\[
 B_{V_{(\mathfrak a,B_{\mathfrak a})}(1,0)}\bigl({V_{(\mathfrak a,B_{\mathfrak a})}(1,0)}_n,{V_{(\mathfrak a,B_{\mathfrak a})}(1,0)}_{n'}\bigr)=0,\qquad n,n'\in\N,\quad n\ne n'.
\]
For all $A_1,A_2\in U(\widehat{\mathfrak a}_{B_{\mathfrak a}})$,
\[
 B_{V_{(\mathfrak a,B_{\mathfrak a})}(1,0)}(A_1\one,A_2\one)=\varepsilon\bigl(\tau(A_1)A_2\one\bigr).
\]
In particular, for $x,y\in\mathfrak a$ and $m,m'\in\Z_{\ge1}$,
\[
 B_{V_{(\mathfrak a,B_{\mathfrak a})}(1,0)}(x(-m)\one,y(-m')\one)=-m\delta_{m,m'}B_{\mathfrak a}(x,y).
\]
\end{proposition}
\begin{proof}
\smallskip\noindent\emph{(1) Existence.}
By \eqref{eq:currentvac},
\[
 {V_{(\mathfrak a,B_{\mathfrak a})}(1,0)}=
 U(\widehat{\mathfrak a}_{B_{\mathfrak a}})\Big/
 \Bigl(\sum_{\substack{x\in\mathfrak a\\m\in\Z_{\ge0}}}
 U(\widehat{\mathfrak a}_{B_{\mathfrak a}})x(m)
 +U(\widehat{\mathfrak a}_{B_{\mathfrak a}})(K-1)\Bigr),
\]
where the left ideal used in the quotient is generated by the elements $x(m)$, $x\in\mathfrak a$, $m\in\Z_{\ge0}$, and $K-1$.
Hence every $u\in{V_{(\mathfrak a,B_{\mathfrak a})}(1,0)}$ can be written as $u=A\one$ with $A\in U(\widehat{\mathfrak a}_{B_{\mathfrak a}})$; such an $A$ is called a \emph{representative} of $u$.
Moreover, for $A,A'\in U(\widehat{\mathfrak a}_{B_{\mathfrak a}})$, we have $A\one=A'\one$ if and only if $A-A'$ lies in the above left ideal.
For $u,v\in{V_{(\mathfrak a,B_{\mathfrak a})}(1,0)}$, choose representatives $A_1,A_2$ with $u=A_1\one$ and $v=A_2\one$, and set
\[
 B_{V_{(\mathfrak a,B_{\mathfrak a})}(1,0)}(u,v)=\varepsilon\bigl(\tau(A_1)A_2\one\bigr).
\]
We must verify that the right-hand side does not depend on the choice of representatives. The right-hand side depends on $A_2$ only through $v=A_2\one$, so it is independent of the choice of $A_2$.
As for $A_1$, if $A_1'$ is another representative of $u$, then $A_1-A_1'$ lies in the above left ideal,
that is, it is a sum of elements of the form $A\,x(m)$, $x\in\mathfrak a$, $m\in\Z_{\ge0}$, and $A(K-1)$, where $A\in U(\widehat{\mathfrak a}_{B_{\mathfrak a}})$, and
\[
 \tau\bigl(A\,x(m)\bigr)=-x(-m)\tau(A),\qquad \tau\bigl(A(K-1)\bigr)=(K-1)\tau(A).
\]
Moreover, $x(m)$ raises the degree for $m<0$, while $x(0)$ preserves the degree and annihilates $\one$, so
\[
 \varepsilon(x(m)w)=0,\quad m\in\Z_{\le0},\qquad
 \varepsilon((K-1)w)=0,\qquad w\in{V_{(\mathfrak a,B_{\mathfrak a})}(1,0)}.
\]
Hence $\varepsilon\bigl(\tau(A_1-A_1')w\bigr)=0$, and the form $B_{V_{(\mathfrak a,B_{\mathfrak a})}(1,0)}$ so defined is also independent of the choice of $A_1$.
Taking $A_1=A_2=1$, we get $B_{V_{(\mathfrak a,B_{\mathfrak a})}(1,0)}(\one,\one)=\varepsilon(\one)=1$ from $\tau(1)=1$.
For $A\in U(\widehat{\mathfrak a}_{B_{\mathfrak a}})$, the elements $AA_1$ and $\tau(A)A_2$ are representatives of $Au$ and $\tau(A)v$, respectively, so by $\tau(AA_1)=\tau(A_1)\tau(A)$,
\[
 B_{V_{(\mathfrak a,B_{\mathfrak a})}(1,0)}(Au,v)=\varepsilon\bigl(\tau(A_1)\tau(A)A_2\one\bigr)=B_{V_{(\mathfrak a,B_{\mathfrak a})}(1,0)}(u,\tau(A)v).
\]
Bilinearity follows from the linearity of $\tau$ and $\varepsilon$. Therefore $B_{V_{(\mathfrak a,B_{\mathfrak a})}(1,0)}$ is a normalized contravariant form,
and by definition it satisfies the formula $B_{V_{(\mathfrak a,B_{\mathfrak a})}(1,0)}(A_1\one,A_2\one)=\varepsilon(\tau(A_1)A_2\one)$ of the proposition.

\smallskip\noindent\emph{(2) Uniqueness.}
Let $B'$ be any normalized contravariant form on ${V_{(\mathfrak a,B_{\mathfrak a})}(1,0)}$.
For $x\in\mathfrak a$, $m\in\Z_{\ge1}$, and $w\in{V_{(\mathfrak a,B_{\mathfrak a})}(1,0)}$, since $x(m)\one=0$,
\[
 B'(\one,x(-m)w)=-B'(x(m)\one,w)=0.
\]
By Proposition~\ref{prop:affine-vacuum-structure}, $\bigoplus_{n\ge1}{V_{(\mathfrak a,B_{\mathfrak a})}(1,0)}_n$ is spanned by vectors of the form $x(-m)w'$ with $x\in\mathfrak a$, $m\in\Z_{\ge1}$, and $w'\in{V_{(\mathfrak a,B_{\mathfrak a})}(1,0)}$,
so by the above identity it is orthogonal to $\one$ with respect to $B'$. Writing $w=\varepsilon(w)\one+w_+$ with $w_+\in\bigoplus_{n\ge1}{V_{(\mathfrak a,B_{\mathfrak a})}(1,0)}_n$,
we then obtain $B'(\one,w)=\varepsilon(w)$ from $B'(\one,\one)=1$.
Now, for $u,v\in{V_{(\mathfrak a,B_{\mathfrak a})}(1,0)}$, choose representatives $A_1,A_2$ with $u=A_1\one$ and $v=A_2\one$. By the contravariance of $B'$ (with $A=A_1$) and the identity $B'(\one,w)=\varepsilon(w)$,
\[
 B'(u,v)=B'(A_1\one,A_2\one)=B'\bigl(\one,\tau(A_1)A_2\one\bigr)=\varepsilon\bigl(\tau(A_1)A_2\one\bigr)=B_{V_{(\mathfrak a,B_{\mathfrak a})}(1,0)}(u,v).
\]
Hence $B'=B_{V_{(\mathfrak a,B_{\mathfrak a})}(1,0)}$, which proves uniqueness.

\smallskip\noindent\emph{(3) Symmetry.}
Let $B_{V_{(\mathfrak a,B_{\mathfrak a})}(1,0)}^{\mathsf T}(u,v)=B_{V_{(\mathfrak a,B_{\mathfrak a})}(1,0)}(v,u)$. By $\tau^2=\id$, for $A\in U(\widehat{\mathfrak a}_{B_{\mathfrak a}})$ we have
\[
 B_{V_{(\mathfrak a,B_{\mathfrak a})}(1,0)}^{\mathsf T}(Au,v)=B_{V_{(\mathfrak a,B_{\mathfrak a})}(1,0)}(v,Au)=B_{V_{(\mathfrak a,B_{\mathfrak a})}(1,0)}(\tau(A)v,u)=B_{V_{(\mathfrak a,B_{\mathfrak a})}(1,0)}^{\mathsf T}(u,\tau(A)v),
\]
and $B_{V_{(\mathfrak a,B_{\mathfrak a})}(1,0)}^{\mathsf T}(\one,\one)=1$, so $B_{V_{(\mathfrak a,B_{\mathfrak a})}(1,0)}^{\mathsf T}$ is also a normalized contravariant form.
By uniqueness, $B_{V_{(\mathfrak a,B_{\mathfrak a})}(1,0)}^{\mathsf T}=B_{V_{(\mathfrak a,B_{\mathfrak a})}(1,0)}$, that is, $B_{V_{(\mathfrak a,B_{\mathfrak a})}(1,0)}$ is symmetric.

\smallskip\noindent\emph{(4) Orthogonality of different degrees.}
Let $P,Q$ be products of operators of the form $x(-m)$, $x\in\mathfrak a$, $m\in\Z_{\ge1}$, with $P\one\in{V_{(\mathfrak a,B_{\mathfrak a})}(1,0)}_n$ and $Q\one\in{V_{(\mathfrak a,B_{\mathfrak a})}(1,0)}_{n'}$, where $n\ne n'$.
Since $\tau(P)$ lowers the degree by $n$,
\[
 B_{V_{(\mathfrak a,B_{\mathfrak a})}(1,0)}(P\one,Q\one)=\varepsilon\bigl(\tau(P)Q\one\bigr)=0.
\]
By Proposition~\ref{prop:affine-vacuum-structure}, each ${V_{(\mathfrak a,B_{\mathfrak a})}(1,0)}_n$ is spanned by such vectors $P\one$, so different degrees are orthogonal.

\smallskip\noindent\emph{(5) The basic pairing.}
Let $x,y\in\mathfrak a$ and $m,m'\in\Z_{\ge1}$. By $\tau(x(-m))=-x(m)$, $x(m)\one=0$, and \eqref{eq:currentbracket},
\[
 \begin{aligned}
 B_{V_{(\mathfrak a,B_{\mathfrak a})}(1,0)}(x(-m)\one,y(-m')\one)
 &=-\varepsilon\bigl(x(m)y(-m')\one\bigr)\\
 &=-\varepsilon\bigl([x,y](m-m')\one\bigr)-m\delta_{m,m'}B_{\mathfrak a}(x,y)\\
 &=-m\delta_{m,m'}B_{\mathfrak a}(x,y),
 \end{aligned}
\]
where the last step uses $\varepsilon([x,y](m-m')\one)=0$: for $m\ge m'$ we have $[x,y](m-m')\one=0$, and for $m<m'$ this vector lies in degree $m'-m\ge1$.
\end{proof}

We now determine the radical
\[
 \rad B_{V_{(\mathfrak a,B_{\mathfrak a})}(1,0)}=\{u\in{V_{(\mathfrak a,B_{\mathfrak a})}(1,0)}:B_{V_{(\mathfrak a,B_{\mathfrak a})}(1,0)}(u,v)=0
              \text{ for all }v\in{V_{(\mathfrak a,B_{\mathfrak a})}(1,0)}\}
\]
of $B_{V_{(\mathfrak a,B_{\mathfrak a})}(1,0)}$.
We prove that it is precisely the maximal proper graded vertex ideal $J_{(\mathfrak a,B_{\mathfrak a})}(1,0)$ of Definition~\ref{def:affine-simple-quotient}. This is analogous to the case of Verma modules for Kac--Moody algebras;
see \cite[Proposition~9.2~(c), Sections~9.3--9.4]{KacIDLA}.

\begin{lemma}\label{lem:radical}
$\rad B_{V_{(\mathfrak a,B_{\mathfrak a})}(1,0)}=J_{(\mathfrak a,B_{\mathfrak a})}(1,0)$. Hence
\[
 {V_{(\mathfrak a,B_{\mathfrak a})}(1,0)}/\rad B_{V_{(\mathfrak a,B_{\mathfrak a})}(1,0)}=L_{(\mathfrak a,B_{\mathfrak a})}(1,0).
\]
\end{lemma}
\begin{proof}
We first prove that $J_{(\mathfrak a,B_{\mathfrak a})}(1,0)\subseteq\rad B_{V_{(\mathfrak a,B_{\mathfrak a})}(1,0)}$. Since $J_{(\mathfrak a,B_{\mathfrak a})}(1,0)$ is a proper graded vertex ideal, it is stable under all $x(m)=x_m$, and $J_{(\mathfrak a,B_{\mathfrak a})}(1,0)\cap{V_{(\mathfrak a,B_{\mathfrak a})}(1,0)}_0=0$,
so $\varepsilon\bigl(J_{(\mathfrak a,B_{\mathfrak a})}(1,0)\bigr)=0$.
For $u\in J_{(\mathfrak a,B_{\mathfrak a})}(1,0)$ and any product $P$ of operators of the form $x(-m)$, $x\in\mathfrak a$, $m\in\Z_{\ge1}$, Proposition~\ref{prop:currentform} gives
\[
 B_{V_{(\mathfrak a,B_{\mathfrak a})}(1,0)}(P\one,u)=\varepsilon(\tau(P)u)=0,
\]
since $\tau(P)u\in J_{(\mathfrak a,B_{\mathfrak a})}(1,0)$. By Proposition~\ref{prop:affine-vacuum-structure}, these vectors $P\one$ span ${V_{(\mathfrak a,B_{\mathfrak a})}(1,0)}$;
by symmetry, $J_{(\mathfrak a,B_{\mathfrak a})}(1,0)\subseteq\rad B_{V_{(\mathfrak a,B_{\mathfrak a})}(1,0)}$.
The above argument uses only the fact that $J_{(\mathfrak a,B_{\mathfrak a})}(1,0)$ is a proper graded subspace stable under all $x(m)$, $x\in\mathfrak a$, $m\in\Z$,
so every such subspace is contained in $\rad B_{V_{(\mathfrak a,B_{\mathfrak a})}(1,0)}$.

It remains to show that $\rad B_{V_{(\mathfrak a,B_{\mathfrak a})}(1,0)}$ is a proper graded vertex ideal; the assertion then follows from the maximality of $J_{(\mathfrak a,B_{\mathfrak a})}(1,0)$ in Definition~\ref{def:affine-simple-quotient}.
By contravariance, $\rad B_{V_{(\mathfrak a,B_{\mathfrak a})}(1,0)}$ is stable under all $x(m)$;
by the orthogonality of different degrees, the same argument as in Remark~\ref{rem:weight-gram} gives
\[
 \rad B_{V_{(\mathfrak a,B_{\mathfrak a})}(1,0)}=\bigoplus_{n\in\N}\rad\Bigl(B_{V_{(\mathfrak a,B_{\mathfrak a})}(1,0)}\big|_{{V_{(\mathfrak a,B_{\mathfrak a})}(1,0)}_n\times{V_{(\mathfrak a,B_{\mathfrak a})}(1,0)}_n}\Bigr),
\]
so it is graded with respect to the natural grading \eqref{eq:affine-grading}; since $B_{V_{(\mathfrak a,B_{\mathfrak a})}(1,0)}(\one,\one)=1$, we have $\one\notin\rad B_{V_{(\mathfrak a,B_{\mathfrak a})}(1,0)}$, so it is proper.
Thus, by the preceding paragraph, it is the largest proper graded subspace stable under all $x(m)$.
By Proposition~\ref{prop:affine-restricted} (see \cite[Proposition~4.12]{Lian}), it is stable under all modes $v_j$, $v\in{V_{(\mathfrak a,B_{\mathfrak a})}(1,0)}$, $j\in\Z$;
by its maximality and \cite[Proposition~3.2, Corollary~3.3]{Lian}, it is also stable under the translation operator $D$, and is therefore a proper graded vertex ideal.
\end{proof}

\Needspace{3\baselineskip}
\subsection{Gram determinants of the affine vertex algebras}\label{subsec:gram-affine}
In this subsection we prove Theorem~\ref{thm:currentsimple}: for a quadratic Lie algebra $(\q,B_{\q})$ without nonzero semisimple ideals,
the normalized contravariant form on ${V_{(\q,B_{\q})}(1,0)}$ is nondegenerate, and hence $J_{(\q,B_{\q})}(1,0)=0$ by Lemma~\ref{lem:radical}.
As in Sections~\ref{subsec:det-criterion}--\ref{subsec:pbw-simplicity}, we pull this form back to $W_n^{\mathfrak{Lie}}$ and compare Gram determinants with an auxiliary bilinear form;
the difference is that here the Gram determinant of the form in each degree remains unchanged along each polynomial limit of Proposition~\ref{prop:degenerate},
while the final Lie bracket is zero, and for it the Gram determinant can be computed directly and is nonzero.
In order to compare these determinants as the Lie bracket varies, we first use PBW symmetrization
to identify the space in each degree with a space $W_n^{\mathfrak{Lie}}$ depending only on the underlying vector space,
and then define on it the bilinear form $G_n^{\mathfrak{Lie}}$ pulled back from this form, together with the bilinear form $H_n^{\mathfrak{Lie}}$ corresponding to the case of zero bracket.
From here through Lemma~\ref{lem:heisenbergpairing}, let $(\mathfrak a,B_{\mathfrak a})$ be a quadratic Lie algebra.

According to the natural grading in Proposition~\ref{prop:affine-vacuum-structure}, we define the following spaces as in Definition~\ref{def:weight-space}.
\begin{definition}\label{def:lie-weight-space}
For $m\in\Z_{\ge1}$, let $\mathfrak a_m$ be a copy of the vector space $\mathfrak a$, assigned degree $m$, and identify the elements of the two spaces via $\mathfrak a_m\simeq\mathfrak a$.
Fix $n\in\N$. A sequence $(r_m)_{m\ge1}$ with $r_m\in\N$ and $\sum_{m\ge1}mr_m=n$ is called \emph{occupancy data} of degree $n$. Define
\[
 W_n^{\mathfrak{Lie}}=\bigoplus_{\sum_{m=1}^{\infty}mr_m=n}\ \bigotimes_{m=1}^{\infty}\Sym^{r_m}\mathfrak a_m;
\]
the summand corresponding to given occupancy data is called an \emph{occupancy summand}.
As in Definition~\ref{def:weight-space}, there are only finitely many summands; omitting the factors $\Sym^0\mathfrak a_m=\C$ and taking the empty tensor product to be $\C$,
we have $W_0^{\mathfrak{Lie}}=\C$. From now on we use the notation
\[
 \prod_{m\ge1}\prod_{i=1}^{r_m}x_{m,i}
 :=\bigotimes_{m\ge1}(x_{m,1}\cdots x_{m,r_m}),\qquad x_{m,i}\in\mathfrak a_m,
\]
where each parenthesized expression is a symmetric product; such a vector is called a \emph{symmetric monomial}.
\end{definition}

As in Definition~\ref{def:auxiliary-pairing}, we now define a bilinear form $H_n^{\mathfrak{Lie}}$ on $W_n^{\mathfrak{Lie}}$;
Lemma~\ref{lem:heisenbergpairing} will show that it corresponds to the case where the Lie bracket is zero.
\begin{definition}\label{def:lie-auxiliary-pairing}
For $m\in\Z_{\ge1}$ and $r\in\N$, define a bilinear form $B_{\Sym^r\mathfrak a_m}$ on $\Sym^r\mathfrak a_m$ as follows: for $x_i,y_i\in\mathfrak a_m$, $1\le i\le r$, set
\[
 B_{\Sym^r\mathfrak a_m}(x_1\cdots x_r,y_1\cdots y_r)
 =(-m)^r\sum_{\sigma\in S_r}\prod_{i=1}^rB_{\mathfrak a}(x_i,y_{\sigma(i)});
\]
for $r=0$, define $B_{\Sym^0\mathfrak a_m}(1,1)=1$.
For $n\in\N$, we define the bilinear form $H_n^{\mathfrak{Lie}}$ on $W_n^{\mathfrak{Lie}}$ as follows: distinct occupancy summands are mutually orthogonal,
and on a given occupancy summand we take the tensor product of the forms $B_{\Sym^{r_m}\mathfrak a_m}$. That is, for $x_{m,i},y_{m,i}\in\mathfrak a_m$,
\begin{equation}\label{eq:currentheisenberg}
 H_n^{\mathfrak{Lie}}\!\left(\prod_{m\ge1}\prod_{i=1}^{r_m}x_{m,i},\prod_{m\ge1}\prod_{i=1}^{r_m}y_{m,i}\right)
 =\prod_{m\ge1}\left((-m)^{r_m}\sum_{\sigma\in S_{r_m}}\prod_{i=1}^{r_m}B_{\mathfrak a}(x_{m,i},y_{m,\sigma(i)})\right).
\end{equation}
This formula is extended bilinearly, and a factor with $r_m=0$ is taken to be $1$; in particular, $H_0^{\mathfrak{Lie}}(1,1)=1$.
\end{definition}

As in Lemma~\ref{lem:symmetric-pairing}, $H_n^{\mathfrak{Lie}}$ is nondegenerate.
\begin{lemma}\label{lem:lie-auxiliary-nondegenerate}
For each $n\in\N$, the form $H_n^{\mathfrak{Lie}}$ is symmetric and nondegenerate.
\end{lemma}
\begin{proof}
Since $B_{\mathfrak a}$ is symmetric and nondegenerate and $-m\ne0$ for $m\in\Z_{\ge1}$,
it suffices to replace $\eta$ and $\alpha_m$ by $B_{\mathfrak a}$ and $-m$, respectively, in the proof of Lemma~\ref{lem:symmetric-pairing}.
\end{proof}

As in Definition~\ref{def:frob-pbw-pairing}, we now use PBW symmetrization to pull back the normalized contravariant form on ${V_{(\mathfrak a,B_{\mathfrak a})}(1,0)}$ to $W_n^{\mathfrak{Lie}}$.
\begin{definition}\label{def:currentgram}
Fix $n\in\N$. Let $\ell\in\N$, $m_j\in\Z_{\ge1}$, and $x_j\in\mathfrak a_{m_j}$ for $1\le j\le\ell$, with $\sum_{j=1}^{\ell}m_j=n$.
As in Definition~\ref{def:frob-pbw-pairing}, grouping the inputs according to the spaces $\mathfrak a_m$ to which they belong, we regard the symmetric product $x_1\cdots x_\ell$ as a symmetric monomial in $W_n^{\mathfrak{Lie}}$. Define the linear map
\[
\begin{aligned}
 \Phi_n^{\mathfrak{Lie}}:W_n^{\mathfrak{Lie}}&\longrightarrow{V_{(\mathfrak a,B_{\mathfrak a})}(1,0)}_n,\\
 x_1\cdots x_\ell&\longmapsto\frac1{\ell!}\sum_{\sigma\in S_\ell}
 x_{\sigma(1)}(-m_{\sigma(1)})\cdots x_{\sigma(\ell)}(-m_{\sigma(\ell)})\one,
\end{aligned}
\]
where each $x_j(-m_j)$ is understood via $\mathfrak a_{m_j}\simeq\mathfrak a$; in particular, $\Phi_0^{\mathfrak{Lie}}(1)=\one$.
Using $\Phi_n^{\mathfrak{Lie}}$, we define a bilinear form on $W_n^{\mathfrak{Lie}}$ by
\[
 G_n^{\mathfrak{Lie}}(u,v)=B_{V_{(\mathfrak a,B_{\mathfrak a})}(1,0)}\bigl(\Phi_n^{\mathfrak{Lie}}(u),\Phi_n^{\mathfrak{Lie}}(v)\bigr),\qquad u,v\in W_n^{\mathfrak{Lie}}.
\]
\end{definition}

As in the discussion following Definition~\ref{def:frob-pbw-pairing}, $\Phi_n^{\mathfrak{Lie}}$ is well defined;
by the PBW symmetrization of \cite[Corollary~13.7]{Et} and Proposition~\ref{prop:affine-vacuum-structure}, each $\Phi_n^{\mathfrak{Lie}}$ is a vector space isomorphism.

When the Lie bracket of $\mathfrak a$ is zero, that is, when $\mathfrak a$ is an abelian Lie algebra, ${V_{(\mathfrak a,B_{\mathfrak a})}(1,0)}$ is a Heisenberg vertex algebra (see \cite[Section~6.3]{LL}), and the pulled-back form is precisely $H_n^{\mathfrak{Lie}}$.
\begin{lemma}\label{lem:heisenbergpairing}
If $\mathfrak a$ is an abelian Lie algebra, then $G_n^{\mathfrak{Lie}}=H_n^{\mathfrak{Lie}}$ for each $n\in\N$.
\end{lemma}
\begin{proof}
In this case, \eqref{eq:currentbracket} becomes
\[
 [x(m),y(m')]=m\delta_{m+m',0}B_{\mathfrak a}(x,y)K,\qquad x,y\in\mathfrak a,\ m,m'\in\Z.
\]
In particular, the operators $x(-m)$, $x\in\mathfrak a$, $m\in\Z_{\ge1}$, commute with one another.
Let
\[
 u=\prod_{m\ge1}\prod_{i=1}^{r_m}x_{m,i},\qquad v=\prod_{m\ge1}\prod_{i=1}^{r'_m}y_{m,i}
\]
be symmetric monomials, where $(r_m)_{m\ge1}$ and $(r'_m)_{m\ge1}$ are occupancy data of degree $n$, and let
\[
 P=\prod_{m\ge1}\prod_{i=1}^{r_m}x_{m,i}(-m),\qquad Q=\prod_{m\ge1}\prod_{i=1}^{r'_m}y_{m,i}(-m).
\]
Then $\Phi_n^{\mathfrak{Lie}}(u)=P\one$ and $\Phi_n^{\mathfrak{Lie}}(v)=Q\one$, so by Proposition~\ref{prop:currentform},
\[
 G_n^{\mathfrak{Lie}}(u,v)=\varepsilon\bigl(\tau(P)Q\one\bigr),\qquad
 \tau(P)=\prod_{m\ge1}\prod_{i=1}^{r_m}\bigl(-x_{m,i}(m)\bigr).
\]
Move the operators $x_{m,i}(m)$ in $\tau(P)$ one by one to the right of $Q$: by the above relation, $x_{m,i}(m)$ has a nonzero commutator only with the operators $y_{m,j}(-m)$ of the same degree, namely
$mB_{\mathfrak a}(x_{m,i},y_{m,j})$, and after being moved to the far right it annihilates $\one$.
Hence only the terms in which the factors of $\tau(P)$ and $Q$ are all paired off according to degree contribute to $\varepsilon(\tau(P)Q\one)$:
if a positive mode remains unpaired, it annihilates $\one$ after being moved to the far right; if a negative mode remains unpaired, the resulting vector has positive degree and vanishes under $\varepsilon$.
Each pair contributes $-mB_{\mathfrak a}(x_{m,i},y_{m,j})$.
If $(r_m)_{m\ge1}\ne(r'_m)_{m\ge1}$, then in some degree the two sides have different numbers of factors, so they cannot all be paired off, and the above expression is zero;
if the two occupancy data coincide, the pairings in each degree $m$ are given exactly by the permutations $\sigma\in S_{r_m}$, which yields the right-hand side of \eqref{eq:currentheisenberg}.
By bilinearity, $G_n^{\mathfrak{Lie}}=H_n^{\mathfrak{Lie}}$.
\end{proof}

For the four-dimensional Nappi--Witten Lie algebra, \cite[Theorem~3.2]{BJP} proves that the corresponding induced vacuum module is irreducible at nonzero levels;
the following theorem proves, for all quadratic Lie algebras without nonzero semisimple ideals, that the corresponding induced vacuum module has no nonzero proper graded vertex ideals.
Fix a basis of $W_n^{\mathfrak{Lie}}$; as in Lemma~\ref{lem:detcriterion}, the determinant of each bilinear form below means its Gram determinant with respect to this basis.

\begin{theorem}\label{thm:currentsimple}
Let $(\q,B_{\q})$ be a quadratic Lie algebra without nonzero semisimple ideals, and take $(\mathfrak a,B_{\mathfrak a})=(\q,B_{\q})$ in the above definitions.
Then for each $n\in\N$,
\[
 \det G_n^{\mathfrak{Lie}}=\det H_n^{\mathfrak{Lie}}\ne0.
\]
Consequently, $\rad B_{V_{(\q,B_{\q})}(1,0)}=0$; hence
\[
 J_{(\q,B_{\q})}(1,0)=0,\quad\text{that is,}\quad L_{(\q,B_{\q})}(1,0)={V_{(\q,B_{\q})}(1,0)}.
\]
\end{theorem}
\begin{proof}
For each $n\in\N$, the proof proceeds in four steps. Step (1) shows that, with respect to a fixed basis, the entries of the Gram matrix of $G_n^{\mathfrak{Lie}}$ are polynomials in the structure constants of $\q$ and the values of $B_{\q}$,
and that these polynomials do not depend on the Lie bracket; step (2) uses this to prove that $\det G_n^{\mathfrak{Lie}}$ remains unchanged under a single polynomial limit;
step (3) applies step (2) repeatedly along the sequence of Lie brackets in Proposition~\ref{prop:degenerate}, and in the final abelian case obtains $\det G_n^{\mathfrak{Lie}}=\det H_n^{\mathfrak{Lie}}\ne0$ from
Lemmas~\ref{lem:heisenbergpairing} and~\ref{lem:lie-auxiliary-nondegenerate};
step (4) deduces from this that $\rad B_{V_{(\q,B_{\q})}(1,0)}=0$ and $J_{(\q,B_{\q})}(1,0)=0$.

\smallskip\noindent\emph{(1) Polynomial expressions for the Gram matrix entries.}
Let $d=\dim\q$, fix a basis $e_1,\ldots,e_d$ of $\q$, and write
\[
 [e_i,e_j]=\sum_{k=1}^d\mathsf c_{ij}^{k}e_k,\qquad
 \mathsf b_{ij}=B_{\q}(e_i,e_j),\qquad 1\le i,j\le d.
\]
For products $P,Q$ of operators of the form $e_i(-m)$, $1\le i\le d$, $m\in\Z_{\ge1}$, Proposition~\ref{prop:currentform} gives
\[
 B_{V_{(\q,B_{\q})}(1,0)}(P\one,Q\one)=\varepsilon\bigl(\tau(P)Q\one\bigr).
\]
By \eqref{eq:currentbracket}, on ${V_{(\q,B_{\q})}(1,0)}$ we have
\[
 e_i(p)e_j(q)
 =e_j(q)e_i(p)
 +\sum_{k=1}^d\mathsf c_{ij}^{k}e_k(p+q)
 +p\delta_{p+q,0}\mathsf b_{ij}\id
 ,\qquad p,q\in\Z.
\]
We reorder $\tau(P)Q$ by the following rule: at each step, take the leftmost pair of adjacent modes $e_i(p)e_j(q)$ with $p\ge0>q$, and replace it using the above identity.
For a product of modes, let $N$ be the number of modes in it, and let $I$ be the number of pairs of modes, not necessarily adjacent, in which a nonnegative mode stands to the left of a negative mode.
Each replacement turns a product into terms whose values of $(N,I)$ are strictly smaller in the lexicographic order, so the reordering terminates after finitely many steps.
When the reordering is complete, terms containing nonnegative modes annihilate $\one$, while a nonempty product of negative modes applied to $\one$ lies in $\bigoplus_{n'\ge1}{V_{(\q,B_{\q})}(1,0)}_{n'}$; both vanish under $\varepsilon$,
and only the scalar terms contribute. The reordering rule depends only on the basis indices and the mode indices of the modes in $P$ and $Q$, and each replacement only multiplies the coefficient by some $\mathsf c_{ij}^{k}$ or some $p\,\mathsf b_{ij}$.
Hence there exists a polynomial $\mathsf f_{P,Q}$ with integer coefficients in $d^3+d^2$ variables, depending only on $P$ and $Q$ and not on the Lie bracket, such that
\[
 \varepsilon\bigl(\tau(P)Q\one\bigr)=\mathsf f_{P,Q}\bigl((\mathsf c_{ij}^{k}),(\mathsf b_{ij})\bigr).
\]

Let $d_n=\dim W_n^{\mathfrak{Lie}}$, and take the basis $w_1,\ldots,w_{d_n}$ consisting of the symmetric monomials formed from the copies of the $e_i$ in the spaces $\q_m$.
A change of basis multiplies $\det G_n^{\mathfrak{Lie}}$ and $\det H_n^{\mathfrak{Lie}}$ by the same nonzero factor, so it suffices to prove the assertion for this basis, which we use from now on.
For $1\le\alpha\le d_n$, write $w_\alpha=e_{i_1}\cdots e_{i_{\ell_\alpha}}$, where $e_{i_j}\in\q_{m_j}$ for $1\le j\le\ell_\alpha$.
By the symmetrization formula for $\Phi_n^{\mathfrak{Lie}}$ in Definition~\ref{def:currentgram},
\[
 \Phi_n^{\mathfrak{Lie}}(w_\alpha)=\frac1{\ell_\alpha!}\sum_{\sigma\in S_{\ell_\alpha}}P_{\alpha,\sigma}\one,\qquad
 P_{\alpha,\sigma}=e_{i_{\sigma(1)}}(-m_{\sigma(1)})\cdots e_{i_{\sigma(\ell_\alpha)}}(-m_{\sigma(\ell_\alpha)}).
\]
Let
\[
 \mathsf F_{\alpha\beta}=\frac1{\ell_\alpha!\,\ell_\beta!}\sum_{\sigma\in S_{\ell_\alpha}}\sum_{\sigma'\in S_{\ell_\beta}}\mathsf f_{P_{\alpha,\sigma},P_{\beta,\sigma'}},
 \qquad 1\le\alpha,\beta\le d_n;
\]
this is a polynomial with rational coefficients in $d^3+d^2$ variables. By bilinearity,
\[
 G_n^{\mathfrak{Lie}}(w_\alpha,w_\beta)
 =\mathsf F_{\alpha\beta}\bigl((\mathsf c_{ij}^{k}),(\mathsf b_{ij})\bigr),\qquad 1\le\alpha,\beta\le d_n.
\]

\smallskip\noindent\emph{(2) Invariance of the Gram determinant under a single polynomial limit.}
Let $[\,,\,]_0$ be the polynomial limit of the bracket of $(\q,[\,,\,],B_{\q})$ with respect to $B_{\q}$-isometries $(\psi_\lambda)_{\lambda\in\C\setminus\{0\}}$, and let $[\,,\,]_\lambda$ be as in \eqref{eq:pullback}.
For $\lambda\in\C$, Lemma~\ref{lem:polynomial-limit} shows that $(\q,[\,,\,]_\lambda,B_{\q})$ is a quadratic Lie algebra.
Denote by $\mathcal V_\lambda$ the induced vacuum module of level $1$ constructed for $(\q,[\,,\,]_\lambda,B_{\q})$ as in Definition~\ref{def:affine-vacuum},
by $B_{\mathcal V_\lambda}$ the normalized contravariant form on it,
and by $\Phi_{n,\lambda}^{\mathfrak{Lie}}$ and $G_{n,\lambda}^{\mathfrak{Lie}}$ the maps $\Phi_n^{\mathfrak{Lie}}$ and the forms $G_n^{\mathfrak{Lie}}$ defined for $(\q,[\,,\,]_\lambda,B_{\q})$.

For $\lambda\in\C$ with $\lambda\ne0$, Lemma~\ref{lem:polynomial-limit} shows that $\psi_\lambda$ gives an isomorphism of quadratic Lie algebras
$(\q,[\,,\,]_\lambda,B_{\q})\to(\q,[\,,\,],B_{\q})$.
It induces an isomorphism $x(m)\mapsto(\psi_\lambda x)(m)$, $K\mapsto K$, between the two affine Lie algebras,
and hence a linear isomorphism $F_\lambda:\mathcal V_\lambda\to{V_{(\q,B_{\q})}(1,0)}$ satisfying
\[
 F_\lambda\one=\one,\qquad
 F_\lambda(x(m)u)=(\psi_\lambda x)(m)F_\lambda u
 ,\qquad x\in\q,\ m\in\Z,\ u\in\mathcal V_\lambda.
\]
Thus $(u,v)\mapsto B_{V_{(\q,B_{\q})}(1,0)}(F_\lambda u,F_\lambda v)$ is a bilinear form on $\mathcal V_\lambda$, and
\[
 \begin{aligned}
 B_{V_{(\q,B_{\q})}(1,0)}(F_\lambda\one,F_\lambda\one)&=1,\\
 B_{V_{(\q,B_{\q})}(1,0)}(F_\lambda(x(m)u),F_\lambda v)
 &=B_{V_{(\q,B_{\q})}(1,0)}((\psi_\lambda x)(m)F_\lambda u,F_\lambda v)\\
 &=-B_{V_{(\q,B_{\q})}(1,0)}(F_\lambda u,(\psi_\lambda x)(-m)F_\lambda v)\\
 &=-B_{V_{(\q,B_{\q})}(1,0)}(F_\lambda u,F_\lambda(x(-m)v)).
 \end{aligned}
\]
By the uniqueness in Proposition~\ref{prop:currentform},
\[
 B_{\mathcal V_\lambda}(u,v)=B_{V_{(\q,B_{\q})}(1,0)}(F_\lambda u,F_\lambda v),\qquad u,v\in\mathcal V_\lambda.
\]
On $W_n^{\mathfrak{Lie}}$, let
\[
 T_n(\psi_\lambda)=
 \bigoplus_{\sum_{m\ge1}mr_m=n}\ \bigotimes_{m\ge1}\Sym^{r_m}(\psi_\lambda),
\]
where $\Sym^{r}(\psi_\lambda)$ denotes the map $x_1\cdots x_r\mapsto(\psi_\lambda x_1)\cdots(\psi_\lambda x_r)$, each $\psi_\lambda x_i$ being understood via $\q_m\simeq\q$.
By the definition of $\Phi_n^{\mathfrak{Lie}}$, replacing each factor $x$ by $\psi_\lambda x$ commutes with symmetrization, so
$F_\lambda\Phi_{n,\lambda}^{\mathfrak{Lie}}=\Phi_n^{\mathfrak{Lie}}T_n(\psi_\lambda)$. Hence
\[
 G_{n,\lambda}^{\mathfrak{Lie}}(u,v)=G_n^{\mathfrak{Lie}}\bigl(T_n(\psi_\lambda)u,T_n(\psi_\lambda)v\bigr),\qquad u,v\in W_n^{\mathfrak{Lie}},
\]
and therefore $\det G_{n,\lambda}^{\mathfrak{Lie}}=(\det T_n(\psi_\lambda))^2\det G_n^{\mathfrak{Lie}}$.
Moreover, since $\psi_\lambda$ preserves $B_{\q}$, \eqref{eq:currentheisenberg} shows that $T_n(\psi_\lambda)$ preserves $H_n^{\mathfrak{Lie}}$, that is,
$H_n^{\mathfrak{Lie}}(T_n(\psi_\lambda)u,T_n(\psi_\lambda)v)=H_n^{\mathfrak{Lie}}(u,v)$;
by Lemma~\ref{lem:lie-auxiliary-nondegenerate}, $\det H_n^{\mathfrak{Lie}}\ne0$, so $(\det T_n(\psi_\lambda))^2=1$.
Therefore $\det G_{n,\lambda}^{\mathfrak{Lie}}=\det G_n^{\mathfrak{Lie}}$ for all nonzero $\lambda\in\C$.

On the other hand, write $[e_i,e_j]_\lambda=\sum_{k=1}^d\mathsf c_{ij}^{k}(\lambda)e_k$. By Definition~\ref{def:polynomial-limit}, $\mathsf c_{ij}^{k}(\lambda)\in\C[\lambda]$;
applying step (1) to $(\q,[\,,\,]_\lambda,B_{\q})$, we get
\[
 G_{n,\lambda}^{\mathfrak{Lie}}(w_\alpha,w_\beta)=\mathsf F_{\alpha\beta}\bigl((\mathsf c_{ij}^{k}(\lambda)),(\mathsf b_{ij})\bigr),\qquad \lambda\in\C,
\]
so $\det G_{n,\lambda}^{\mathfrak{Lie}}$ is a polynomial in $\lambda$. It takes the constant value $\det G_n^{\mathfrak{Lie}}$ at all nonzero $\lambda\in\C$, so by the continuity of polynomials at $\lambda=0$,
\[
 \det G_{n,0}^{\mathfrak{Lie}}=\lim_{\lambda\to0}\det G_{n,\lambda}^{\mathfrak{Lie}}=\det G_n^{\mathfrak{Lie}}.
\]

\smallskip\noindent\emph{(3) Reduction to the abelian case.}
Take the Lie brackets $[\,,\,]^{(0)}=[\,,\,],\ldots,[\,,\,]^{(s)}=0$ of Proposition~\ref{prop:degenerate},
and, for $0\le j\le s$, denote by $G_n^{\mathfrak{Lie},(j)}$ the form $G_n^{\mathfrak{Lie}}$ defined for $(\q,[\,,\,]^{(j)},B_{\q})$.
Since each $[\,,\,]^{(j+1)}$ is a polynomial limit of $[\,,\,]^{(j)}$, since the space $W_n^{\mathfrak{Lie}}$, its basis, and $H_n^{\mathfrak{Lie}}$ do not depend on the Lie bracket,
and since steps (1) and (2) do not use the assumption that there is no nonzero semisimple ideal, applying step (2) successively to each $(\q,[\,,\,]^{(j)},B_{\q})$ gives
\[
 \det G_n^{\mathfrak{Lie},(0)}=\det G_n^{\mathfrak{Lie},(1)}=\cdots=\det G_n^{\mathfrak{Lie},(s)}.
\]
Since $[\,,\,]^{(s)}=0$, that is, $(\q,[\,,\,]^{(s)})$ is an abelian Lie algebra, Lemma~\ref{lem:heisenbergpairing} gives $G_n^{\mathfrak{Lie},(s)}=H_n^{\mathfrak{Lie}}$;
then Lemma~\ref{lem:lie-auxiliary-nondegenerate} gives $\det G_n^{\mathfrak{Lie}}=\det H_n^{\mathfrak{Lie}}\ne0$.

\smallskip\noindent\emph{(4) Vanishing of the radical.}
Since $\det G_n^{\mathfrak{Lie}}\ne0$ and $\Phi_n^{\mathfrak{Lie}}$ is a vector space isomorphism, the restriction of $B_{V_{(\q,B_{\q})}(1,0)}$ to ${V_{(\q,B_{\q})}(1,0)}_n\times{V_{(\q,B_{\q})}(1,0)}_n$ is nondegenerate for every $n\in\N$.
By the orthogonality of different degrees, the same argument as in Remark~\ref{rem:weight-gram} gives
\[
 \rad B_{V_{(\q,B_{\q})}(1,0)}=\bigoplus_{n\in\N}\rad\Bigl(B_{V_{(\q,B_{\q})}(1,0)}\big|_{{V_{(\q,B_{\q})}(1,0)}_n\times{V_{(\q,B_{\q})}(1,0)}_n}\Bigr)=0.
\]
Then, by Lemma~\ref{lem:radical}, $J_{(\q,B_{\q})}(1,0)=0$ and $L_{(\q,B_{\q})}(1,0)={V_{(\q,B_{\q})}(1,0)}$.
\end{proof}

\section{Zhu algebras and the classification in conformal weight one}\label{sec:affine}
In this section we prove Theorem~\ref{thm:one}. Let $V$ satisfy its hypotheses. By Proposition~\ref{prop:lian}, $V$ is isomorphic, as a vertex algebra, to
the simple graded quotient of $V_{(\g,B_{\g})}(1,0)$; moreover, for the part $\q$ of $\g$ without nonzero semisimple ideals,
$V_{(\q,B_{\q})}(1,0)$ has no nonzero proper graded vertex ideals (Theorem~\ref{thm:currentsimple}).
The proof rests on two key points: first, the rationality of $V$ ensures that the Zhu algebra $A(V)$ is finite-dimensional, which rules out $\q$;
second, the conformal vector of $V$ determines the Sugawara structure of each simple affine factor.
Only after the conformal vector has been determined and the critical level has been excluded does the tensor product decomposition of $V$ yield a vertex operator algebra isomorphism,
and only then does the rationality of $V$ pass to the factors.

Section~\ref{subsec:zhu-decomposition} decomposes $\g$ into simple ideals and $\q$, and obtains a grading-preserving tensor product decomposition of $V$ (Proposition~\ref{prop:curfactors});
it then uses the top-level action of \cite{Zhu,DLMtw} to prove that $U(\q)\to A(V)$ is injective, and combines this with rationality to obtain $\q=0$; see Corollary~\ref{cor:qzero}.
Section~\ref{subsec:casimir} introduces the Casimir elements of the simple ideals, determines the conformal vector of $V$, and proves that the level determined by $B_{\g}$ on each simple ideal is not the critical level;
this lifts the above decomposition to a vertex operator algebra isomorphism; see Proposition~\ref{prop:sugawara}.
Section~\ref{subsec:proof-weight-one} uses the rationality of tensor products and the classification of simple affine vertex operator algebras (Proposition~\ref{prop:affineclass})
to determine the levels, completing the proof of Theorem~\ref{thm:one}.

\Needspace{3\baselineskip}
\subsection{Tensor product decomposition and the Zhu algebra}\label{subsec:zhu-decomposition}
Throughout this subsection, $V$ satisfies the hypotheses of Theorem~\ref{thm:one}, that is, $V$ is a rational vertex operator algebra of CFT type that is generated by $V_1$ as a vertex algebra; we keep the notation $\g=V_1$ and $B_{\g}$.
We first decompose $\g$ into simple ideals and a part $\q$ without nonzero semisimple ideals; see Lemma~\ref{lem:onesplit}.
Then Proposition~\ref{prop:lian}, Lemma~\ref{lem:radical}, and Theorem~\ref{thm:currentsimple}
give a grading-preserving tensor product decomposition of $V$; see Proposition~\ref{prop:curfactors}.
Next, working on the top level of an admissible $V$-module, we prove that $U(\q)\to A(V)$ is injective
(Proposition~\ref{prop:injection}), and deduce $\q=0$ from rationality; see Corollary~\ref{cor:qzero}.

By Lemma~\ref{lem:simple}, $V$ is simple; then, by Proposition~\ref{prop:lian}~(1), $(\g,B_{\g})$ is a quadratic Lie algebra.
We now apply the decomposition of Lemma~\ref{lem:split} to $\g$, and further decompose its semisimple part into simple ideals.

\begin{lemma}\label{lem:onesplit}
There exist $r\in\N$, simple ideals $\s_1,\ldots,\s_r$, and an ideal $\q$ such that
\[
 \g=\s_1\perp\cdots\perp\s_r\perp\q
\]
and $\q$ has no nonzero semisimple ideal.
Write $B_{\s_i}=B_{\g}|_{\s_i\times\s_i}$ and $B_{\q}=B_{\g}|_{\q\times\q}$,
and let $(\, ,\,)_i$ be the normalized Killing form of $\s_i$ (see Definition~\ref{def:affine-sugawara}) for $1\le i\le r$.
Then there exist $k_i\in\C$ with $k_i\ne0$ such that $B_{\s_i}=k_i(\, ,\,)_i$,
and $(\s_i,B_{\s_i})$ and $(\q,B_{\q})$ are quadratic Lie algebras.
\end{lemma}
\begin{proof}
By Lemma~\ref{lem:split}, $\g=\s\perp\q$, where $\s$ is a semisimple ideal and $\q$ has no nonzero semisimple ideal.
By Proposition~\ref{prop:lie-basic}~(4), $\s=\bigoplus_{i=1}^r\s_i$, where $r\in\N$ and the $\s_i$ are simple ideals of $\s$;
in what follows, $1\le i,j\le r$. Since $[\q,\s_i]\subseteq[\q,\s]=0$, each $\s_i$ is also an ideal of $\g$.
Since each $\s_i$ is simple, $[\s_i,\s_i]=\s_i$; for $i\ne j$, $[\s_i,\s_j]\subseteq\s_i\cap\s_j=0$.
By the invariance of $B_{\g}$,
\[
 B_{\g}([u,v],y)=B_{\g}(u,[v,y])=0,\qquad u,v\in\s_i,\ y\in\s_j.
\]
Since $[\s_i,\s_i]=\s_i$, we get $B_{\g}(\s_i,\s_j)=0$; moreover, $B_{\g}(\s,\q)=0$. Hence $\g=\s_1\perp\cdots\perp\s_r\perp\q$.
The above summands are pairwise orthogonal and $B_{\g}$ is nondegenerate, so $B_{\s_i}$ and $B_{\q}$ are nondegenerate.
Since $B_{\s_i}$ and $(\, ,\,)_i$ are both nondegenerate symmetric invariant bilinear forms on $\s_i$ (for the nondegeneracy of the latter, see \cite[Section~5.1, Theorem]{Hum}),
\cite[Section~22.1, Lemma]{Hum} gives $k_i\in\C$ with $k_i\ne0$ such that $B_{\s_i}=k_i(\, ,\,)_i$.
Hence $(\s_i,B_{\s_i})$ and $(\q,B_{\q})$ are quadratic Lie algebras.
\end{proof}

By Lemma~\ref{lem:onesplit}, the induced vacuum modules $V_{(\s_i,B_{\s_i})}(1,0)$ and $V_{(\q,B_{\q})}(1,0)$ corresponding to these summands
are given by Definition~\ref{def:affine-vacuum}.
For $1\le i\le r$, by \eqref{eq:affine-simple-quotient} and Remark~\ref{rem:affine-level-scaling},
$L_{(\s_i,(\, ,\,)_i)}(k_i,0)=L_{(\s_i,B_{\s_i})}(1,0)$ is
the simple graded quotient of $V_{(\s_i,B_{\s_i})}(1,0)$.
Let $h_i^\vee$ be the dual Coxeter number of $\s_i$. Since we have not yet proved that $k_i\ne-h_i^\vee$,
$L_{(\s_i,(\, ,\,)_i)}(k_i,0)$ is used only as a vertex algebra, with the grading induced by the natural grading \eqref{eq:affine-grading}.
Combining Proposition~\ref{prop:lian}, Lemma~\ref{lem:radical}, and Theorem~\ref{thm:currentsimple},
we now prove that $V$ has the following grading-preserving decomposition as a vertex algebra.

\begin{proposition}\label{prop:curfactors}
Let $\g=\s_1\perp\cdots\perp\s_r\perp\q$ be the decomposition of Lemma~\ref{lem:onesplit}, where the $\s_i$, $1\le i\le r$, are simple ideals,
$\q$ has no nonzero semisimple ideal, and $B_{\s_i}=k_i(\, ,\,)_i$ with $k_i\ne0$. Then there is a grading-preserving vertex algebra isomorphism
\begin{equation}\label{eq:curfactors}
 V\cong\left(\bigotimes_{i=1}^r L_{(\s_i,(\, ,\,)_i)}(k_i,0)\right)\otimes{V_{(\q,B_{\q})}(1,0)},
\end{equation}
where the right-hand side is equipped with the grading induced by the natural gradings \eqref{eq:affine-grading} of the tensor factors,
and the left-hand side with the conformal weight grading of $V$.
\end{proposition}
\begin{proof}
Write
\[
 T=\left(\bigotimes_{i=1}^r{V_{(\s_i,B_{\s_i})}(1,0)}\right)\otimes{V_{(\q,B_{\q})}(1,0)},
\]
and denote its vacuum by $\one_T$.
As in the first paragraph of the proof of Proposition~\ref{prop:frobfactors}, for $x\in\s_i$, $1\le i\le r$, or $x\in\q$,
let $\iota(x)\in T$ be the tensor product of $x(-1)\one$ in the corresponding tensor factor with the vacuum vectors of the other factors,
and extend linearly to $\iota:\g\to T$.
By \eqref{eq:tensor-vertex}, the modes of $\iota(x)$ act only on the corresponding tensor factor, where they act as $x(m)$.
Let $x(m)$, $x\in\g$, $m\in\Z$, act on $T$ as $\iota(x)_m$, and let $K$ act as $\id_T$.
If $x,y$ lie in the same summand of the decomposition of Lemma~\ref{lem:onesplit}, then \eqref{eq:currentbracket} is the commutation relation in the induced vacuum module
corresponding to that summand; if $x,y$ lie in different summands, then $[x,y]=0$ and $B_{\g}(x,y)=0$,
while $\iota(x)_m$ and $\iota(y)_n$ act on different tensor factors and commute, so both sides of \eqref{eq:currentbracket} are zero.
By bilinearity, $T$ becomes an $\widehat{\g}_{B_{\g}}$-module of level $1$;
moreover, since $x(m)\one=0$ for $m\in\N$ in each factor, we have $x(m)\one_T=0$ for $x\in\g$ and $m\in\N$.
By the universal property of the induced module in \eqref{eq:currentvac} and \cite[Proposition~5.7.9]{LL}, we obtain a
vertex algebra homomorphism
\[
\begin{aligned}
 \Xi:{V_{(\g,B_{\g})}(1,0)}&\longrightarrow T,\\
 x(-1)\one&\longmapsto\iota(x),\qquad x\in\g.
\end{aligned}
\]
Take the ordered basis of $\g$ obtained by concatenating bases of the summands, and order the negative modes in the PBW monomials in the order $\s_1,\ldots,\s_r,\q$.
Then $\Xi$ maps the PBW basis of Proposition~\ref{prop:affine-vacuum-structure} bijectively onto the tensor product basis of $T$
and preserves degrees; hence $\Xi$ is a grading-preserving vertex algebra isomorphism.

Denote by ${B_{V_{(\g,B_{\g})}(1,0)}}$, $B_{V_{(\s_i,B_{\s_i})}(1,0)}$, and $B_{V_{(\q,B_{\q})}(1,0)}$ the normalized contravariant forms on the corresponding induced vacuum modules
given by Proposition~\ref{prop:currentform}; by that proposition, they are symmetric.
As in \eqref{eq:tensor-form}, let
\[
 B_T=B_{V_{(\s_1,B_{\s_1})}(1,0)}\otimes\cdots\otimes B_{V_{(\s_r,B_{\s_r})}(1,0)}\otimes B_{V_{(\q,B_{\q})}(1,0)}.
\]
By the normalization of the factors, $B_T(\one_T,\one_T)=1$.
For $x\in\s_i$, $1\le i\le r$, and $m\in\Z$, the operator $x(m)$ acts only on the $i$-th tensor factor,
so the contravariance of $B_{V_{(\s_i,B_{\s_i})}(1,0)}$ gives
\[
 B_T(x(m)u,v)=-B_T(u,x(-m)v),\qquad u,v\in T;
\]
the same holds for $x\in\q$, and hence, by linearity, for all $x\in\g$. Therefore $(u,v)\mapsto B_T(\Xi u,\Xi v)$ satisfies \eqref{eq:contracurrent},
and the uniqueness in Proposition~\ref{prop:currentform} gives
\[
 B_{V_{(\g,B_{\g})}(1,0)}(u,v)=B_T(\Xi u,\Xi v),\qquad u,v\in{V_{(\g,B_{\g})}(1,0)}.
\]
Since $\Xi$ is bijective, we have, in particular, $\Xi\bigl(\rad B_{V_{(\g,B_{\g})}(1,0)}\bigr)=\rad B_T$.

By Propositions~\ref{prop:affine-vacuum-structure} and~\ref{prop:currentform}, the homogeneous components of each factor are finite-dimensional,
and different degrees are orthogonal with respect to the corresponding normalized contravariant forms; by Lemma~\ref{lem:radical}, each radical is a vertex ideal.
Hence, as in Proposition~\ref{prop:tensorforms}, \eqref{eq:tensor-quotient} gives a grading-preserving vertex algebra isomorphism
\[
 T/\rad B_T\cong
 \left(\bigotimes_{i=1}^r{V_{(\s_i,B_{\s_i})}(1,0)}/\rad B_{V_{(\s_i,B_{\s_i})}(1,0)}\right)
 \otimes{V_{(\q,B_{\q})}(1,0)}/\rad B_{V_{(\q,B_{\q})}(1,0)}.
\]
By Lemma~\ref{lem:radical}, the $i$-th factor on the right-hand side is $L_{(\s_i,B_{\s_i})}(1,0)=L_{(\s_i,(\, ,\,)_i)}(k_i,0)$;
by Theorem~\ref{thm:currentsimple}, the last factor is ${V_{(\q,B_{\q})}(1,0)}$.
On the other hand, by Proposition~\ref{prop:lian}~(3), Lemma~\ref{lem:radical}, and $\Xi\bigl(\rad B_{V_{(\g,B_{\g})}(1,0)}\bigr)=\rad B_T$,
\[
 V\cong L_{(\g,B_{\g})}(1,0)={V_{(\g,B_{\g})}(1,0)}/\rad B_{V_{(\g,B_{\g})}(1,0)}\cong T/\rad B_T,
\]
and all these isomorphisms preserve the gradings. Combining the last two displays gives \eqref{eq:curfactors}.
\end{proof}

By Remark~\ref{rem:top-level-action}, the part $\q\subseteq V_1$ without nonzero semisimple ideals in Lemma~\ref{lem:onesplit}
gives an algebra homomorphism $U(\q)\to A(V)$.
For the four-dimensional Nappi--Witten Lie algebra $H_4$, \cite[Remark~4.9]{BJP} points out that
the Zhu algebra of the corresponding affine vertex operator algebra at a nonzero level is isomorphic to $U(H_4)$.
We now prove that the above homomorphism is injective.

\begin{proposition}\label{prop:injection}
Let $\g=\s_1\perp\cdots\perp\s_r\perp\q$ be the decomposition of Lemma~\ref{lem:onesplit}, where $\q$ has no nonzero semisimple ideal.
Then the algebra homomorphism given by Remark~\ref{rem:top-level-action},
\[
\begin{aligned}
 \rho:U(\q)&\longrightarrow A(V),\\
 x&\longmapsto[x],\qquad x\in\q,
\end{aligned}
\]
is injective. Consequently, $\dim A(V)<\infty$ implies $\q=0$.
\end{proposition}
\begin{proof}
Take the left regular module $U(\q)$, and let $\q[t]\oplus\C K$ act on it by
\[
 x(0)u=xu,\qquad x(m)u=0,\quad m\in\Z_{\ge1},\qquad Ku=u,
 \qquad x\in\q,\ u\in U(\q).
\]
By \eqref{eq:currentbracket}, this gives a $\q[t]\oplus\C K$-module structure. Define the induced module
\[
 \mathcal M_{(\q,B_{\q})}(1,U(\q))=U(\widehat{\q}_{B_{\q}})
 \otimes_{U(\q[t]\oplus\C K)}U(\q),
\]
which we abbreviate to $\mathcal M$. By the PBW theorem, there is a vector space isomorphism
\[
 \mathcal M\cong
 U(\q\otimes t^{-1}\C[t^{-1}])\otimes U(\q).
\]
Assigning degree zero to $U(\q)$ and degree $m$ to $x(-m)$, $x\in\q$, $m\in\Z_{\geq1}$,
we have $\mathcal M=\bigoplus_{d\in\N}\mathcal M(d)$ with
$\mathcal M(0)=U(\q)$. The relation \eqref{eq:currentbracket} is compatible with this grading, so
\[
 x(m)\mathcal M(d)\subseteq\mathcal M(d-m)
 ,\qquad x\in\q,\ m\in\Z,\ d\in\N.
\]
Here the components of negative degree are zero, so the left-hand side is zero for $m>d$.
Hence $\mathcal M$ is a restricted $\widehat{\q}_{B_{\q}}$-module of level $1$.
By Proposition~\ref{prop:affine-restricted} (see also \cite[Proposition~4.12]{Lian}), it becomes a weak ${V_{(\q,B_{\q})}(1,0)}$-module with
\[
 Y_{\mathcal M}(x(-1)\one,z)
 =\sum_{m\in\Z}x(m)z^{-m-1},\qquad x\in\q.
\]

Write $S=\bigotimes_{i=1}^rL_{(\s_i,(\, ,\,)_i)}(k_i,0)$, with vacuum $\one_S$ and natural grading $S=\bigoplus_{d_1\in\N}S_{d_1}$.
By \eqref{eq:curfactors}, the tensor product $\mathcal N=S\otimes\mathcal M$ of the adjoint module of $S$ with $\mathcal M$
is a weak $V$-module (for tensor product modules, see \cite[Section~4.6]{LL}); we give it the grading
\[
 \mathcal N(d)=\bigoplus_{\substack{d_1,d_2\in\N\\d_1+d_2=d}}S_{d_1}\otimes\mathcal M(d_2)
 ,\quad d\in\N,\qquad \mathcal N(d)=0,\quad d\in\Z_{<0}.
\]
For $x\in\g$ and $m\in\Z$, we again write $x(m)$ for the $m$-th mode of $x$ on $\mathcal N$.
If $x\in\s_i$, $1\le i\le r$, or $x\in\q$, then $x(m)$ lowers the degree by $m$ on the corresponding tensor factor only. By linearity,
\[
 x(m)\mathcal N(d)\subseteq\mathcal N(d-m)
 ,\qquad x\in \g,\ m\in\Z,\ d\in\N.
\]
We now extend this relation to all homogeneous vectors of $V$.
By Remark~\ref{rem:strong-generation}~(2), $V$ is spanned by $\one$ and the vectors
$x^1_{-n_1}\cdots x^\ell_{-n_\ell}\one$ with $\ell\in\Z_{\ge1}$, $x^1,\ldots,x^\ell\in\g$, and $n_1,\ldots,n_\ell\in\Z_{\ge1}$,
and by \eqref{eq:weights} these vectors are homogeneous. We prove by induction on $\ell$ that, for every nonzero such vector $v$,
\[
 v_j\mathcal N(d)\subseteq\mathcal N(d+\wt v-j-1),
 \qquad j\in\Z,\quad d\in\N.
\]
For $\ell=0$, we have $v=\one$ and $\one_j=\delta_{j,-1}\id$, so the assertion holds.
Suppose that the assertion holds for $v'$, and let $v=x_{-n}v'\ne0$, where $x\in\g$ and $n\in\Z_{\ge1}$. By \eqref{eq:negative-mode-module},
\[
 Y_{\mathcal N}(x_{-n}v',z)
 =\frac1{(n-1)!}:
 \bigl(\partial_z^{\,n-1}Y_{\mathcal N}(x,z)\bigr)Y_{\mathcal N}(v',z):.
\]
Comparing the coefficients of $z^{-j-1}$, we see that each resulting term is a scalar multiple of $x(p)v'_{p'}$ or $v'_{p'}x(p)$,
where $p,p'\in\Z$ and $p+p'=j-n$; on any fixed vector, only finitely many terms are nonzero by truncation.
Each term changes the degree by
\[
 -p+\wt v'-p'-1
 =\wt v'+n-j-1
 =\wt v-j-1,
\]
so the assertion holds for $v$. Since each homogeneous subspace of $V$ is spanned by the above homogeneous vectors, linearity gives
\[
 v_j\mathcal N(d)\subseteq\mathcal N(d+\wt v-j-1)
 ,\qquad v\in V\text{ homogeneous},\ j\in\Z,\quad d\in\N.
\]
This is exactly the admissible grading condition of Definition~\ref{def:admissible} for $g=\id_V$,
so $\mathcal N$ is an admissible $V$-module, and $\mathcal N(0)=\C\one_S\otimes U(\q)$.

By Proposition~\ref{prop:zhu}, $A(V)$ acts on $\mathcal N(0)$ via
$[v]\mapsto o(v)$.
For $x\in\q\subseteq \g$, we have $\rho(x)=[x]$ and $o(x)=x_0=x(0)$, so
\[
 \rho(x)\cdot(\one_S\otimes u)
 =x(0)(\one_S\otimes u)
 =\one_S\otimes xu
 ,\qquad u\in U(\q).
\]
Since $\q$ generates $U(\q)$ and $\rho$ is an algebra homomorphism, the above identity gives
\[
 \rho(u)\cdot(\one_S\otimes u')=\one_S\otimes uu',\qquad u,u'\in U(\q).
\]
Hence, for $u\in\ker\rho$,
\[
 0=\rho(u)\cdot(\one_S\otimes1)=\one_S\otimes u,
\]
so $u=0$, that is, $\rho$ is injective.
If $\q\ne0$, the PBW theorem gives $\dim U(\q)=\infty$,
so $\dim A(V)<\infty$ implies $\q=0$.
\end{proof}

Since $V$ is rational, Proposition~\ref{prop:zhu} shows that the Zhu algebra $A(V)$ of $V$ is finite-dimensional;
hence, by Proposition~\ref{prop:injection}, the ideal $\q$ in Lemma~\ref{lem:onesplit} must be zero,
and therefore, by Proposition~\ref{prop:lie-basic}~(4), $\g=\s_1\perp\cdots\perp\s_r$ is a semisimple Lie algebra.
In this case $V_{(\q,B_{\q})}(1,0)=\C\one$, so Proposition~\ref{prop:curfactors} reduces to an isomorphism between $V$ and a tensor product of simple affine vertex algebras.
We record this as follows.

\begin{corollary}\label{cor:qzero}
In the decomposition $\g=\s_1\perp\cdots\perp\s_r\perp\q$ of Lemma~\ref{lem:onesplit}, we have $\q=0$.
Hence $\g=\s_1\perp\cdots\perp\s_r$ is a semisimple Lie algebra,
and \eqref{eq:curfactors} reduces to a grading-preserving vertex algebra isomorphism
\[
 V\cong\bigotimes_{i=1}^rL_{(\s_i,(\, ,\,)_i)}(k_i,0).
\]
\end{corollary}

The isomorphism in Corollary~\ref{cor:qzero} is only a grading-preserving vertex algebra isomorphism: at this stage we have not yet proved that $k_i\ne-h_i^\vee$ for $1\le i\le r$,
the factors on the right-hand side have not yet been equipped with Sugawara conformal vectors, and the image of the conformal vector of $V$ under this isomorphism is not yet known.
In the next subsection we prove that $k_i\ne-h_i^\vee$, and show that the conformal vector of $V$ is the sum of the
Sugawara conformal vectors of the simple ideals, thereby lifting this isomorphism to a vertex operator algebra isomorphism;
see Proposition~\ref{prop:sugawara}.

\Needspace{3\baselineskip}
\subsection{Casimir elements and the Sugawara conformal vector}\label{subsec:casimir}
In this subsection we keep the hypotheses and notation of Section~\ref{subsec:zhu-decomposition}: $V$ satisfies the hypotheses of Theorem~\ref{thm:one}, that is, $V$ is a rational vertex operator algebra of CFT type
that is generated by $V_1$ as a vertex algebra, and $\g=V_1$ is equipped with the Lie bracket and the bilinear form $B_{\g}$ given by \eqref{eq:currentproducts}.
By Corollary~\ref{cor:qzero}, the decomposition of Lemma~\ref{lem:onesplit} reduces to
\[
 \g=\s_1\perp\cdots\perp\s_r,
\]
where $r\in\N$ and the $\s_i$, $1\le i\le r$, are simple ideals; hence $\g$ is semisimple.
Here $(\, ,\,)_i$ is the normalized Killing form of $\s_i$ (see Definition~\ref{def:affine-sugawara}),
$h_i^\vee$ is the dual Coxeter number of $\s_i$, and $B_{\s_i}=k_i(\, ,\,)_i$ with $k_i\in\C$ and $k_i\ne0$.
Corollary~\ref{cor:qzero} also gives a grading-preserving vertex algebra isomorphism
\[
 V\cong\bigotimes_{i=1}^rL_{(\s_i,(\, ,\,)_i)}(k_i,0).
\]
We prove that $k_i\ne-h_i^\vee$, $1\le i\le r$, and that the conformal vector $\omega$ of $V$ is the sum of the
Sugawara conformal vectors of the simple ideals, so that this isomorphism is a vertex operator algebra isomorphism; see Proposition~\ref{prop:sugawara}.
To this end, we first introduce in $V$ the Casimir elements $\Omega_i$ corresponding to the $\s_i$, and determine the tensors in $\g\otimes\g$ annihilated by $\g$;
see Lemma~\ref{lem:sym-invariants}.

For $1\le i\le r$, choose dual bases $\{u_{i,p}\}_{p=1}^{\dim\s_i}$ and $\{u_i^p\}_{p=1}^{\dim\s_i}$ of $\s_i$ with respect to $(\, ,\,)_i$,
and set
\[
 \Omega_i=\sum_{p=1}^{\dim\s_i}(u_{i,p})_{-1}u_i^p\in V_2.
\]
As in \cite[Remark~6.2.14]{LL}, $\Omega_i$ is called the \emph{Casimir element} corresponding to $\s_i$.
By the same computation as in \cite[(6.2.37), Remark~6.2.14, (6.2.53)]{LL} (without dividing by $k_i+h_i^\vee$), $\Omega_i$ is independent of the choice of dual bases, and
\begin{equation}\label{eq:casimir-action}
 x_1\Omega_j=2\delta_{ij}(k_i+h_i^\vee)x,\qquad x\in\s_i,\ 1\le i,j\le r.
\end{equation}

Let $\g$ act on the adjoint module $\g$, and give $\g\otimes\g$ the tensor product module structure as in \cite[Section~6.1]{Hum}:
\[
 a\cdot x=[a,x],\qquad a\cdot(y\otimes z)=[a,y]\otimes z+y\otimes[a,z],\qquad a,x,y,z\in\g.
\]
We now determine the elements of $\g\otimes\g$ annihilated by $\g$.

\begin{lemma}\label{lem:sym-invariants}
For $1\le i\le r$, let $\widetilde\Omega_i=\sum_{p=1}^{\dim\s_i}u_{i,p}\otimes u_i^p\in\g\otimes\g$. Then
\begin{equation}\label{eq:tensor-invariants}
 \{\xi\in\g\otimes\g:a\cdot \xi=0\text{ for all }a\in\g\}=\Span\{\widetilde\Omega_1,\ldots,\widetilde\Omega_r\}.
\end{equation}
\end{lemma}
\begin{proof}
By the nondegeneracy and invariance of $B_{\g}$, the map $v\mapsto B_{\g}(v,\cdot)$ is a $\g$-module isomorphism $\g\to\g^*$.
Combining this with the $\g$-module isomorphism $\g^*\otimes\g\cong\End(\g)$ of \cite[Section~6.1]{Hum}, we obtain a $\g$-module isomorphism
\[
\begin{aligned}
 \Upsilon:\g\otimes\g&\longrightarrow\End(\g),\\
 y\otimes z&\longmapsto\bigl(v\mapsto B_{\g}(v,y)z\bigr),\qquad y,z\in\g,
\end{aligned}
\]
where the action on $\End(\g)$ is $a\cdot\phi=[\operatorname{ad}a,\phi]$, $a\in\g$, $\phi\in\End(\g)$.
Hence, for $\xi\in\g\otimes\g$, we have $a\cdot \xi=0$ for all $a\in\g$ if and only if $\Upsilon(\xi)$ commutes with all $\operatorname{ad}a$, $a\in\g$,
that is, $\Upsilon(\xi)$ lies in
\[
 \End_{\g}(\g)=\{\phi\in\End(\g):\phi\circ\operatorname{ad}a=\operatorname{ad}a\circ\phi\text{ for all }a\in\g\}.
\]
For $1\le i\le r$, let $\mathsf p_i\in\End(\g)$ be the projection onto $\s_i$ with respect to $\g=\s_1\perp\cdots\perp\s_r$,
that is, $\mathsf p_i|_{\s_i}=\id_{\s_i}$ and $\mathsf p_i|_{\s_j}=0$ for $j\ne i$. Since the $\s_j$ are ideals, $\mathsf p_i\in\End_{\g}(\g)$.

Since $\s_i$ is simple, it is an irreducible $\g$-module. For $j\ne i$, $\s_j$ acts on $\s_i$ as zero,
while $\s_j$ acts nontrivially on $\s_j$ because $[\s_j,\s_j]=\s_j$; hence $\s_i$ and $\s_j$ are not isomorphic.
Therefore, for $\phi\in\End_{\g}(\g)$ and $1\le i,j\le r$, Schur's lemma (see \cite[Section~6.1]{Hum}) shows that
$\mathsf p_i\circ\phi|_{\s_j}:\s_j\to\s_i$ is zero for $j\ne i$, and, for $j=i$, equals $\theta_i\,\id_{\s_i}$ for some $\theta_i\in\C$.
Hence $\phi=\sum_{i=1}^r\theta_i\mathsf p_i$, that is,
\[
 \End_{\g}(\g)=\Span\{\mathsf p_1,\ldots,\mathsf p_r\}.
\]

On the other hand, for $1\le i,j\le r$ and $v\in\s_j$, using $B_{\g}(\s_j,\s_i)=0$ for $j\ne i$, $B_{\s_i}=k_i(\, ,\,)_i$,
and the dual basis relation
\[
 w=\sum_{p=1}^{\dim\s_i}(w,u_{i,p})_i\,u_i^p,\qquad w\in\s_i,
\]
we obtain
\[
 \Upsilon(\widetilde\Omega_i)(v)=\sum_{p=1}^{\dim\s_i}B_{\g}(v,u_{i,p})\,u_i^p=\delta_{ij}k_iv,
\]
so $\Upsilon(\widetilde\Omega_i)=k_i\mathsf p_i$. Since $k_i\ne0$ and $\Upsilon$ is a linear isomorphism, the space of tensors annihilated by $\g$ is
\[
 \Upsilon^{-1}\bigl(\End_{\g}(\g)\bigr)=\Span\{\Upsilon^{-1}(\mathsf p_1),\ldots,\Upsilon^{-1}(\mathsf p_r)\}=\Span\{\widetilde\Omega_1,\ldots,\widetilde\Omega_r\}.
\]
Hence \eqref{eq:tensor-invariants} holds.
\end{proof}

Using \eqref{eq:casimir-action} and Lemma~\ref{lem:sym-invariants}, we now determine the conformal vector $\omega$ of $V$,
and thereby show that the grading-preserving vertex algebra isomorphism of Corollary~\ref{cor:qzero} is a vertex operator algebra isomorphism.

\begin{proposition}\label{prop:sugawara}
Let $\g=\s_1\perp\cdots\perp\s_r$ be the decomposition of Corollary~\ref{cor:qzero}, where the $\s_i$, $1\le i\le r$, are simple ideals,
$B_{\s_i}=k_i(\, ,\,)_i$ with $k_i\ne0$, and $h_i^\vee$ is the dual Coxeter number of $\s_i$; let $\Omega_i$ be as defined before \eqref{eq:casimir-action}.
Then $k_i\ne-h_i^\vee$ for each $1\le i\le r$, and the conformal vector $\omega$ of $V$ is
\begin{equation}\label{eq:sugawara}
 \omega=\sum_{i=1}^r\frac{\Omega_i}{2(k_i+h_i^\vee)}.
\end{equation}
Consequently, when each $L_{(\s_i,(\, ,\,)_i)}(k_i,0)$ is equipped with the image of the Sugawara conformal vector $\omega_{\s_i,k_i}$ (see Definition~\ref{def:affine-sugawara}),
the grading-preserving vertex algebra isomorphism of Corollary~\ref{cor:qzero}
\[
 V\cong\bigotimes_{i=1}^rL_{(\s_i,(\, ,\,)_i)}(k_i,0)
\]
is a vertex operator algebra isomorphism.
\end{proposition}
\begin{proof}
\smallskip\noindent\emph{(1) The action of $\g$ on $\omega$.}
Let $x\in\g$. Since $V$ is of CFT type, \eqref{eq:weights} gives $L(0)x=x$, $L(1)x\in\C\one$, and $L(n)x=0$ for $n\in\Z_{\ge2}$.
Taking $a=x$ and $b=\omega$ in \eqref{eq:skew-symmetry},
comparing the coefficients of $z^{-1}$ and $z^{-2}$, and using $L(-1)\one=0$, we obtain
\[
 \begin{aligned}
 x_0\omega&=-L(-1)x+L(-1)L(0)x-\tfrac12L(-1)^2L(1)x=0,\\
 x_1\omega&=L(0)x-L(-1)L(1)x=x.
 \end{aligned}
\]

\smallskip\noindent\emph{(2) $\omega$ is a linear combination of the $\Omega_i$.}
By Remark~\ref{rem:strong-generation}~(1), (2), $V_2=\Span\{x_{-2}\one,\ y_{-1}z:x,y,z\in\g\}$, and $x_{-2}\one=L(-1)x$. Define the linear map
\[
 \begin{aligned}
 F:\g\oplus(\g\otimes\g)&\longrightarrow V_2,\\
 (x,y\otimes z)&\longmapsto L(-1)x+y_{-1}z,\qquad x,y,z\in\g,
 \end{aligned}
\]
extended linearly; then $F$ is surjective.
Let $\g$ act on $\g\oplus(\g\otimes\g)$ componentwise. For $a\in\g$, by $L(-1)=\omega_0$ and \eqref{eq:mode-commutator},
$[L(-1),a_0]=(L(-1)a)_0$; moreover, by \eqref{eq:voa-derivative}, $(L(-1)a)_n=-na_{n-1}$, $n\in\Z$, so $[L(-1),a_0]=0$.
Also, by \eqref{eq:mode-commutator} and \eqref{eq:currentproducts}, $[a_0,y_{-1}]=(a_0y)_{-1}=[a,y]_{-1}$ for $y\in\g$.
Hence, for $a,x,y,z\in\g$,
\[
 \begin{aligned}
 a_0F(x,0)&=a_0L(-1)x=L(-1)[a,x]=F\bigl(a\cdot(x,0)\bigr),\\
 a_0F(0,y\otimes z)&=a_0(y_{-1}z)=[a,y]_{-1}z+y_{-1}[a,z]=F\bigl(a\cdot(0,y\otimes z)\bigr).
 \end{aligned}
\]
By linearity, $F\bigl(a\cdot(x,\xi)\bigr)=a_0F(x,\xi)$ for all $a,x\in\g$ and $\xi\in\g\otimes\g$. In particular, $\ker F$ is a $\g$-submodule of $\g\oplus(\g\otimes\g)$.
By Proposition~\ref{prop:lie-basic}~(3), there is a $\g$-stable complement $\mathcal E$ such that
\[
 \g\oplus(\g\otimes\g)=\ker F\oplus\mathcal E,\qquad F|_{\mathcal E}:\mathcal E\xrightarrow{\sim}V_2.
\]
Choose $(\zeta_\omega,\xi_\omega)\in\mathcal E$ with $F(\zeta_\omega,\xi_\omega)=\omega$, where $\zeta_\omega\in\g$ and $\xi_\omega\in\g\otimes\g$.
By step (1), for every $a\in\g$, $F\bigl(a\cdot(\zeta_\omega,\xi_\omega)\bigr)=a_0\omega=0$, so $a\cdot(\zeta_\omega,\xi_\omega)\in\mathcal E\cap\ker F=0$.
Since the action is componentwise, $a\cdot(\zeta_\omega,\xi_\omega)=([a,\zeta_\omega],a\cdot\xi_\omega)$,
so $[a,\zeta_\omega]=0$ and $a\cdot\xi_\omega=0$ for all $a\in\g$. The former shows that $\zeta_\omega\in Z(\g)$;
since $\g$ is semisimple, $Z(\g)=0$ by \cite[Section~5.1]{Hum}, so $\zeta_\omega=0$.
By Lemma~\ref{lem:sym-invariants}, the latter shows that $\xi_\omega$ is a linear combination of the $\widetilde\Omega_i$.
Moreover, by the definition of $\Omega_i$, $F(0,\widetilde\Omega_i)=\Omega_i$ for $1\le i\le r$. Hence
\[
 \omega=F(\zeta_\omega,\xi_\omega)=F(0,\xi_\omega)=\sum_{i=1}^r\gamma_i\Omega_i,\qquad \gamma_i\in\C.
\]

\smallskip\noindent\emph{(3) Determination of the coefficients and the vertex operator algebra isomorphism.}
Let $1\le i\le r$ and let $x\in\s_i$ be nonzero. By step (1) and \eqref{eq:casimir-action},
\[
 x=x_1\omega=\sum_{j=1}^r\gamma_jx_1\Omega_j=2\gamma_i(k_i+h_i^\vee)x,
\]
so $2\gamma_i(k_i+h_i^\vee)=1$. Hence $k_i\ne-h_i^\vee$ for each $1\le i\le r$, and \eqref{eq:sugawara} holds.

Denote the grading-preserving vertex algebra isomorphism of Corollary~\ref{cor:qzero} by
\[
 \varphi:V\xrightarrow{\sim}\bigotimes_{i=1}^rL_{(\s_i,(\, ,\,)_i)}(k_i,0).
\]
By the proof of Proposition~\ref{prop:curfactors}, for $x\in\s_i$, $1\le i\le r$, $\varphi(x)$ is the tensor product of the image of $x(-1)\one$ in the $i$-th tensor factor
with the vacuum vectors of the other factors. Since $\varphi$ is a vertex algebra isomorphism,
$\varphi(\Omega_i)$, $1\le i\le r$, is the tensor product of the image of $\sum_{p=1}^{\dim\s_i}u_{i,p}(-1)u_i^p(-1)\one$ in the $i$-th factor with the vacuum vectors of the other factors;
by Definition~\ref{def:affine-sugawara}, the former equals $2(k_i+h_i^\vee)$ times the image of $\omega_{\s_i,k_i}$.
Hence, by \eqref{eq:sugawara}, $\varphi(\omega)$ is the sum of the images of the $\omega_{\s_i,k_i}$ in the respective factors (with the vacuum vectors in the other factors),
which by \eqref{eq:tensor-vectors} is the conformal vector of the tensor product.
By Definition~\ref{def:affine-sugawara}, each $L_{(\s_i,(\, ,\,)_i)}(k_i,0)$, $1\le i\le r$, equipped with the image of $\omega_{\s_i,k_i}$, is a
simple affine vertex operator algebra, so $\varphi$ is a vertex operator algebra isomorphism.
\end{proof}

\Needspace{3\baselineskip}
\subsection{Proof of the theorem in conformal weight one}\label{subsec:proof-weight-one}
In this subsection we prove Theorem~\ref{thm:one} by combining the tensor product decomposition of Corollary~\ref{cor:qzero}, the determination of the conformal vector in Proposition~\ref{prop:sugawara},
the rationality and the $C_2$ quotients of tensor products in Proposition~\ref{prop:tensor}, and the classification in Proposition~\ref{prop:affineclass}.
The key point of the proof is that, by Proposition~\ref{prop:sugawara}, the isomorphism in Corollary~\ref{cor:qzero} is a
vertex operator algebra isomorphism, so that the rationality of $V$ passes to the simple affine factors,
whose levels are then determined by Proposition~\ref{prop:affineclass}.

\begin{proof}[Proof of Theorem~\ref{thm:one}]
Let $V$ satisfy the hypotheses of Theorem~\ref{thm:one}, that is, $V$ is a rational vertex operator algebra of CFT type that is generated by $V_1$ as a vertex algebra;
write $\g=V_1$, equipped with the Lie bracket and the bilinear form $B_{\g}$ given by \eqref{eq:currentproducts}.
By Corollary~\ref{cor:qzero} and Proposition~\ref{prop:sugawara}, $\g$ is a semisimple Lie algebra and $\g=\s_1\perp\cdots\perp\s_r$,
where $r\in\N$ and the $\s_i$, $1\le i\le r$, are simple ideals of $\g$; by Proposition~\ref{prop:lie-basic}~(4), they are precisely the simple ideals of $\g$.
If $\g=0$, then $r=0$ and $V=\C\one$.
For $1\le i\le r$, let $(\, ,\,)_i$ be the normalized Killing form of $\s_i$ (see Definition~\ref{def:affine-sugawara})
and $h_i^\vee$ the dual Coxeter number of $\s_i$; then $B_{\s_i}=k_i(\, ,\,)_i$, where $k_i\in\C$, $k_i\ne0$, and $k_i\ne-h_i^\vee$.
Furthermore, when each $L_{(\s_i,(\, ,\,)_i)}(k_i,0)$ is equipped with the image of the Sugawara conformal vector $\omega_{\s_i,k_i}$,
there is a vertex operator algebra isomorphism
\begin{equation}\label{eq:onevoa}
 V\cong\bigotimes_{i=1}^rL_{(\s_i,(\, ,\,)_i)}(k_i,0).
\end{equation}

We now determine the levels. Since $V$ is rational, \eqref{eq:onevoa} and Proposition~\ref{prop:tensor}~(1) show that
each $L_{(\s_i,(\, ,\,)_i)}(k_i,0)$, $1\le i\le r$, is rational.
Since $k_i\ne-h_i^\vee$, Proposition~\ref{prop:affineclass} gives $k_i\in\N$; together with $k_i\ne0$, this yields $k_i\in\Z_{>0}$, $1\le i\le r$.
Hence \eqref{eq:onevoa} is exactly the decomposition in Theorem~\ref{thm:one}.

Finally, we prove $C_2$-cofiniteness and regularity. Each $L_{(\s_i,(\, ,\,)_i)}(k_i,0)$ is $C_2$-cofinite by Proposition~\ref{prop:affineclass};
by \eqref{eq:onevoa} and Proposition~\ref{prop:tensor}~(2),
\[
 \dim R_V=\prod_{i=1}^r\dim R_{L_{(\s_i,(\, ,\,)_i)}(k_i,0)}<\infty,
\]
so $V$ is $C_2$-cofinite. Moreover, since $V$ is of CFT type and rational, Proposition~\ref{prop:knownC2} shows that $V$ is regular.
\end{proof}

We end the paper with a discussion of the relations among the conformal weight one analogues of the four equivalent conditions in Theorem~\ref{thm:two}.

\begin{remark}\label{rem:weight-one-conditions}
Let $V$ be a simple vertex operator algebra of CFT type generated by $V_1$ as a vertex algebra,
and let $g\in\Aut(V)$ be an automorphism of finite order.
Consider the following conditions, analogous to those in Theorem~\ref{thm:two}:
\begin{enumerate}
\item there are only finitely many isomorphism classes of irreducible ordinary $g$-twisted $V$-modules;
\item $V$ is $g$-rational;
\item $V$ is rational;
\item there is a vertex operator algebra isomorphism
\[
 V\cong\bigotimes_{i=1}^{r}L_{(\s_i,(\, ,\,)_i)}(k_i,0),
\]
where $r\in\N$, and, for $1\le i\le r$, $\s_i$ is a finite-dimensional simple Lie algebra,
$(\, ,\,)_i$ is its normalized Killing form, and $k_i\in\Z_{>0}$;
for $r=0$, the right-hand side is taken to be $\C\one$.
\end{enumerate}
Then (iii) and (iv) are equivalent, and
\[
 \text{(iii)}\ \Longrightarrow\ \text{(ii)}\ \Longrightarrow\ \text{(i)},
\]
but neither of these implications can be reversed in general. We justify these assertions as follows.
\begin{enumerate}[label=(\arabic*)]
\item The implication (iii)$\Rightarrow$(iv) follows from Theorem~\ref{thm:one},
and the converse follows from Propositions~\ref{prop:tensor}~(1) and~\ref{prop:affineclass}.

\item If $V$ is rational, then Theorem~\ref{thm:one} also gives regularity,
so $V$ is $g$-rational by \cite[Theorem~2.8]{DR}.
This proves (iii)$\Rightarrow$(ii).
The implication (ii)$\Rightarrow$(i) follows from Remark~\ref{rem:module-relations}~(3).

\item Condition (ii) does not imply (iii).
Indeed, the rank-one Heisenberg vertex operator algebra $M(1)$ is simple, of CFT type,
and generated by $M(1)_1$ as a vertex algebra;
see \cite[Section~2.3 and Theorem~2.3]{DN1} and \cite[Section~6.3]{LL}.
Let $\theta$ be the automorphism of $M(1)$ acting as $-1$ on $M(1)_1$.
Then $M(1)$ is $\theta$-rational by \cite[Theorem~3.3 and Corollary~3.9]{DN2}.
On the other hand, \cite[Theorem~2.3]{DN1} gives infinitely many pairwise non-isomorphic
irreducible ordinary $M(1)$-modules $M(1,\lambda)$, $\lambda\in\C$.
Thus $M(1)$ is irrational by Remark~\ref{rem:module-relations}~(3), applied with $g=\id$.

\item Condition (i) does not imply (ii) or (iii), already for $g=\id$.
Let $V=L_{(\mathfrak{sl}_2,(\, ,\,))}(-\frac12,0)$ be the simple affine vertex operator algebra,
where $(\, ,\,)$ is the normalized Killing form.
With its Sugawara conformal vector, $V$ is of CFT type and generated by $V_1$ as a vertex algebra;
see \cite[Section~2]{CHY}.
The level $-\frac12=-2+\frac32$ is admissible for $\widehat{\mathfrak{sl}}_2$,
as recalled in \cite[Section~1]{CHY}.
Hence, by \cite[Corollary~3.3~(1)]{CHY}, there are only finitely many isomorphism classes of irreducible ordinary $V$-modules.
However, since $-\frac12\notin\N$, Proposition~\ref{prop:affineclass} shows that $V$ is irrational.
Thus (i) holds for $g=\id$, while neither (ii) nor (iii) holds.
\end{enumerate}
\end{remark}

\medskip\noindent\textbf{Use of AI tools.} ChatGPT (OpenAI) was used for mathematical exploration, literature searches, preparation of preliminary drafts, proof checking, language editing, and \LaTeX{} typesetting. The research problems, the main ideas and strategies of the proofs, and all final mathematical decisions are due to the authors. Every statement, proof, computation, and citation in this paper has been checked by the authors against the original sources. The authors assume full responsibility for the mathematical arguments, references, and final text, and confirm that the paper contains no fabricated results or references.


\begin{thebibliography}{DHVW}
\addcontentsline{toc}{section}{\refname}
\setlength{\itemsep}{1pt}
\bibitem[Abe1]{AbeSF}  T. Abe,
A $\Z_2$-orbifold model of the symplectic fermionic vertex operator superalgebra,
\emph{Math. Z.} \textbf{255} (2007), no.~4, 755--792.

\bibitem[Abe2]{AbeT}  T. Abe,
$C_2$-cofiniteness of 2-cyclic permutation orbifold models,
\emph{Comm. Math. Phys.} \textbf{317} (2013), no.~2, 425--445.

\bibitem[ABD]{ABD}  T. Abe, G. Buhl and C.-Y. Dong,
Rationality, regularity, and $C_2$-cofiniteness,
\emph{Trans. Amer. Math. Soc.} \textbf{356} (2004), no.~8, 3391--3402.

\bibitem[Abr]{AbE}  L. Abrams,
The quantum Euler class and the quantum cohomology of the Grassmannians,
\emph{Israel J. Math.} \textbf{117} (2000), no.~1, 335--352.

\bibitem[AL]{AL}  C.-R. Ai and X.-J. Lin,
On the unitary structures of vertex operator superalgebras,
\emph{J. Algebra} \textbf{487} (2017), 217--243.

\bibitem[AlM]{AleM}  A. Alekseev and E. Meinrenken,
Lie theory and the Chern--Weil homomorphism,
\emph{Ann. Sci. \'Ecole Norm. Sup.} (4) \textbf{38} (2005), no.~2, 303--338.

\bibitem[AK]{AKl}  A.-B. Altman and S.-L. Kleiman,
\emph{A Term of Commutative Algebra},
Worldwide Center of Mathematics, Cambridge, MA, 2013; revised version of April 11, 2021.

\bibitem[AE]{AvE}  T. Arakawa and J. van Ekeren,
Modularity of relatively rational vertex algebras and fusion rules of principal affine $W$-algebras,
\emph{Comm. Math. Phys.} \textbf{370} (2019), no.~1, 205--247.

\bibitem[ALY1]{ALY1}  T. Arakawa, C.-H. Lam and H. Yamada,
Zhu's algebra, $C_2$-algebra and $C_2$-cofiniteness of parafermion vertex operator algebras,
\emph{Adv. Math.} \textbf{264} (2014), 261--295.

\bibitem[ALY2]{ALY2}  T. Arakawa, C.-H. Lam and H. Yamada,
Parafermion vertex operator algebras and $W$-algebras,
\emph{Trans. Amer. Math. Soc.} \textbf{371} (2019), no.~6, 4277--4301.

\bibitem[ArM]{AM2}  T. Arakawa and A. Moreau,
\emph{Arc Spaces and Vertex Algebras},
book-version lecture notes, 2024.

\bibitem[BJP]{BJP}  Y.-X. Bao, C.-P. Jiang and Y.-F. Pei,
Representations of affine Nappi--Witten algebras,
\emph{J. Algebra} \textbf{342} (2011), 111--133.

\bibitem[BDM]{BDM}  K. Barron, C.-Y. Dong and G. Mason,
Twisted sectors for tensor product vertex operator algebras associated to permutation groups,
\emph{Comm. Math. Phys.} \textbf{227} (2002), no.~2, 349--384.

\bibitem[B]{B}  N. Bourbaki,
\emph{Algebra I. Chapters 1--3},
Elem. Math. (Berlin), Springer-Verlag, Berlin, 1989.

\bibitem[CS]{CS}  S.-J. Cao and J.-C. Sun,
Decreasing filtrations, $C_2$-algebra and twisted modules,
{\hypersetup{urlcolor=black}\href{https://arxiv.org/abs/2511.01474v1}{arXiv:2511.01474v1}}.

\bibitem[CM]{CM}  S. Carnahan and M. Miyamoto,
Regularity of fixed-point vertex operator subalgebras,
{\hypersetup{urlcolor=black}\href{https://arxiv.org/abs/1603.05645}{arXiv:1603.05645}}.

\bibitem[CKLW]{CKLW}  S. Carpi, Y. Kawahigashi, R. Longo and M. Weiner,
From vertex operator algebras to conformal nets and back,
\emph{Mem. Amer. Math. Soc.} \textbf{254} (2018), no.~1213.

\bibitem[CDW]{CDW}  F.-L. Chen, L.-G. Ding and Q. Wang,
Vertex algebras and TKK algebras,
\emph{J. Algebra} \textbf{640} (2024), 147--173.

\bibitem[C]{Con}  K. Conrad,
\emph{Bilinear Forms}, expository notes, 2008.

\bibitem[CHY]{CHY}  T. Creutzig, Y.-Z. Huang and J.-W. Yang,
Braided tensor categories of admissible modules for affine Lie algebras,
\emph{Comm. Math. Phys.} \textbf{362} (2018), no.~3, 827--854.

\bibitem[DVVV]{DVVV}  R. Dijkgraaf, C. Vafa, E. Verlinde and H. Verlinde,
The operator algebra of orbifold models,
\emph{Comm. Math. Phys.} \textbf{123} (1989), no.~3, 485--526.

\bibitem[DHVW]{DHVW}  L. Dixon, J.-A. Harvey, C. Vafa and E. Witten,
Strings on orbifolds,
\emph{Nuclear Phys. B} \textbf{261} (1985), 678--686.

\bibitem[DGL]{DGL}  C.-Y. Dong, R.-L. Griess, Jr. and C.-H. Lam,
Uniqueness results for the moonshine vertex operator algebra,
\emph{Amer. J. Math.} \textbf{129} (2007), no.~2, 583--609.

\bibitem[DLWY]{DLWY}  C.-Y. Dong, C.-H. Lam, Q. Wang and H. Yamada,
The structure of parafermion vertex operator algebras,
\emph{J. Algebra} \textbf{323} (2010), no.~2, 371--381.

\bibitem[DLM1]{DLMreg}  C.-Y. Dong, H.-S. Li and G. Mason,
Regularity of rational vertex operator algebras,
\emph{Adv. Math.} \textbf{132} (1997), no.~1, 148--166.

\bibitem[DLM2]{DLMtw}  C.-Y. Dong, H.-S. Li and G. Mason,
Twisted representations of vertex operator algebras,
\emph{Math. Ann.} \textbf{310} (1998), no.~3, 571--600.

\bibitem[DLM3]{DLMmod}  C.-Y. Dong, H.-S. Li and G. Mason,
Modular-invariance of trace functions in orbifold theory and generalized moonshine,
\emph{Comm. Math. Phys.} \textbf{214} (2000), no.~1, 1--56.

\bibitem[DL]{DL}  C.-Y. Dong and X.-J. Lin,
Unitary vertex operator algebras,
\emph{J. Algebra} \textbf{397} (2014), 252--277.

\bibitem[DM1]{DMec}  C.-Y. Dong and G. Mason,
Rational vertex operator algebras and the effective central charge,
\emph{Int. Math. Res. Not.} \textbf{2004} (2004), no.~56, 2989--3008.

\bibitem[DM2]{DMinteg}  C.-Y. Dong and G. Mason,
Integrability of $C_2$-cofinite vertex operator algebras,
\emph{Int. Math. Res. Not.} \textbf{2006} (2006), Paper No.~80468, 15 pp.

\bibitem[DMZ]{DMZ}  C.-Y. Dong, G. Mason and Y.-C. Zhu,
Discrete series of the Virasoro algebra and the moonshine module,
in: \emph{Algebraic Groups and Their Generalizations: Quantum and Infinite-Dimensional Methods}
(University Park, PA, 1991), Proc. Sympos. Pure Math. \textbf{56}, Part~2,
Amer. Math. Soc., Providence, RI, 1994, \mbox{295--316}.

\bibitem[DN1]{DN1}  C.-Y. Dong and K. Nagatomo,
Classification of irreducible modules for the vertex operator algebra $M(1)^+$,
\emph{J. Algebra} \textbf{216} (1999), no.~1, 384--404.

\bibitem[DN2]{DN2}  C.-Y. Dong and K. Nagatomo,
Automorphism groups and twisted modules for lattice vertex operator algebras,
in: \emph{Recent Developments in Quantum Affine Algebras and Related Topics}
(Raleigh, NC, 1998), Contemp. Math. \textbf{248},
Amer. Math. Soc., Providence, RI, 1999, \mbox{117--133}.

\bibitem[DR]{DR}  C.-Y. Dong and L. Ren,
Congruence property in orbifold theory,
\emph{Proc. Amer. Math. Soc.} \textbf{146} (2018), no.~2, 497--506.

\bibitem[DRX]{DRX}  C.-Y. Dong, L. Ren and F. Xu,
On orbifold theory,
\emph{Adv. Math.} \textbf{321} (2017), 1--30.

\bibitem[DY]{DY}  C.-Y. Dong and N. Yu,
$\Z$-graded weak modules and regularity,
\emph{Comm. Math. Phys.} \textbf{316} (2012), no.~1, 269--277.

\bibitem[DuM]{DMu}  O. Dumitrescu and M. Mulase,
Edge-contraction on dual ribbon graphs, 2D TQFT, and the mirror of
orbifold Hurwitz numbers,
{\hypersetup{urlcolor=black}\href{https://arxiv.org/abs/1508.05922v1}{arXiv:1508.05922v1}}.

\bibitem[EKM]{EKM}  R. Elman, N. Karpenko and A. Merkurjev,
\emph{The Algebraic and Geometric Theory of Quadratic Forms},
Amer. Math. Soc. Colloq. Publ. \textbf{56},
Amer. Math. Soc., Providence, RI, 2008.

\bibitem[E]{Et}  P. Etingof,
\emph{Lie Groups and Lie Algebras},
lecture notes, MIT OpenCourseWare, 2024.

\bibitem[EGH]{EGH}  P. Etingof, O. Golberg, S. Hensel, T.-K. Liu,
A. Schwendner, D. Vaintrob and E.~Yudovina,
\emph{Introduction to Representation Theory},
Stud. Math. Libr. \textbf{59},
Amer. Math. Soc., Providence, RI, 2011.

\bibitem[FS]{FS}  J.-M. Figueroa-O'Farrill and S. Stanciu,
On the structure of symmetric self-dual Lie algebras,
\emph{J. Math. Phys.} \textbf{37} (1996), no.~8, 4121--4134.

\bibitem[FM]{FM}  D. Fiorenza and R. Murri,
Feynman diagrams via graphical calculus,
\emph{J. Knot Theory Ramifications} \textbf{11} (2002), no.~7, 1095--1131.

\bibitem[FHL]{FHL}  I.-B. Frenkel, Y.-Z. Huang and J. Lepowsky,
On axiomatic approaches to vertex operator algebras and modules,
\emph{Mem. Amer. Math. Soc.} \textbf{104} (1993), no.~494.

\bibitem[FLM]{FLM}  I.-B. Frenkel, J. Lepowsky and A. Meurman,
\emph{Vertex Operator Algebras and the Monster},
Pure Appl. Math. \textbf{134}, Academic Press, Boston, 1988.

\bibitem[FZ]{FZ}  I.-B. Frenkel and Y.-C. Zhu,
Vertex operator algebras associated to representations of affine and Virasoro algebras,
\emph{Duke Math. J.} \textbf{66} (1992), no.~1, 123--168.

\bibitem[G1]{G}  R.-L. Griess, Jr.,
The friendly giant,
\emph{Invent. Math.} \textbf{69} (1982), no.~1, 1--102.

\bibitem[G2]{GNAVOA}  R.-L. Griess, Jr.,
GNAVOA, I. Studies in groups, nonassociative algebras and vertex operator algebras,
in: \emph{Vertex Operator Algebras in Mathematics and Physics} (Toronto, ON, 2000),
Fields Inst. Commun. \textbf{39}, Amer. Math. Soc., Providence, RI, 2003, 71--88.

\bibitem[Hua]{HuangHaag}  T.-C. Huang,
Haagerup conformal field theory at central charge eight,\newline
{\hypersetup{urlcolor=black}\href{https://arxiv.org/abs/2610.14441}{arXiv:2610.14441}}.

\bibitem[Hum]{Hum}  J.-E. Humphreys,
\emph{Introduction to Lie Algebras and Representation Theory},
Grad. Texts in Math. \textbf{9}, Springer-Verlag, New York, 1972.

\bibitem[J]{Jacobson}  N. Jacobson,
\emph{Lie Algebras},
Dover Publications, New York, 1979, corrected reprint of the 1962 edition.

\bibitem[JLY]{JLY}  C.-P. Jiang, C.-H. Lam and H. Yamauchi,
The classification of vertex operator algebras of OZ-type generated by
Ising vectors of $\sigma$-type,
\emph{Adv. Math.} \textbf{465} (2025), Paper No.~110146.

\bibitem[JS]{JS}  C.-P. Jiang and J.-T. Song,
Associated varieties of simple affine VOAs $L_k(\mathfrak{sl}_3)$ and $W$-algebras $W_k(\mathfrak{sl}_3,f)$,
\emph{Comm. Math. Phys.} \textbf{406} (2025), no.~5, Paper No.~104.

\bibitem[JZ]{JZ}  W. Jiang and W. Zhang,
Vertex algebras associated with truncated Virasoro algebras and free field realizations,
\emph{Comm. Algebra} \textbf{54} (2026), no.~2, 722--735.

\bibitem[K1]{KacIDLA}  V.-G. Kac,
\emph{Infinite-Dimensional Lie Algebras}, third edition,
Cambridge University Press, Cambridge, 1990.

\bibitem[K2]{Kac}  V.-G. Kac,
\emph{Vertex Algebras for Beginners}, second edition,
Univ. Lecture Ser. \textbf{10}, Amer. Math. Soc., Providence, RI, 1998.

\bibitem[KL]{KL}  M. Karel and H.-S. Li,
Certain generating subspaces for vertex operator algebras,
\emph{J. Algebra} \textbf{217} (1999), no.~2, 393--421.

\bibitem[KKEG]{KKEG}  M. Kliesch, R. Kueng, J. Eisert and D. Gross,
Guaranteed recovery of quantum processes from few measurements,
\emph{Quantum} \textbf{3} (2019), Paper No.~171.

\bibitem[Lam1]{Lam}  C.-H. Lam,
Construction of vertex operator algebras from commutative associative algebras,
\emph{Comm. Algebra} \textbf{24} (1996), no.~14, 4339--4360.

\bibitem[Lam2]{LamT}  C.-H. Lam,
\emph{On the Structure of Vertex Operator Algebras and Their Weight Two Subspaces},
Ph.D. dissertation, The Ohio State University, 1996.

\bibitem[Lamt]{Lamt}  T.-Y. Lam,
\emph{A First Course in Noncommutative Rings},
second edition, Grad. Texts in Math. \textbf{131},
Springer-Verlag, New York, 2001.

\bibitem[LL]{LL}  J. Lepowsky and H.-S. Li,
\emph{Introduction to Vertex Operator Algebras and Their Representations},
Progr. Math. \textbf{227}, Birkh\"auser, Boston, 2004.

\bibitem[L1]{LiForm}  H.-S. Li,
Symmetric invariant bilinear forms on vertex operator algebras,
\emph{J. Pure Appl. Algebra} \textbf{96} (1994), no.~3, 279--297.

\bibitem[L2]{LiLocal}  H.-S. Li,
Local systems of vertex operators, vertex superalgebras and modules,
\emph{J. Pure Appl. Algebra} \textbf{109} (1996), no.~2, 143--195.

\bibitem[L3]{Li}  H.-S. Li,
Some finiteness properties of regular vertex operator algebras,
\emph{J. Algebra} \textbf{212} (1999), no.~2, 495--514.

\bibitem[LY]{LYu}  H.-S. Li and N. Yu,
On a family of vertex operator superalgebras,
\emph{J. Algebra} \textbf{608} (2022), 290--324.

\bibitem[Lia]{Lian}  B.-H. Lian,
On the classification of simple vertex operator algebras,
\emph{Comm. Math. Phys.} \textbf{163} (1994), no.~2, 307--357.

\bibitem[MR]{MR}  A. Medina and Ph. Revoy,
Alg\`ebres de Lie et produit scalaire invariant,
\emph{Ann. Sci. \'Ecole Norm. Sup.} (4) \textbf{18} (1985), no.~3, 553--561.

\bibitem[Mil]{Milas}  A. Milas,
Tensor product of vertex operator algebras,
{\hypersetup{urlcolor=black}\href{https://arxiv.org/abs/q-alg/9602026}{arXiv:q-alg/9602026}}.

\bibitem[Miy]{Miy}  M. Miyamoto,
$C_2$-cofiniteness of cyclic-orbifold models,
\emph{Comm. Math. Phys.} \textbf{335} (2015), no.~3, 1279--1286.

\bibitem[MS]{MS}  G.-W. Moore and G. Segal,
D-branes and K-theory in 2D topological field theory,
{\hypersetup{urlcolor=black}\href{https://arxiv.org/abs/hep-th/0609042}{arXiv:hep-th/0609042}}.

\bibitem[PP]{PP}  K.-B. Petersen and M.-S. Pedersen,
\emph{The Matrix Cookbook}, version of November~15, 2012.

\bibitem[P]{Primc}  M. Primc,
Vertex algebras generated by Lie algebras,
\emph{J. Pure Appl. Algebra} \textbf{135} (1999), no.~3, 253--293.

\bibitem[Ri]{Riesen}  A. Riesen,
The Haagerup VOA at $c=8$,
{\hypersetup{urlcolor=black}\href{https://arxiv.org/abs/2610.13933}{arXiv:2610.13933}}.

\bibitem[Roi]{Roi}  M. Roitman,
Invariant bilinear forms on a vertex algebra,
\emph{J. Pure Appl. Algebra} \textbf{194} (2004), no.~3, 329--345.

\bibitem[Rot]{Rotman}  J.-J. Rotman,
\emph{An Introduction to the Theory of Groups},
fourth edition, Grad. Texts in Math. \textbf{148}, Springer-Verlag, New York, 1995.

\bibitem[Sa]{Saw}  S. Sawin,
Direct sum decompositions and indecomposable TQFTs,
\emph{J. Math. Phys.} \textbf{36} (1995), no.~12, 6673--6680.

\bibitem[Sc]{Sch}  A. Schrijver,
On Lie algebra weight systems for 3-graphs,
\emph{J. Pure Appl. Algebra} \textbf{219} (2015), no.~10, 4597--4606.

\bibitem[SY]{SY}  J.-C. Sun and R.-C. Yang,
$g$-Rationality, $g$-regularity and $C_2$-cofiniteness of unitary
vertex operator algebras of CFT type,
submitted, 2026.

\bibitem[TW]{TW}  S.-T. Tsou and A.-G. Walker,
Metrisable Lie groups and algebras,
\emph{Proc. Roy. Soc. Edinburgh Sect. A} \textbf{64} (1957), no.~3, 290--304.

\bibitem[W]{Wang}  W.-Q. Wang,
Rationality of Virasoro vertex operator algebras,
\emph{Int. Math. Res. Not.} \textbf{1993} (1993), no.~7, 197--211.

\bibitem[Y]{Yamauchi}  H. Yamauchi,
Modularity on vertex operator algebras arising from semisimple primary vectors,
\emph{Internat. J. Math.} \textbf{15} (2004), no.~1, 87--109.

\bibitem[YS]{YS}  R.-C. Yang and J.-C. Sun,
$g$-Rationality, $g$-regularity and $C_2$-cofiniteness of vertex operator algebras,
submitted, 2025.

\bibitem[ZD]{ZD}  W. Zhang and C.-Y. Dong,
$W$-algebra $W(2,2)$ and the vertex operator algebra
$L(\frac12,0)\otimes L(\frac12,0)$,
\emph{Comm. Math. Phys.} \textbf{285} (2009), no.~3, 991--1004.

\bibitem[Z]{Zhu}  Y.-C. Zhu,
Modular invariance of characters of vertex operator algebras,
\emph{J. Amer. Math. Soc.} \textbf{9} (1996), no.~1, 237--302.

\end{thebibliography}
\end{document}